\documentclass[12pt,leqno]{article}
\usepackage{amssymb,amsmath,amsfonts,esint}
\usepackage{relsize}
\newtheorem{theorem}{Theorem}[section]
\newtheorem{problem}[theorem]{Problem}
\newtheorem{definition}[theorem]{Definition}
\newtheorem{corollary}[theorem]{Corollary}
\newtheorem{proposition}[theorem]{Proposition}
\newtheorem{remark}[theorem]{Remark}
\newtheorem{lemma}[theorem]{Lemma}
\newtheorem{claim}[theorem]{Claim}
\newtheorem{statement}[theorem]{Statement}
\newcommand {\Kc}      {{\mathcal K}}
\newcommand {\Dc}      {{\mathcal D}}
\newcommand {\Gc}      {{\mathcal G}}
\newcommand {\Mc}      {{\mathcal M}}
\newcommand {\Uc}      {{\mathcal U}}
\newcommand {\Vc}      {{\mathcal V}}
\newcommand {\Ac}      {{\mathcal A}}
\newcommand {\Tc}      {{\mathcal T}}
\newcommand {\Bc}      {{\mathcal B}}
\newcommand {\Cc}      {{\mathcal C}}
\newcommand {\Lc}      {{\mathcal L}}
\newcommand {\Ic}      {{\mathcal I}}
\newcommand {\Pc}      {{\mathcal P}}
\newcommand {\Qc}      {{\mathcal Q}}
\newcommand {\Rc}      {{\mathcal R}}
\newcommand {\Sc}      {{\mathcal S}}
\newcommand {\Yc}      {{\mathcal Y}}
\newcommand {\Xc}      {{\mathcal X}}
\newcommand {\mZ}      {{\mathbb Z}}
\newcommand {\LOP}     {L^1_p(\RT)}
\newcommand {\LOPN}    {L^1_p(\RN)}
\newcommand {\VLOP}    {\vec{L}^1_p(\RT)}
\newcommand {\LTP}     {L^2_p(\RT)}
\newcommand {\LMP}     {L^m_p(\RN)}
\newcommand {\LMPR}    {L^m_p(\R)}
\newcommand {\LMIR}    {L^m_\infty(\R)}
\newcommand {\LMIRN}   {L^m_\infty(\RN)}
\newcommand {\PMRN}    {\Pc_{m-1}(\RN)}
\newcommand {\PO}      {\Pc_1}
\newcommand {\LPRN}    {L_p(\RN)}
\newcommand {\LPRT}    {L_p(\RT)}
\newcommand {\LTIN}    {L^2_\infty(\RT)}
\newcommand {\LPM}     {L_p(\RT;\mu)}
\newcommand {\LPMN}    {L_p(\RN;\mu)}
\newcommand {\VLPM}    {\vec{L}_p(\RT;\mu)}
\newcommand {\VLPME}   {\vec{L}_p(\RT;\mu_E)}
\newcommand {\VSUM}    {\VLOP+\VLPME}
\newcommand {\VF}      {\vec{F}}
\newcommand {\VS}      {\vec{\Sigma}}
\newcommand {\CORT}    {C^{1,1}(\RT)}
\newcommand {\tK}      {\widetilde{K}}
\newcommand {\tT}      {\widetilde{T}}
\newcommand {\tQ}      {\widetilde{Q}}
\newcommand {\tQc}     {\widetilde{\Qc}}
\newcommand {\tSc}     {\widetilde{\Sc}}
\newcommand {\tL}      {\widetilde{{\mathcal L}}}
\newcommand {\tTc}     {\widetilde{{\mathcal T}}}
\newcommand {\hL}      {\hat{{\mathcal L}}}
\newcommand {\tg}      {\tilde{g}}
\newcommand {\tth}     {\tilde{h}}
\newcommand {\tdl}     {\tilde{\delta}}
\newcommand {\R}       {{\bf R}}
\newcommand {\N}       {{\bf N}}
\newcommand {\RN}      {\R^n}
\newcommand {\RT}      {\R^2}
\newcommand {\ve}      {\varepsilon}
\newcommand {\emp}     {\emptyset}
\newcommand {\intl}    {\int\limits}
\newcommand {\sq}      {square }
\newcommand {\sqs}     {squares }
\newcommand {\br}      {bridge }
\newcommand {\brs}     {bridges }
\newcommand {\ip}[1]   {\langle{#1}\rangle}
\newcommand {\q}       {90}
\newcommand {\lr}      {\leftrightarrow}
\newcommand {\LE}      {\Lc_E}
\newcommand {\GE}      {\Gc\LE}
\newcommand {\BRE}     {\Bc\Rc_E}
\newcommand {\GBE}     {\Gc\Bc_E}
\newcommand {\bcn}     {\leftrightsquigarrow}
\newcommand {\Br}      {T}
\newcommand {\QL}      {Q^{(L)}}
\newcommand {\ART}     {\Aff(\RT)}
\newcommand {\dg}      {degenerate }
\newcommand {\KE}      {\Kc_E}
\newcommand {\TF}      {\Tc(f)}
\newcommand {\hQ}      {\widehat{Q}}
\newcommand {\KFL}     {{\Kc_d}}
\newcommand {\SP}      {Sobolev-Poincar\'e }
\newcommand {\HLM}     {Hardy-Littlewood maximal }
\newcommand {\PQ }     {P^{(Q)}}
\newcommand {\A}[1]    {A^{[#1]}}
\newcommand {\B}[1]    {B^{[#1]}}
\newcommand {\WKE}     {\widetilde{\Kc}_E}
\newcommand {\TT}      {\tTc_1}
\newcommand {\smed}    {\mathlarger{\sum}}
\newcommand {\sbig}    {\mathlarger{\mathlarger{\sum}}}
\newcommand {\shuge}    {\mathlarger{\mathlarger{\mathlarger{\sum}}}}
\DeclareMathOperator{\esssup}{ess\,sup}
\DeclareMathOperator{\essinf}{ess\,inf}
\newcommand {\Id}      {\operatorname{Id}}
\newcommand {\cu}      {\mathlarger{\mathbf c}}
\newcommand {\card}    {\#}
\newcommand {\supp}    {\operatorname{supp}}
\newcommand {\diam}    {\operatorname{diam}}
\newcommand {\dist}    {\operatorname{dist}}
\newcommand {\PR}      {\operatorname{Pr}}
\newcommand {\Lip}     {\operatorname{Lip}}
\newcommand {\Aff}     {\operatorname{Aff}}
\newcommand {\Tri}     {\operatorname{Triangle}}
\newcommand {\PRL}     {\mathcal{PR}}
\newcommand {\bx}      {\hfill$\blacksquare$}
\newcommand {\rbx}     {\hfill$\vartriangleleft$}
\newcommand {\BX}      {\hspace{10mm}\blacksquare}
\newcommand {\nn}      {\nonumber}
\newcommand {\rf}[1]    {(\ref{#1})}      %no references
\newcommand {\reff}[1] {\ref{#1}}        %no references
\newcommand{\lbl}[1]    {\label{#1}}       %no ref
\newcommand{\be}        {\begin{eqnarray}}
\newcommand{\bel}[1]   {\begin{eqnarray} \label{#1}}
\newcommand{\ee}           {\end{eqnarray}}
\newcommand {\SECT}[2] {\section*{\centerline{\normalsize
{\bf #1}}} \setcounter{section}{#2}
\setcounter{theorem}{0}\setcounter{equation}{0}}
\newcommand {\SECTLONG}[3]
{\section*{\centerline{\normalsize {\bf #1}}
\centerline{\normalsize {\bf #2}}} \setcounter{section}{#3}
\setcounter{theorem}{0}\setcounter{equation}{0}}
\begin{document}
%----------------------------------------------------------
\parindent 1em
%----------------------------------------------------------
\parskip 0mm
%----------------------------------------------------------
\medskip
%@@@@@@@@@@@@@@@@@@@@@@@@@@@@@@@@@@@@@@@@@@@@@@@@@@@@@@@@@@
%@@@@@@@@@@@@@@@@@@@@@@@@@@@@@@@@@@@@@@@@@@@@@@@@@@@@@@@@@@
%----------------------------------------------------------
\centerline{\large{\bf Sobolev $L^2_p$-functions on closed subsets of $\R^2$}}\vspace*{5mm}
%----------------------------------------------------------
\vspace*{10mm} \centerline{By~  {\it Pavel Shvartsman}}\vspace*{5 mm}
%----------------------------------------------------------
\centerline {\it Department of Mathematics, Technion - Israel Institute of Technology}\vspace*{2 mm}
%----------------------------------------------------------
\centerline{\it 32000 Haifa, Israel}\vspace*{2 mm}
%----------------------------------------------------------
\centerline{\it e-mail: pshv@tx.technion.ac.il}
%----------------------------------------------------------
\vspace*{10 mm}
%----------------------------------------------------------
\renewcommand{\thefootnote}{ }
%----------------------------------------------------------
\footnotetext[1]{{\it\hspace{-6mm}Math Subject
Classification} 46E35}
%@@@@@@@@@@@@@@@@@@@@@@@@@@@@@@@@@@@@@@@@@@@@@@@@@@@@@@@@@@
%@@@@@@@@@@@@@@@@@@@@@@@@@@@@@@@@@@@@@@@@@@@@@@@@@@@@@@@@@@
%----------------------------------------------------------
\begin{abstract} For each $p>2$ we give intrinsic characterizations of the restriction of the homogeneous Sobolev space $\LTP$ to an arbitrary finite subset $E$ of $\RT$. The trace criterion is expressed in terms of certain weighted oscillations of the second order with respect to a measure generated by the Menger curvature of triangles with vertices in $E$.
%----------------------------------------------------------
\end{abstract}
%----------------------------------------------------------
\renewcommand{\contentsname}{ }
\tableofcontents
%----------------------------------------------------------
\addtocontents{toc}{{\centerline{\sc{Contents}}}
\vspace*{10mm}\par}
%----------------------------------------------------------
%@@@@@@@@@@@@@@@@@@@@@@@@@@@@@@@@@@@@@@@@@@@@@@@@@@@@@@@@@@
%@@@@@@@@@@@@@@@@@@@@@@@@@@@@@@@@@@@@@@@@@@@@@@@@@@@@@@@@@@
%@@@@@@@@@@@@@@@@@@@@@@@@@@@@@@@@@@@@@@@@@@@@@@@@@@@@@@@@@@
%@@@@@@@@@@@@@@@@@@@@@@@@@      @@@@@@@@@@@@@@@@@@@@@@@@@@@
%@@@@@@@@@@@@@@@@@@@@@@@          @@@@@@@@@@@@@@@@@@@@@@@@@
%@@@@@@@@@@@@@@@@@@@@@              @@@@@@@@@@@@@@@@@@@@@@@
%@@@@@@@@@@@@@@@@@@@     SECTION 1    @@@@@@@@@@@@@@@@@@@@@
%@@@@@@@@@@@@@@@@@@@@@              @@@@@@@@@@@@@@@@@@@@@@@
%@@@@@@@@@@@@@@@@@@@@@@@          @@@@@@@@@@@@@@@@@@@@@@@@@
%@@@@@@@@@@@@@@@@@@@@@@@@@      @@@@@@@@@@@@@@@@@@@@@@@@@@@
%@@@@@@@@@@@@@@@@@@@@@@@@@@@@@@@@@@@@@@@@@@@@@@@@@@@@@@@@@@
%@@@@@@@@@@@@@@@@@@@@@@@@@@@@@@@@@@@@@@@@@@@@@@@@@@@@@@@@@@
%@@@@@@@@@@@@@@@@@@@@@@@@@@@@@@@@@@@@@@@@@@@@@@@@@@@@@@@@@@
%----------------------------------------------------------
\SECT{1. Introduction.}{1}
%----------------------------------------------------------
\addtocontents{toc}{~~~1. Introduction. \hfill \thepage\\\par}
%----------------------------------------------------------
\indent
%@@@@@@@@@@@@@@@@@@@@@@@@@@@@@@@@@@@@@@@@@@@@@@@@@@@@@@@@@@
%----------------------------------------------------------
\par Let $p\in[1,\infty]$ and let $m$ be a positive integer. By $\LMP$ we denote the homogeneous Sobolev space
consisting of all (equivalence classes of) real valued
functions $F\in L_{p,loc}(\RN)$ whose distributional partial derivatives of order $m$ belong to the space $\LPRN$. The space $\LMP$ is equipped with the seminorm
%----------------------------------------------------------
$$
\|F\|_{\LMP}:=\|\nabla^mF\|_{\LPRN}
$$
%----------------------------------------------------------
where
%----------------------------------------------------------
$$
\nabla^mF(x):=\left(\sum_{|\alpha|= m}(D^\alpha F(x))^2\right)^{\frac{1}{2}},~~~~x\in\RN.
$$
%----------------------------------------------------------
\par When $p>n$, it follows from the Sobolev embedding theorem that every function $F\in\LMP$ coincides almost everywhere with a $C^{m-1}$-function.  This fact enables us {\it to identify each element $F\in \LMP$, $p>n$, with its unique $C^{m-1}$-representative}. In particular, this identification implies that $F$ {\it has a well defined restriction to any given subset of $\RN$.} It also enables us to identify the Sobolev space $\LMIRN$ with the
space $C^{m-1,1}(\RN)$ of $C^{m-1}$-functions whose partial derivatives of order $m-1$ are Lipschitz continuous on $\RN$. The space $C^{m-1,1}(\RN)$ is equipped with the seminorm
%----------------------------------------------------------
$$
\|F\|_{C^{m-1,1}(\RN)}=
\sum_{|\alpha|=m-1}\|D^\alpha F\|_{\Lip(\RN)}.
$$
%----------------------------------------------------------
\par In this paper we study the following
%----------------------------------------------------------
%@@@@@@@@@@@@@@@@@@@@@@@@@@@@@@@@@@@@@@@@@@@@@@@@@@@@@@@@@@
%@@@@@@@@@@@@@@@@@@@@@@@@@@@@@@@@@@@@@@@@@@@@@@@@@@@@@@@@@@
%@@@@@@@@@@@@@@@@@@@@@@@@@@@@@@@@@@@@@@@@@@@@@@@@@@@@@@@@@@
%@@@@@@@@@@@@@@@@@@@@@@@@@@@@@@@@@@@@@@@@@@@@@@@@@@@@@@@@@@
%----------------------------------------------------------
\begin{problem}\lbl{MAIN-PR} {\em Given a finite set $E\subset\mathbf{R}^{2}$ and a function $f:E\to\mathbf{R}$, we consider the $\LTP$-norms of all $C^{1}$-functions $F:\mathbf{R}^{2}\to\mathbf{R}$ which coincide with $f$ on $E$. How small can they be?}
%----------------------------------------------------------
\end{problem}
%----------------------------------------------------------
%@@@@@@@@@@@@@@@@@@@@@@@@@@@@@@@@@@@@@@@@@@@@@@@@@@@@@@@@@@
%@@@@@@@@@@@@@@@@@@@@@@@@@@@@@@@@@@@@@@@@@@@@@@@@@@@@@@@@@@
%@@@@@@@@@@@@@@@@@@@@@@@@@@@@@@@@@@@@@@@@@@@@@@@@@@@@@@@@@@
%----------------------------------------------------------
\par We denote the infimum of all these norms by $\|f\|_{\LTP|_E}$; thus
%----------------------------------------------------------
$$
\|f\|_{\LTP|_E}
:=\inf \{\|F\|_{\LTP}:F\in\LTP\cap C^1(\RT), F|_{E}=f\}.
$$
%----------------------------------------------------------
As is customary, we refer to $\|f\|_{\LTP|_E}$ as the {\it trace norm of the function $f$} (in $\LTP$). This quantity provides the standard quotient space seminorm in {\it the trace space} $\LTP|_{E}$ of all restrictions of $\LTP$-functions to $E$, i.e., in the space %@@@@@@@@@@@@@@@@@@@@@@@@@@@@@@@@@@@@@@@@@@@@@@@@@@@@@@@@@@
$$
\LTP|_{E}:=\{f:E\to\R:\text{there exists}~~F\in
\LTP\cap C^1(\RT)\ \ \text{such that}\ \ F|_{E}=f\}.
$$
%@@@@@@@@@@@@@@@@@@@@@@@@@@@@@@@@@@@@@@@@@@@@@@@@@@@@@@@@@@
\par Problem \reff{MAIN-PR} is a variant of a classical extension problem posed by H. Whitney in 1934 in his pioneering papers \cite{W1,W2}, namely: {\it How can one tell whether a given function $f$ defined on an arbitrary subset $E\subset\RN$ extends to a $C^m$-function on all of $\RN$?} Over the years since 1934 this problem, often called the Whitney Extension Problem, has attracted a lot of attention, and there is an extensive literature devoted to different aspects of this problem and its analogues for various spaces of smooth functions. Among the multitude of results we mention the papers by G. Glaeser \cite{Gl},
Y. Brudnyi, P. Shvartsman \cite{BS1,BS2,BS3,S-D,S1,S-Tr,S-GAN}, E. Bierstone, P. Milman and W. Pawlucki \cite{BM,BMP1}, and N. Zobin \cite{Z1,Z2}. In the last decade C. Fefferman \cite{F2,F3,F4,F-J,F-IO,F-Bl} made several important
breakthroughs in this area and developed, partially in collaboration with B. Klartag, a series of new directions of investigation related to computational and algorithmic
aspects of the Whitney Extension Problem, see \cite{F6,F7,F8,F-Bl}. We refer the reader to all of the above-mentioned papers, and references therein, for numerous results and techniques concerning this topic.
%----------------------------------------------------------
\par All of these papers deal with functions whose partial derivatives are bounded and continuous or in $L_\infty$.
It is natural to ask analogous questions when $L_\infty$ is replaced by $L_p$, i.e., to study the extension and restriction properties of functions whose partial derivatives belong to the space $L_p$ with $p<\infty$. In the papers \cite{S-W1,S4} we have considered such questions and presented several constructive descriptions of the trace space $\LOPN|_E$ for an arbitrary subset $E\subset\RN$ provided $p>n$. (See Remark \reff{L-1} for more details.)\medskip
%----------------------------------------------------------
\par In this paper for $n=2$ we extend these results from the case of first order derivatives to second order derivatives. More explicitly, we solve Problem \reff{MAIN-PR} by presenting a constructive formula for calculation of the order of magnitude of the trace norm $\|f\|_{\LTP|_E}$. This formula is expressed only in terms of the values of the function $f$ on $E$ and certain geometric characteristics of the set $E$.
%----------------------------------------------------------
\par We also prove the existence of a {\it continuous linear extension operator} from ${\LTP|_E}$ to ${\LTP}$.
Note that the first results related to the existence of continuous linear extension operators acting on traces of
Sobolev space to {\it arbitrary} closed subsets were recently obtained by A. Israel \cite{Is} and by C. Fefferman, A. Israel and G. K. Luli \cite{FIL}. (These results are discussed after Theorem \reff{LINEXT}.)\smallskip
%----------------------------------------------------------
\par Before we formulate the main result of the paper
we need to define several notions and fix some notation:
%----------------------------------------------------------
\par Throughout this paper, the word ``square'' will mean a closed \sq in $\RT$ whose sides are parallel to the
coordinate axes. For each \sq $Q$ we let $c_Q$ denote its center. Given $\lambda>0$ we let $\lambda\, Q$ denote the dilation of $Q$ with respect to its center by a factor of $\lambda $. The Lebesgue measure of a measurable set $A\subset \RT$ will be denoted by $\left|A\right|$. By $\# A$ we denote the number of elements of a finite set $A$.
Let $\Ac$ be a family of sets in $\RN$. By $M(\Ac)$ we denote its covering multiplicity, i.e., the minimal positive integer $M$ such that every point $x\in\RN$ is covered by at most $M$ sets from $\Ac$.
%----------------------------------------------------------
\par By $\|\cdot\|$ and $\|\cdot\|_2$ we denote, respectively, the uniform and the Euclidean norms in $\RT$. For each pair of points $z_1$ and $z_2$ in $\RT$ we let $(z_1,z_2)$ denote the open line segment joining them.
%----------------------------------------------------------
\par The word ``triangle" will always mean a subset $\Delta=\{z_1,z_2,z_3\}\subset\RT$ consisting of three non-collinear points. By $\Tri(E)$ we denote the family of all triangles with vertices in $E$.
%----------------------------------------------------------
\par We let $R_\Delta$ denote the radius of the circle passing through the points $z_1,z_2,z_3.$
The reciprocal of $R_\Delta$, the quantity
%----------------------------------------------------------
\bel{M-CURV}
\cu_\Delta:=\frac{1}{R_\Delta},
\ee
%----------------------------------------------------------
is called {\it the Menger curvature} of the ``triangle" $\Delta$ (as discussed in detail, e.g. in \cite{P}).
%----------------------------------------------------------
\par Given a function $f:E\to\R$ and a triangle $\Delta=\{z_1,z_2,z_3\}\subset E$ we let $P_\Delta[f]$ denote  the affine polynomial which interpolates $f$ on $\Delta$. Thus
%----------------------------------------------------------
$$
P_\Delta[f]\in\PO~~~~\text{and}~~~~
P_\Delta[f](z_i)=f(z_i),~i=1,2,3.
$$
%----------------------------------------------------------
Here $\PO$ denotes the space of polynomials on $\RT$ of degree at most one.
%----------------------------------------------------------
%@@@@@@@@@@@@@@@@@@@@@@@@@@@@@@@@@@@@@@@@@@@@@@@@@@@@@@@@@@
%@@@@@@@@@@@@@@@@@@@@@@@@@@@@@@@@@@@@@@@@@@@@@@@@@@@@@@@@@@
%@@@@@@@@@@@@@@@@@@@@@@@@@@@@@@@@@@@@@@@@@@@@@@@@@@@@@@@@@@
%@@@@@@@@@@@@@@@@@@@@@@@@@@@@@@@@@@@@@@@@@@@@@@@@@@@@@@@@@@
%@@@@@@@@@@@@@@@@@@@@@@@@@@@@@@@@@@@@@@@@@@@@@@@@@@@@@@@@@@
%@@@@@@@@@@@@@@@@@@@@@@@@@@@@@@@@@@@@@@@@@@@@@@@@@@@@@@@@@@
%@@@@@@@@@@@@@@@@@@@@@@@@@@@@@@@@@@@@@@@@@@@@@@@@@@@@@@@@@@
%----------------------------------------------------------
\par Here now is the main result of our paper:
%----------------------------------------------------------
\begin{theorem}\lbl{MAIN} Let $2<p<\infty$. Let $E$ be a finite subset of $\RT$ and let $f$ be a function defined on $E$. Then
%----------------------------------------------------------
\bel{M-EQN}
\|f\|_{\LTP|_E}\sim \inf \lambda^{\frac{1}{p}}
\ee
%----------------------------------------------------------
where the infimum is taken over all positive constants $\lambda$ which satisfy all of the following conditions for a certain absolute positive constant $\gamma$:
%----------------------------------------------------------
\par (i). For every finite family $\{Q_i: i=1,...,m\}$ of pairwise disjoint \sqs and every choice of
collinear points $z^{(i)}_1,z^{(i)}_2,z^{(i)}_3\in E$ such that $z^{(i)}_2\in(z^{(i)}_1,z^{(i)}_3)$ and
%----------------------------------------------------------
$$
z^{(i)}_j\in \gamma Q_i,~~i=1,...,m,~j=1,2,3,
$$
%----------------------------------------------------------
the following inequality
%----------------------------------------------------------
\bel{DDIF}
\sbig_{i=1}^m \left|
\frac{f(z^{(i)}_1)-f(z^{(i)}_2)}{\|z^{(i)}_1-z^{(i)}_2\|_2}
-\frac{f(z^{(i)}_2)-f(z^{(i)}_3)}{\|z^{(i)}_2-z^{(i)}_3\|_2}
\right|^p
(\diam Q_i)^{2-p}\le \lambda
\ee
%----------------------------------------------------------
holds.\medskip
%----------------------------------------------------------
%$$
%\hspace*{-15mm}\shuge_{i=1}^\ell
%\left|\frac{\mathlarger{
%\frac{f(z^{(i)}_1)-f(z^{(i)}_2)}{\|z^{(i)}_1-z^{(i)}_2\|_2}
%-\frac{f(z^{(i)}_2)-f(z^{(i)}_3)}{\|z^{(i)}_2-z^{(i)}_3\|_2}}}
%{\diam Q_i}
%\right|^p|Q_i|\le \lambda
%$$
%----------------------------------------------------------
%@@@@@@@@@@@@@@@@@@@@@@@@@@@@@@@@@@@@@@@@@@@@@@@@@@@@@@@@@@
%----------------------------------------------------------
\par (ii). Let $\Qc$ and $\Kc$ be arbitrary finite families of pairwise disjoint squares. Suppose that to each \sq $K\in\Kc$ we have arbitrarily assigned a triangle $\Delta(K)$ in $E$ such that
%----------------------------------------------------------
\bel{GSEP}
\Delta(K)\subset\gamma K~~~\text{and}~~~\diam K\le\gamma\diam \Delta(K).
\ee
%----------------------------------------------------------
\par Suppose that to each \sq $Q\in\Qc$ we have arbitrarily assigned two \sqs $Q',Q''\in\Qc$ such that\,\,
$Q'\cup Q''\subset \gamma Q.$
%----------------------------------------------------------
\par Let
%----------------------------------------------------------
\bel{SGM}
\sigma_p(Q;\Kc):=
\sum\left\{\cu_{\Delta(K)}^p|K|: K\in\Kc, c_K\in\,Q\right\},~~~~~Q\in\Qc,
\ee
%----------------------------------------------------------
%@@@@@@@@@@@@@@@@@@@@@@@@@@@@@@@@@@@@@@@@@@@@@@@@@@@@@@@@@@
%----------------------------------------------------------
and
%----------------------------------------------------------
$$
S_p(f:Q',Q'';\Kc):=
\sum_{\substack {K'\in\,\Kc\\c_{K'}\in Q'}}\,\,
\sum_{\substack {K''\in\,\Kc\\c_{K''}\in Q''}}
\|\nabla P_{\Delta(K')}[f]-
\nabla P_{\Delta(K'')}[f]\|^p\,
\cu_{\Delta(K')}^p |K'|\,\cu_{\Delta(K'')}^p |K''|.
$$
%----------------------------------------------------------
\par Then the following inequality
%----------------------------------------------------------
$$
\sbig\limits_{Q\in\Qc}\,
\frac{(\diam Q)^{2-p}\,\,S_p(f:Q',Q'';\Kc)}
{\left\{(\diam Q')^{2-p} +\sigma_p(Q';\Kc)\right\}
\left\{(\diam Q'')^{2-p} +\sigma_p(Q'';\Kc)\right\}}
\le\lambda
$$
%----------------------------------------------------------
holds.
%----------------------------------------------------------
\par The constants of equivalence in \rf{M-EQN} depend only on $p$.
%----------------------------------------------------------
\end{theorem}
%----------------------------------------------------------
\medskip
%----------------------------------------------------------
%@@@@@@@@@@@@@@@@@@@@@@@@@@@@@@@@@@@@@@@@@@@@@@@@@@@@@@@@@@
%@@@@@@@@@@@@@@@@@@@@@@@@@@@@@@@@@@@@@@@@@@@@@@@@@@@@@@@@@@
%@@@@@@@@@@@@@@@@@@@@@@@@@@@@@@@@@@@@@@@@@@@@@@@@@@@@@@@@@@
%\par  PICTURE 2: ILLUSTRATION TO THE THEOREM                                   %@@@@@@@@@@@@@@@@@@@@@@@@@@@@@@@@@@@@@@@@@@@@@@@@@@@@@@@@@@
%@@@@@@@@@@@@@@@@@@@@@@@@@@@@@@@@@@@@@@@@@@@@@@@@@@@@@@@@@@
%@@@@@@@@@@@@@@@@@@@@@@@@@@@@@@@@@@@@@@@@@@@@@@@@@@@@@@@@@@
%@@@@@@@@@@@@@@@@@@@@@@@@@@@@@@@@@@@@@@@@@@@@@@@@@@@@@@@@@@
%----------------------------------------------------------
%@@@@@@@@@@@@@@@@@@@@@@@@@@@@@@@@@@@@@@@@@@@@@@@@@@@@@@@@@@
%----------------------------------------------------------
\begin{remark} {\em For the peace of mind of any particularly pedantic reader, we mention that ultimately we can choose the constant $\gamma$ to equal, for example, $2^{22}$.}\rbx
\end{remark}
%----------------------------------------------------------
\par The trace criterion given in this theorem describes  the structure of the trace space $\LTP|_E$ and shows which
properties of a function $f$ on $E$ control its almost optimal extension to a function from $\LTP$. At the same time it is not quite clear how one could check the conditions of part (i) and part (ii) of Theorem \reff{MAIN} for a given function $f$ on $E$. In fact, these conditions depend on an infinite number of families of \sqs ($\Qc,\Kc$ etc.), triangles and triples of collinear points.
%----------------------------------------------------------
\par Nevertheless a careful examination of our proof shows that it constructs {\it four particular families}  of squares, triangles and points, etc., depending only on the set $E$ and that it is enough to examine the behavior of given functions $f$ only on these particular families.
%%----------------------------------------------------------
\par It is convenient to express this fact by the following theorem, which refines one part of Theorem \reff{MAIN}.
We are very grateful to Charles Fefferman for conjecturing and motivating us to seek a result along these lines.
%----------------------------------------------------------
%@@@@@@@@@@@@@@@@@@@@@@@@@@@@@@@@@@@@@@@@@@@@@@@@@@@@@@@@@@
%@@@@@@@@@@@@@@@@@@@@@@@@@@@@@@@@@@@@@@@@@@@@@@@@@@@@@@@@@@
%@@@@@@@@@@@@@@@@@@@@@@@@@@@@@@@@@@@@@@@@@@@@@@@@@@@@@@@@@@
%@@@@@@@@@@@@@@@@@@@@@@@@@@@@@@@@@@@@@@@@@@@@@@@@@@@@@@@@@@
%@@@@@@@@@@@@@@@@@@@@@@@@@@@@@@@@@@@@@@@@@@@@@@@@@@@@@@@@@@
%@@@@@@@@@@@@@@@@@@@@@@@@@@@@@@@@@@@@@@@@@@@@@@@@@@@@@@@@@@
%@@@@@@@@@@@@@@@@@@@@@@@@@@@@@@@@@@@@@@@@@@@@@@@@@@@@@@@@@@
%----------------------------------------------------------
%@@@@@@@@@@@@@@@@@@@@@@@@@@@@@@@@@@@@@@@@@@@@@@@@@@@@@@@@@@
%----------------------------------------------------------
\begin{theorem}\lbl{REF-MAIN} Let $2<p<\infty$ and let $E$ be a finite subset of $\RT$. There exist absolute constants $\gamma >0$, $C>0$ and $N \in \N$ and\,:
\medskip
%----------------------------------------------------------
\par (i) A family $\{Q_i:i=1,...,m\}$ with $m\le C \,\#E$ of pairwise disjoint \sqs and a family
%----------------------------------------------------------
$$
\{z^{(i)}_{1},z^{(i)}_{2},z^{(i)}_{3}\in E\cap (\gamma Q_i): z^{(i)}_{2}\in(z^{(i)}_{1},z^{(i)}_{3}),\,\, i=1,...,m\}
$$
%---------------------------------------------------------
of triples of collinear points; \medskip
%----------------------------------------------------------
\par (ii) A family $\Kc$ of pairwise disjoint \sqs with $\#\Kc\le C\,\# E$, and a mapping
%----------------------------------------------------------
$$
\Kc\ni K\mapsto \Delta(K)\in \Tri(E)
$$
%----------------------------------------------------------
such that $\Delta(K)\subset \gamma K$ and $\diam K\le\gamma \diam \Delta(K)$ for every $K\in\Kc$;\medskip
%----------------------------------------------------------
%@@@@@@@@@@@@@@@@@@@@@@@@@@@@@@@@@@@@@@@@@@@@@@@@@@@@@@@@@@
%----------------------------------------------------------
\par (iii) A family of \sqs $\Qc$ with covering multiplicity $M(\Qc)\le N$ and $\#\Qc\le C\,\# E$, and mappings $\Qc\ni Q\mapsto Q'\in \Qc$ and $\Qc\ni Q\mapsto Q''\in \Qc$ satisfying the condition $Q'\cup Q''\subset \gamma Q$ for all $Q\in\Qc$, \medskip
%----------------------------------------------------------
\par such that for every function $f:E\to\R$ the following equivalence
%----------------------------------------------------------
\be
\|f\|_{\LTP|_E}&\sim&
\left(\sbig_{i=1}^{m} \left|
\frac{f(z^{(i)}_{1})-f(z^{(i)}_{2})}
{\|z^{(i)}_{1}-z^{(i)}_{2}\|_2}
-\frac{f(z^{(i)}_{2})-f(z^{(i)}_{3})}
{\|z^{(i)}_{2}-z^{(i)}_{3}\|_2}
\right|^p
(\diam Q_i)^{2-p}\right)^{\frac1p}\nn\\
&+&
\left(\sbig\limits_{Q\in\Qc}\,
\frac{(\diam Q)^{2-p}\,\,S_p(f:Q',Q'';\Kc)}
{\left\{(\diam Q')^{2-p} +\sigma_p(Q';\Kc)\right\}
\left\{(\diam Q'')^{2-p} +\sigma_p(Q'';\Kc)\right\}}\right)^{\frac1p}
\nn
\ee
%----------------------------------------------------------
holds. The constants of this equivalence depend only on $p$.
%----------------------------------------------------------
\end{theorem}
%----------------------------------------------------------
%@@@@@@@@@@@@@@@@@@@@@@@@@@@@@@@@@@@@@@@@@@@@@@@@@@@@@@@@@@
%@@@@@@@@@@@@@@@@@@@@@@@@@@@@@@@@@@@@@@@@@@@@@@@@@@@@@@@@@@
%@@@@@@@@@@@@@@@@@@@@@@@@@@@@@@@@@@@@@@@@@@@@@@@@@@@@@@@@@@
%@@@@@@@@@@@@@@@@@@@@@@@@@@@@@@@@@@@@@@@@@@@@@@@@@@@@@@@@@@
%----------------------------------------------------------
\par A structure of the trace norms of functions defined on a finite set $E\subset\RN$ in the space $L^m_p(\RN)|_E$ has been discussed and studied by A. Israel in \cite{Is} and by C. Fefferman, A. Israel and G. K. Luli in \cite{FIL}. In these papers the authors introduce a notion of a linear functional with the so-called ``assisted bounded depth'' and characterize the trace norm of a function on $E$ using this notion. They also formulate several open problems related to the structure of the trace norm. Let us formulate one of them.
%---------------------------------------------------------
\begin{problem}\lbl{RT-DEPTH}(A. Israel \cite{Is}) Given finite set $E\subset\R^2$ do there exist linear functionals
%----------------------------------------------------------
$$
\{\lambda_i\}_{i=1}^L\subset (\LTP|_E)^*,~~~~
L\le C\,\# E,
$$
%----------------------------------------------------------
each depending only on $N$ values such that
%----------------------------------------------------------
$$
\|f\|_{\LTP|_E}^p\sim \sum_{i=1}^L\,\,|\lambda_i(f)|^p
$$
%----------------------------------------------------------
for every function $f:E\to\R$\,?
\par Here $N>0$ is a certain absolute constant.
%---------------------------------------------------------
\end{problem}
%---------------------------------------------------------
\par Our next result provides a solution to this problem with $N=6$ whenever $2<p<\infty$.
%---------------------------------------------------------
\begin{theorem}\lbl{REF-SP} Let $2<p<\infty$ and let $E\subset\RT$ be a finite set. There exist linear functionals
%---------------------------------------------------------
$$
\{\lambda_1,\lambda_2,...,\lambda_L\},~~~~ L\le C\card E,
$$
%---------------------------------------------------------
each depending on at most six values of a function on $E$ such that for every function $f:E\to\R$ the following equivalence
%----------------------------------------------------------
$$
\|f\|_{\LTP|_E}^p\sim \mathlarger{\sum}_{i=1}^L\,\,|\lambda_i(f)|^p
$$
%----------------------------------------------------------
holds.
%----------------------------------------------------------
\par The constants of this equivalence depend only on $p$.
%---------------------------------------------------------
\end{theorem}
%---------------------------------------------------------
\par The proof of this theorem and Theorem \reff{REF-MAIN}, including an explicit construction of the objects which they mention and use, are given in Section 11. Note that the objects described in parts (i), (ii) and (iii) of Theorem \reff{REF-MAIN} and the linear functionals $\{\lambda_i\}$ from Theorem \reff{REF-SP} depend only on $E$ and $p$.
%----------------------------------------------------------
%@@@@@@@@@@@@@@@@@@@@@@@@@@@@@@@@@@@@@@@@@@@@@@@@@@@@@@@@@@
%@@@@@@@@@@@@@@@@@@@@@@@@@@@@@@@@@@@@@@@@@@@@@@@@@@@@@@@@@@
%@@@@@@@@@@@@@@@@@@@@@@@@@@@@@@@@@@@@@@@@@@@@@@@@@@@@@@@@@@
%@@@@@@@@@@@@@@@@@@@@@@@@@@@@@@@@@@@@@@@@@@@@@@@@@@@@@@@@@@
%----------------------------------------------------------
\par In the next four remarks we will briefly review several previous results about Sobolev extensions which are related to Theorem \reff{MAIN}.
%----------------------------------------------------------
%@@@@@@@@@@@@@@@@@@@@@@@@@@@@@@@@@@@@@@@@@@@@@@@@@@@@@@@@@@
%@@@@@@@@@@@@@@@@@@@@@@@@@@@@@@@@@@@@@@@@@@@@@@@@@@@@@@@@@@
%@@@@@@@@@@@@@@@@@@@@@@@@@@@@@@@@@@@@@@@@@@@@@@@@@@@@@@@@@@
%@@@@@@@@@@@@@@@@@@@@@@@@@@@@@@@@@@@@@@@@@@@@@@@@@@@@@@@@@@
%----------------------------------------------------------
\begin{remark}\lbl{P-INF} {\em  The case $p=\infty$ was studied in the author's papers \cite{S-82,S-D,S1,S-GAN}. (Recall that $\LTIN$ coincides with the space $\CORT$ of all $C^1$-functions whose partial derivatives of the first order are Lipschitz continuous on $\RT$.) In those papers we proved that a function $f$ defined on an arbitrary closed set $E\subset\RT$ extends to a function $F\in\LTIN$ if and only if there exists a constant  $\lambda>0$ such that:\medskip
%----------------------------------------------------------
\par (a). For every three collinear points $z_1,z_2,z_3\in E$ such that $z_2\in(z_1,z_3)$ the following inequality
%----------------------------------------------------------
$$
\left|
\frac{f(z_1)-f(_2)}{\|z_1-z_2\|_2}
-\frac{f(z_2)-f(z_3)}{\|z_2-z_3\|_2}
\right|\le \lambda \|z_1-z_3\|_2
$$
%----------------------------------------------------------
holds;
%----------------------------------------------------------
\par (b). For every two triangles $\Delta_1$ and $\Delta_2$ with vertices in $E$ we have
%------------------------------------------------------------
$$
\|\nabla P_{\Delta_1}[f]-\nabla P_{\Delta_2}[f]\|\le
\lambda\{R_{\Delta_1}+R_{\Delta_2}+
\diam (\Delta_1\cup \Delta_2)\}.
$$
%------------------------------------------------------------
\par Furthermore, $\|f\|_{\LTIN|_E}\sim \inf \lambda$.\bigskip
%------------------------------------------------------------
\par Note that in parts (a) and (b) of this criterion we use the values of the function $f$ at no more than $6$ points of the set $E$. This phenomenon is a consequence of the so-called ``finiteness principle'' for the space $\CORT$ proven in \cite{S-D,S1}. (See also \cite{BS3}.) This principle provides a criterion for calculation of the trace norm $\|f\|_{\CORT|_E}$ in the following form: Let $N=6$. Then
%------------------------------------------------------------
$$
\|f\|_{\CORT|_E}\sim \sup\{\|f|_{E'}\|_{\CORT|_{E'}}: E'\subset E, \# E'\le N\}.
$$
%------------------------------------------------------------
\par The same formula is true for the space $L^2_\infty(\RN) =C^{1,1}(\RN)$ with $N=3\cdot 2^{n-1}$ (and this value is sharp), see \cite{S1}.
%------------------------------------------------------------
\par The finiteness principle for the space $L^m_\infty(\RN) =C^{m-1,1}(\RN)$ with arbitrary $m$ and
$n$ and a certain constant $N=N(m,n)$ has been proved by C. Fefferman \cite{F2}.\rbx
%------------------------------------------------------------
\medskip}
%------------------------------------------------------------
%@@@@@@@@@@@@@@@@@@@@@@@@@@@@@@@@@@@@@@@@@@@@@@@@@@@@@@@@@@@@
%------------------------------------------------------------
\end{remark}
%------------------------------------------------------------
%@@@@@@@@@@@@@@@@@@@@@@@@@@@@@@@@@@@@@@@@@@@@@@@@@@@@@@@@@@@@
%@@@@@@@@@@@@@@@@@@@@@@@@@@@@@@@@@@@@@@@@@@@@@@@@@@@@@@@@@@@@
%@@@@@@@@@@@@@@@@@@@@@@@@@@@@@@@@@@@@@@@@@@@@@@@@@@@@@@@@@@@@
%@@@@@@@@@@@@@@@@@@@@@@@@@@@@@@@@@@@@@@@@@@@@@@@@@@@@@@@@@@@@
%------------------------------------------------------------
\begin{remark}\lbl{L-1} {\em Whitney-type extension problems for the Sobolev space $\LOPN$ have been studied in the author's papers \cite{S-W1,S4}. We have proved that {\it the classical extension operator  constructed by Whitney in \cite{W1} for the space $C^1(\RN)$ provides an almost optimal extension operator for the trace space $\LOPN|_E$ whenever $n<p<\infty$} and $E$ is an arbitrary closed subset of $\RN$. This enables us to give several constructive trace criteria. Let us mention two of them here.
%----------------------------------------------------------
\par In particular, a function $f:E\to\R$ can be extended to a (continuous) function $F\in\LOPN$, $n<p<\infty,$ if and only if there exists a constant $\lambda>0$ such that for every finite family $\{Q_{i}:i=1,...,m\}$ of pairwise disjoint cubes in $\RN$ and every choice of points $x_{i},y_{i}\in (11Q_{i})\cap E$ the following inequality
%----------------------------------------------------------
$$
\smed_{i=1}^{m}\frac{|f(x_{i})-f(y_{i})|^{p}}
{(\operatorname{diam}Q_{i})^{p-n}}\le \lambda
$$
%----------------------------------------------------------
holds. Furthermore, $\|f\|_{\LOPN|_E}\sim \inf\lambda^{\frac1p}$.\medskip
%----------------------------------------------------------
\par The second criterion is expressed as an explicit formula for the order of magnitude of the trace norm of a function $f:E\to\R$ in the space $\LOPN$ whenever $p\in(n,\infty)$:
%----------------------------------------------------------
\bel{CR-LOP}
\|f\|_{L^1_p(\RN)|_E}\sim \left\{\,\,\intl_{\RN}\left(\,\sup_{y,z\in E}\frac{|f(y)-f(z)|}{\|x-y\|+\|x-z\|}\right)^p\,dx
\right\}^{\frac{1}{p}}\,.
\ee
%---------------------------------------------------------
See \cite{S-W1}, Theorem 1.4, and \cite{S4}, Theorem 6.1.
%---------------------------------------------------------
\par In both of these criteria the constants of the equivalences depend only on $n$ and $p$.\medskip
%----------------------------------------------------------
\par We remark that in \cite{S-W1} we also present a series of trace criteria expressed in terms of local oscillations of functions with respect to certain {\it doubling measures} supported on $E$. These results have been inspired by a paper of  A. Jonsson \cite{J}. In particular, in \cite{J} A. Jonsson presents a trace criterion of a similar kind which provides a constructive characterization of the restrictions of Besov spaces $B^\alpha_{p,q}(\RN)$ to an arbitrary closed set $E\subset\RN$ whenever $0<\alpha<1$ and $p\,\alpha>n$. \rbx
%------------------------------------------------------------
\medskip}
%------------------------------------------------------------
%@@@@@@@@@@@@@@@@@@@@@@@@@@@@@@@@@@@@@@@@@@@@@@@@@@@@@@@@@@@@
%------------------------------------------------------------
\end{remark}
%------------------------------------------------------------
%@@@@@@@@@@@@@@@@@@@@@@@@@@@@@@@@@@@@@@@@@@@@@@@@@@@@@@@@@@@@
%@@@@@@@@@@@@@@@@@@@@@@@@@@@@@@@@@@@@@@@@@@@@@@@@@@@@@@@@@@@@
%@@@@@@@@@@@@@@@@@@@@@@@@@@@@@@@@@@@@@@@@@@@@@@@@@@@@@@@@@@@@
%@@@@@@@@@@@@@@@@@@@@@@@@@@@@@@@@@@@@@@@@@@@@@@@@@@@@@@@@@@@@
%@@@@@@@@@@@@@@@@@@@@@@@@@@@@@@@@@@@@@@@@@@@@@@@@@@@@@@@@@@@@
%@@@@@@@@@@@@@@@@@@@@@@@@@@@@@@@@@@@@@@@@@@@@@@@@@@@@@@@@@@@@
%@@@@@@@@@@@@@@@@@@@@@@@@@@@@@@@@@@@@@@@@@@@@@@@@@@@@@@@@@@@@
%------------------------------------------------------------
\begin{remark} {\em Let $\{P_x\in\PMRN: x\in E\}$ be a  polynomial field, i.e., a mapping which to every point $x\in E$ assigns a polynomial $P_x$ on $\RN$ of degree at most $m-1$. Given an $m$-times differentiable function $F$ and a point $x\in \RN$,	 we let $T^m_x[F]$ denote the Taylor polynomial of $F$ at $x$ of order $m$. We recall one of Whitney classical extension theorems \cite{W1} which states that
%---------------------------------------------------------
\be
&&\inf\left\{\|F\|_{L^m_\infty(\RN)}:F\in L^m_\infty(\RN), T^{m-1}_x[F]=P_x~~\text{for every}~~x\in E\right\}\nn\\\nn\\
&\sim&\sum_{|\alpha|\le\, m-1}\,\,
\sup_{y,z\in E} |D^{\alpha}(P_y-P_z)(y)|/\|y-z\|^{m-|\alpha|}\,\,.
\nn
\ee
%---------------------------------------------------------
See also Glaeser \cite{G}.
%---------------------------------------------------------
\par In \cite{S4} we generalize this theorem to the case of the Sobolev space $\LMP$, $p>n$, proving that for every polynomial field $\{P_x\in\PMRN: x\in E\}$ the following equivalence
%---------------------------------------------------------
\be
&&\inf\left\{\|F\|_{L^m_p(\RN)}:F\in \LMP, T^{m-1}_x[F]=P_x~~\text{for every}~~x\in E\right\}\nn\\\nn\\
&\sim& \,\,\smed_{|\alpha|\le m-1}\,\,
\left(\,\,\int\limits_{\RN}
\left(\,\,\,
\sup_{y,z\in E}\,\, \frac{|D^{\alpha}(P_y-P_z)(y)|}{\|x-y\|^{m-|\alpha|}+
\|x-z\|^{m-|\alpha|}}\right)^pdx\right)^{\frac{1}{p}}\nn
\ee
%---------------------------------------------------------
holds with constants depending only on $m,p,$ and $n$. Clearly, for $m=1$ we obtain \rf{CR-LOP}.
%----------------------------------------------------------
\par As in the case $m=1$ (see Remark \reff{L-1}) we also show in \cite{S4} that the Whitney extension operator constructed in \cite{W1} for the space $C^m(\RN)$ provides an almost optimal extension for every polynomial field $\{P_x\in\PMRN: x\in E\}$ on $E$  whenever $n<p<\infty$ and $E$ is an arbitrary closed subset of $\RN$.\rbx
\medskip}
%------------------------------------------------------------
%@@@@@@@@@@@@@@@@@@@@@@@@@@@@@@@@@@@@@@@@@@@@@@@@@@@@@@@@@@@@
%------------------------------------------------------------
\end{remark}
%------------------------------------------------------------
%@@@@@@@@@@@@@@@@@@@@@@@@@@@@@@@@@@@@@@@@@@@@@@@@@@@@@@@@@@@@
%@@@@@@@@@@@@@@@@@@@@@@@@@@@@@@@@@@@@@@@@@@@@@@@@@@@@@@@@@@@@
%@@@@@@@@@@@@@@@@@@@@@@@@@@@@@@@@@@@@@@@@@@@@@@@@@@@@@@@@@@@@
%@@@@@@@@@@@@@@@@@@@@@@@@@@@@@@@@@@@@@@@@@@@@@@@@@@@@@@@@@@@@
%------------------------------------------------------------
%@@@@@@@@@@@@@@@@@@@@@@@@@@@@@@@@@@@@@@@@@@@@@@@@@@@@@@@@@@@@
%@@@@@@@@@@@@@@@@@@@@@@@@@@@@@@@@@@@@@@@@@@@@@@@@@@@@@@@@@@@@
%@@@@@@@@@@@@@@@@@@@@@@@@@@@@@@@@@@@@@@@@@@@@@@@@@@@@@@@@@@@@
%@@@@@@@@@@@@@@@@@@@@@@@@@@@@@@@@@@@@@@@@@@@@@@@@@@@@@@@@@@@@
%------------------------------------------------------------
\begin{remark}\lbl{L-R1} {\em In \cite{S4} we study a Whitney-type extension problem for the space $\LMPR$, $p>1$. Let $E\subset\R$ be an arbitrary closed set and let $T_E$ be the extension operator  constructed by Whitney in \cite{W2} for the trace space $\LMIR|_E=C^{m-1,1}(\R)|_E$. We show that this very same Whitney extension operator $T_E$ also provides an almost optimal extension of each function $f\in\LMPR|_E$ whenever $1<p<\infty.$
%------------------------------------------------------------
\par This leads us to an analogue of Whitney trace criterion \cite{W2} for the space $L^m_\infty(\R)$. (See also J. Merrien \cite{Mer}.) We recall that this classical result  states that, for arbitrary positive integer $m$ and every function $f$ on $E$,
%----------------------------------------------------------
$$
\|f\|_{L^m_\infty(\R)|_E}\sim
\sup_{S\subset E,\,\,\card S=m+1}
|\Delta^mf[S]|.
$$
%------------------------------------------------------------
In the preceding formula $\Delta^mf[S]$ denotes {\it the divided difference} of $f$ on $S$, i.e., for
each finite set $S=\{x_0,x_1,...,x_m\}\subset \R$, it is given by
%---------------------------------------------------------
$$
\Delta^mf[S]=\Delta^mf[x_0,x_1,...,x_{m}]
=\sum_{i=0}^m~\frac{f(x_i)}
{\prod\limits_{0\le k\le m,\,\,k\neq i}(x_i-x_k)}\,\,.
$$
%---------------------------------------------------------
\par In \cite{S4} we find a counterpart of the preceding result, in which we replace $L^m_\infty$ by $L^m_p$. We
prove that for arbitrary $1<p<\infty$ and every function $f$ on $E$
%----------------------------------------------------------
\be
\|f\|_{L^m_p(\R)|_E}&\sim& \left\{ \int\limits_{\R}\sup_{S\subset E,\,\,\card S=m+1}\left(\frac{|\Delta^mf[S]|\diam S}{\diam (\{x\}\cup S)}\right)^pdx\right\}^{\frac{1}{p}}\nn\\\nn\\
&\sim& \left\{ \int\limits_{\R}\sup_{x_0<x_1<...<x_m,\,\,x_i\in E}\,\,\frac{|\,\Delta^{m-1}f[x_0,...,x_{m-1}]-
\Delta^{m-1}f[x_1,...,x_{m}]|^{\,p}}{|x-x_0|^p+|x-x_m|^{p}}
\,\,dx\right\}^{\frac{1}{p}}\nn
\ee
%---------------------------------------------------------
with constants of equivalence depending only on $p$ and $m$.\rbx}
%------------------------------------------------------------
\end{remark}
%----------------------------------------------------------
%@@@@@@@@@@@@@@@@@@@@@@@@@@@@@@@@@@@@@@@@@@@@@@@@@@@@@@@@@@
%@@@@@@@@@@@@@@@@@@@@@@@@@@@@@@@@@@@@@@@@@@@@@@@@@@@@@@@@@@
%@@@@@@@@@@@@@@@@@@@@@@@@@@@@@@@@@@@@@@@@@@@@@@@@@@@@@@@@@@
%@@@@@@@@@@@@@@@@@@@@@@@@@@@@@@@@@@@@@@@@@@@@@@@@@@@@@@@@@@
%@@@@@@@@@@@@@@@@@@@@@@@@@@@@@@@@@@@@@@@@@@@@@@@@@@@@@@@@@@
%@@@@@@@@@@@@@@@@@@@@@@@@@@@@@@@@@@@@@@@@@@@@@@@@@@@@@@@@@@
%----------------------------------------------------------
\par Our next result, Theorem \reff{NEW-MAIN}, provides a slight modification of the trace criterion given in Theorem \reff{MAIN}. More specifically, part (ii) of the new criterion is expressed in terms of families of {\it disks} with certain geometric constraints on the curvature of their boundaries.
%----------------------------------------------------------
\par Before we formulate Theorem \reff{NEW-MAIN} let us
define several additional notions and fix additional notation. For a disk $D$ we let $c_D$ and $r_D$ denote its center and radius respectively. By $\cu_D$ we denote the curvature of the boundary of the disk $D$, i.e., the reciprocal of its radius:
%----------------------------------------------------------
\bel{C-DS}
\cu_D=\frac{1}{r_D}\,.
\ee
%----------------------------------------------------------
We refer to $\cu_D$ as {\it the curvature of the disk $D$}.
%----------------------------------------------------------
\par Let $p\in[1,\infty)$. Let $\Dc$ be a family of pairwise disjoint disks in $\RT$ and let $\Delta$ be a mapping which to every disk $D\in\Dc$ assigns a triangle $\Delta(D)\subset\RT$.
%----------------------------------------------------------
\par Given a disk $B\subset\RT$ we define a quantity $\cu_p(B:\Dc,\Delta)$ by the formula:
%----------------------------------------------------------
$$
\cu_p(B:\Dc,\Delta):=
\left\{\frac{1}{|B|}
\sum\limits_{c_D\in B}\cu^p_{\Delta(D)}|D|\right\}^{\frac1p}\,.
$$
%----------------------------------------------------------
We refer to $\cu_p(B:\Dc,\Delta)$ as {\it the $p$-average of the Menger curvature on $B$} with respect to the family $\Dc$ and the mapping $\Delta$.
%----------------------------------------------------------
%@@@@@@@@@@@@@@@@@@@@@@@@@@@@@@@@@@@@@@@@@@@@@@@@@@@@@@@@@@
%@@@@@@@@@@@@@@@@@@@@@@@@@@@@@@@@@@@@@@@@@@@@@@@@@@@@@@@@@@
%@@@@@@@@@@@@@@@@@@@@@@@@@@@@@@@@@@@@@@@@@@@@@@@@@@@@@@@@@@
%@@@@@@@@@@@@@@@@@@@@@@@@@@@@@@@@@@@@@@@@@@@@@@@@@@@@@@@@@@
%@@@@@@@@@@@@@@@@@@@@@@@@@@@@@@@@@@@@@@@@@@@@@@@@@@@@@@@@@@
%@@@@@@@@@@@@@@@@@@@@@@@@@@@@@@@@@@@@@@@@@@@@@@@@@@@@@@@@@@
%@@@@@@@@@@@@@@@@@@@@@@@@@@@@@@@@@@@@@@@@@@@@@@@@@@@@@@@@@@
%----------------------------------------------------------
\begin{theorem}\lbl{NEW-MAIN} Let $2<p<\infty$. Let $E$ be a finite subset of $\RT$ and let $f$ be a function defined on $E$. Then
%----------------------------------------------------------
\bel{N-E}
\|f\|_{\LTP|_E}\sim \inf \lambda^{\frac{1}{p}}
\ee
%----------------------------------------------------------
where the infimum is taken over all positive constants $\lambda$ which satisfy all of the following conditions for a certain absolute positive constant $\gamma$:\smallskip
%----------------------------------------------------------
\par (i). The condition (i) of Theorem \reff{MAIN} holds;\medskip
%----------------------------------------------------------
\par (ii). Let $\Bc$ and $\Dc$ be arbitrary finite families of pairwise disjoint disks. Suppose that to each disk $D\in\Dc$ we have arbitrarily assigned a triangle $\Delta(D)$ in $E$ such that
%----------------------------------------------------------
$$
\Delta(D)\subset\gamma D~~~\text{and}~~~\diam D\le\gamma\diam \Delta(D).
$$
%----------------------------------------------------------
\par Suppose that to each disk $B\in\Bc$ we have arbitrarily assigned two disks $B',B''\in\Bc$ such that\,\, $B'\cup B''\subset \gamma B$ and
%----------------------------------------------------------
\bel{B-CV}
\cu_p(B':\Dc,\Delta)\le \cu_{B'}~~~~\text{and}~~~~\cu_p(B'':\Dc,\Delta)\le \cu_{B''}.
\ee
%----------------------------------------------------------
%@@@@@@@@@@@@@@@@@@@@@@@@@@@@@@@@@@@@@@@@@@@@@@@@@@@@@@@@@@
%----------------------------------------------------------
Let
%----------------------------------------------------------
$$
S_p(f:B',B'';\Dc):=
\sum_{\substack {D'\in\,\Dc\\c_{D'}\in B'}}\,\,
\sum_{\substack {D''\in\,\Dc\\c_{D''}\in B''}}
\|\nabla P_{\Delta(D')}[f]-
\nabla P_{\Delta(D'')}[f]\|^p\,
\cu_{\Delta(D')}^p |D'|\,\cu_{\Delta(D'')}^p |D''|.
$$
%----------------------------------------------------------
\par Then the following inequality
%----------------------------------------------------------
$$
\smed\limits_{B\in\Bc}\,
\left(\frac{\diam B' \diam B''}{\diam B} \right)^{p-2} S_p(f:B',B'';\Dc)
\le\lambda
$$
%----------------------------------------------------------
holds.
%----------------------------------------------------------
\par The constants of equivalence in \rf{N-E} depend only on $p$.
%----------------------------------------------------------
\end{theorem}
%----------------------------------------------------------
%@@@@@@@@@@@@@@@@@@@@@@@@@@@@@@@@@@@@@@@@@@@@@@@@@@@@@@@@@@@@
%@@@@@@@@@@@@@@@@@@@@@@@@@@@@@@@@@@@@@@@@@@@@@@@@@@@@@@@@@@@@
%@@@@@@@@@@@@@@@@@@@@@@@@@@@@@@@@@@@@@@@@@@@@@@@@@@@@@@@@@@@@
%@@@@@@@@@@@@@@@@@@@@@@@@@@@@@@@@@@@@@@@@@@@@@@@@@@@@@@@@@@@@
%------------------------------------------------------------
\par We prove this theorem in Section 11.
%------------------------------------------------------------
\bigskip
%------------------------------------------------------------
\par The proof of Theorem \reff{MAIN} uses the extension methods developed for the Sobolev space $L^2_\infty(\RN)=C^{1,1}(\RN)$, see \cite{S1,S-GAN}.
It also uses Theorem \reff{SL-DEC} about sums of Sobolev and $L_p$-weighted spaces. (See also the discussion in Section 9.)
%------------------------------------------------------------
\par The methods of the proof of Theorem \reff{MAIN} and some related ideas enable us to construct an almost optimal extension algorithm which {\it linearly} depends on a function from the trace space $\LTP|_E$. Thus we can give a new and different proof of the following recent theorem of Arie Israel \cite{Is}.
%----------------------------------------------------------
%@@@@@@@@@@@@@@@@@@@@@@@@@@@@@@@@@@@@@@@@@@@@@@@@@@@@@@@@@@
%@@@@@@@@@@@@@@@@@@@@@@@@@@@@@@@@@@@@@@@@@@@@@@@@@@@@@@@@@@
%----------------------------------------------------------
\begin{theorem} \lbl{LINEXT} For every finite subset
$E\subset \RT$ and every $p>2$ there exists a linear extension operator which maps the trace space $\LTP|_E$ continuously into $\LTP$. Its operator norm is bounded by a constant depending only on $p$.
%----------------------------------------------------------
\end{theorem}
%@@@@@@@@@@@@@@@@@@@@@@@@@@@@@@@@@@@@@@@@@@@@@@@@@@@@@@@@@@
%@@@@@@@@@@@@@@@@@@@@@@@@@@@@@@@@@@@@@@@@@@@@@@@@@@@@@@@@@@
%@@@@@@@@@@@@@@@@@@@@@@@@@@@@@@@@@@@@@@@@@@@@@@@@@@@@@@@@@@
%@@@@@@@@@@@@@@@@@@@@@@@@@@@@@@@@@@@@@@@@@@@@@@@@@@@@@@@@@@
%@@@@@@@@@@@@@@@@@@@@@@@@@@@@@@@@@@@@@@@@@@@@@@@@@@@@@@@@@@
%@@@@@@@@@@@@@@@@@@@@@@@@@@@@@@@@@@@@@@@@@@@@@@@@@@@@@@@@@@
%----------------------------------------------------------
\par It seems worthwhile at this point to briefly describe and compare the methods used in these two different proofs. (Both A. Israel and the author initially presented details of their respective approaches to this result at the workshop ``Differentiable structures on finite sets'' organized by American Institute of Mathematics, Palo Alto, CA, in 2010.)
%----------------------------------------------------------
\par In \cite{Is} A. Israel generalizes and introduces new elements into an approach developed earlier by C. Fefferman \cite{F2,F-COM} in his proof of the ``finiteness principle'' for the trace space $\LMIRN|_E=C^{m-1,1}(\RN)|_E$. (See also \cite{F-IL,F3,F-IO}.) One of the main ingredients of this approach is an elaborate Calder\'on-Zygmund type decomposition of $\RT$ into appropriate
``Calder\'on-Zygmund" cubes each containing a portion of $E$ which is in a certain sense ``convenient'' for extension. In Fefferman's earlier proof, which deals with the case $p=\infty$, each of the ``Calder\'on-Zygmund" cubes are treated on an ``equal footing"; no cube has a more important role than any of the other cubes. However, in \cite{Is}, for the case $p<\infty$, a variant of the above-mentioned decomposition yields cubes, in fact squares since $n=2$, which are no longer all on an equal footing. A certain special subfamily of them, which are referred to as Keystone squares, have an important role, and all the others can be ignored.  The Keystone squares concentrate all information  which is necessary for an almost optimal extension of a function defined on $E$ to a function from the Sobolev space $\LTP$.\smallskip
%----------------------------------------------------------
\par Let us now turn to the author's different approach to the proof of Theorem \reff{LINEXT} whose details appear below in Section 2. The reader will be able to note some features which are common to both approaches. However, the starting point here comes from a method developed and used earlier in the papers \cite{S-82,S1,S-GAN} where the author proved the above-mentioned finiteness principle for the space $\LTIN=C^{1,1}(\RT)$. (Cf. Remark \ref{P-INF}.) After modifying this method to cope with the case $p<\infty$ we are led to a new problem of constructing an almost optimal decomposition of a given function into a sum of a Sobolev and a weighted $L_p$-function. The particular weight function which must be used here is generated by the Menger curvature of certain triangles with vertices in $E$.
The solution of a rather more general version of this problem (for an arbitrary weight function or measure) is presented separately in \cite{S3}. The main tool there is a filtering procedure which determines a special family of \sqs in $\RT$ (we call them ``important squares''). These \sqs possess certain measure concentration properties and together contain all necessary information about extension properties of a function defined on $E$.
%----------------------------------------------------------
\par At this stage we do not know if there is any relation between the Keystone squares introduced in \cite{Is} and the ``important \sqs'' with respect to the Menger curvature measure which are treated in this paper. Perhaps, for a given set $E$ they even coincide exactly with each other. We share the opinion expressed in \cite{Is} and \cite{FIL} that it would be interesting and useful to investigate possible connections between the Keystone \sqs  and the ``important'' squares.
%----------------------------------------------------------
\par Let us mention some other results which are related to Theorem \reff{LINEXT}. We note that the existence of a continuous linear extension operator for the space $L^2_\infty(\RN)=C^{1,1}(\RN)$ was shown by Yu. Brudnyi and the author \cite{BS2}. C. Fefferman \cite{F-IL} proved that such an operator exists for the space $L^m_\infty(\RN)$ for arbitrary $m$.
%----------------------------------------------------------
\par Quite recently  C. Fefferman, A. Israel and G. K. Luli \cite{FIL} proved that a continuous linear extension operator exists for the trace space $\LMP|_E$  whenever $n<p<\infty$ and $E\subset\RN$ is an arbitrary closed set. This remarkable result also relies on ideas and techniques related to the above mentioned Calder\'on-Zygmund type decomposition and certain properties of Keystone cubes.
%----------------------------------------------------------
\par Finally we note that G. K. Luli \cite{L} proved the existence of a continuous linear extension operator for the space $L^m_p(\R)|_E,$ $p>1,$ for every {\it finite} set $E\subset\R$. As we have noted in Remark \reff{L-R1} the extension operator $T_E$ constructed by H. Whitney \cite{W2} for the space $C^m(\R)|_E$ provides an almost optimal extension for every trace space $L^m_p(\R)|_E$ and {\it every} closed set $E\subset\R$ whenever $p>1$. (See \cite{S4}.) Since $T_E$ is {\it linear}, this proves the existence of a continuous linear extension operator $L^m_p(\R)|_E,$ $p>1,$ for an {\it arbitrary} closed set $E\subset\R$.\smallskip
%----------------------------------------------------------
\medskip
%@@@@@@@@@@@@@@@@@@@@@@@@@@@@@@@@@@@@@@@@@@@@@@@@@@@@@@@@@@
\par {\bf Acknowledgements.} I am very thankful to M. Cwikel for useful suggestions and remarks. I am also very grateful to C. Fefferman, N. Zobin, A. Israel, G. K. Luli, and all participants of the ``Whitney Problems Workshops" in Palo-Alto, August 2010, and in Williamsburg, August 2011, for stimulating discussions and valuable advice. In particular, I thank C. Fefferman for some very helpful remarks (already mentioned above) in connection with Theorem \reff{REF-MAIN}. I am pleased to thank B. Klartag for very useful discussions related to his proof of Theorem \reff{SPR-P} and various ideas of sparsification.
%----------------------------------------------------------
%@@@@@@@@@@@@@@@@@@@@@@@@@@@@@@@@@@@@@@@@@@@@@@@@@@@@@@@@@@
%------------------------------------------------------------
%@@@@@@@@@@@@@@@@@@@@@@@@@@@@@@@@@@@@@@@@@@@@@@@@@@@@@@@@@@@@
%@@@@@@@@@@@@@@@@@@@@@@@@@@@@@@@@@@@@@@@@@@@@@@@@@@@@@@@@@@@@
%@@@@@@@@@@@@@@@@@@@@@@@@@@@@@@@@@@@@@@@@@@@@@@@@@@@@@@@@@@@@
%@@@@@@@@@@@@@@@@@@@@@@@@@@@@@@@@@@@@@@@@@@@@@@@@@@@@@@@@@@@@
%----------------------------------------------------------
%@@@@@@@@@@@@@@@@@@@@@@@@@@@@@@@@@@@@@@@@@@@@@@@@@@@@@@@@@@
%@@@@@@@@@@@@@@@@@@@@@@@@@@@@@@@@@@@@@@@@@@@@@@@@@@@@@@@@@@
%@@@@@@@@@@@@@@@@@@@@@@@@@@@@@@@@@@@@@@@@@@@@@@@@@@@@@@@@@@
%@@@@@@@@@@@@@@@@@@@@@@@@@      @@@@@@@@@@@@@@@@@@@@@@@@@@@
%@@@@@@@@@@@@@@@@@@@@@@@          @@@@@@@@@@@@@@@@@@@@@@@@@
%@@@@@@@@@@@@@@@@@@@@@              @@@@@@@@@@@@@@@@@@@@@@@
%@@@@@@@@@@@@@@@@@@@     SECTION 2    @@@@@@@@@@@@@@@@@@@@@
%@@@@@@@@@@@@@@@@@@@@@              @@@@@@@@@@@@@@@@@@@@@@@
%@@@@@@@@@@@@@@@@@@@@@@@          @@@@@@@@@@@@@@@@@@@@@@@@@
%@@@@@@@@@@@@@@@@@@@@@@@@@      @@@@@@@@@@@@@@@@@@@@@@@@@@@
%@@@@@@@@@@@@@@@@@@@@@@@@@@@@@@@@@@@@@@@@@@@@@@@@@@@@@@@@@@
%@@@@@@@@@@@@@@@@@@@@@@@@@@@@@@@@@@@@@@@@@@@@@@@@@@@@@@@@@@
%@@@@@@@@@@@@@@@@@@@@@@@@@@@@@@@@@@@@@@@@@@@@@@@@@@@@@@@@@@
%----------------------------------------------------------
\SECT{2. Plan of the proof of the Main Theorem \reff{MAIN}.}{2}
%----------------------------------------------------------
\addtocontents{toc}{2. Plan of the proof of the Main Theorem \reff{MAIN}. \hfill \thepage\\\par}
%----------------------------------------------------------
\indent
%@@@@@@@@@@@@@@@@@@@@@@@@@@@@@@@@@@@@@@@@@@@@@@@@@@@@@@@@@@
%----------------------------------------------------------
\par Let us briefly describe the main stages of the proof of Theorem \reff{MAIN}.\bigskip
%-----------------------------------------------------------
\par In {\bf Section 3} we prove that for every function $f$ on $E$ the conditions (i) and (ii) of this theorem are satisfied with $\lambda=C\|f\|_{\LTP|_E}^p$ where $C$ is a constant depending only on $p$. We refer to this part of the proof of Theorem \reff{MAIN} as the ``necessity'' part.
%-----------------------------------------------------------
\par  More specifically, we fix a constant $\gamma\in[1,\infty)$ and a function $F\in\LTP, p>2,$ such that $F|_E=f$ and show that the statements (i) and (ii) of Theorem \reff{MAIN} hold with $\lambda=C(p,\gamma)\|F\|_{\LTP}^p$. The proof of these statements is based on the classical \SP inequality for Sobolev functions which we recall in Proposition \reff{SP-IN}. Let us give some remarks related to the proof of part (ii); the part (i) is much more simpler and its proof follows from the \SP inequality and the \HLM theorem.
%-----------------------------------------------------------
\par The point of departure of our proof  of part (ii) is Proposition \reff{TR1}. Given a \sq $Q$, a triangle
$\Delta\subset Q$ and a point $x\in Q$ this proposition provides a bound of the distance between $\nabla F(x)$ and $\nabla P_{\Delta}[F]$. Recall that $P_{\Delta}[F]$ is the affine polynomial which interpolates $F$ at the vertices of $\Delta$. This bound depends on the radius $R_\Delta$ of the circle passing through the vertices of $\Delta.$ Recall that the Menger curvature $\cu_\Delta=1/R_\Delta$ so that $\cu_\Delta$ is implicitly involved in this bound.
%-----------------------------------------------------------
\par Let $\Kc$ be a family of pairwise disjoint \sqs and let $\Cc_{\Kc}=\{c_K:~K\in\Kc\}$ be the family of its centers. Given $K\in\Kc$ let $\Delta(K)$ be a triangle satisfying condition \rf{GSEP} of Theorem \reff{MAIN}. Proposition \reff{TR1} motivates us to introduces two important objects - a mapping $\Tc:\RT\to\RT$ and a non-negative Borel measure $\mu$ on $\RT$. The mapping $\Tc$ is supported on the set $\Cc_{\Kc}$ and for each
\sq $K\in\Kc$ it takes the value $\Tc(c_K):=\nabla P_{\Delta(K)}[F]$ at the center of $K$. The measure $\mu$ is a discrete measure supported on $\Cc_{\Kc}$ such that $\mu(\{c_K\})=\cu_{\Delta(K)}^p\,|K|$ for each $K\in\Kc$
%-----------------------------------------------------------
\par Using Proposition \reff{TR1} and the \HLM theorem we prove that the mapping $\Tc$ belongs to the sum $\VS$ of the vector Sobolev space $\VLOP$ and the vector $L_p$-space $\VLPM$ with respect to the measure $\mu$. Furthermore, we show that $\|\Tc\|_{\VS}\le C(p)\|F\|_{\LTP}.$ This is a crucial point of the proof of the necessity.
%-----------------------------------------------------------
\par In \cite{S3} we present a general formula for calculation of the norm of a function in the sum of the Sobolev space $\LOP$, $p>n$, and a space $\LPMN$ where $\mu$ is an arbitrary non-negative Borel measure on $\RN$, see Theorem \reff{SL-DEC}. Applying this theorem to the mapping $\Tc$ we immediately obtain the required statement of part (ii) of Theorem \reff{MAIN}.\medskip
%-----------------------------------------------------------
\par For the reader's convenience, we also give a direct proof of this statement which does not use results of the work \cite{S3}. This proof relies only on the classical \HLM theorem and Proposition \reff{2TR} which provides a Sobolev-Poincar\'{e} type inequality for pairs of interpolating affine polynomials in $\RT$. \bigskip
%------------------------------------------------------------
%@@@@@@@@@@@@@@@@@@@@@@@@@@@@@@@@@@@@@@@@@@@@@@@@@@@@@@@@@@@@
%@@@@@@@@@@@@@@@@@@@@@@@@@@@@@@@@@@@@@@@@@@@@@@@@@@@@@@@@@@@@
%@@@@@@@@@@@@@@@@@@@@@@@@@@@@@@@@@@@@@@@@@@@@@@@@@@@@@@@@@@@@
%@@@@@@@@@@@@@@@@@@@@@@@@@@@@@@@@@@@@@@@@@@@@@@@@@@@@@@@@@@@@
%------------------------------------------------------------
%@@@@@@@@@@@@@@@@@@@@@@@@@@@@@@@@@@@@@@@@@@@@@@@@@@@@@@@@@@@@
%@@@@@@@@@@@@@@@@@@@@@@@@@@@@@@@@@@@@@@@@@@@@@@@@@@@@@@@@@@@@
%-----------------------------------------------------------
\par In Sections 4-8 we prove that if a function $f:E\to\R$ and a constant $\lambda>0$ satisfy conditions (i) and (ii) of Theorem \reff{MAIN}, then $\|f\|_{\LTP|_E}^p\le C(p)\lambda$. We refer to this part of the proof of Theorem \reff{MAIN} as the ``sufficiency'' part. Let us describe its main stages.\bigskip
%-----------------------------------------------------------
\par {\bf Sections 4-6} of the paper deal with some geometric preparations which are basic elements of our construction of an extension operator for the Sobolev space $\LTP$. The main goal of these preparations is to define a certain kind of geometric ``structure'' within the set $E$ which will enable us to organize and characterize certain ``holes'' in the complement $\RT\setminus E$ and to characterize certain kinds of ``contacts'' between them. We refer to these ``holes'' as ``lacunae''.
Although in this paper we will only be dealing with subsets of $\R^2$, it turns out to be just as easy to present these geometric preparations in the case of subsets of $\R^n$ for arbitrary $n$. Therefore the material of Sections 4-6 deals with the general $n$-dimensional case, with a view towards possible future applications of this approach. Later on, in {\bf Section 7}, we go back to dealing only with the two dimensional case.\medskip
%-----------------------------------------------------------
\par Let $E\subset\RN$ be an arbitrary closed subset. We let $W_E$ denote a {\it Whitney covering} of the open set $\RN\setminus E$, i.e., a family of non-overlapping cubes (so-called ``Whitney cubes") such that for each $Q\in W_E$ we have $\diam Q\sim \dist(Q,E)$. (See Theorem \reff{Wcov}.) In {Section 4} we introduce and study one of the main ingredients of our approach - the notion of a {\it lacuna} of Whitney cubes with respect to $E$.
%-----------------------------------------------------------
\par Roughly speaking a lacuna $L$ is a ``hole'' in the set $\RN\setminus E$, a collection of cubes in $W_E$ which surrounds a certain small subset $V_L$ of the set $E$.
$L$ is defined by the two requirements that
%------------------------------------------------------------
$$
(10Q)\cap E=(90Q)\cap E ~~\text{for each}~~Q\in L
$$
%------------------------------------------------------------
and
%------------------------------------------------------------
$$
(10Q)\cap E=(10Q')\cap E=V_L~~\text{for each}~~Q~~\text{and}~~Q'~~\text{in}~~L.
$$
%------------------------------------------------------------
\par We let $\Lc_E$ denote the set of all lacunae with respect to $E$.
%------------------------------------------------------------
\bigskip
%------------------------------------------------------------
%@@@@@@@@@@@@@@@@@@@@@@@@@@@@@@@@@@@@@@@@@@@@@@@@@@@@@@@@@@@@
%@@@@@@@@@@@@@@@@@@@@@@@@@@@@@@@@@@@@@@@@@@@@@@@@@@@@@@@@@@@@
%@@@@@@@@@@@@@@@@@@@@@@@@@@@@@@@@@@@@@@@@@@@@@@@@@@@@@@@@@@@@
%@@@@@@@@@@@@@@@@@@@@@@@@@@@@@@@@@@@@@@@@@@@@@@@@@@@@@@@@@@@@
%------------------------------------------------------------
%@@@@@@@@@@@@@@@@@@@@@@@@@@@@@@@@@@@@@@@@@@@@@@@@@@@@@@@@@@@@
%@@@@@@@@@@@@@@@@@@@@@@@@@@@@@@@@@@@@@@@@@@@@@@@@@@@@@@@@@@@@
%-----------------------------------------------------------
\par In {Section 4} we establish several important geometric properties of lacunae. In particular, we show that, for each lacuna $L$, the diameter of the associated set $V_L$ is equivalent to the diameter of the minimal cube of $L$. (See Proposition \reff{DL-M}.) We also show that the diameter of the maximal cube of $L$ is equivalent to the diameter of the union of all the cubes of $L$. (See Proposition \reff{D-MAX}.)\bigskip
%------------------------------------------------------------
%@@@@@@@@@@@@@@@@@@@@@@@@@@@@@@@@@@@@@@@@@@@@@@@@@@@@@@@@@@@@
%@@@@@@@@@@@@@@@@@@@@@@@@@@@@@@@@@@@@@@@@@@@@@@@@@@@@@@@@@@@@
%@@@@@@@@@@@@@@@@@@@@@@@@@@@@@@@@@@@@@@@@@@@@@@@@@@@@@@@@@@@@
%@@@@@@@@@@@@@@@@@@@@@@@@@@@@@@@@@@@@@@@@@@@@@@@@@@@@@@@@@@@@
%------------------------------------------------------------
%@@@@@@@@@@@@@@@@@@@@@@@@@@@@@@@@@@@@@@@@@@@@@@@@@@@@@@@@@@@@
%@@@@@@@@@@@@@@@@@@@@@@@@@@@@@@@@@@@@@@@@@@@@@@@@@@@@@@@@@@@@
%-----------------------------------------------------------
\par In {\bf Section 5} we continue the study of geometric properties of lacunae. In particular, in Proposition \reff{L-PE}, we show that there exists a mapping $\Lc_E\ni L\mapsto \PRL(L)\in E$ such that:
%-----------------------------------------------------------
\par (i) $\PRL(L)$ lies in a fixed dilation of the minimal cube of the lacuna, and
%-----------------------------------------------------------
\par (ii) every point $A\in E$ has at most $C(n)$ ``sources", i.e., lacunae $L'\in\Lc_E$ such that $\PRL(L')=A$.\smallskip
%-----------------------------------------------------------
\par We refer to the mapping $\PRL$ as a ``projection'' of $\Lc_E$ into the set $E$. In particular, the existence of a projection enables us to show that for each finite set $E$ the number of its lacunae does not exceed $C(n)$ times the number of points in $E$. (See Corollary \reff{CR-PE}.)
%-----------------------------------------------------------
\par In this section we also introduce a special graph $\Gc$ whose vertices are all the lacunae of the set $E$. We define the edges of $\Gc$ by saying that two lacunae $L,L'\in\Lc_E$ are joined by an edge  if there exist cubes $Q\in L, Q'\in L'$ such that $Q\cap Q'\ne\emp.$ In this case we call the cubes $Q,Q'$ {\it contacting cubes} and the lacunae $L,L'$ {\it contacting lacunae}. We use the notation $L\lr L'$ to denote that $L$ and $L'$ are contacting.
%-----------------------------------------------------------
\par The above-mentioned properties of lacunae imply certain special properties of the graph $\Gc$. In particular, we prove that the degree of every vertex of $\Gc$ is bounded by a constant $C(n)$; furthermore, the number of vertices does not exceed $C(n)$ times the number of points in $E$. We also show that if a cube $Q\in L$
intersects some cube $Q'$ of some other lacuna, (i.e. if $Q$ and $Q'$ are contacting cubes) then either $\diam Q$ is almost minimal or almost maximal in the lacuna $L$.
%-----------------------------------------------------------
\par For each lacuna $L$ we need to choose two points $A_L$ and $B_L$ in $E$ which will play an important role in the sequel. If $V_L$ is a single point, we choose $A_L$ to be that point, and we may choose $B_L$ to be some point in $E\setminus \{A_L\}$ whose distance from $A_L$ is minimal. Otherwise we may choose $A_L$ and $B_L$ to be any pair of distinct points in $V_L$ which satisfy $\|A_L-B_L\|=\diam V_L$. We refer to the ordered pair $(A_L,B_L)$ as an {\it interior bridge} of the lacuna $L$. We call the points $A_L$ and $B_L$ {\it the ends of the interior bridge}.
%-----------------------------------------------------------
\par We note that the fact that two lacunae $L$ and $L'$ are contacting, does not guarantee any connection between their respective interior bridges.\bigskip
%------------------------------------------------------------
%@@@@@@@@@@@@@@@@@@@@@@@@@@@@@@@@@@@@@@@@@@@@@@@@@@@@@@@@@@@@
%@@@@@@@@@@@@@@@@@@@@@@@@@@@@@@@@@@@@@@@@@@@@@@@@@@@@@@@@@@@@
%@@@@@@@@@@@@@@@@@@@@@@@@@@@@@@@@@@@@@@@@@@@@@@@@@@@@@@@@@@@@
%@@@@@@@@@@@@@@@@@@@@@@@@@@@@@@@@@@@@@@@@@@@@@@@@@@@@@@@@@@@@
%------------------------------------------------------------
%@@@@@@@@@@@@@@@@@@@@@@@@@@@@@@@@@@@@@@@@@@@@@@@@@@@@@@@@@@@@
%@@@@@@@@@@@@@@@@@@@@@@@@@@@@@@@@@@@@@@@@@@@@@@@@@@@@@@@@@@@@
%-----------------------------------------------------------
\par In {\bf Section 6} we construct additional bridges  which join interior bridges between con\-tac\-ting lacunae. We refer to these new bridges as {\it exterior bridges}. More specifically, given contacting lacunae $L$ and $L'$ we choose points $C(L,L')\in\{A_L,B_L\}$ and $C(L',L)\in\{A_{L'},B_{L'}\}$ in such a way that the following two conditions are satisfied: \medskip %----------------------------------------------------------
\par {\it In the triangle $\Delta\{A_L,B_L,C(L',L)\}$ the side of the triangle which is opposite to the vertex $C(L,L')$ is not the smallest side of the triangle.
%----------------------------------------------------------
\par A similar condition holds for the triangle $\Delta\{A_{L'},B_{L'},C(L,L')\}$: the side of the triangle which is opposite to the vertex $C(L',L)$ is not the smallest side of the triangle.}\medskip
%----------------------------------------------------------
\par We refer to the (non-ordered) pair $\{C(L,L'),C(L',L)\}$ as an {\it exterior bridge connecting the interior bridges $(A_L,B_L)$ and $(A_{L'},B_{L'})$}. We prove that for every pair of contacting lacunae $L$ and $L'$ there always exists an exterior bridge which connects $L$ to $L'$. We say that the interior bridge $T=(A_L,B_L)$ and the exterior bridge $\tT=\{C(L,L'),C(L',L)\}$ are {\it connected bridges}. In this case we write $T\bcn \tT$.\smallskip
%----------------------------------------------------------
\par This special choice of the ends $C(L,L')$ and $C(L',L)$ of the exterior bridge implies a certain special  property of the triangle $\Delta=\Delta\{A_L,B_L,C(L',L)\}$. We know that $C(L,L')$ is one of the vertices of $\Delta$. This point is a common point of the interior bridge $(A_L,B_L)$ and the exterior bridge $\{C(L,L'),C(L',L)\}$. Let $\alpha$ be the angle of the triangle $\Delta$ corresponding to the vertex $C(L,L')$. Then
%----------------------------------------------------------
$$
|\sin \alpha|\sim\tfrac{1}{R_\Delta}\,\diam\Delta
=\cu_\Delta\,\diam\Delta.
$$
%----------------------------------------------------------
We use this specific property of the triangle $\Delta$ later in our estimates of the Sobolev norm of the  extension operator. In particular, it explains the appearance of the Menger curvature $\cu_\Delta$ in the criterion of Theorem \reff{MAIN}.
%-----------------------------------------------------------
\par By $\BRE$ we denote the family of all bridges (interior and exterior). For every bridge $T\in\BRE$ we let $\A{T}$ and $\B{T}$ denote its ends. We know that some pairs $T,T'\in\BRE$ of the bridges are connected to each other ($T\bcn T'$). We also know that one of these connected bridges is an interior bridge of a lacuna, and the second one is an exterior bridge joining $L$ with another contacting lacuna.
%-----------------------------------------------------------
\par We finish our geometric preparations with the following statement proven in Proposition \reff{BR-TO-Q}:
\smallskip
%-----------------------------------------------------------
\par {\it There exists a family $\KE$ of pairwise disjoint cubes and a one-to-one mapping which to every pair of connected bridges $T,T'\in\BRE$, $T\bcn T'$, assigns a cube $K(T,T')\in\KE$ such that
$\diam \{\A{T},\B{T},\A{T'},\B{T'}\} \sim \diam K(\Br,\Br')$ and
$\{\A{T},\B{T},\A{T'},\B{T'}\}\subset \gamma K(\Br,\Br').$}
%------------------------------------------------------------
%@@@@@@@@@@@@@@@@@@@@@@@@@@@@@@@@@@@@@@@@@@@@@@@@@@@@@@@@@@@@
%@@@@@@@@@@@@@@@@@@@@@@@@@@@@@@@@@@@@@@@@@@@@@@@@@@@@@@@@@@@@
%@@@@@@@@@@@@@@@@@@@@@@@@@@@@@@@@@@@@@@@@@@@@@@@@@@@@@@@@@@@@
%@@@@@@@@@@@@@@@@@@@@@@@@@@@@@@@@@@@@@@@@@@@@@@@@@@@@@@@@@@@@
%------------------------------------------------------------
%@@@@@@@@@@@@@@@@@@@@@@@@@@@@@@@@@@@@@@@@@@@@@@@@@@@@@@@@@@@@
%@@@@@@@@@@@@@@@@@@@@@@@@@@@@@@@@@@@@@@@@@@@@@@@@@@@@@@@@@@@@
\smallskip
%------------------------------------------------------------
\par In {\bf Section 7}, as already mentioned, we return to dealing only with the case where the dimension $n=2$. We turn to the construction of an extension operator for the trace space $\LTP|_E$ whenever $E\subset\RT$ is an arbitrary {\it finite} subset and $p>2$.\smallskip
%------------------------------------------------------------
\par We fix a function $f:E\to\R$. To every lacuna $L\in\Lc_E$ we are going to assign an affine polynomial $P_L\in\PO$ which {\it interpolates $f$ at the points $A_L$ and $B_L$}. Note that we have some flexibility (one degree of freedom) in the choice of this polynomial.
A collection of several results, in Section 7 and also in Section 8, ultimately enable us to develop an appropriate strategy for choosing $P_L$ in an ``almost optimal" way.
%------------------------------------------------------------
\par Once we have obtained a suitable family of affine polynomials $\{P_L\in\PO:L\in\Lc_E\}$ by this procedure, we can use it to generate a Whitney-type extension of $f$ to all of $\RT$ in the following way: First, we assign a polynomial $\PQ$ to every Whitney \sq $Q\in W_E$. We do this by setting $\PQ:=P_L$ for all the squares $Q$ in the lacuna $L$, for each $L\in \Lc_E$. Then we define an extension $F$ of $f$ by the Whitney formula
%----------------------------------------------------------
\bel{D-W}
F(x):=
\sum\limits_{Q\in W_E}
\varphi_Q(x)\,\PQ(x),~~~\text{whenever}~~x\in\RT\setminus E,
\ee
%----------------------------------------------------------
and $F(x)=f(x)$, if $x\in E$. Here $\{\varphi_Q:Q\in W_E\}$ is a smooth partition of unity subordinated to the Whitney decomposition $W_E$ with standard properties (as described in Lemma \reff{P-U}).
%@@@@@@@@@@@@@@@@@@@@@@@@@@@@@@@@@@@@@@@@@@@@@@@@@@@@@@@@@@
%----------------------------------------------------------
\par Now let us begin sketching the (rather long) series of results which we need to subsequently fix our strategy for choosing the affine polynomials $P_L$. Our goal will be to choose them in a way that ensures that the function $F$ defined by \rf{D-W} satisfies
%------------------------------------------------------------
\bel{F-LM}
\|F\|_{\LOP}\le C(p)\lambda^{\frac1p}
\ee
%------------------------------------------------------------
whenever the constant $\lambda$ is related to the function $f$ by the conditions (i), (ii) of Theorem \reff{MAIN}. This will of course complete the proof of the sufficiency part of Theorem \reff{MAIN}.
%----------------------------------------------------------
\par Our first step is to show that, no matter how we fix the remaining degree of freedom in our choice of the polynomials $P_L$, the function $F$ of \rf{D-W} will always satisfy
%----------------------------------------------------------
$$
\|F\|^p_{\LTP}\le C\,\sum d(L,L')^{2-2p}
\max\limits_{Q(L,L')}|P_L-P_{L'}|^p
$$
%----------------------------------------------------------
where the sum is taken over all pairs of {\it contacting} lacunae $L,L'\in\LE$ ($L\lr L'$). Here $d(L,L'):=\diam\{A_L,B_L,A_{L'},B_{L'}\}$, and $Q(L,L'):=Q(A_L,d(L,L'))$. (See Proposition \reff{F-INAF}).
%----------------------------------------------------------
\par In the next step, in addition to the polynomials $P_L$ (which interpolate $f$ at the ends of the interior bridges $(A_L,B_L)$) we also bring into play other affine polynomials which interpolate $f$ at the ends of {\it exterior} bridges. Our aim is to obtain an estimate for the norm $\|F\|^p_{\LTP}$ of the function $F$ of \rf{D-W}, in terms of the {\it gradients} of all these affine polynomials. As in the previous step, the result here will hold no matter how we fix the remaining degree of freedom in the choice of the polynomials.
%------------------------------------------------------------
\par It is convenient to consider ``gradient" mappings $g:\BRE\to\RT$ which assign a vector in $\RT$
to each (interior or exterior) bridge $T\in\BRE$. We note that such a mapping $g$ satisfies
%----------------------------------------------------------
\bel{END-I}
\ip{g(T),\A{T}-\B{T}}=f(\A{T})-f(\B{T})
\ee
%----------------------------------------------------------
if and only if it is the gradient of an affine polynomial which interpolates $f$ at the ends $\A{T}$ and $\B{T}$ of the bridge $T$. (Here $\ip{\cdot,\cdot}$ denotes the standard inner product in $\RT$.)
%------------------------------------------------------------
\par Suppose that the mapping $g$ satisfies \rf{END-I} for {\it every} \br $T\in\BRE$ with ends at points $\A{T},\B{T}\in E$. Then we prove (see Proposition \reff{F-GN}) that, if we take the polynomials $P_L$ to be
$P_L(x)=f(A_L)+\ip{g(T_L),x-A_L}$ where $T_L=(A_L,B_L)$
is the interior bridge of the lacuna $L$, then the function $F$ of \rf{D-W} satisfies
%------------------------------------------------------------
\par
$$
\|F\|^p_{\LTP}\le C\,\sum \left\{D(T,T')^{2-p}
\|g(T)-g(T')\|^p:~T,T'\in\BRE, T\bcn T'\right\}.
$$
%------------------------------------------------------------
Here $D(T,T'):=\diam \{\A{T},\B{T},\A{T'},\B{T'}\}$ for each pair $T,T'\in\BRE$ of connected \brs ($T\bcn T'$) with the ends at points $\{\A{T},\B{T}\}$ and $\{\A{T'},\B{T'}\}$ respectively.
%------------------------------------------------------------
\par  We turn to the next step. Here we introduce a new ``parameter", a non-negative $L_p$-function $h:\RT\to\R_+$. We fix a constant $q\in(2,p)$ (for instance we may take
$q=(p+2)/2$) and again consider a mapping $g:\BRE\to\RT$ satisfying condition \rf{END-I} for  every \br $T\in\BRE$ with ends at points $\A{T},\B{T}\in E$. We assume that $g$ also satisfies the following condition: for every pair of connected bridges $T,T'\in\BRE,$ $T\bcn T'$, %----------------------------------------------------------
\bel{S-SEL1}
\|g(T)-g(T')\|\le \diam Q(T,T')\,\left(\frac{1}{|\gamma Q(T,T')|}
\intl_{\gamma Q(T,T')}h^q(z)dz\right)^{\frac1q}.
\ee
%----------------------------------------------------------
where $Q(T,T')=Q(\A{T},D(T,T'))$. Here $\gamma\ge 1$ is an absolute constant. We prove (see Proposition \reff{A-H})
that the existence of a mapping $g$ and a function $h$ satisfying all the above conditions implies the existence of an extension $F\in\LTP$ of the function $f$ such that
%----------------------------------------------------------
\bel{F-NE}
\|F\|_{\LTP}\le C(p)\|h\|_{\LPRT}.
\ee
%----------------------------------------------------------
\smallskip
%----------------------------------------------------------
%@@@@@@@@@@@@@@@@@@@@@@@@@@@@@@@@@@@@@@@@@@@@@@@@@@@@@@@@@@
%@@@@@@@@@@@@@@@@@@@@@@@@@@@@@@@@@@@@@@@@@@@@@@@@@@@@@@@@@@
%@@@@@@@@@@@@@@@@@@@@@@@@@@@@@@@@@@@@@@@@@@@@@@@@@@@@@@@@@@
%----------------------------------------------------------
%@@@@@@@@@@@@@@@@@@@@@@@@@@@@@@@@@@@@@@@@@@@@@@@@@@@@@@@@@@
%@@@@@@@@@@@@@@@@@@@@@@@@@@@@@@@@@@@@@@@@@@@@@@@@@@@@@@@@@@
%@@@@@@@@@@@@@@@@@@@@@@@@@@@@@@@@@@@@@@@@@@@@@@@@@@@@@@@@@@
%----------------------------------------------------------
\par This is the last result of Section 7. We begin {\bf Section 8} by introducing the notions of of {\it set-valued mappings} and their {\it Sobolev-type selections.} More specifically, motivated by the equality \rf{END-I}, we introduce a set-valued mapping $G_f$ which to every bridge $T\in\BRE$ assigns the straight line consisting of all points $z\in\RT$ which satisfy
%----------------------------------------------------------
$$
\ip{z,\A{T}-\B{T}}=f(\A{T})-f(\B{T}).
$$
%----------------------------------------------------------
We say that a (vector valued) mapping $g:\BRE\to \RT$ is a  {\it selection of the set-valued mapping $G_f$} if $g(T)\in G_f(T)$ for every bridge $T\in\BRE$. So any mapping $g$ which satisfies \rf{END-I} is a selection of $G_f$.
%----------------------------------------------------------
\par We observe that the inequality \rf{S-SEL1} resembles the classical \SP inequality (which is recalled in Proposition \reff{SP-IN}). This motivates us to introduce the terminology {\it Sobolev-type selection of the set-valued mapping $G_f$ with respect to $h$} to designate any mapping $g:\BRE\to\RT$ which satisfies conditions \rf{END-I} and \rf{S-SEL1} for some non-negative function $h\in L_p(\RT)$.
%----------------------------------------------------------
\par The remaining part of Section 8, which we shall presently describe in some detail, is a step by step construction of a non-negative function $h\in \LPRT$
and a Sobolev-type selection $g$ with respect to $h$ of $G_f$ for which
%------------------------------------------------------------
\bel{H-LM}
\|h\|_{\LPRT}\le C(p)\lambda^{\frac1p}.
\ee
%------------------------------------------------------------
Here $\lambda$ and $f$ are related as already mentioned above. (See \rf{F-LM}.) This construction finally achieves the goal mentioned just before \rf{F-LM}, since of course \rf{H-LM} and \rf{F-NE} together imply \rf{F-LM}.
%------------------------------------------------------------
\par Many of the steps in Section 8 generalize ideas and methods for the case $p=q=\infty$
presented and used in \cite{S-GAN} to the case of Sobolev-type selections, i.e., to the case $p<\infty$.
%----------------------------------------------------------
\begin{remark} {\em Note that for $p=q=\infty$ the problem of the existence of a Sobolev-type selection is a special case of the following general {\it Lipschitz selection problem}:
%----------------------------------------------------------
\par Let $G$ be a set-valued mapping which to every point of a metric space $(\Mc,\rho)$ assigns a convex closed subset of $\RN$.
%----------------------------------------------------------
\par {\it How do we know whether there exists a selection $g:\Mc\to\RN$ of the mapping $G$ satisfying the Lipschitz condition with respect to the metric $\rho$? What is the order of magnitude of the minimal Lipschitz constant of such a selection of $G$?}
%----------------------------------------------------------
\par These and other problems related to Lipschitz selections of set-valued mapping have been studied in the author's papers \cite{S-GAN,S-GF,S-JF}. A solution to one them given in \cite{S-GAN}, Theorem 3.14, leads us to the trace criterion presented in Remark \reff{P-INF}.\rbx}
%----------------------------------------------------------
\end{remark}
\medskip
%----------------------------------------------------------
\par We can now begin our description of the construction of the selection $g$ and the function $h$ mentioned above. It can be divided into two main parts. In the first and rather longer part (a process of ``pre-selection") we obtain the required function $h$ and a mapping $\tg$ such that $\tg(T)$ is in some sense ``close" to the straight line $G_f(T)$ for every bridge $T\in\BRE$. In the second part a relatively simple procedure enables us to obtain $g$ from $\tg$.
%----------------------------------------------------------
\par We begin the first part by associating a triangle
$\Delta(T,T')$ to each pair of connected bridges $T,T'\in\BRE$, ($T\bcn T'$).
The vertices of $\Delta(T,T')$
are the ends of the bridges $T$ and $T'$.
It is relevant to recall here that in {Section 6,} Proposition \reff{BR-TO-Q}, we construct a family $\KE$ of pairwise disjoint \sqs
%----------------------------------------------------------
$$
\{K(\Br,\Br'): T,T'\in\BRE, T\bcn T'\}
$$
%----------------------------------------------------------
with the following properties:
%----------------------------------------------------------
\par (i) for each pair of connected bridges $T,T'\in\BRE, T\bcn T',$ the size of the triangle $\Delta(T,T')$ is equivalent to the size of the \sq $K(\Br,\Br')$;
%----------------------------------------------------------
\par (ii)  $\Delta(T,T')\subset \gamma K(\Br,\Br')$ where  $\gamma>0$ is an absolute constant.\medskip
%----------------------------------------------------------
\par Then we introduce two basic elements of our construction. The first of them is a mapping $\TF:\RT\to\RT$ which we define as follows:
%----------------------------------------------------------
\par At the center of every \sq $K(\Br,\Br')\in \KE$ this mapping takes  the value of {\it the gradient  $\nabla P_{\Delta(T,T')}[f]$ of the affine polynomial $P_{\Delta(T,T')}[f]$ which interpolates the function $f$ at the vertices of the triangle $\Delta(T,T')$.}
At all other points we set $\TF\equiv 0$.
%----------------------------------------------------------
\par Note that the vector $\nabla P_{\Delta(T,T')}[f]\in\RT$ has a simple geometrical description: {\it it coincides with the point of intersection of the straight lines $G_f(T)$ and $G_f(T')$.}
%----------------------------------------------------------
\par The second element is a discrete Borel measure on $\RT$ which we denote by $\mu_E$. This measure is supported on the set of centers of \sqs of the family $\KE$. For every \sq $K=K(\Br,\Br')\in \KE$ the $\mu_E$-measure of its center equals $\cu_\Delta\,|K|$ where {\it $\cu_\Delta$ is the Menger curvatures of the triangle $\Delta=\Delta(T,T')$}. (See \rf{M-CURV}.)
%----------------------------------------------------------
\par We prove that the mapping $\TF$ belongs to the space $\VS=\VLOP+\VLPME$ and that its norm in this space
satisfies $\|\TF\|_{\VS}\le C\lambda^{\frac1p}$. The proof of this fact uses a description of the norm in the space $\VS$ given in \cite{S3}. This description enables us to decompose $\TF$ into the sum $\TF=\Tc_1(f)+\Tc_2(f)$ of a mapping $\Tc_1(f)\in\VLOP$ and a mapping $\Tc_2(f)\in \VLPME$ provided the constant $\lambda$ is related to the function $f$ by the conditions (i), (ii) of Theorem \reff{MAIN}. The norms of these mappings in the corresponding spaces are bounded by a constant $C(p)\lambda^{\frac1p}$.\bigskip
%----------------------------------------------------------
\par We can next introduce a mapping $\tg(f):\BRE\to\RT$ which, to every bridge $T\in\BRE$
whose ends are the points $\A{T},\B{T}\in E$, assigns the vector
%----------------------------------------------------------
$$
\tg(T;f)=\Tc_1(\A{T};f).
$$
%----------------------------------------------------------
\par We show that the mapping $\tg(f)$ satisfies the condition \rf{S-SEL1} when we choose
$h_f=\|\nabla \Tc_1(f)\|$. Since $\Tc_1(f)\in\VLOP$ and $\|\Tc_1(f)\|_{\VLOP}\le C\lambda^{\frac1p}$, the function $h_f$ belongs to the space $\LPRT$ and $\|h_f\|_{\LPRT}\le C\lambda^{\frac1p}$. Thus {\it $\tg(f)$ is a Sobolev-type mapping with respect to $h_f$.}\medskip
%----------------------------------------------------------
\par However we cannot guarantee that $\tg(T;f)\in G_f(T)$ for every bridge $T\in\BRE$.
This means  that {\it in general $\tg(f)$ is not a selection of $G_f$}. We refer to the mapping $\tg(f)$ as {\it a Sobolev-type pre-selection} of the set-valued mapping $G_f$ (with respect to the function $h_f$).
%----------------------------------------------------------
\par It turns out that the pre-selection $\tg(f)$ possesses various additional properties which will enable us to transform $\tg(f)$ into a true  selection $g$ of $G_f$ satisfying the Sobolev-type condition \rf{S-SEL1} with respect to a certain function $h\in \LPRT$. In particular, we know that the $L_p(\RT;\mu_E)$-norm of the component $\Tc_2(f)$ is at most $C\lambda^{\frac1p}$. Using this property we show that,
for each bridge $T\in\BRE$,
in spite of the fact that the pre-selection $\tg(T;f)$ does not necessarily lie on the straight line $G_f(T)$,
it lies ``rather close'' to $G_f(T)$.\bigskip
%----------------------------------------------------------
\par
This concludes the first part of our construction.
We turn to the second and final part: Given a bridge $T\in\BRE$ we define the required selection $g(T;f)$
as {\it the orthogonal projection of the pre-selection $\tg(T;f)$ onto the straight line $G_f(T)$}.
%----------------------------------------------------------
\par Using the ``closeness'' of the pre-selection $\tg(T;f)$ to $G_f(T)$ we prove that the mapping $g(f):\BRE\to\RT$ is a Sobolev-type selection of the set-valued mapping $G_f$ with respect to a
certain
non-negative function $h\in\LPRT$ whose norm satisfies 	 $\|h\|_{\LPRT}\le C\lambda^{\frac1p}$. Hence, by \rf{F-NE}, there exists a function $F\in\LTP$ such that $F|_E=f$ and
%----------------------------------------------------------
$$
\|F\|_{\LTP}\le C(p)\|h\|_{\LPRT}\le C(p)\lambda^{\frac1p}.
$$
%----------------------------------------------------------
This completes the proof of the Main Theorem \reff{MAIN}.
\bigskip
%----------------------------------------------------------
\par {\bf Section 9} further discusses some aspects
of one of the main ingredients of the proof of Theorem \reff{MAIN}, namely
an almost optimal decomposition of the mapping  $\TF$ into the sum $\TF=\Tc_1(f)+\Tc_2(f)$ of a mapping $\Tc_1(f)\in\VLOP$ and a mapping $\Tc_2(f)\in \VLPME$.
As already mentioned above, criteria for the existence of such a decomposition are established in \cite{S3}.
Here, using methods of \cite{S3} we describe a constructive algorithm which provides such a decomposition.
%----------------------------------------------------------
\medskip
%----------------------------------------------------------
\par In {\bf Section 10} we re-examine all the methods and ideas which we have used in the proof of Theorem \reff{MAIN} and present a constructive algorithm which, to every function $f$ defined on a finite set $E$, assigns its almost optimal extension to a function $F(f)\in\LTP$. At every step of this algorithm we verify that its corresponding components which we construct at this step depend on $f$ linearly. Consequently, at the final step we obtain an almost optimal extension of $f$ which {\it depends linearly on $f$.} This gives a constructive proof of Theorem \reff{LINEXT}, i.e., the existence of a continuous linear extension operator $Ext_E:\LTP|_E\to\LTP$ whose operator norm satisfies $\|Ext_E\|\le C(p)$.
%----------------------------------------------------------
\medskip
%----------------------------------------------------------
\par In {\bf Section 11} we prove Theorems \reff{REF-MAIN}, \reff{REF-SP} and \reff{NEW-MAIN}. The proof of Theorem \reff{REF-MAIN} is a slight modification of the proof of Theorem \reff{MAIN}. More specifically, we replace in this proof the criterion of the norm in the space $\VS=\VLOP+\VLPM$ given in Theorem \reff{S-CMS} with its refinement presented in Theorem \reff{REF-S1}.
%----------------------------------------------------------
\par We prove Theorem \reff{REF-SP} using another refinement of Theorem \reff{S-CMS} given in Theorem \reff{REF-S2}. We obtain a new criterion for the norm of a function in the trace space $\LTP|_E$ which we present in Theorem \reff{REF-ML}. We ``sparse'' this new criterion using the $p$-sparsification Theorem \reff{SPR-P} proven by B. Klartag. (In turn this result uses a sparsification theorem by J. D. Batson, D. A. Spielman
and N. Srivastava \cite{BSS} as a black box.) This leads us to a proof of Theorem \reff{REF-SP}.
%----------------------------------------------------------
%@@@@@@@@@@@@@@@@@@@@@@@@@@@@@@@@@@@@@@@@@@@@@@@@@@@@@@@@@@
%@@@@@@@@@@@@@@@@@@@@@@@@@@@@@@@@@@@@@@@@@@@@@@@@@@@@@@@@@@
%@@@@@@@@@@@@@@@@@@@@@@@@@@@@@@@@@@@@@@@@@@@@@@@@@@@@@@@@@@
%@@@@@@@@@@@@@@@@@@@@@@@@@@@@@@@@@@@@@@@@@@@@@@@@@@@@@@@@@@
%@@@@@@@@@@@@@@@@@@@@@@@@@@@@@@@@@@@@@@@@@@@@@@@@@@@@@@@@@@
%@@@@@@@@@@@@@@@@@@@@@@@@@@@@@@@@@@@@@@@@@@@@@@@@@@@@@@@@@@
%----------------------------------------------------------
%@@@@@@@@@@@@@@@@@@@@@@@@@@@@@@@@@@@@@@@@@@@@@@@@@@@@@@@@@@
%@@@@@@@@@@@@@@@@@@@@@@@@@@@@@@@@@@@@@@@@@@@@@@@@@@@@@@@@@@
%@@@@@@@@@@@@@@@@@@@@@@@@@@@@@@@@@@@@@@@@@@@@@@@@@@@@@@@@@@
%@@@@@@@@@@@@@@@@@@@@@@@@@      %@@@@@@@@@@@@@@@@@@@@@@@@@@@
%@@@@@@@@@@@@@@@@@@@@@@@          @@@@@@@@@@@@@@@@@@@@@@@@@
%@@@@@@@@@@@@@@@@@@@@@              %@@@@@@@@@@@@@@@@@@@@@@@
%@@@@@@@@@@@@@@@@@@@     SECTION 3
%@@@@@@@@@@@@@@@@@@@@@
%@@@@@@@@@@@@@@@@@@@@@              %@@@@@@@@@@@@@@@@@@@@@@@
%@@@@@@@@@@@@@@@@@@@@@@@          %@@@@@@@@@@@@@@@@@@@@@@@@@
%@@@@@@@@@@@@@@@@@@@@@@@@@      %@@@@@@@@@@@@@@@@@@@@@@@@@@@
%@@@@@@@@@@@@@@@@@@@@@@@@@@@@@@@@@@@@@@@@@@@@@@@@@@@@@@@@@@
%@@@@@@@@@@@@@@@@@@@@@@@@@@@@@@@@@@@@@@@@@@@@@@@@@@@@@@@@@@
%@@@@@@@@@@@@@@@@@@@@@@@@@@@@@@@@@@@@@@@@@@@@@@@@@@@@@@@@@@
%----------------------------------------------------------
\SECT{3. Main Theorem \reff{MAIN}: Necessity.}{3}
%----------------------------------------------------------
\addtocontents{toc}{3. Main Theorem \reff{MAIN}: Necessity. \hfill \thepage\par}
%----------------------------------------------------------
\indent
%@@@@@@@@@@@@@@@@@@@@@@@@@@@@@@@@@@@@@@@@@@@@@@@@@@@@@@@@@@
%----------------------------------------------------------
\par Let us fix several additional notation. Throughout the paper $\gamma,\gamma_1,\gamma_2...,$ and $C,C_1,C_2,...$ will be generic positive constants which depend only on parameters determining function spaces ($p,q,$ etc). These constants can change even in a single string of estimates. The dependence of a constant on certain parameters is expressed, for example, by the notation $C=C(p)$. We write $A\sim B$ if there is a constant $C\ge 1$ such that $A/C\le B\le CA$.
%----------------------------------------------------------
\par We let $Q(c,r)$ denote the \sq in $\RT$ centered at $c$ with side length $2r$. Given a \sq $Q$ by $c_Q$ we denote its center and by $r_Q$ half of its side length. Thus $Q=Q(c_Q,r_Q)$ and $\lambda Q=Q(c_Q,\lambda r_Q)$ for every constant $\lambda>0$.
%----------------------------------------------------------
\par Given $x=(x_1,x_2)\in\RT$ by $\|x\|:=\max \{|x_1|,|x_2|\}$ and by $\|x\|_2:=(|x_1|^2+|x_2|^2)^{\frac12}$ we denote the uniform and the Euclidean norms in $\RT$ respectively.
%----------------------------------------------------------
\par Let $A,B\subset \RT$. We put
$\diam A:=\sup\{\|a-a'\|:~a,a'\in A\}$ and
%----------------------------------------------------------
$$
\dist(A,B):=\inf\{\|a-b\|:~a\in A, b\in B\}.
$$
%----------------------------------------------------------
For $x\in \RT$ we also set $\dist(x,A):=\dist(\{x\},A)$.
We put  $\dist(A,B)=+\infty$ and $\dist({x},B)=+\infty$ whenever $B=\emp.$
%----------------------------------------------------------
\par By $\ip{\cdot,\cdot}$ we denote the standard inner product in $\RT$. Given a finite set $A$ by $\card A$ we denote the number of its elements. Finally, we let $\PO$  denote the space of all polynomials of degree at most $1$ on $\RT$.\bigskip
%----------------------------------------------------------
\par In this section we prove that
%----------------------------------------------------------
$$
\inf \lambda^{\frac{1}{p}}\le C(p)\|f\|_{\LTP|_E}
$$
%----------------------------------------------------------
where $\lambda$ is the parameter from Theorem \reff{MAIN} and $E$ is an {\it arbitrary} (not necessarily finite) closed subset of $\RT$. This inequality follows from the next
%---------------------------------------------------------- %@@@@@@@@@@@@@@@@@@@@@@@@@@@@@@@@@@@@@@@@@@@@@@@@@@@@@@@@@@
%@@@@@@@@@@@@@@@@@@@@@@@@@@@@@@@@@@@@@@@@@@@@@@@@@@@@@@@@@@
%@@@@@@@@@@@@@@@@@@@@@@@@@@@@@@@@@@@@@@@@@@@@@@@@@@@@@@@@@@
%@@@@@@@@@@@@@@@@@@@@@@@@@@@@@@@@@@@@@@@@@@@@@@@@@@@@@@@@@@
%----------------------------------------------------------
\begin{statement}\lbl{ST-M} Let $E$ be a closed subset of $\RT$ and let $F\in\LTP$ be a $C^1$-function on $\RT$ such that $F|_E=f$. Then the conditions (i) and (ii) of Theorem \reff{MAIN} are satisfied with $\lambda=C(p,\gamma)\|F\|_{\LTP}^p$ and arbitrary $\gamma\in [1,\infty)$.
%----------------------------------------------------------
\end{statement}
\bigskip
%---------------------------------------------------------- %@@@@@@@@@@@@@@@@@@@@@@@@@@@@@@@@@@@@@@@@@@@@@@@@@@@@@@@@@@
%@@@@@@@@@@@@@@@@@@@@@@@@@@@@@@@@@@@@@@@@@@@@@@@@@@@@@@@@@@
%@@@@@@@@@@@@@@@@@@@@@@@@@@@@@@@@@@@@@@@@@@@@@@@@@@@@@@@@@@
%@@@@@@@@@@@@@@@@@@@@@@@@@@@@@@@@@@@@@@@@@@@@@@@@@@@@@@@@@@
%@@@@@@@@@@@@@@@@@@@@@@@@@@@@@@@@@@@@@@@@@@@@@@@@@@@@@@@@@@
%@@@@@@@@@@@@@@@@@@@@@@@@@@@@@@@@@@@@@@@@@@@@@@@@@@@@@@@@@@
%@@@@@@@@@@@@@@@@@@@@@@@@@@@@@@@@@@@@@@@@@@@@@@@@@@@@@@@@@@
%@@@@@@@@@@@@@@@@@@@@@@@@@@@@@@@@@@@@@@@@@@@@@@@@@@@@@@@@@@
%@@@@@@@@@@@@@@@@@@@@@@@@@@@@@@@@@@@@@@@@@@@@@@@@@@@@@@@@@@
\medskip
\par {\bf 3.1. Sobolev-Poincar\'e type inequalities for interpolating polynomials.}
%----------------------------------------------------------
\addtocontents{toc}{~~~~3.1. Sobolev-Poincar\'e type inequalities for interpolating polynomials. \hfill \thepage\par}
%---------------------------------------------------------- %@@@@@@@@@@@@@@@@@@@@@@@@@@@@@@@@@@@@@@@@@@@@@@@@@@@@@@@@@@
%@@@@@@@@@@@@@@@@@@@@@@@@@@@@@@@@@@@@@@@@@@@@@@@@@@@@@@@@@@
%@@@@@@@@@@@@@@@@@@@@@@@@@@@@@@@@@@@@@@@@@@@@@@@@@@@@@@@@@@
%@@@@@@@@@@@@@@@@@@@@@@@@@@@@@@@@@@@@@@@@@@@@@@@@@@@@@@@@@@
%----------------------------------------------------------
Our proof of Statement \reff{ST-M} uses Proposition \reff{SP-IN} which presents the classical Sobolev imbedding inequality for the space $\LTP$ whenever $p>2$, see, e.g. \cite{M}, p. 61, or \cite{MP}, p. 55. (This inequality is also known in the literature as Sobolev-Poincar\'{e} inequality for Sobolev $L^2_p$-functions.)
%@@@@@@@@@@@@@@@@@@@@@@@@@@@@@@@@@@@@@@@@@@@@@@@@@@@@@@@@@@
%@@@@@@@@@@@@@@@@@@@@@@@@@@@@@@@@@@@@@@@@@@@@@@@@@@@@@@@@@@
\begin{proposition} \lbl{SP-IN} Let $F\in\LTP$ be a $C^1$-function on $\RT$ and let $2<q\le p<\infty$. Then for every \sq $Q\subset\RT$ and every $x,y\in Q$ the following inequalities
%----------------------------------------------------------
\bel{S-1}
~~~~~~|F(x)-F(y)-\ip{\nabla F(y),x-y}|
\le C(q)\|x-y\|\diam Q \left(\frac{1}{|Q|} \intl_Q(\nabla^2F)^q\,dz\right)^{\frac{1}{q}},
\ee
%----------------------------------------------------------
%@@@@@@@@@@@@@@@@@@@@@@@@@@@@@@@@@@@@@@@@@@@@@@@@@@@@@@@@@@
%----------------------------------------------------------
\bel{S-2}
\|\nabla F(x)-\nabla F(y)\|
\le C(q) \diam Q \left(\frac{1}{|Q|} \intl_Q(\nabla^2F)^q\,dz\right)^{\frac{1}{q}},
\ee
%----------------------------------------------------------
hold.
%----------------------------------------------------------
\end{proposition}
%----------------------------------------------------------
%@@@@@@@@@@@@@@@@@@@@@@@@@@@@@@@@@@@@@@@@@@@@@@@@@@@@@@@@@@
%@@@@@@@@@@@@@@@@@@@@@@@@@@@@@@@@@@@@@@@@@@@@@@@@@@@@@@@@@@
%@@@@@@@@@@@@@@@@@@@@@@@@@@@@@@@@@@@@@@@@@@@@@@@@@@@@@@@@@@
%@@@@@@@@@@@@@@@@@@@@@@@@@@@@@@@@@@@@@@@@@@@@@@@@@@@@@@@@@@
%----------------------------------------------------------
\begin{remark} {\em Note that inequalities \rf{S-1} and \rf{S-2} are equivalent to the inequality
%----------------------------------------------------------
$$
\max_{Q_{xy}}|T_x[F]-T_y[F]|
\le C(q) \|x-y\|^2 \left(\frac{1}{|Q_{xy}|} \intl_{Q_{xy}} (\nabla^2F)^q\,dz\right)^{\frac{1}{q}}.
$$
%----------------------------------------------------------
Here given $a\in\RT$
%----------------------------------------------------------
$$
T_a[F](z)=F(a)+\ip{\nabla F(a),z-a}
$$
%----------------------------------------------------------
is the first order Taylor polynomial of the function $F\in C^1(\RT)\cap \LTP$ at the point $a$. By  $Q_{xy}$ we denote a \sq $Q(x,\|x-y\|)$.\rbx}
%----------------------------------------------------------
\end{remark}
%----------------------------------------------------------
%@@@@@@@@@@@@@@@@@@@@@@@@@@@@@@@@@@@@@@@@@@@@@@@@@@@@@@@@@@
%@@@@@@@@@@@@@@@@@@@@@@@@@@@@@@@@@@@@@@@@@@@@@@@@@@@@@@@@@@
%@@@@@@@@@@@@@@@@@@@@@@@@@@@@@@@@@@@@@@@@@@@@@@@@@@@@@@@@@@
%@@@@@@@@@@@@@@@@@@@@@@@@@@@@@@@@@@@@@@@@@@@@@@@@@@@@@@@@@@
%----------------------------------------------------------
\medskip
%----------------------------------------------------------
\par First let us prove that condition (i) of Theorem \reff{MAIN} holds with $\lambda=C(p)\|F\|_{\LTP}^p$ provided  $F\in\LTP\cap C^1(\RT)$ and $F|_E=f$.
%----------------------------------------------------------
%@@@@@@@@@@@@@@@@@@@@@@@@@@@@@@@@@@@@@@@@@@@@@@@@@@@@@@@@@@
%@@@@@@@@@@@@@@@@@@@@@@@@@@@@@@@@@@@@@@@@@@@@@@@@@@@@@@@@@@
%@@@@@@@@@@@@@@@@@@@@@@@@@@@@@@@@@@@@@@@@@@@@@@@@@@@@@@@@@@
%@@@@@@@@@@@@@@@@@@@@@@@@@@@@@@@@@@@@@@@@@@@@@@@@@@@@@@@@@@
%----------------------------------------------------------
\begin{lemma} \lbl{3POINTS} Let $F\in\LTP$ be a $C^1$-function on $\RT$, and let $2<q\le p<\infty$. Let $Q$ be a \sq in $\RT$ which contains three collinear points $z_1,z_2,z_3\in\RT$ such that $z_2\in(z_1,z_3).$ Then
%----------------------------------------------------------
$$
\left|
\frac{F(z_1)-F(z_2)}{\|z_1-z_2\|_2}
-\frac{F(z_2)-F(z_3)}{\|z_2-z_3\|_2}
\right|\le
C(q)\diam Q \left(\frac{1}{|Q|} \intl_Q(\nabla^2F)^q\,dz\right)^{\frac{1}{q}}.
$$
%----------------------------------------------------------
%@@@@@@@@@@@@@@@@@@@@@@@@@@@@@@@@@@@@@@@@@@@@@@@@@@@@@@@@@@
%----------------------------------------------------------
\end{lemma}
%@@@@@@@@@@@@@@@@@@@@@@@@@@@@@@@@@@@@@@@@@@@@@@@@@@@@@@@@@@
%@@@@@@@@@@@@@@@@@@@@@@@@@@@@@@@@@@@@@@@@@@@@@@@@@@@@@@@@@@
\par {\it Proof.} Since $F$ is a $C^1$-function, there exist points $z',z''\in (z_1,z_3)$ such that
%----------------------------------------------------------
$$
F(z_1)-F(z_2)=\ip{\nabla F(z'),z_1-z_2}
 ~~~~\text{and}~~~~
F(z_2)-F(z_3)=\ip{\nabla F(z''),z_2-z_3}.
$$
%----------------------------------------------------------
Hence
%----------------------------------------------------------
\be
I&:=&\left|
\frac{F(z_1)-F(z_2)}{\|z_1-z_2\|_2}
-\frac{F(z_2)-F(z_3)}{\|z_2-z_3\|_2}
\right|\nn\\&=&
\left|\left\langle\nabla F(z'), \frac{z_1-z_2}{\|z_1-z_2\|_2} \right\rangle-\left\langle\nabla F(z''), \frac{z_2-z_3}{\|z_2-z_3\|_2}\right\rangle\right|.\nn
\ee
%----------------------------------------------------------
Since $z_2\in(z_1,z_3)$,
%----------------------------------------------------------
$$
\frac{z_1-z_2}{\|z_1-z_2\|_2}= \frac{z_2-z_3}{\|z_2-z_3\|_2}=\vec{n}
$$
%----------------------------------------------------------
where $\vec{n}\in\RT$ is a vector with $\|\vec{n}\|_2=1$. Hence, by \rf{S-2},
%----------------------------------------------------------
\be
I&=&\left|\ip{\nabla F(z')-\nabla F(z''),\vec{n}}
\right|\nn\\&\le&
\|\nabla F(z')-\nabla F(z'')\|_2\le
C(q)\diam Q \left(\frac{1}{|Q|} \intl_Q(\nabla^2F)^q\,dz\right)^{\frac{1}{q}}
\nn
\ee
%----------------------------------------------------------
proving the lemma.\bx\medskip
%@@@@@@@@@@@@@@@@@@@@@@@@@@@@@@@@@@@@@@@@@@@@@@@@@@@@@@@@@@
%@@@@@@@@@@@@@@@@@@@@@@@@@@@@@@@@@@@@@@@@@@@@@@@@@@@@@@@@@@
%@@@@@@@@@@@@@@@@@@@@@@@@@@@@@@@@@@@@@@@@@@@@@@@@@@@@@@@@@@
%@@@@@@@@@@@@@@@@@@@@@@@@@@@@@@@@@@@@@@@@@@@@@@@@@@@@@@@@@@
%----------------------------------------------------------
\par Let $\gamma\ge 1$ and let $q:=(p+2)/2$. Clearly $2<q<p$. Let $Q$ be a \sq in $\RT$ and let $Z=\{z_1,z_2,z_3\}$ be a subset of $E$ such that $ Z\subset\gamma Q$ and $z_2\in(z_1,z_3)$. Since $F|_E=f$, by Lemma \reff{3POINTS},
%----------------------------------------------------------
$$
I(Z;f):=\left|
\frac{f(z_1)-f(z_2)}{\|z_1-z_2\|_2}
-\frac{f(z_2)-f(z_3)}{\|z_2-z_3\|_2}
\right|\le
C\diam (\gamma Q) \left(\frac{1}{|\gamma Q|} \intl_{\gamma Q}(\nabla^2F)^q\,dz\right)^{\frac{1}{q}}
$$
%----------------------------------------------------------
so that
%----------------------------------------------------------
\bel{IZ}
I(Z;f)^p(\diam Q)^{2-p}\le
C|Q|\left(\frac{1}{|\gamma Q|} \intl_{\gamma Q}(\nabla^2F)^q\,dz\right)^{\frac{p}{q}}
\ee
%----------------------------------------------------------
where $C=C(p,\gamma)$.
%----------------------------------------------------------
\par Let $y\in Q$. As usual given a function $g\in L_{1,loc}(\RT)$ by $\Mc[g](y)$ we denote the \HLM function  %----------------------------------------------------------
$$
\Mc[g](y):=\sup_{K\ni y}\frac{1}{|K|} \intl_{K}|g(z)|\,dz.
$$
%----------------------------------------------------------
Here the supremum is taken over all \sqs $K$ containing $y$. Hence
%----------------------------------------------------------
$$
\left(\frac{1}{|\gamma Q|} \intl_{\gamma Q}(\nabla^2F)^q\,dz\right)^{\frac{p}{q}}\le
\left(\Mc[(\nabla^2F)^q](y)\right)^{\frac{p}{q}}.
$$
%----------------------------------------------------------
Integrating this inequality over \sq $Q$ we obtain
%----------------------------------------------------------
\bel{IQ-P}
|Q|\left(\frac{1}{|\gamma Q|} \intl_{\gamma Q}(\nabla^2F)^q\,dz\right)^{\frac{p}{q}}\le\intl_Q
\Mc[(\nabla^2F)^q]^{\frac{p}{q}}\,dz.
\ee
%----------------------------------------------------------
Combining this inequality with \rf{IZ}, we obtain
%----------------------------------------------------------
\bel{EMF}
I(Z;f)^p(\diam Q)^{2-p}\le
C\intl_Q
\Mc[(\nabla^2F)^q]^{\frac{p}{q}}\,dz.
\ee
%----------------------------------------------------------
\par Let $\{Q_i:i=1,...,m\}$ be a finite family of pairwise disjoint squares. Consider three arbitrary collinear points $z^{(i)}_1, z^{(i)}_2, z^{(i)}_3\in E\cap(\gamma Q_i)$ such that $z^{(i)}_2\in(z^{(i)}_1, z^{(i)}_3), i=1,...,m.$ We have
%----------------------------------------------------------
$$
J:=\smed_{i=1}^m\left|
\frac{f(z_1^{(i)})-f(z_2^{(i)})}{\|z_1^{(i)}-z_2^{(i)}\|}
-\frac{f(z_2^{(i)})-f(z_3^{(i)})}{\|z_2^{(i)}-z_3^{(i)}\|}
\right|^p(\diam Q_i)^{2-p}=
\sum_{i=1}^m I(Z_i;f)^p(\diam Q_i)^{2-p}
$$
%----------------------------------------------------------
so that, by \rf{EMF},
%----------------------------------------------------------
$$
J\le C\smed_{i=1}^m\, \intl_{Q_i}\Mc[(\nabla^2F)^q]^{\frac{p}{q}}\,dz\le C \intl_{\RT}\Mc[(\nabla^2F)^q]^{\frac{p}{q}}\,dz
$$
%----------------------------------------------------------
where $C=C(p,\gamma)$. Since $p/q>1$, by the Hardy-Littlewood maximal theorem,
%----------------------------------------------------------
$$
J\le C \intl_{\RT}((\nabla^2F)^q)^{\frac{p}{q}}\,dz
=C \intl_{\RT}(\nabla^2F)^p\,dz=C\|F\|_{\LTP}^p.
$$
%----------------------------------------------------------
This inequality proves that for every function $F\in\LTP\cap C^1(\RT)$ such that $F|_E=f$ the condition (i) of Theorem \reff{MAIN} holds with the constant $\lambda=C(p,\gamma)\|F\|_{\LTP}^p$.
%----------------------------------------------------------
%@@@@@@@@@@@@@@@@@@@@@@@@@@@@@@@@@@@@@@@@@@@@@@@@@@@@@@@@@@
%@@@@@@@@@@@@@@@@@@@@@@@@@@@@@@@@@@@@@@@@@@@@@@@@@@@@@@@@@@
%@@@@@@@@@@@@@@@@@@@@@@@@@@@@@@@@@@@@@@@@@@@@@@@@@@@@@@@@@@
%@@@@@@@@@@@@@@@@@@@@@@@@@@@@@@@@@@@@@@@@@@@@@@@@@@@@@@@@@@
%----------------------------------------------------------
\bigskip
%----------------------------------------------------------
\par Now let us prove that for every $\gamma\in [1,\infty)$ and every function $F\in\LTP\cap C^1(\RT)$ such that $F|_E=f$ the condition (ii) of Theorem \reff{MAIN} holds with $\lambda=C(p,\gamma)\|F\|_{\LTP}^p$.
%----------------------------------------------------------
\par We will be needed two auxiliary results. The first of them is the following
%----------------------------------------------------------
%@@@@@@@@@@@@@@@@@@@@@@@@@@@@@@@@@@@@@@@@@@@@@@@@@@@@@@@@@@
%@@@@@@@@@@@@@@@@@@@@@@@@@@@@@@@@@@@@@@@@@@@@@@@@@@@@@@@@@@
%@@@@@@@@@@@@@@@@@@@@@@@@@@@@@@@@@@@@@@@@@@@@@@@@@@@@@@@@@@
%@@@@@@@@@@@@@@@@@@@@@@@@@@@@@@@@@@@@@@@@@@@@@@@@@@@@@@@@@@
%----------------------------------------------------------
\begin{proposition} \lbl{TR1} Let $F\in\LTP$ be a $C^1$-function on $\RT$, and let $2<q\le p<\infty$. Let $\Delta=\{z_1,z_2,z_3\}$ be a triangle.
Then for every \sq $Q\supset \Delta$ and every $x\in Q$ we have
%----------------------------------------------------------
\bel{DI-1}
\|\nabla P_\Delta[F]-\nabla F(x)\|\le C(q) \frac{R_\Delta}
{\diam \Delta} \diam Q\left(\frac{1}{|Q|} \intl_Q(\nabla^2F)^q\,dz\right)^{\frac{1}{q}}.
\ee
%----------------------------------------------------------
Recall that $R_\Delta$ denotes the radius of the circle passing through the points $z_1,z_2,z_3,$ and $P_\Delta[F]$ denotes the affine polynomial interpolating $F$ at the points $z_1,z_2$ and $z_3$.
%@@@@@@@@@@@@@@@@@@@@@@@@@@@@@@@@@@@@@@@@@@@@@@@@@@@@@@@@@@
%----------------------------------------------------------
\end{proposition}
%@@@@@@@@@@@@@@@@@@@@@@@@@@@@@@@@@@@@@@@@@@@@@@@@@@@@@@@@@@
%@@@@@@@@@@@@@@@@@@@@@@@@@@@@@@@@@@@@@@@@@@@@@@@@@@@@@@@@@@
\par {\it Proof.} Without loss of generality we may assume that
%----------------------------------------------------------
\bel{TS}
\|z_1-z_2\|_2\le\|z_2-z_3\|_2\le\|z_1-z_3\|_2.
\ee
%----------------------------------------------------------
Note that, by \rf{S-2},
%----------------------------------------------------------
$$
\|\nabla F(x)-\nabla F(z_1)\|\le C(q) \diam Q \left(\frac{1}{|Q|} \intl_Q(\nabla^2F)^q\,dz\right)^{\frac{1}{q}},
$$
%----------------------------------------------------------
so that it suffice to prove the proposition for the case $x=z_1$.
%@@@@@@@@@@@@@@@@@@@@@@@@@@@@@@@@@@@@@@@@@@@@@@@@@@@@@@@@@@
\par Note that $P_\Delta[F]\in\PO$ so that
%----------------------------------------------------------
$$
P_\Delta[F](u)=P_\Delta[F](z_1)+\ip{\nabla P_\Delta[F], u-z_1}, ~~~u\in\RT.
$$
%----------------------------------------------------------
Also recall that $P_\Delta[F]$ interpolates $F$ on the vertices of $\Delta$ so that
%----------------------------------------------------------
$$
P_\Delta[F](z_i)=F(z_i),~~~~i=1,2,3.
$$
%----------------------------------------------------------
%@@@@@@@@@@@@@@@@@@@@@@@@@@@@@@@@@@@@@@@@@@@@@@@@@@@@@@@@@@
\par Let
%----------------------------------------------------------
$$
T_{z_1}[F](u)=F(z_1)+\ip{\nabla F(z_1),u-z_1}
=P_\Delta[F](z_1)+\ip{\nabla F(z_1),u-z_1},~~~u\in\RT,
$$
%----------------------------------------------------------
be the first order Taylor polynomial of $F$ at the point $z_1$. Then, by \rf{S-1},
%----------------------------------------------------------
$$
|P_\Delta[F](z_2)-T_{z_1}[F](z_2)|=|F(z_2)-T_{z_1}[F](z_2)|\le
C\|z_1-z_2\|_2\diam Q \left(\frac{1}{|Q|} \intl_Q(\nabla^2F)^q\,dz\right)^{\frac{1}{q}}.
$$
%----------------------------------------------------------
On the hand,
%----------------------------------------------------------
$$
P_\Delta[F](z_2)-T_{z_1}[F](z_2)=\ip{\nabla P_\Delta[F]-\nabla F(z_1),z_2-z_1}
$$
%----------------------------------------------------------
so that
%----------------------------------------------------------
\bel{D1}
\left|\left\langle\nabla P_\Delta[F]-\nabla F(z_1),\frac{z_2-z_1}{\|z_2-z_1\|_2}\right\rangle\right|
\le
C\diam Q \left(\frac{1}{|Q|} \intl_Q(\nabla^2F)^q\,dz\right)^{\frac{1}{q}}.
\ee
%----------------------------------------------------------
In a similar way we prove that
%----------------------------------------------------------
\bel{D2}
\left|\left\langle\nabla P_\Delta[F]-\nabla F(z_1),\frac{z_3-z_1}{\|z_3-z_1\|_2}\right\rangle\right|
\le
C\diam Q \left(\frac{1}{|Q|} \intl_Q(\nabla^2F)^q\,dz\right)^{\frac{1}{q}}.
\ee
%----------------------------------------------------------
\par Now let
%----------------------------------------------------------
$$
\vec{a}:=\nabla P_\Delta[F]-\nabla F(z_1), ~~\vec{n}:=\frac{z_2-z_1}{\|z_2-z_1\|_2},~~
\vec{m}:=\frac{z_3-z_1}{\|z_3-z_1\|_2},
$$
%----------------------------------------------------------
and
%----------------------------------------------------------
$$
A:=C\diam Q \left(\frac{1}{|Q|} \intl_Q(\nabla^2F)^q\,dz\right)^{\frac{1}{q}}.
$$
%----------------------------------------------------------
Then, by \rf{D1} and \rf{D2},
%----------------------------------------------------------
\bel{EAN}
|\ip{\vec{a},\vec{n}}|\le \,A ~~~\text{and}~~~|\ip{\vec{a},\vec{m}}|\le \,A.
\ee
%----------------------------------------------------------
Let $\alpha\in(0,\pi)$ be the angle between the sides $[z_1,z_2]$ and $[z_1,z_3]$ of the triangle $\Delta$. Thus $\alpha$ is also the angle between $\vec{n}$ and $\vec{m}$. %----------------------------------------------------------
\par Let $\vec{n}_1$ be a unit vector in $\RT$ which is orthogonal to $\vec{n}$, i.e., $\|\vec{n}_1\|_2=1$ and $\ip{\vec{n},\vec{n}_1}=0.$ Then $\vec{m}=(\cos\alpha)\vec{n}+\ve(\sin\alpha)\vec{n_1}$ where $\ve\in\{1,-1\}$. Hence
%----------------------------------------------------------
$$
|\ip{\vec{a},\vec{n_1}}|=\frac{1}{\sin\alpha} |\ip{\vec{a},\vec{m}-(\cos\alpha)\vec{n}}|\le
\frac{1}{\sin\alpha}(|\ip{\vec{a},\vec{m}}| +|\ip{\vec{a},\vec{n}}|)
$$
%----------------------------------------------------------
so that, by \rf{EAN},
$|\ip{\vec{a},\vec{n_1}}|\le 2A/\sin\alpha.$ We obtain
%----------------------------------------------------------
$$
\|a\|_2=(\ip{\vec{a},\vec{n}}^2 +\ip{\vec{a},\vec{n_1}}^2)^{\frac12}\le \left(A^2+\left(\frac{2A}{\sin\alpha} \right)^2\right)^{\frac12}\le C_1A/\sin\alpha.
$$
%----------------------------------------------------------
Hence
%----------------------------------------------------------
\bel{FEP}
\|\nabla P_\Delta[F]-\nabla F(z_1)\|\le C\,\frac{\diam Q}{\sin\alpha} \left(\frac{1}{|Q|} \intl_Q(\nabla^2F)^q\,dz\right)^{\frac{1}{q}}.
\ee
%----------------------------------------------------------
\par Recall that $\alpha$ is the angle in the triangle $\Delta$ which is opposite to the side $[z_2,z_3]$  so that %----------------------------------------------------------
$$
R_\Delta=\frac{\|z_2-z_3\|_2}{2\sin\alpha}\,\,,
$$
%----------------------------------------------------------
see, e.g. \cite{P}, p. 29.  But, by \rf{TS} and the triangle inequality,
%----------------------------------------------------------
$$
\diam \Delta=\|z_1-z_3\|_2\le 2\|z_2-z_3\|_2
$$
%----------------------------------------------------------
so that
%----------------------------------------------------------
$$
\frac{1}{\sin\alpha}=\frac{2R_\Delta}{\|z_2-z_3\|_2}\le \frac{4R_\Delta}{\diam \Delta}\,\,.
$$
%----------------------------------------------------------
Combining this inequality with \rf{FEP} we obtain the required inequality \rf{DI-1}.\bx
%----------------------------------------------------------
%@@@@@@@@@@@@@@@@@@@@@@@@@@@@@@@@@@@@@@@@@@@@@@@@@@@@@@@@@@
%@@@@@@@@@@@@@@@@@@@@@@@@@@@@@@@@@@@@@@@@@@@@@@@@@@@@@@@@@@
%@@@@@@@@@@@@@@@@@@@@@@@@@@@@@@@@@@@@@@@@@@@@@@@@@@@@@@@@@@
%@@@@@@@@@@@@@@@@@@@@@@@@@@@@@@@@@@@@@@@@@@@@@@@@@@@@@@@@@@
%----------------------------------------------------------
\bigskip
%----------------------------------------------------------
\par Proposition \reff{TR1} and inequality \rf{S-1} imply the following
%@@@@@@@@@@@@@@@@@@@@@@@@@@@@@@@@@@@@@@@@@@@@@@@@@@@@@@@@@@
%@@@@@@@@@@@@@@@@@@@@@@@@@@@@@@@@@@@@@@@@@@@@@@@@@@@@@@@@@@
%@@@@@@@@@@@@@@@@@@@@@@@@@@@@@@@@@@@@@@@@@@@@@@@@@@@@@@@@@@
%@@@@@@@@@@@@@@@@@@@@@@@@@@@@@@@@@@@@@@@@@@@@@@@@@@@@@@@@@@
\begin{proposition} \lbl{2TR} Let $2<q\le p<\infty$ and let $\Delta=\{z_1,z_2,z_3\},\Delta'=\{z_1',z_2',z_3'\}$ be two triangles in $\RT$. Let $Q,Q',\tQ$ be squares in $\RT$ such that
%----------------------------------------------------------
$$
\Delta\subset Q,~~\Delta'\subset Q',~~Q\cup Q'\subset\tQ.
$$
%----------------------------------------------------------
\par Then for every $C^1$-function $F\in\LTP$ the following inequality
%----------------------------------------------------------
\be
&&\|\nabla P_\Delta[F]-\nabla P_{\Delta'}[F]\|\le C(q)\left\{\frac{R_\Delta}{\diam\Delta}\diam Q \left(\frac{1}{|Q|} \intl_Q(\nabla^2F)^q\,dz\right)^{\frac{1}{q}}\right.\nn\\
&+&
\left.\frac{R_{\Delta'}}{\diam\Delta'} \diam Q' \left(\frac{1}{|Q'|} \intl_{Q'}(\nabla^2F)^q\,dz\right)^{\frac{1}{q}}
+\diam \tQ \left(\frac{1}{|\tQ|} \intl_{\tQ}(\nabla^2F)^q\,dz\right)^{\frac{1}{q}}
\right\}\nn
\ee
%----------------------------------------------------------
holds.
%@@@@@@@@@@@@@@@@@@@@@@@@@@@@@@@@@@@@@@@@@@@@@@@@@@@@@@@@@@
%----------------------------------------------------------
\end{proposition}
%----------------------------------------------------------
%@@@@@@@@@@@@@@@@@@@@@@@@@@@@@@@@@@@@@@@@@@@@@@@@@@@@@@@@@@
%@@@@@@@@@@@@@@@@@@@@@@@@@@@@@@@@@@@@@@@@@@@@@@@@@@@@@@@@@@
%@@@@@@@@@@@@@@@@@@@@@@@@@@@@@@@@@@@@@@@@@@@@@@@@@@@@@@@@@@
%@@@@@@@@@@@@@@@@@@@@@@@@@@@@@@@@@@@@@@@@@@@@@@@@@@@@@@@@@@
%@@@@@@@@@@@@@@@@@@@@@@@@@@@@@@@@@@@@@@@@@@@@@@@@@@@@@@@@@@
%@@@@@@@@@@@@@@@@@@@@@@@@@@@@@@@@@@@@@@@@@@@@@@@@@@@@@@@@@@
%---------------------------------------------------------- %@@@@@@@@@@@@@@@@@@@@@@@@@@@@@@@@@@@@@@@@@@@@@@@@@@@@@@@@@@
%@@@@@@@@@@@@@@@@@@@@@@@@@@@@@@@@@@@@@@@@@@@@@@@@@@@@@@@@@@
%@@@@@@@@@@@@@@@@@@@@@@@@@@@@@@@@@@@@@@@@@@@@@@@@@@@@@@@@@@
%@@@@@@@@@@@@@@@@@@@@@@@@@@@@@@@@@@@@@@@@@@@@@@@@@@@@@@@@@@
%@@@@@@@@@@@@@@@@@@@@@@@@@@@@@@@@@@@@@@@@@@@@@@@@@@@@@@@@@@
%@@@@@@@@@@@@@@@@@@@@@@@@@@@@@@@@@@@@@@@@@@@@@@@@@@@@@@@@@@
%@@@@@@@@@@@@@@@@@@@@@@@@@@@@@@@@@@@@@@@@@@@@@@@@@@@@@@@@@@
%@@@@@@@@@@@@@@@@@@@@@@@@@@@@@@@@@@@@@@@@@@@@@@@@@@@@@@@@@@
%@@@@@@@@@@@@@@@@@@@@@@@@@@@@@@@@@@@@@@@@@@@@@@@@@@@@@@@@@@
%----------------------------------------------------------
\bigskip
%----------------------------------------------------------
\par {\bf 3.2. Part (ii) of the necessity and the space
$\VS:=\VLOP+\VLPM$.}
%----------------------------------------------------------
\addtocontents{toc}{~~~~3.2. Part (ii) of the necessity and the space $\VS:=\VLOP+\VLPM$. \hfill \thepage\par}
%---------------------------------------------------------- %@@@@@@@@@@@@@@@@@@@@@@@@@@@@@@@@@@@@@@@@@@@@@@@@@@@@@@@@@@
%@@@@@@@@@@@@@@@@@@@@@@@@@@@@@@@@@@@@@@@@@@@@@@@@@@@@@@@@@@
%@@@@@@@@@@@@@@@@@@@@@@@@@@@@@@@@@@@@@@@@@@@@@@@@@@@@@@@@@@
%@@@@@@@@@@@@@@@@@@@@@@@@@@@@@@@@@@@@@@@@@@@@@@@@@@@@@@@@@@
%----------------------------------------------------------
The next auxiliary result relates to an optimal decomposition of a function into a sum of a Sobolev function and an $L_p$-weighted function.
%----------------------------------------------------------
\par Let $\mu $ be a non-trivial non-negative Borel measure on $\RT$. Let $2<p<\infty $ and let $\LPM$ be the $L_p-$space on $\RT$ with respect to the measure $\mu$. We norm this space by
%----------------------------------------------------------
$$ \|f\|_{\LPM}=\left(\,\,\int\limits_{\RT}|f|^{p}d\mu\right)^
{\frac{1}{p}}.
$$
%----------------------------------------------------------
%@@@@@@@@@@@@@@@@@@@@@@@@@@@@@@@@@@@@@@@@@@@@@@@@@@@@@@@@@@
%----------------------------------------------------------
\par By $\sum$ we denote the space $\LOP+\LPM$ equipped with the norm:
%----------------------------------------------------------
$$
\|f\|_{\sum} =\inf\{\|f_1\|_{\LOP}+\|f_2\|_{\LPM}:
f=f_1+f_2 , f_1 \in \LOP ,f_2\in\LPM\}.
$$
%----------------------------------------------------------
\par Let us also define vector versions of the above spaces. We let $\VLOP$ denote the space of Sobolev mappings $\VF=(F_1,F_2):\RT\to\RT$ whose components $F_1,F_2\in \LOP.$ This space is normed by
%----------------------------------------------------------
\bel{N-VLOP}
\|\VF\|_{\VLOP}:=\|F_1\|_{\LOP}+\|F_2\|_{\LOP}.
\ee
%----------------------------------------------------------
In turn, by  $\VLPM$ we denote the space of all mappings  $\VF=(F_1,F_2):\RT\to\RT$ whose components $F_1,F_2\in \LPM$. We norm $\VLPM$ by
%----------------------------------------------------------
$$
\|\VF\|_{\VLPM}:=\|F_1\|_{\LPM}+\|F_2\|_{\LPM}.
$$
%----------------------------------------------------------
By $\VS=\VLOP+\VLPM$ we denote the sum of these spaces. The space  $\VS$ is normed by
%----------------------------------------------------------
$$
\|\VF\|_{\VS} =\inf\{\|\VF_1\|_{\VLOP}+\|\VF_2\|_{\VLPM}:
\VF=\VF_1+\VF_2 ,\VF_1 \in \VLOP ,\VF_2\in\VLPM\}.
$$
%----------------------------------------------------------
%@@@@@@@@@@@@@@@@@@@@@@@@@@@@@@@@@@@@@@@@@@@@@@@@@@@@@@@@@@
%@@@@@@@@@@@@@@@@@@@@@@@@@@@@@@@@@@@@@@@@@@@@@@@@@@@@@@@@@@
%@@@@@@@@@@@@@@@@@@@@@@@@@@@@@@@@@@@@@@@@@@@@@@@@@@@@@@@@@@
%----------------------------------------------------------
\par In \cite{S3}, Theorem 2.4, we present
a necessary condition for a function to belong to the space   $\sum=\LOP+\LPM$ whenever $p>2$. Applying this result to every component of a mapping $\Tc:\RT\to\RT$ we obtain the following
%----------------------------------------------------------
%@@@@@@@@@@@@@@@@@@@@@@@@@@@@@@@@@@@@@@@@@@@@@@@@@@@@@@@@@@
%@@@@@@@@@@@@@@@@@@@@@@@@@@@@@@@@@@@@@@@@@@@@@@@@@@@@@@@@@@
%@@@@@@@@@@@@@@@@@@@@@@@@@@@@@@@@@@@@@@@@@@@@@@@@@@@@@@@@@@
%@@@@@@@@@@@@@@@@@@@@@@@@@@@@@@@@@@@@@@@@@@@@@@@@@@@@@@@@@@
%@@@@@@@@@@@@@@@@@@@@@@@@@@@@@@@@@@@@@@@@@@@@@@@@@@@@@@@@@@
%----------------------------------------------------------
\begin{theorem}\lbl{SL-DEC} Let $2<p<\infty, 1\le\gamma<\infty,$ and let $\mu$ be a non-trivial non-negative Borel measure on $\RT$. Suppose that a mapping $\Tc\in\VS=\VLOP+\VLPM$. Then the following statement is true:
%----------------------------------------------------------
\par Let $\Ac$ be a finite family of  \sqs in $\RT$, and let $\Sc$ be a finite family of closed subsets in $\RT$ with covering multiplicity $M(\Ac),M(\Sc)\le N$. Suppose that to every \sq $Q\in\Ac$ we have assigned  two subsets $S_Q',S_Q''\in \Sc$ such that
%----------------------------------------------------------
\bel{InG}
S_Q'\cup S_Q''\subset \gamma Q.
\ee
%----------------------------------------------------------
\par Then
%----------------------------------------------------------
\bel{F-H}
\shuge_{Q\in\Ac}\,\,
\frac{(\diam Q)^{2-p}\iint \limits_{S'_Q\times S''_Q}
\|\Tc(x)-\Tc(y)\|^p\, d\mu(x)d\mu(y)}
{ \{(\diam S'_Q)^{2-p}+\mu(S'_Q)\} \{(\diam S''_Q)^{2-p}+\mu(S''_Q)\}}
\le C(\gamma,N)\|\Tc\|^p_{\VS}\,.
\ee
%----------------------------------------------------------
\end{theorem}
%----------------------------------------------------------
%@@@@@@@@@@@@@@@@@@@@@@@@@@@@@@@@@@@@@@@@@@@@@@@@@@@@@@@@@@
%@@@@@@@@@@@@@@@@@@@@@@@@@@@@@@@@@@@@@@@@@@@@@@@@@@@@@@@@@@
%@@@@@@@@@@@@@@@@@@@@@@@@@@@@@@@@@@@@@@@@@@@@@@@@@@@@@@@@@@
%@@@@@@@@@@@@@@@@@@@@@@@@@@@@@@@@@@@@@@@@@@@@@@@@@@@@@@@@@@
%----------------------------------------------------------
\par We are in a position to prove the condition (ii) of Theorem \reff{MAIN} with $\lambda=C(p,\gamma)\|F\|_{\LTP}^p$ for every $\gamma\in[1,\infty)$ and every $C^1$-function $F\in\LTP$ such that $F|_E=f$.
%----------------------------------------------------------
\par Let $\Kc$ be a finite family of pairwise disjoint squares. Given $K\in\Kc$ let $\Delta(K)$ be a triangle satisfying conditions \rf{GSEP} with some constant $\gamma\in[1,\infty)$. Let
%----------------------------------------------------------
$$
\Cc:=\{c_K: K\in \Kc\}
$$
%----------------------------------------------------------
be the family of centers of all \sqs from $\Kc$.
%----------------------------------------------------------
\par We introduce a discrete Borel measure $\mu$ on $\RT$ with $\supp \mu:=\Cc$ as follows: for every \sq $K\in\Kc$ we put
%----------------------------------------------------------
\bel{D-MU}
\mu(\{c_K\}):=\cu_{\Delta(K)}^p\,|K|.
\ee
%----------------------------------------------------------
Thus for every set $S\subset\R^2$ we have
%----------------------------------------------------------
$$
\mu(S)=\sum\left\{\cu_{\Delta(K)}^p\,|K|: K\in\Kc,
c_K\in S\right\}.
$$
%----------------------------------------------------------
%@@@@@@@@@@@@@@@@@@@@@@@@@@@@@@@@@@@@@@@@@@@@@@@@@@@@@@@@@@
%@@@@@@@@@@@@@@@@@@@@@@@@@@@@@@@@@@@@@@@@@@@@@@@@@@@@@@@@@@
%@@@@@@@@@@@@@@@@@@@@@@@@@@@@@@@@@@@@@@@@@@@@@@@@@@@@@@@@@@
%@@@@@@@@@@@@@@@@@@@@@@@@@@@@@@@@@@@@@@@@@@@@@@@@@@@@@@@@@@
\begin{proposition} \lbl{M-PL} Let $p\in(2,\infty)$ and let $F\in\LTP$. Let $\Tc:\RT\to\RT$ be a mapping such that
%----------------------------------------------------------
\bel{DGV-1}
\Tc(c_K):=\nabla P_{\Delta(K)}[F]~~~\text{for every}~~~K\in\Kc,
\ee
%----------------------------------------------------------
and
%----------------------------------------------------------
\bel{DGV-2}
\Tc(x):=0~~~\text{whenever}~~~x\in \RT\setminus \Cc.
\ee
%----------------------------------------------------------
\par Then $\Tc\in \VS=\VLOP+\VLPM$. Furthermore,
%----------------------------------------------------------
\bel{COORD}
\|\Tc\|_{\VS}\le C(p,\gamma)\|F\|_{\LTP}.
\ee
%----------------------------------------------------------
%@@@@@@@@@@@@@@@@@@@@@@@@@@@@@@@@@@@@@@@@@@@@@@@@@@@@@@@@@@
%----------------------------------------------------------
\end{proposition}
%@@@@@@@@@@@@@@@@@@@@@@@@@@@@@@@@@@@@@@@@@@@@@@@@@@@@@@@@@@
%@@@@@@@@@@@@@@@@@@@@@@@@@@@@@@@@@@@@@@@@@@@@@@@@@@@@@@@@@@
\par {\it Proof.} Fix $q\in(2,p)$ (say $q=(p+2)/2$) and a \sq $K\in \Kc$. Applying Proposition \reff{TR1} to the \sq $Q=\gamma K$ and arbitrary  $x=c_K$ we obtain
%----------------------------------------------------------
$$
\|\nabla P_{\Delta(K)}[F]-\nabla F(c_K)\|\le C R_{\Delta(K)}
\frac{\diam K}{\diam \Delta(K)} \left(\frac{1}{|\gamma K|} \intl_{\gamma K}(\nabla^2F)^q\,dz\right)^{\frac{1}{q}}.
$$
%----------------------------------------------------------
where $C=C(p,\gamma)$. Since $\diam \Delta(K)\sim \diam K$, we have
%----------------------------------------------------------
\be
\|\nabla P_{\Delta(K)}[F]-\nabla F(c_K)\|&\le& C R_{\Delta(K)}
\left(\frac{1}{|\gamma K|} \intl_{\gamma K}(\nabla^2F)^q\,dz\right)^{\frac{1}{q}}\nn\\&\le& C R_{\Delta(K)}\, \Mc[(\nabla^2F)^q]^{\frac{1}{q}}(y), ~~~y\in K.\nn
\ee
%----------------------------------------------------------
(Recall that $\Mc$ denotes the Hardy-Littlewood maximal function.) Hence
%----------------------------------------------------------
$$
\|\nabla P_{\Delta(K)}[F]-\nabla F(c_K)\|^p \cu_{\Delta(K)}^p
\le C\Mc[(\nabla^2F)^q]^{\frac{p}{q}}(y), ~~~y\in K.
$$
%----------------------------------------------------------
Integrating this inequality over \sq $K$ we obtain
%----------------------------------------------------------
$$
\|\nabla P_{\Delta(K)}[F]-\nabla F(x)\|^p \cu_{\Delta(K)}^p|K|
\le C\intl_K\Mc[(\nabla^2F)^q]^{\frac{p}{q}}(y)dy.
$$
%----------------------------------------------------------
Hence, by \rf{DGV-1} and \rf{DGV-2},
%----------------------------------------------------------
\be
\intl_{\RT}\|\Tc(x)-\nabla F(x)\|^p d\mu(x)&=&\sum_{K\in \Kc}\|\nabla P_{\Delta(K)}[F]-\nabla F(x)\|^p\, \cu_{\Delta(K)}^p|K|\nn\\&\le& C\sum_{K\in \Kc}\intl_K \Mc[(\nabla^2F)^q]^{\frac{p}{q}}(y)dy\nn\\&\le& C\intl_{\RT} \Mc[(\nabla^2F)^q]^{\frac{p}{q}}(y)dy.\nn
\ee
%----------------------------------------------------------
Since $p/q>1$, by the Hardy-Littlewood maximal theorem,
%----------------------------------------------------------
$$
\intl_{\RT}\|\Tc(x)-\nabla F(c_K)\|^p d\mu(x)\le C\intl_{\RT}
((\nabla^2F)^q)^{\frac{p}{q}}(y)dy= C\intl_{\RT}(\nabla^2F)^p(y)dy
$$
%----------------------------------------------------------
proving that
%----------------------------------------------------------
$$
\|\Tc-\nabla F\|_{\VLPM}^p=\intl_{\RT}\|\Tc(x)-\nabla F(x)\|^p d\mu(x)\le C\|F\|_{\LTP}^p.
$$
%----------------------------------------------------------
\par Finally,
%----------------------------------------------------------
$$
\|\Tc\|_{\VS}\le \|\nabla F\|_{\VLOP}+\|\Tc-
\nabla F\|_{\VLPM}\le C\|F\|_{\LTP}.
$$
%----------------------------------------------------------
\par The proposition is proved.\bx\medskip
%----------------------------------------------------------
%@@@@@@@@@@@@@@@@@@@@@@@@@@@@@@@@@@@@@@@@@@@@@@@@@@@@@@@@@@
%@@@@@@@@@@@@@@@@@@@@@@@@@@@@@@@@@@@@@@@@@@@@@@@@@@@@@@@@@@
%@@@@@@@@@@@@@@@@@@@@@@@@@@@@@@@@@@@@@@@@@@@@@@@@@@@@@@@@@@
%@@@@@@@@@@@@@@@@@@@@@@@@@@@@@@@@@@@@@@@@@@@@@@@@@@@@@@@@@@
%@@@@@@@@@@@@@@@@@@@@@@@@@@@@@@@@@@@@@@@@@@@@@@@@@@@@@@@@@@
%@@@@@@@@@@@@@@@@@@@@@@@@@@@@@@@@@@@@@@@@@@@@@@@@@@@@@@@@@@
%----------------------------------------------------------
\par Let $\Qc$ be a finite family of squares with covering multiplicity $M(\Qc)\le N$. Given $Q\in\Qc$ let $Q',Q''\in\Qc$ be squares such that $Q'\cup Q''\subset \gamma Q$.
%----------------------------------------------------------
\par Let us apply Theorem \reff{SL-DEC} to the mapping $\Tc$ whenever
%----------------------------------------------------------
\bel{ASQ}
\Ac=\Sc=\Qc
\ee
%----------------------------------------------------------
and $S'_Q=Q', S''_Q=Q''$ for arbitrary $Q\in\Ac$. Then inclusion \rf{InG} is equivalent to inclusion $Q'\cup Q''\subset\gamma Q$. Furthermore, $\diam S'_Q=\diam Q',\,\diam S''_Q=\diam Q''$ and, by \rf{SGM},
%----------------------------------------------------------
$$
\mu(S'_Q)=\sigma_p(Q';\Kc)~~~\text{and}~~~
\mu(S''_Q)=\sigma_p(Q'';\Kc).
$$
%----------------------------------------------------------
Finally,
%----------------------------------------------------------
\be
&&\iint \limits_{S'_Q\times S''_Q}
\|\Tc(x)-\Tc(y)\|^p\, d\mu(x)d\mu(y)\nn\\
&=&
\sum_{\substack {K'\in\,\Kc\\c_{K'}\in Q'}}\,\,
\sum_{\substack {K''\in\,\Kc\\c_{K''}\in Q''}}
\|\nabla P_{\Delta(K')}[F]-
\nabla P_{\Delta(K'')}[F]\|^p\,
\cu_{\Delta(K')}^p |K'|\,\cu_{\Delta(K'')}^p |K''|\nn\\
&=&
S_p(F:Q',Q'';\Kc).\nn
\ee
%----------------------------------------------------------
But $F|_E=f$ so that $S_p(F:Q',Q'';\Kc)=S_p(f:Q',Q'';\Kc)$ proving that
%----------------------------------------------------------
$$
\iint \limits_{S'_Q\times S''_Q}
\|\Tc(x)-\Tc(y)\|^p\, d\mu(x)d\mu(y)=S_p(f:Q',Q'';\Kc).
$$
%----------------------------------------------------------
\par Hence, by inequalities \rf{F-H} and \rf{COORD},
%----------------------------------------------------------
\be
&&\sbig\limits_{Q\in\Qc}\,
\frac{(\diam Q)^{2-p}\,\,S_p(f:Q',Q'';\Kc)}
{\left\{(\diam Q')^{2-p} +\sigma_p(Q';\Kc)\right\}
\left\{(\diam Q'')^{2-p} +\sigma_p(Q'';\Kc)\right\}}\nn\\
&=&
\sbig_{Q\in\Ac}\,\,
\frac{(\diam Q)^{2-p}\iint \limits_{S'_Q\times S''_Q}
\|\Tc(x)-\Tc(y)\|^p\, d\mu(x)d\mu(y)}
{ \{(\diam S'_Q)^{2-p}+\mu(S'_Q)\} \{(\diam S''_Q)^{2-p}+\mu(S''_Q)\}}\nn\\\nn\\&\le&
C\|\Tc\|^p_{\VS}\le C(p,\gamma)\|F\|_{\LTP}^p\,.\nn
%----------------------------------------------------------
\ee
%----------------------------------------------------------
\par This proves that for every function $F\in\LTP\cap C^1(\RT)$ such that  $F|_E=f$ the condition (ii) of Theorem \reff{MAIN} holds with $\lambda=C(p,\gamma)\|F\|_{\LTP}^p$.
%----------------------------------------------------------
%@@@@@@@@@@@@@@@@@@@@@@@@@@@@@@@@@@@@@@@@@@@@@@@@@@@@@@@@@@
%@@@@@@@@@@@@@@@@@@@@@@@@@@@@@@@@@@@@@@@@@@@@@@@@@@@@@@@@@@
%@@@@@@@@@@@@@@@@@@@@@@@@@@@@@@@@@@@@@@@@@@@@@@@@@@@@@@@@@@
%@@@@@@@@@@@@@@@@@@@@@@@@@@@@@@@@@@@@@@@@@@@@@@@@@@@@@@@@@@
%@@@@@@@@@@@@@@@@@@@@@@@@@@@@@@@@@@@@@@@@@@@@@@@@@@@@@@@@@@
%@@@@@@@@@@@@@@@@@@@@@@@@@@@@@@@@@@@@@@@@@@@@@@@@@@@@@@@@@@
%---------------------------------------------------------- %@@@@@@@@@@@@@@@@@@@@@@@@@@@@@@@@@@@@@@@@@@@@@@@@@@@@@@@@@@
%@@@@@@@@@@@@@@@@@@@@@@@@@@@@@@@@@@@@@@@@@@@@@@@@@@@@@@@@@@
%@@@@@@@@@@@@@@@@@@@@@@@@@@@@@@@@@@@@@@@@@@@@@@@@@@@@@@@@@@
%@@@@@@@@@@@@@@@@@@@@@@@@@@@@@@@@@@@@@@@@@@@@@@@@@@@@@@@@@@
%@@@@@@@@@@@@@@@@@@@@@@@@@@@@@@@@@@@@@@@@@@@@@@@@@@@@@@@@@@
%@@@@@@@@@@@@@@@@@@@@@@@@@@@@@@@@@@@@@@@@@@@@@@@@@@@@@@@@@@
%@@@@@@@@@@@@@@@@@@@@@@@@@@@@@@@@@@@@@@@@@@@@@@@@@@@@@@@@@@
%@@@@@@@@@@@@@@@@@@@@@@@@@@@@@@@@@@@@@@@@@@@@@@@@@@@@@@@@@@
%@@@@@@@@@@@@@@@@@@@@@@@@@@@@@@@@@@@@@@@@@@@@@@@@@@@@@@@@@@
%----------------------------------------------------------
\bigskip
%----------------------------------------------------------
\par {\bf 3.3. A direct proof of part (ii) of the necessity.}
%----------------------------------------------------------
\addtocontents{toc}{~~~~3.3. A direct proof of part (ii) of the necessity. \hfill \thepage\\\par}
%---------------------------------------------------------- %@@@@@@@@@@@@@@@@@@@@@@@@@@@@@@@@@@@@@@@@@@@@@@@@@@@@@@@@@@
%@@@@@@@@@@@@@@@@@@@@@@@@@@@@@@@@@@@@@@@@@@@@@@@@@@@@@@@@@@
%@@@@@@@@@@@@@@@@@@@@@@@@@@@@@@@@@@@@@@@@@@@@@@@@@@@@@@@@@@
%@@@@@@@@@@@@@@@@@@@@@@@@@@@@@@@@@@@@@@@@@@@@@@@@@@@@@@@@@@
%----------------------------------------------------------
As we have mentioned in Section 2, for the reader's convenience, we also give a direct proof of part (ii) of Theorem \reff{MAIN} which does not use any results of \cite{S3}.
%----------------------------------------------------------
%@@@@@@@@@@@@@@@@@@@@@@@@@@@@@@@@@@@@@@@@@@@@@@@@@@@@@@@@@@
\par We will prove a general result which implies the statement of part (ii).
%----------------------------------------------------------
\par Let $\gamma\ge 1$ be a constant. Let  $\Qc$ be a finite families of pairwise disjoint squares and let $\Sc$ be a finite families of pairwise disjoint sets in $\RT$. Given a \sq $Q\in\Qc$ let $S'_Q,S''_Q$ be a pair of sets from $\Sc$ such that $S'_Q\cup S''_Q\subset\gamma Q$.
%----------------------------------------------------------
\par Let $\Kc$ be a finite family of pairwise disjoint squares. Given $K\in\Kc$ let $\Delta(K)$ be a triangle such that
%----------------------------------------------------------
$$
\Delta(K)\subset\gamma K~~~\text{and}~~~\diam K\le\gamma\diam \Delta(K).
$$
%----------------------------------------------------------
\par Let $S\in\Sc$ and let
%----------------------------------------------------------
\bel{SG-S}
\sigma_p(S;\Kc):=
\sum\left\{\cu_{\Delta(K)}^p|K|: K\in\Kc, c_K\in\,S\right\}.
\ee
%----------------------------------------------------------
%@@@@@@@@@@@@@@@@@@@@@@@@@@@@@@@@@@@@@@@@@@@@@@@@@@@@@@@@@@
%----------------------------------------------------------
Given a function $F:\RT\to\RT$ and sets $S'_Q,S''_Q\in\Sc$ let
%----------------------------------------------------------
$$
S_p(F:S'_Q,S''_Q;\Kc):=
\sum_{\substack {K'\in\,\Kc\\c_{K'}\in S'_Q}}\,\,
\sum_{\substack {K''\in\,\Kc\\c_{K''}\in S''_Q}}
\|\nabla P_{\Delta(K')}[F]-
\nabla P_{\Delta(K'')}[F]\|^p\,
\cu_{\Delta(K')}^p |K'|\,\cu_{\Delta(K'')}^p |K''|.
$$
%---------------------------------------------------------
\par Finally we put
%----------------------------------------------------------
$$
U\{\Qc\}:=\bigcup_{Q\in\Qc} Q~~~~~\text{and}~~~~
U\{\Kc\}:=\bigcup_{K\in\Kc} K.
$$
%----------------------------------------------------------
%@@@@@@@@@@@@@@@@@@@@@@@@@@@@@@@@@@@@@@@@@@@@@@@@@@@@@@@@@@
%@@@@@@@@@@@@@@@@@@@@@@@@@@@@@@@@@@@@@@@@@@@@@@@@@@@@@@@@@@
%@@@@@@@@@@@@@@@@@@@@@@@@@@@@@@@@@@@@@@@@@@@@@@@@@@@@@@@@@@
%@@@@@@@@@@@@@@@@@@@@@@@@@@@@@@@@@@@@@@@@@@@@@@@@@@@@@@@@@@
%----------------------------------------------------------
\begin{proposition} \lbl{MF-N} Let $2<q<p$. Then for every smooth function $F\in\LTP$ the following inequality
%----------------------------------------------------------
\be
&&\sbig\limits_{Q\in\Qc}\,
\frac{(\diam Q)^{2-p}\,\,S_p(F:S'_Q,S''_Q;\Kc)}
{\left\{(\diam S'_Q)^{2-p} +\sigma_p(S'_Q;\Kc)\right\}
\left\{(\diam S''_Q)^{2-p} +\sigma_p(S''_Q;\Kc)\right\}}\nn\\
&\le&
C\, \left\{\intl_{U\{\Kc\}}
\Mc[(\nabla^2F)^q]^{\frac{p}{q}}\,dz+
\intl_{U\{\Qc\}}
\Mc[(\nabla^2F)^q]^{\frac{p}{q}}\,dz\right\}.
\nn
\ee
%----------------------------------------------------------
holds. Here $C$ is a constant depending only on $q,p,$ and $\gamma$.
%----------------------------------------------------------
%@@@@@@@@@@@@@@@@@@@@@@@@@@@@@@@@@@@@@@@@@@@@@@@@@@@@@@@@@@
%----------------------------------------------------------
\end{proposition}
%@@@@@@@@@@@@@@@@@@@@@@@@@@@@@@@@@@@@@@@@@@@@@@@@@@@@@@@@@@
%@@@@@@@@@@@@@@@@@@@@@@@@@@@@@@@@@@@@@@@@@@@@@@@@@@@@@@@@@@
\par {\it Proof.} Consider two \sqs $K',K''\in\Kc$ such that $c_{K'}\in S'_Q$ and $c_{K''}\in S''_Q$. Let us apply Proposition \reff{2TR} to \sqs $\gamma K',\gamma K'',\gamma Q,$ and triangles $\Delta(K')\subset\gamma K'$ and $\Delta(K'')\subset\gamma K''$. By this proposition,
%----------------------------------------------------------
\be
\|\nabla P_{\Delta(K')}[F]-\nabla P_{{\Delta(K'')}}[F]\|&\le& C\left\{\frac{R_{\Delta(K')}}{\diam\Delta(K')}\diam (\gamma K') \left(\frac{1}{|\gamma K'|} \intl_{\gamma K'}(\nabla^2F)^q\,dz\right)^{\frac{1}{q}}\right.\nn\\
&+&
\frac{R_{\Delta(K'')}}{\diam\Delta(K'')}\diam (\gamma K'') \left(\frac{1}{|\gamma K''|} \intl_{\gamma K''}(\nabla^2F)^q\,dz\right)^{\frac{1}{q}}
\nn\\
&+&\left.\diam (\gamma Q) \left(\frac{1}{|\gamma Q|} \intl_{\gamma Q}(\nabla^2F)^q\,dz\right)^{\frac{1}{q}}
\right\}.\nn
\ee
%----------------------------------------------------------
Recall that $\diam \Delta(K')\sim \diam K',
\diam \Delta(K'')\sim \diam K'',$. We also recall that for every triangle $\Delta\subset\RT$ we have $\cu_\Delta=1/R_\Delta$. Hence
%----------------------------------------------------------
\be
&&\|\nabla P_{\Delta(K')}[F]-\nabla P_{{\Delta(K'')}}[F]\|^p\cu_{\Delta(K')}^p
\cu_{\Delta(K'')}^p|K'||K''|\nn\\&\le& C\left\{\cu_{\Delta(K'')}^p|K''||K'|
\left(\frac{1}{|\gamma K'|}\intl_{\gamma K'}(\nabla^2F)^q\,dz\right)^{\frac{p}{q}}
\right.\nn\\&+&
\cu_{\Delta(K')}^p|K'||K''|
\left(\frac{1}{|\gamma K''|}\intl_{\gamma K''}(\nabla^2F)^q\,dz\right)^{\frac{p}{q}}\nn\\
&+&\left.\cu_{\Delta(K')}^p\cu_{\Delta(K'')}^p
|K'||K''|(\diam Q)^p \left(\frac{1}{|\gamma Q|} \intl_{\gamma Q}(\nabla^2F)^q\,dz\right)^{\frac{p}{q}}
\right\}.\nn
\ee
%----------------------------------------------------------
\par We obtain
%----------------------------------------------------------
\be
&&S_p(F:S'_Q,S''_Q;\Kc)\nn\\
&\le& C\left
\{\left(\sum_{\substack {K''\in\,\Kc\\c_{K''}\in S''_Q}}
\cu_{\Delta(K'')}^p |K''|\right)\left(\sum_{\substack {K'\in\,\Kc\\c_{K'}\in S'_Q}}|K'|
\left(\frac{1}{|\gamma K'|}\intl_{\gamma K'}(\nabla^2F)^q\,dz\right)^{\frac{p}{q}}\right)\right.
\nn\\
&+&
\left(\sum_{\substack {K'\in\,\Kc\\c_{K'}\in S'_Q}}
\cu_{\Delta(K')}^p |K'|\right)\left(\sum_{\substack {K''\in\,\Kc\\c_{K''}\in S''_Q}}|K''|
\left(\frac{1}{|\gamma K''|}\intl_{\gamma K''}(\nabla^2F)^q\,dz\right)^{\frac{p}{q}}\right)\nn\\
&+&
\left.\left(\sum_{\substack {K'\in\,\Kc\\c_{K'}\in S'_Q}}
\cu_{\Delta(K')}^p |K'|\right)
\left(\sum_{\substack {K''\in\,\Kc\\c_{K''}\in S''_Q}}
\cu_{\Delta(K'')}^p |K''|\right)
(\diam Q)^p \left(\frac{1}{|\gamma Q|} \intl_{\gamma Q}(\nabla^2F)^q\,dz\right)^{\frac{p}{q}}\right\}
\nn
\ee
%----------------------------------------------------------
\par Let $S\in\Sc$ and let
%----------------------------------------------------------
$$
V\{S\}:=\bigcup\{K\in\Kc: c_K\in S\}.
$$
%----------------------------------------------------------
For the sake of brevity we put
%----------------------------------------------------------
$$
H(x)=\Mc[(\nabla^2F)^q]^{\frac{p}{q}}(x),~~~~x\in\RT.
$$
%----------------------------------------------------------
By \rf{IQ-P},
%----------------------------------------------------------
$$
|K'|\left(\frac{1}{|\gamma K'|} \intl_{\gamma K'}(\nabla^2F)^q\,dz\right)^{\frac{p}{q}}\le\intl_{K'}
H(z)\,dz,
$$
%----------------------------------------------------------
so that
%----------------------------------------------------------
$$
\sum_{\substack {K'\in\,\Kc\\c_{K'}\in S'_Q}}|K'|
\left(\frac{1}{|\gamma K'|}\intl_{\gamma K'}
(\nabla^2F)^q\,dz\right)^{\frac{p}{q}}
\le
\sum_{\substack {K'\in\,\Kc\\c_{K'}\in S'_Q}}\,\,
\intl_{K'}
H(z)\,dz\le \intl_{V\{S'_Q\}}
H(z)\,dz.
$$
%----------------------------------------------------------
In the same way we prove that
%----------------------------------------------------------
$$
\sum_{\substack {K''\in\,\Kc\\c_{K''}\in S''_Q}}|K''|
\left(\frac{1}{|\gamma K''|}\intl_{\gamma K''}
(\nabla^2F)^q\,dz\right)^{\frac{p}{q}}
\le \intl_{V\{S''_Q\}}
H(z)\,dz.
$$
%----------------------------------------------------------
Also, by \rf{IQ-P},
%----------------------------------------------------------
$$
|Q|\left(\frac{1}{|\gamma Q|} \intl_{\gamma Q}(\nabla^2F)^q\,dz\right)^{\frac{p}{q}}\le\intl_Q
H(z)\,dz.
$$
%----------------------------------------------------------
\par Combining these inequalities with the above estimate of $S_p(F:S'_Q,S''_Q;\Kc)$ and definition \rf{SG-S}, we obtain
%----------------------------------------------------------
\be
S_p(F:S'_Q,S''_Q;\Kc)&\le&
C\left\{\sigma_p(S''_Q;\Kc)\intl_{V\{S'_Q\}}
H(z)\,dz+\sigma_p(S'_Q;\Kc)
\intl_{V\{S''_Q\}} H(z)\,dz\right.\nn\\
&+&
\left.\sigma_p(S'_Q;\Kc)\sigma_p(S''_Q;\Kc)
(\diam Q)^{p-2}\intl_Q
H(z)\,dz.\right\}
\nn
\ee
%----------------------------------------------------------
\par Let
%----------------------------------------------------------
$$
A:=\smed\limits_{Q\in\Qc}\,
\frac{(\diam Q)^{2-p}\,\,S_p(F:S'_Q,S''_Q;\Kc)}
{\left\{(\diam S'_Q)^{2-p} +\sigma_p(S'_Q;\Kc)\right\}
\left\{(\diam S''_Q)^{2-p} +\sigma_p(S''_Q;\Kc)\right\}}\,.
$$
%----------------------------------------------------------
Then the above estimate of $S_p(F:S'_Q,S''_Q;\Kc)$ implies the following inequality:
%----------------------------------------------------------
\bel{A-123}
A\le C(I_1+I_2+I_3).
\ee
%----------------------------------------------------------
Here
%----------------------------------------------------------
$$
I_1:=\smed\limits_{Q\in\Qc}\,
\frac{(\diam Q)^{2-p} \,\sigma_p(S''_Q;\Kc)}
{\left\{(\diam S'_Q)^{2-p} +\sigma_p(S'_Q;\Kc)\right\}
\left\{(\diam S''_Q)^{2-p} +\sigma_p(S''_Q;\Kc)\right\}}
\intl_{V\{S'_Q\}}H(z)\,dz\,,
$$
%----------------------------------------------------------
$$
I_2:=\smed\limits_{Q\in\Qc}\,
\frac{(\diam Q)^{2-p} \,\sigma_p(S'_Q;\Kc)}
{\left\{(\diam S'_Q)^{2-p} +\sigma_p(S'_Q;\Kc)\right\}
\left\{(\diam S''_Q)^{2-p} +\sigma_p(S''_Q;\Kc)\right\}}
\intl_{V\{S''_Q\}}H(z)\,dz\,,
$$
%----------------------------------------------------------
and
%----------------------------------------------------------
$$
I_3:=\smed\limits_{Q\in\Qc}\,
\frac{\sigma_p(S'_Q;\Kc)\sigma_p(S''_Q;\Kc)}
{\left\{(\diam S'_Q)^{2-p} +\sigma_p(S'_Q;\Kc)\right\}
\left\{(\diam S''_Q)^{2-p} +\sigma_p(S''_Q;\Kc)\right\}}
\intl_Q H(z)\,dz.
$$
%----------------------------------------------------------
\par Clearly,
%----------------------------------------------------------
$$
I_3\le\smed\limits_{Q\in\Qc}\, \intl_Q H(z)\,dz.
$$
%----------------------------------------------------------
Since the squares of the family $\Qc$ are pairwise disjoint, we obtain
%----------------------------------------------------------
\bel{I3-F}
I_3\le \intl_{U\{\Qc\}} H(z)\,dz.
\ee
%----------------------------------------------------------
\par Prove that
%----------------------------------------------------------
$$
I_1\le \intl_{U\{\Kc\}} H(z)\,dz.
$$
%----------------------------------------------------------
\par Obviously,
%----------------------------------------------------------
$$
I_1\le\sbig\limits_{Q\in\Qc}\,
\left(\frac{\diam S'_Q}{\diam Q}\right)^{p-2}
\intl_{V\{S'_Q\}}H(z)\,dz
$$
%----------------------------------------------------------
so that
%----------------------------------------------------------
\bel{I1-J}
I_1\le\smed\limits_{S\in\tSc}\,
J(S)(\diam S)^{p-2}
\intl_{V\{S\}}H(z)\,dz.
\ee
%----------------------------------------------------------
Here
%----------------------------------------------------------
$$
\tSc:=\{S\in\Sc:\diam S>0\}
$$
%----------------------------------------------------------
and
%----------------------------------------------------------
$$
J(S):=\sum\,\left\{(\diam Q)^{2-p}:Q\in\Qc,\, S'_Q=S\right\}\,.
$$
%----------------------------------------------------------
\par Prove that
%----------------------------------------------------------
\bel{JS}
J(S)\le C(\diam S)^{2-p}
\ee
%----------------------------------------------------------
for every $S\in\tSc$. Let
%----------------------------------------------------------
$$
G(S):=\{Q\in\Qc: S'_Q=S\}.
$$
%----------------------------------------------------------
Since $S\subset \gamma Q$ whenever $Q\in G(S)$, we have
%----------------------------------------------------------
\bel{DI-G}
\diam S\le \gamma\diam Q~~~\text{for every}~~~Q\in G(S).
\ee
%----------------------------------------------------------
\par Fix a point $x_0\in S$ and introduce two family of squares:
%----------------------------------------------------------
$$
G_1(S):=\{Q\in\Qc: S'_Q=S,  x_0\in Q\}\,,
$$
%----------------------------------------------------------
and
%----------------------------------------------------------
$$
G_2(S):=\{Q\in\Qc: S'_Q=S,   x_0\notin Q\}.
$$
%----------------------------------------------------------
Clearly, $G_1(S)$ and $G_2(S)$ is a partition of $G(S)$.
%----------------------------------------------------------
\par Since the \sqs of $G(S)$ are pairwise disjoint, the family $G_1(S)$ consists of at most one element. Hence, by \rf{DI-G},
%----------------------------------------------------------
\bel{G1-E}
\sum\,\{(\diam Q)^{2-p}:Q\in G_1(S)\}\le \gamma^{p-2}\,(\diam S)^{2-p}.
\ee
%----------------------------------------------------------
\par Prove that
%----------------------------------------------------------
$$
\sum\,\{(\diam Q)^{2-p}:Q\in G_2(S)\}\le C(\gamma,p)(\diam S)^{2-p}.
$$
%----------------------------------------------------------
\par In fact, let $Q\in G_2(S)$ and let $\tQ=\frac12 Q$.
Since $S\subset \gamma Q$ and $\gamma\ge 1$,  we have  $y,x_0\in \gamma Q$ for every $y\in\tQ$. Hence,  %----------------------------------------------------------
$$
\|x_0-y\|\le \diam (\gamma Q)=\gamma\diam Q
$$
%----------------------------------------------------------
so that
%----------------------------------------------------------
$$
(\diam Q)^{-p}\le \gamma^p\|x_0-y\|^{-p},~~~y\in\tQ.
$$
%----------------------------------------------------------
Integrating this inequality over $\tQ$, we obtain
%----------------------------------------------------------
$$
(\diam Q)^{2-p}\le 4\gamma^p\intl_{\tQ}\|x_0-y\|^{-p}\,dy
$$
%----------------------------------------------------------
so that
%----------------------------------------------------------
$$
\sum\,\{(\diam Q)^{2-p}:Q\in G_2(S)\}\le C(\gamma,p)\intl_{\tQ}\|x_0-y\|^{-p}\,dy.
$$
%----------------------------------------------------------
\par Let $U:=\cup\{\tQ: Q\in G_2(S)\}$. Since the \sqs of the family $\{\tQ: Q\in G_2(S)\}$ are pairwise disjoint, we have
%----------------------------------------------------------
\bel{INT-L}
\sum\,\{(\diam Q)^{2-p}:Q\in G_2(S)\}\le C(\gamma,p)\intl_{U}\|x_0-y\|^{-p}\,dy.
\ee
%----------------------------------------------------------
\par Let $Q\in G_2(S)$ and let $y\in \tQ=\frac12 Q$. Since $x_0\notin Q$, we have
%----------------------------------------------------------
$$
\|x_0-y\|\ge \frac14\diam Q\ge \frac{1}{4\gamma}\diam S,
$$
%----------------------------------------------------------
see inequality \rf{DI-G}. Hence
%----------------------------------------------------------
\bel{IR}
U\subset \{y\in\RT:\|x_0-y\|\ge \tfrac{1}{4\gamma}\diam S\}.
\ee
%----------------------------------------------------------
\par Let $R:=\frac{1}{4\gamma}\diam S$. Then, by \rf{INT-L} and \rf{IR},
%----------------------------------------------------------
$$
\sum\,\{(\diam Q)^{2-p}:Q\in G_2(S)\}\le C(\gamma,p)\intl_{\|x_0-y\|\ge R}\|x_0-y\|^{-p}\,dy
\le C(\gamma,p)R^{2-p}.
$$
%----------------------------------------------------------
\par We obtain
%----------------------------------------------------------
$$
\sum\,\{(\diam Q)^{2-p}:Q\in G_2(S)\}\le C(\gamma,p)(\diam S)^{2-p}.
$$
%----------------------------------------------------------
\par This estimate and inequality \rf{G1-E} imply the required inequality \rf{JS}.
%----------------------------------------------------------
\par In turn \rf{I1-J} and \rf{JS} imply the following inequality:
%----------------------------------------------------------
$$
I_1\le C(\gamma,p)\,\sum\limits_{S\in\tSc}\,\,\,
\intl_{V\{S\}}H(z)\,dz.
$$
%----------------------------------------------------------
Recall that
%----------------------------------------------------------
$$
V\{S\}:=\bigcup\{K\in\Kc: c_K\in S\}.
$$
%----------------------------------------------------------
Since the families $\Sc$ and $\Kc$ consist of pairwise disjoint sets, the sets of the family
$\{V\{S\}: S\in\Sc\}$ are pairwise disjoint as well. Hence %----------------------------------------------------------
$$
I_1\le C(\gamma,p)\,
\intl_{U\{\Kc\}}H(z)\,dz.
$$
%----------------------------------------------------------
\par In the same way we prove that
%----------------------------------------------------------
$$
I_2\le C(\gamma,p)\,
\intl_{U\{\Kc\}}H(z)\,dz.
$$
%----------------------------------------------------------
\par Combining these inequalities with \rf{I3-F} and \rf{A-123} we obtain the statement of the proposition.\bx\medskip
%----------------------------------------------------------
%@@@@@@@@@@@@@@@@@@@@@@@@@@@@@@@@@@@@@@@@@@@@@@@@@@@@@@@@@@
%@@@@@@@@@@@@@@@@@@@@@@@@@@@@@@@@@@@@@@@@@@@@@@@@@@@@@@@@@@
%@@@@@@@@@@@@@@@@@@@@@@@@@@@@@@@@@@@@@@@@@@@@@@@@@@@@@@@@@@
%----------------------------------------------------------
\par Let us finish the proof of the necessity of part (ii) of Theorem \reff{MAIN}.
%----------------------------------------------------------
\par Let $F\in\LTP$ and let $F|_E=f$. Let us estimate the quantity
%----------------------------------------------------------
$$
A(f)=\sbig\limits_{Q\in\Qc}\,
\frac{(\diam Q)^{2-p}\,\,S_p(f:S'_Q,S''_Q;\Kc)}
{\left\{(\diam S'_Q)^{2-p} +\sigma_p(S'_Q;\Kc)\right\}
\left\{(\diam S''_Q)^{2-p} +\sigma_p(S''_Q;\Kc)\right\}}.
$$
%----------------------------------------------------------
\par Since $F|_E=f$, we have $S_p(f:Q',Q'';\Kc)=S_p(F:Q',Q'';\Kc)$ so that
$A(f)=A(F)$. Let $q:=(p+2)/2$. Then, by Proposition \reff{MF-N},
%----------------------------------------------------------
$$
A(f)\le
C\, \left\{\intl_{U\{\Kc\}}
\Mc[(\nabla^2F)^q]^{\frac{p}{q}}\,dz+
\intl_{U\{\Qc\}}
\Mc[(\nabla^2F)^q]^{\frac{p}{q}}\,dz\right\}\le C\,\intl_{\RT}
\Mc[(\nabla^2F)^q]^{\frac{p}{q}}\,dz
$$
%----------------------------------------------------------
where $C$ is a constant depending only on $p$ and $\gamma$.
%----------------------------------------------------------
\par Since $p/q>1$, by the \HLM theorem,
%----------------------------------------------------------
$$
A(f)\le C\,\intl_{\RT}\,[(\nabla^2F)^q]^{\frac{p}{q}}\,dz=
C\,\intl_{\RT}\,(\nabla^2F)^p\,dz=C\|F\|_{\LTP}^p
$$
%----------------------------------------------------------
proving again that for every $F\in\LTP\cap C^1(\RT)$ such that  $F|_E=f$ the condition (ii) of Theorem \reff{MAIN} holds with $\lambda=C(p,\gamma)\|F\|_{\LTP}^p$.\medskip
%----------------------------------------------------------
\par The necessity part of Theorem \reff{MAIN} is completely proved. \bx
%----------------------------------------------------------
%@@@@@@@@@@@@@@@@@@@@@@@@@@@@@@@@@@@@@@@@@@@@@@@@@@@@@@@@@@
%@@@@@@@@@@@@@@@@@@@@@@@@@@@@@@@@@@@@@@@@@@@@@@@@@@@@@@@@@@
%@@@@@@@@@@@@@@@@@@@@@@@@@@@@@@@@@@@@@@@@@@@@@@@@@@@@@@@@@@
%----------------------------------------------------------
%@@@@@@@@@@@@@@@@@@@@@@@@@@@@@@@@@@@@@@@@@@@@@@@@@@@@@@@@@@
%@@@@@@@@@@@@@@@@@@@@@@@@@@@@@@@@@@@@@@@@@@@@@@@@@@@@@@@@@@
%@@@@@@@@@@@@@@@@@@@@@@@@@@@@@@@@@@@@@@@@@@@@@@@@@@@@@@@@@@
%@@@@@@@@@@@@@@@@@@@@@@@@@@@@@@@@@@@@@@@@@@@@@@@@@@@@@@@@@@
%----------------------------------------------------------
%@@@@@@@@@@@@@@@@@@@@@@@@@@@@@@@@@@@@@@@@@@@@@@@@@@@@@@@@@@
%@@@@@@@@@@@@@@@@@@@@@@@@@@@@@@@@@@@@@@@@@@@@@@@@@@@@@@@@@@
%@@@@@@@@@@@@@@@@@@@@@@@@@@@@@@@@@@@@@@@@@@@@@@@@@@@@@@@@@@
%@@@@@@@@@@@@@@@@@@@@@@@@@      @@@@@@@@@@@@@@@@@@@@@@@@@@@
%@@@@@@@@@@@@@@@@@@@@@@@          @@@@@@@@@@@@@@@@@@@@@@@@@
%@@@@@@@@@@@@@@@@@@@@@              @@@@@@@@@@@@@@@@@@@@@@@
%@@@@@@@@@@@@@@@@@@@     SECTION 4    @@@@@@@@@@@@@@@@@@@@@
%@@@@@@@@@@@@@@@@@@@@@              @@@@@@@@@@@@@@@@@@@@@@@
%@@@@@@@@@@@@@@@@@@@@@@@          @@@@@@@@@@@@@@@@@@@@@@@@@
%@@@@@@@@@@@@@@@@@@@@@@@@@      @@@@@@@@@@@@@@@@@@@@@@@@@@@
%@@@@@@@@@@@@@@@@@@@@@@@@@@@@@@@@@@@@@@@@@@@@@@@@@@@@@@@@@@
%@@@@@@@@@@@@@@@@@@@@@@@@@@@@@@@@@@@@@@@@@@@@@@@@@@@@@@@@@@
%@@@@@@@@@@@@@@@@@@@@@@@@@@@@@@@@@@@@@@@@@@@@@@@@@@@@@@@@@@
%----------------------------------------------------------
\SECT{4. Lacunae of Whitney cubes.}{4}
%----------------------------------------------------------
\addtocontents{toc}{4. Lacunae of Whitney cubes. \hfill \thepage\par}
%----------------------------------------------------------
\indent
%@@@@@@@@@@@@@@@@@@@@@@@@@@@@@@@@@@@@@@@@@@@@@@@@@@@@@@@@@@
%----------------------------------------------------------
\par We turn to the proof of the sufficiency part of Main Theorem \reff{MAIN}.
%----------------------------------------------------------
\par We prove the sufficiency in several steps. At the first step we present a modification of the Whitney extension method which is based on the notion of a {\it ``lacuna of Whitney \sqs''}. We have briefly described this object in Section 2. This notion enables us to identify and characterize all possible ``holes'' in the set $E$. In this section we give main definitions and describe several main properties of lacunae.
%----------------------------------------------------------
\par As we have mentioned in Section 2, {\it in this and  the next two sections $E$ is an arbitrary closed subset of $\RN$.} We equip $\RN$ with the uniform norm $\|x\|:=\max\{|x_i|:i=1,...,n\}$.
%----------------------------------------------------------
\par  As usual the word ``cube" will mean a closed cube in $\RN$ whose sides are parallel to the coordinate axes. We will use the same notation for cubes and distances between sets as in two dimensional case.
%---------------------------------------------------------- %@@@@@@@@@@@@@@@@@@@@@@@@@@@@@@@@@@@@@@@@@@@@@@@@@@@@@@@@@@
%@@@@@@@@@@@@@@@@@@@@@@@@@@@@@@@@@@@@@@@@@@@@@@@@@@@@@@@@@@
%@@@@@@@@@@@@@@@@@@@@@@@@@@@@@@@@@@@@@@@@@@@@@@@@@@@@@@@@@@
%@@@@@@@@@@@@@@@@@@@@@@@@@@@@@@@@@@@@@@@@@@@@@@@@@@@@@@@@@@
%@@@@@@@@@@@@@@@@@@@@@@@@@@@@@@@@@@@@@@@@@@@@@@@@@@@@@@@@@@
%@@@@@@@@@@@@@@@@@@@@@@@@@@@@@@@@@@@@@@@@@@@@@@@@@@@@@@@@@@
%@@@@@@@@@@@@@@@@@@@@@@@@@@@@@@@@@@@@@@@@@@@@@@@@@@@@@@@@@@
%@@@@@@@@@@@@@@@@@@@@@@@@@@@@@@@@@@@@@@@@@@@@@@@@@@@@@@@@@@
%@@@@@@@@@@@@@@@@@@@@@@@@@@@@@@@@@@@@@@@@@@@@@@@@@@@@@@@@@@
\medskip
\par {\bf 4.1. Whitney cubes and lacunae of Whitney cubes.}
%----------------------------------------------------------
\addtocontents{toc}{~~~~4.1. Whitney cubes and lacunae of Whitney cubes: main definitions. \hfill \thepage\par}
%----------------------------------------------------------
First let us recall the notion of a Whitney cube. Since $E$ is a closed set, the set $\RN\setminus E$ is open so that it admits a Whitney decomposition $W_E$ into a family of non-overlapping cubes. In the next theorem we present the main properties of this decomposition, see, e.g. \cite{St} or \cite{G}.
%----------------------------------------------------------
%@@@@@@@@@@@@@@@@@@@@@@@@@@@@@@@@@@@@@@@@@@@@@@@@@@@@@@@@@@
%@@@@@@@@@@@@@@@@@@@@@@@@@@@@@@@@@@@@@@@@@@@@@@@@@@@@@@@@@@
%@@@@@@@@@@@@@@@@@@@@@@@@@@@@@@@@@@@@@@@@@@@@@@@@@@@@@@@@@@
%@@@@@@@@@@@@@@@@@@@@@@@@@@@@@@@@@@@@@@@@@@@@@@@@@@@@@@@@@@
%----------------------------------------------------------
\begin{theorem}\lbl{Wcov} $W_E=\{Q_k\}$ is a countable family of non-overlapping cubes  such that
%----------------------------------------------------------
\par (i). $\RN\setminus E=\cup\{Q:Q\in W_E\}$;
%----------------------------------------------------------
\par (ii). For every cube $Q\in W_E$ we have
%----------------------------------------------------------
\bel{DQ-E}
\diam Q\le \dist(Q,E)\le 4\diam Q.
\ee
%----------------------------------------------------------
%@@@@@@@@@@@@@@@@@@@@@@@@@@@@@@@@@@@@@@@@@@@@@@@@@@@@@@@@@@
%----------------------------------------------------------
\end{theorem}
%@@@@@@@@@@@@@@@@@@@@@@@@@@@@@@@@@@@@@@@@@@@@@@@@@@@@@@@@@@
%@@@@@@@@@@@@@@@@@@@@@@@@@@@@@@@@@@@@@@@@@@@@@@@@@@@@@@@@@@
%@@@@@@@@@@@@@@@@@@@@@@@@@@@@@@@@@@@@@@@@@@@@@@@@@@@@@@@@@@
%@@@@@@@@@@@@@@@@@@@@@@@@@@@@@@@@@@@@@@@@@@@@@@@@@@@@@@@@@@
\par We are also needed certain additional properties of
Whitney cubes which we present in the next lemma. These
properties easily follow from constructions of Whitney decomposition presented in \cite{St} and \cite{G}.
%----------------------------------------------------------
%@@@@@@@@@@@@@@@@@@@@@@@@@@@@@@@@@@@@@@@@@@@@@@@@@@@@@@@@@@
\par Given a cube $Q\subset\RN$ let $Q^*:=\frac{9}{8}Q$.
%@@@@@@@@@@@@@@@@@@@@@@@@@@@@@@@@@@@@@@@@@@@@@@@@@@@@@@@@@@
\begin{lemma}\lbl{Wadd}
%----------------------------------------------------------
%@@@@@@@@@@@@@@@@@@@@@@@@@@@@@@@@@@@@@@@@@@@@@@@@@@@@@@@@@@
%----------------------------------------------------------
(1). If $Q,K\in W_E$ and $Q^*\cap K^*\ne\emptyset$, then
%----------------------------------------------------------
$$
\frac{1}{4}\diam Q\le \diam K\le 4\diam Q.
$$
%----------------------------------------------------------
%@@@@@@@@@@@@@@@@@@@@@@@@@@@@@@@@@@@@@@@@@@@@@@@@@@@@@@@@@@
\par (2). For every cube $K\in W_E$ there are at most
$N=N(n)$ cubes from the family
%----------------------------------------------------------
$$W_E^*:=\{Q^*:Q\in W_E\}$$
%----------------------------------------------------------
which intersect $K^*$.
%----------------------------------------------------------
\par (3). If $Q,K\in W_E$, then $Q^*\cap K^*\ne\emptyset$
if and only if  $Q\cap K\ne\emptyset$.
%----------------------------------------------------------
%@@@@@@@@@@@@@@@@@@@@@@@@@@@@@@@@@@@@@@@@@@@@@@@@@@@@@@@@@@
\end{lemma}
%@@@@@@@@@@@@@@@@@@@@@@@@@@@@@@@@@@@@@@@@@@@@@@@@@@@@@@@@@@
\par Note that inequality \rf{DQ-E} implies the following property of Whitney cubes:
%----------------------------------------------------------
\bel{INT-W}
(9Q)\cap E\ne\emp ~~~~\text{for every}~~~Q\in W_E.
\ee
%----------------------------------------------------------
\par By $LW_E$ we denote a subfamily of Whitney cubes satisfying the following condition:
%----------------------------------------------------------
\bel{L-PR}
(10Q)\cap E=(\q Q)\cap E.
\ee
%----------------------------------------------------------
\par Then we introduce a binary relation $\sim$ on $LW_E$: for every $Q_1,Q_2\in LW_E$
%----------------------------------------------------------
$$
Q_1\sim Q_2 ~~\Longleftrightarrow~ (10Q_1)\cap E= (10Q_2)\cap E.
$$
%----------------------------------------------------------
\par It can be easily seen that $\sim$ satisfies the axioms of equivalence relations, i.e., it is reflexive, symmetric and transitive. Given a cube $Q\in LW_E$ by
%----------------------------------------------------------
$$
[Q]:=\{K\in LW_E: K\sim Q\}
$$
%----------------------------------------------------------
we denote the equivalence class of $Q$. We refer to this equivalence class as {\it a true lacuna} with respect to the set $E$.
%----------------------------------------------------------
\par Let
%----------------------------------------------------------
$$
\tL_E=LW_E\backslash\sim\,=\{[Q]: Q\in LW_E\}
$$
%----------------------------------------------------------
be the corresponding quotient set of $LW_E$ by $\sim$\,, i.e., the set of all possible equivalence classes (lacunae) of $LW_E$ by $\sim$\,.
%----------------------------------------------------------
\par Thus for every pair of Whitney cubes $Q_1,Q_2\in W_E$ which belong to a true lacuna $L\in\tL_E$ we have
%----------------------------------------------------------
\bel{I-L}
(10Q_1)\cap E=(\q Q_1)\cap E=(10Q_2)\cap E=(\q Q_2)\cap E.
\ee
%----------------------------------------------------------
By $V_L$ we denote the associated set of the lacuna $L$
%----------------------------------------------------------
\bel{D-VL}
V_L:=(\q Q)\cap E.
\ee
%----------------------------------------------------------
Here $Q$ is an arbitrary cube from $L$. By \rf{I-L}, any choice of a cube $Q\in L$ provides the same set $V_L$ so that $V_L$ is well-defined. Also note that for each cube $Q$ which belong to a true lacuna $L$ we have $V_L=(10Q)\cap E.$\medskip
%----------------------------------------------------------
\par We extend the family $\tL_E$ of true lacunae to a family $\LE$ of {\it all lacunae} in the following way. Suppose that $Q\in W_E\setminus LW_E$, see \rf{L-PR}, i.e.,
%----------------------------------------------------------
\bel{A-L}
(10Q)\cap E\ne(\q Q)\cap E.
\ee
%----------------------------------------------------------
In this case to the cube $Q$ we assign a lacuna $L:=\{Q\}$ consisting of a unique cube - the cube $Q$.  We also put $V_L:=(\q Q)\cap E$ as in \rf{D-VL}.
%----------------------------------------------------------
\par We refer to such a lacuna $L:=\{Q\}$ as an {\it elementary lacuna} with respect to the set $E$. By $\hL_E$ we denote the family of all elementary lacunae with respect to $E$:
%----------------------------------------------------------
$$
\hL_E:=\{L=\{Q\}:Q\in W_E\setminus LW_E\}
$$
%----------------------------------------------------------
\par We note that property \rf{A-L} implies the existence of a point
%----------------------------------------------------------
$$
a\in (E\setminus (10Q))\cap (\q Q).
$$
%----------------------------------------------------------
On the other hand, by \rf{INT-W}, there exists a point $$b\in(9Q)\cap E.$$ Hence
%----------------------------------------------------------
$$
\|a-b\|\ge r_Q=(1/2)\diam Q
$$
%----------------------------------------------------------
so that
%----------------------------------------------------------
\bel{A-DL}
\diam V_L=\diam ((\q Q)\cap E)\ge \tfrac12\diam Q
\ee
%----------------------------------------------------------
provided
%----------------------------------------------------------
$$
L=\{Q\}\in \hL_E
$$
%----------------------------------------------------------
is an elementary lacuna.
%----------------------------------------------------------
\par Finally, by $\LE$ we denote the family of all lacunae with respect to $E$:
%----------------------------------------------------------
$$
\LE=\tL_E\cup \hL_E.
$$
%----------------------------------------------------------
\smallskip
%@@@@@@@@@@@@@@@@@@@@@@@@@@@@@@@@@@@@@@@@@@@@@@@@@@@@@@@@@@
%@@@@@@@@@@@@@@@@@@@@@@@@@@@@@@@@@@@@@@@@@@@@@@@@@@@@@@@@@@
%@@@@@@@@@@@@@@@@@@@@@@@@@@@@@@@@@@@@@@@@@@@@@@@@@@@@@@@@@@
%@@@@@@@@@@@@@@@@@@@@@@@@@@@@@@@@@@@@@@@@@@@@@@@@@@@@@@@@@@
%@@@@@@@@@@@@@@@@@@@@@@@@@@@@@@@@@@@@@@@@@@@@@@@@@@@@@@@@@@
%@@@@@@@@@@@@@@@@@@@@@@@@@@@@@@@@@@@@@@@@@@@@@@@@@@@@@@@@@@
%@@@@@@@@@@@@@@@@@@@@@@@@@@@@@@@@@@@@@@@@@@@@@@@@@@@@@@@@@@
%----------------------------------------------------------
\par {\bf 4.2. Main properties of the lacunae.}
%----------------------------------------------------------
\addtocontents{toc}{~~~~4.2. Main properties of the lacunae. \hfill \thepage\\\par}
%----------------------------------------------------------
%@@@@@@@@@@@@@@@@@@@@@@@@@@@@@@@@@@@@@@@@@@@@@@@@@@@@@@@@@@
%@@@@@@@@@@@@@@@@@@@@@@@@@@@@@@@@@@@@@@@@@@@@@@@@@@@@@@@@@@
%@@@@@@@@@@@@@@@@@@@@@@@@@@@@@@@@@@@@@@@@@@@@@@@@@@@@@@@@@@
%@@@@@@@@@@@@@@@@@@@@@@@@@@@@@@@@@@@@@@@@@@@@@@@@@@@@@@@@@@
%@@@@@@@@@@@@@@@@@@@@@@@@@@@@@@@@@@@@@@@@@@@@@@@@@@@@@@@@@@
%----------------------------------------------------------
\begin{lemma}\lbl{CL-P1} Let $L\in\LE$ be a lacuna and let $Q\in L$. Then
%----------------------------------------------------------
$$
\dist(Q,V_L)=\dist(Q,E).
$$
%----------------------------------------------------------
\end{lemma}
%----------------------------------------------------------
%@@@@@@@@@@@@@@@@@@@@@@@@@@@@@@@@@@@@@@@@@@@@@@@@@@@@@@@@@@
%@@@@@@@@@@@@@@@@@@@@@@@@@@@@@@@@@@@@@@@@@@@@@@@@@@@@@@@@@@
%@@@@@@@@@@@@@@@@@@@@@@@@@@@@@@@@@@@@@@@@@@@@@@@@@@@@@@@@@@
%----------------------------------------------------------
\par {\it Proof.} Recall that for each Whitney cube $Q\in W_E$ we have $(9Q)\cap E\ne\emp.$ Hence
%----------------------------------------------------------
$$
\dist(Q,E)=\dist(Q,(9Q)\cap E)=\dist(Q,(90Q)\cap E).
$$
%----------------------------------------------------------
But $V_L=(90Q)\cap E$ for each $Q\in L$, see \rf{D-VL}, and the proof is finished.\bx\medskip
%@@@@@@@@@@@@@@@@@@@@@@@@@@@@@@@@@@@@@@@@@@@@@@@@@@@@@@@@@@
%@@@@@@@@@@@@@@@@@@@@@@@@@@@@@@@@@@@@@@@@@@@@@@@@@@@@@@@@@@
%@@@@@@@@@@@@@@@@@@@@@@@@@@@@@@@@@@@@@@@@@@@@@@@@@@@@@@@@@@
%@@@@@@@@@@@@@@@@@@@@@@@@@@@@@@@@@@@@@@@@@@@@@@@@@@@@@@@@@@
%----------------------------------------------------------
\begin{proposition}\lbl{DL-M} Let $L\in\LE$ be a lacuna.
Then
%----------------------------------------------------------
\bel{VLQ}
\tfrac{1}{\q}\diam V_L\le \inf_{Q\in L} \diam Q\le \gamma \diam V_L
\ee
%----------------------------------------------------------
where $\gamma$ is an absolute constant.
%@@@@@@@@@@@@@@@@@@@@@@@@@@@@@@@@@@@@@@@@@@@@@@@@@@@@@@@@@@
%----------------------------------------------------------
\end{proposition}
%----------------------------------------------------------
%@@@@@@@@@@@@@@@@@@@@@@@@@@@@@@@@@@@@@@@@@@@@@@@@@@@@@@@@@@
%@@@@@@@@@@@@@@@@@@@@@@@@@@@@@@@@@@@@@@@@@@@@@@@@@@@@@@@@@@
%@@@@@@@@@@@@@@@@@@@@@@@@@@@@@@@@@@@@@@@@@@@@@@@@@@@@@@@@@@
%@@@@@@@@@@@@@@@@@@@@@@@@@@@@@@@@@@@@@@@@@@@@@@@@@@@@@@@@@@
\par {\it Proof.} Suppose that $\diam V_L=0$ so that $V_L=\{A_L\}$ where $A_L$ is a point in $E$. Then for each $\ve>0$ the $\ve$-neighborhood of $A_L$  contains a cube from $L$. This proves \rf{VLQ} in the case under consideration.
%----------------------------------------------------------
\par Let us assume that $\diam V_L>0$. If $L\in \hL$ is an elementary lacuna, then $L=\{Q\}$ where $Q\in W_E\setminus LW_E$. In this case the statement of the proposition follows from \rf{D-VL} and inequality \rf{A-DL} with $\gamma=2$.
%@@@@@@@@@@@@@@@@@@@@@@@@@@@@@@@@@@@@@@@@@@@@@@@@@@@@@@@@@@
\par Let $L\in\tL_E$ be a true lacuna. We recall that in this case every cube $Q\in L$ satisfies the condition
%----------------------------------------------------------
\bel{IN-L}
V_L=(10Q)\cap E=(\q Q)\cap E.
\ee
%----------------------------------------------------------
Hence $\diam V_L\le 10\diam Q$ for every $Q\in L$ proving the first inequality in \rf{VLQ}.
%----------------------------------------------------------
\par Prove the second inequality in \rf{VLQ} with $\gamma:=600$. Suppose that it is not true, i.e.,
%----------------------------------------------------------
\bel{A-G}
\diam Q>\gamma\diam V_L~~~\text{for every}~~~Q\in L.
\ee
%----------------------------------------------------------
\par Let $m=30$. Fix a cube $Q\in L$ and prove that there exists a cube $K\in L$ such that $\diam K\le \frac1m\diam Q$.
%----------------------------------------------------------
\par  Note that, by Lemma \reff{CL-P1},
$\dist(Q,V_L)=\dist(Q,E)$
and, by Theorem \reff{Wcov}, $\diam Q\le \dist(Q,E).$
%----------------------------------------------------------
Hence
%----------------------------------------------------------
$$
\diam Q\le\dist(Q,V_L)\le\dist(x,V_L)~~~
\text{for every}~~~ x\in Q.
$$
%----------------------------------------------------------
Thus $\dist(\cdot,V_L)$ is a continuous function on $\RN$ which takes the value $0$ on $V_L$. Since the values of this function on $Q$ are at least $\diam Q$, there exists a point $x_0\in\RN\setminus E$ such that
%----------------------------------------------------------
\bel{D-X}
\dist(x_0,V_L)=\tfrac1m \diam Q.
\ee
%----------------------------------------------------------
\par Prove that
%----------------------------------------------------------
\bel{X-ME}
\dist(x_0,V_L)=\dist(x_0,E).
\ee
%----------------------------------------------------------
In fact, by \rf{IN-L} and \rf{D-X},
%----------------------------------------------------------
$$
x_0\in (10+2/m)Q.
$$
%----------------------------------------------------------
This and \rf{IN-L} imply the following inequality
%----------------------------------------------------------
$$
\dist(x_0,E\setminus V_L)=\dist(x_0,E\setminus(\q Q))
\ge (\q -(10+2/m))r_Q=(80-2/m)r_Q.
$$
%----------------------------------------------------------
Hence
%----------------------------------------------------------
\bel{DX-Z}
\dist(x_0,E\setminus V_L)\ge(40-1/m)\diam Q.
\ee
%----------------------------------------------------------
Since $1/m<40-1/m$, this inequality and equality \rf{D-X} imply that
%----------------------------------------------------------
$$
\dist(x_0,V_L)<\dist(x_0,E\setminus V_L)
$$
%----------------------------------------------------------
proving \rf{X-ME}.
%----------------------------------------------------------
\par Let $K\in W_E$ be a Whitney cube containing $x_0$. Then, by Theorem \reff{Wcov},
%----------------------------------------------------------
\bel{K-WQ}
\diam K\le \dist(K,E)\le 4\diam K.
\ee
%----------------------------------------------------------
Prove that
%----------------------------------------------------------
\bel{K-TE}
\tfrac{1}{5m}\diam Q\le \diam K\le \tfrac1m \diam Q.
\ee
%----------------------------------------------------------
Since $x_0\in K$, we have
%----------------------------------------------------------
$$
\dist(x_0,V_L)=\dist(x_0,E)\le\dist(K,E)+\diam K\le 5\diam K
$$
%----------------------------------------------------------
so that, by \rf{D-X}, $(1/5m)\diam Q\le \diam K$. Also,
%----------------------------------------------------------
$$
\diam K\le \dist(K,E)\le \dist(x_0,E)=\dist(x_0,V_L)
$$
%----------------------------------------------------------
so that, by \rf{D-X},  $\diam K\le (1/m)\diam Q$, and  \rf{K-TE} is proved.
%----------------------------------------------------------
\par Prove that $K\in L$, i.e.,
%----------------------------------------------------------
\bel{KINL}
V_L=(10K)\cap E=(\q K)\cap E.
\ee
%----------------------------------------------------------
Since $\dist(\cdot,V_L)$ is a Lipschitz function, by \rf{D-X}, for every $y\in K$ we have
%----------------------------------------------------------
$$
\dist(y,V_L)\le\|x_0-y\|+\dist(x_0,V_L)\le\diam K+\tfrac1m\diam Q
$$
%----------------------------------------------------------
so that, by \rf{K-TE},
%----------------------------------------------------------
$$
\dist(y,V_L)\le\tfrac2m\diam Q.
$$
%----------------------------------------------------------
On the other hand, by \rf{DX-Z} and \rf{K-TE},
%----------------------------------------------------------
\be
\dist(y,E\setminus V_L)&\ge&
\dist(x_0,E\setminus V_L)-\diam K\nn\\&\ge&(40-\tfrac1m)\diam Q
-\tfrac1m\diam Q=(40-\tfrac2m)\diam Q.
\nn
\ee
%----------------------------------------------------------
Since $\tfrac2m\le 40-\tfrac2m$ for $m\ge 1$, we obtain
%----------------------------------------------------------
$$
\dist(y,V_L)<\dist(y,E\setminus V_L)
$$
%----------------------------------------------------------
proving that
%----------------------------------------------------------
$$
\dist(y,V_L)=\dist(y,E)~~~\text{for every}~~~y\in K.
$$
%----------------------------------------------------------
Hence
%----------------------------------------------------------
\bel{DKE}
\dist(K,E)=\dist(K,V_L).
\ee
%----------------------------------------------------------
\par By \rf{DKE}, there exists a point $b\in V_L$ which is  a point nearest to $K$ on $E$. Then, by \rf{K-WQ},
%----------------------------------------------------------
$$
\dist(b,K)=\dist(K,E)\le 4\diam K.
$$
%----------------------------------------------------------
Let $K=Q(c_K,r_K)$, i.e., $c_K$ is the center of $K$ and $r_K=\tfrac12\diam K$. Then for each $z\in V_L$
%----------------------------------------------------------
$$
\|z-c_K\|\le\|z-b\|+\dist(b,K)+r_K\le
\diam V_L+4\diam K+r_K
$$
%----------------------------------------------------------
so that
%----------------------------------------------------------
\bel{P-2}
\|z-c_K\|\le\diam V_L+9r_K.
\ee
%----------------------------------------------------------
But, by assumption \rf{A-G}, $\diam V_L<(1/\gamma)\diam Q$, so that, by \rf{K-TE},
%----------------------------------------------------------
$$
\|z-c_K\|\le(1/\gamma)\diam Q+9r_K
\le (5m/\gamma)\diam K+9r_K\le (9+(10m/\gamma))r_K
$$
%---------------------------------------------------------- %@@@@@@@@@@@@@@@@@@@@@@@@@@@@@@@@@@@@@@@@@@@@@@@@@@@@@@@@@@
proving that
%----------------------------------------------------------
$$
V_L\subset (9+(10m/\gamma))K.
$$
%----------------------------------------------------------
Hence
%----------------------------------------------------------
\bel{10KL}
V_L\subset 10K
\ee
%----------------------------------------------------------
provided $10m<\gamma$. Since $V_L\subset E$, this shows that $V_L\subset (10K)\cap E$.
%----------------------------------------------------------
\par Prove that
%----------------------------------------------------------
\bel{10K-QL}
10K\subset 10Q.
\ee
%----------------------------------------------------------
In fact, since $Q\in W_E$ we have $(9Q)\cap E\ne\emp$. Combining this with \rf{IN-L} we obtain
%----------------------------------------------------------
$$
(9Q)\cap E=(9Q)\cap V_L\ne\emp
$$
%----------------------------------------------------------
so that, by \rf{10KL},
%----------------------------------------------------------
$$
(10K)\cap(9Q)\ne\emp.
$$
%----------------------------------------------------------
But, by \rf{K-TE}, $\diam K\le\tfrac1m\diam Q$, so that
%----------------------------------------------------------
$$
10K\subset(9+20/m)Q.
$$
%----------------------------------------------------------
Since $m\ge 20$, we obtain the required imbedding \rf{10K-QL}.
%----------------------------------------------------------
\par By dilation, from \rf{10K-QL} we have
$\q K\subset \q Q.$
%----------------------------------------------------------
\par Finally, by \rf{IN-L},
%----------------------------------------------------------
$$
V_L\subset(10K)\cap E\subset (10Q)\cap E=V_L
$$
%----------------------------------------------------------
and
%----------------------------------------------------------
$$
V_L\subset(\q K)\cap E\subset (\q Q)\cap E=V_L
$$
%----------------------------------------------------------
proving \rf{KINL}. Hence $K\in L$.
%----------------------------------------------------------
\par Thus we have proved that for every $Q\in L$ there exists a cube $K\in L$ such that $\diam K\le \frac1m\diam Q.$ Hence $$\inf_{Q\in L}\diam Q=0.$$ But, by the assumption \rf{A-G}, $\diam V_L<(1/\gamma)\diam Q$ so that %----------------------------------------------------------
$$
\inf_{Q\in L}\diam Q\ge \gamma\diam V_L
$$
%----------------------------------------------------------
which implies that $\diam V_L=0$. But $\diam V_L>0$, a contradiction. Thus the second inequality in \rf{VLQ} holds.
%----------------------------------------------------------
\par The proposition is completely proved.\bx
%@@@@@@@@@@@@@@@@@@@@@@@@@@@@@@@@@@@@@@@@@@@@@@@@@@@@@@@@@@
%@@@@@@@@@@@@@@@@@@@@@@@@@@@@@@@@@@@@@@@@@@@@@@@@@@@@@@@@@@
%@@@@@@@@@@@@@@@@@@@@@@@@@@@@@@@@@@@@@@@@@@@@@@@@@@@@@@@@@@
%@@@@@@@@@@@@@@@@@@@@@@@@@@@@@@@@@@@@@@@@@@@@@@@@@@@@@@@@@@
%----------------------------------------------------------
\bigskip
%----------------------------------------------------------
\par Given a lacuna $L\in\LE$ we introduce a pair of two cubes characterizing the ``size'' of this lacuna - a cube $Q_L$ of the minimal diameter, and the cube $\QL$ of the maximal diameter. These cubes not always exist. Let us describe conditions for their existence and present main properties of these cubes.
%----------------------------------------------------------
\begin{proposition} Let $L\in\LE$ be a lacuna.
If $\diam V_L>0$, then there exists a cube $Q_L\in L$ such that
%----------------------------------------------------------
$$
\diam Q_L=\min\{\diam Q: Q\in L\}.
$$
%----------------------------------------------------------
Furthermore,
%----------------------------------------------------------
\bel{VLT}
\tfrac{1}{\q}\diam V_L\le \diam Q_L\le \gamma \diam V_L
\ee
%----------------------------------------------------------
where $\gamma$ is an absolute constant.
%@@@@@@@@@@@@@@@@@@@@@@@@@@@@@@@@@@@@@@@@@@@@@@@@@@@@@@@@@@
%----------------------------------------------------------
\end{proposition}
%----------------------------------------------------------
%@@@@@@@@@@@@@@@@@@@@@@@@@@@@@@@@@@@@@@@@@@@@@@@@@@@@@@@@@@
%@@@@@@@@@@@@@@@@@@@@@@@@@@@@@@@@@@@@@@@@@@@@@@@@@@@@@@@@@@
%@@@@@@@@@@@@@@@@@@@@@@@@@@@@@@@@@@@@@@@@@@@@@@@@@@@@@@@@@@
%@@@@@@@@@@@@@@@@@@@@@@@@@@@@@@@@@@@@@@@@@@@@@@@@@@@@@@@@@@
%----------------------------------------------------------
\par {\it Proof.} Inequality \rf{VLT} follows from Proposition \reff{DL-M} so we turn to the proof of the existence of the cube $Q_L$.
%----------------------------------------------------------
\par Of course, its trivial for an elementary lacuna so we may assume that $L\in\tL$ is a true lacuna. Fix a cube $\tQ\in L$ and consider a family of cubes
%----------------------------------------------------------
$$
J_{\tQ}=\{K\in L: \diam K\le\diam \tQ\}.
$$
%----------------------------------------------------------
Prove that $J_{\tQ}$ is a {\it finite} family of cubes. In fact, since $L$ is a true lacuna, $V_L=(10\tQ)\cap E$ and
%----------------------------------------------------------
$$
V_L=(10K)\cap E~~~\text{for every}~~~K\in J_{\tQ}
$$
%----------------------------------------------------------
so that $(10K)\cap(10\tQ)\ne\emp$. But $\diam K\le\diam \tQ$ which implies that $K\subset 10K\subset 30\tQ$.
%----------------------------------------------------------
\par Thus $\cup\{K: K\in J_{\tQ}\}\subset 30\tQ$. By Proposition \reff{DL-M}, $\diam K\ge \diam V_L/90$ for every $K\in J_{\tQ}$. Since $J_{\tQ}\subset W_E$, the cubes of this family are non-overlapping. These properties of $J_{\tQ}$ immediately implies that $\# J_{\tQ}<\infty$.
%----------------------------------------------------------
\par Clearly,
%----------------------------------------------------------
$$
\inf_{Q\in L}\diam Q=\inf_{K\in J_{\tQ}}\diam K=\min_{K\in J_{\tQ}}\diam K
$$
%----------------------------------------------------------
so that the minimum of $\diam Q$ is attained on $L$ proving the proposition.\bx\medskip
%----------------------------------------------------------
%@@@@@@@@@@@@@@@@@@@@@@@@@@@@@@@@@@@@@@@@@@@@@@@@@@@@@@@@@@
%@@@@@@@@@@@@@@@@@@@@@@@@@@@@@@@@@@@@@@@@@@@@@@@@@@@@@@@@@@
%@@@@@@@@@@@@@@@@@@@@@@@@@@@@@@@@@@@@@@@@@@@@@@@@@@@@@@@@@@
%----------------------------------------------------------
\par Summarizing inequality \rf{VLT} and inequality \rf{A-DL}, we obtaining the following
%----------------------------------------------------------
%@@@@@@@@@@@@@@@@@@@@@@@@@@@@@@@@@@@@@@@@@@@@@@@@@@@@@@@@@@
%@@@@@@@@@@@@@@@@@@@@@@@@@@@@@@@@@@@@@@@@@@@@@@@@@@@@@@@@@@
%@@@@@@@@@@@@@@@@@@@@@@@@@@@@@@@@@@@@@@@@@@@@@@@@@@@@@@@@@@
%----------------------------------------------------------
\begin{corollary}\lbl{D-QL} For every $L\in\LE$ such that $\diam V_L>0$ we have
%----------------------------------------------------------
$$
\diam Q_L\le \eta\diam ((\theta Q_L)\cap E)
$$
%----------------------------------------------------------
where $\theta=90$ and $\eta$ is an absolute constant.
%@@@@@@@@@@@@@@@@@@@@@@@@@@@@@@@@@@@@@@@@@@@@@@@@@@@@@@@@@@
%----------------------------------------------------------
\end{corollary}
%@@@@@@@@@@@@@@@@@@@@@@@@@@@@@@@@@@@@@@@@@@@@@@@@@@@@@@@@@@
%@@@@@@@@@@@@@@@@@@@@@@@@@@@@@@@@@@@@@@@@@@@@@@@@@@@@@@@@@@
%@@@@@@@@@@@@@@@@@@@@@@@@@@@@@@@@@@@@@@@@@@@@@@@@@@@@@@@@@@
%@@@@@@@@@@@@@@@@@@@@@@@@@@@@@@@@@@@@@@@@@@@@@@@@@@@@@@@@@@
%----------------------------------------------------------
\begin{proposition}\lbl{D-MAX} For every lacuna $L\in\LE$ we have
%----------------------------------------------------------
\bel{D-MAX1}
\tfrac1\gamma\dist(V_L,E\setminus V_L)\le\sup_{Q\in L}\diam Q\le 40\,\dist(V_L,E\setminus V_L)
\ee
%----------------------------------------------------------
where $\gamma$ is an absolute constant.
%@@@@@@@@@@@@@@@@@@@@@@@@@@@@@@@@@@@@@@@@@@@@@@@@@@@@@@@@@@
%----------------------------------------------------------
\end{proposition}
%@@@@@@@@@@@@@@@@@@@@@@@@@@@@@@@@@@@@@@@@@@@@@@@@@@@@@@@@@@
%@@@@@@@@@@@@@@@@@@@@@@@@@@@@@@@@@@@@@@@@@@@@@@@@@@@@@@@@@@
%@@@@@@@@@@@@@@@@@@@@@@@@@@@@@@@@@@@@@@@@@@@@@@@@@@@@@@@@@@
%@@@@@@@@@@@@@@@@@@@@@@@@@@@@@@@@@@@@@@@@@@@@@@@@@@@@@@@@@@
%----------------------------------------------------------
\par {\it Proof.} For an elementary lacuna $L$ the statement of the proposition follows from inequality \rf{A-DL}.
%----------------------------------------------------------
\par Let $L\in\tL$ be a true lacuna so that for every cube   $Q\in L$ we have
%----------------------------------------------------------
\bel{IN-L2}
V_L=(10Q)\cap E=(\q Q)\cap E.
\ee
%----------------------------------------------------------
Hence
%----------------------------------------------------------
$$
\dist(V_L,E\setminus V_L)\ge\dist(10Q,\RN\setminus \q Q)=80 r_Q=40\diam Q
$$
%----------------------------------------------------------
proving the second inequality in \rf{D-MAX1}.
%----------------------------------------------------------
\par Prove the first inequality. This inequality is trivial whenever $\sup_{Q\in L}\diam Q=\infty$. Suppose that
%----------------------------------------------------------
\bel{S-INF}
\sup_{Q\in L}\,\,\diam Q<\,\infty\,.
\ee
%----------------------------------------------------------
\par Let us assume that
%----------------------------------------------------------
\bel{A-GM}
\dist(V_L,E\setminus V_L)>\gamma\sup_{Q\in L}\diam Q
\ee
%----------------------------------------------------------
with
%----------------------------------------------------------
$$
\gamma:=10^4.
$$
%----------------------------------------------------------
Prove that this assumption contradicts \rf{S-INF}. Let
%----------------------------------------------------------
$$
m:=100.
$$
%----------------------------------------------------------
Fix a cube $Q\in L$ and construct a cube $K\in L$ such that %----------------------------------------------------------
$$
\diam K\ge\, \tfrac{m}{5}\,\diam Q.
$$
%----------------------------------------------------------
\par By the assumption \rf{A-GM},
%----------------------------------------------------------
\bel{Q-GA}
\dist(V_L,E\setminus V_L)>\gamma\diam Q.
\ee
%----------------------------------------------------------
\par Note that the function $\dist(\cdot,V_L)$ is continuous and equals $0$ on $V_L$. Besides it is strictly grater than $\gamma\diam Q$ on $E\setminus V_L$. Since $0<m<\gamma$, there exists a point $x_0\in\RN\setminus E$ such that
%----------------------------------------------------------
\bel{X0}
\dist(x_0,V_L)=m\diam Q.
\ee
%----------------------------------------------------------
\par Let $K\in W_E$ be a Whitney cube which contains $x_0$. Prove that
%----------------------------------------------------------
$$
\diam K\ge \tfrac{m}{5}\diam Q~~~\text{and}~~~K\in L.
$$
%----------------------------------------------------------
\par First let us show that %----------------------------------------------------------
\bel{DX0}
\dist(x_0,V_L)=\dist(x_0,E).
\ee
%----------------------------------------------------------
Note that, by \rf{IN-L2}, $V_L\subset 10Q$ so that
$\diam V_L\le 10\diam Q.$
%----------------------------------------------------------
\par By \rf{Q-GA} and \rf{X0},
%----------------------------------------------------------
\be
\dist(x_0,E\setminus V_L)&\ge&\dist(V_L,E\setminus V_L)-
\diam V_L-\dist(x_0,V_L)\nn\\&>&\gamma \diam Q-10\diam Q-m\diam Q.\nn
\ee
%----------------------------------------------------------
Hence,
%----------------------------------------------------------
\bel{X0-G}
\dist(x_0,E\setminus V_L)>(\gamma-m-10)\diam Q>m\diam Q.
\ee
%----------------------------------------------------------
Combining this inequality with \rf{X0} we conclude that  $\dist(x_0,V_L)<\dist(x_0,E\setminus V_L)$ which proves \rf{DX0}.
%----------------------------------------------------------
\par Since $K\in W_E$, by Theorem \reff{Wcov},
%----------------------------------------------------------
$$
\diam K\le \dist(K,E)\le 4\diam K.
$$
%----------------------------------------------------------
Prove that %----------------------------------------------------------
\bel{E-M5}
\tfrac{m}{5}\diam Q\le \diam K\le m \diam Q.
\ee
%----------------------------------------------------------
Since $x_0\in K$, we have
%----------------------------------------------------------
$$
\dist(x_0,V_L)=\dist(x_0,E)\le\dist(K,E)+\diam K\le 5\diam K
$$
%----------------------------------------------------------
so that, by \rf{X0}, $m\diam Q\le 5\diam K$. Conversely,
%----------------------------------------------------------
$$
\diam K\le \dist(K,E)\le \dist(x_0,E)=\dist(x_0,V_L)
$$
%----------------------------------------------------------
so that, by \rf{X0},  $\diam K\le m\diam Q$, and  \rf{E-M5} is proved.
%----------------------------------------------------------
\par Prove that $K\in L$, i.e.,
%----------------------------------------------------------
\bel{LCN}
V_L=(10K)\cap E=(\q K)\cap E.
\ee
%----------------------------------------------------------
First let us show that
%----------------------------------------------------------
\bel{DKE-M}
\dist(K,E)=\dist(K,V_L).
\ee
%----------------------------------------------------------
\par Let $y\in K$. Then, by \rf{E-M5},
%----------------------------------------------------------
$$
\dist(y,V_L)\le\|x_0-y\|+\dist(x_0,V_L)\le\diam K+m\diam Q
$$
%----------------------------------------------------------
so that, by \rf{E-M5},
%----------------------------------------------------------
$$
\dist(y,V_L)\le 2m\diam Q.
$$
%----------------------------------------------------------
On the other hand, by \rf{X0-G} and \rf{E-M5},
%----------------------------------------------------------
\be
\dist(y,E\setminus V_L)&\ge&
\dist(x_0,E\setminus V_L)-\diam K\nn\\&\ge&(\gamma-m-20)\diam Q-m\diam Q\nn\\&=&(\gamma-2m-20)\diam Q>2m\diam Q.
\nn
\ee
%----------------------------------------------------------
Hence
%----------------------------------------------------------
$$
\dist(y,V_L)<\dist(y,E\setminus V_L)~~~\text{for every}~~~y\in K
$$
%----------------------------------------------------------
proving \rf{DKE-M}.
%----------------------------------------------------------
\par Repeating the proof of inequality \rf{P-2} from  Proposition \reff{DL-M} we conclude that for every $z\in V_L$
%----------------------------------------------------------
$$
\|z-c_K\|\le\diam V_L+9r_K.
$$
%----------------------------------------------------------
(Recall that $K=Q(c_K,r_K)$.)
%----------------------------------------------------------
\par Since $V_L\subset 10 Q$, we have $\diam V_L\le 10\diam Q$ so that
%----------------------------------------------------------
$$
\|z-c_K\|\le10\diam Q+9r_K.
$$
%----------------------------------------------------------
Combining this inequality with \rf{E-M5} we obtain
%----------------------------------------------------------
$$
\|z-c_K\|\le 10(5\diam K/m)+9r_K=(9+50/m)r_K<10r_K.
$$
%----------------------------------------------------------
Hence $V_L\subset 10Q.$
%----------------------------------------------------------
\par Prove that
%----------------------------------------------------------
\bel{90-E}
(90K)\cap(E\setminus V_L)=\emp.
\ee
%----------------------------------------------------------
In fact, since $90K\supset V_L$, we have
%----------------------------------------------------------
$$
\dist(90K,E\setminus V_L)\ge \dist(V_L,E\setminus V_L)-
\diam(90K)
$$
%----------------------------------------------------------
so that, by the assumption \rf{Q-GA},
%----------------------------------------------------------
$$
\dist(90K,E\setminus V_L)\ge \gamma\diam Q-
90\diam K.
$$
%----------------------------------------------------------
Since $\diam Q\ge\diam K/m$, see \rf{E-M5}, we obtain
%----------------------------------------------------------
$$
\dist(90K,E\setminus V_L)\ge (\gamma/m-90)\diam K.
$$
%----------------------------------------------------------
Since $\gamma>90m$, we have $\dist(90K,E\setminus V_L)>0$
proving \rf{90-E}.
%----------------------------------------------------------
\par Hence $(90K)\cap E=(90K)\cap V_L=V_L.$ Finally,
%----------------------------------------------------------
$$
V_L\subset (10K)\cap E\subset(90K)\cap E=V_L
$$
%----------------------------------------------------------
proving \rf{LCN}. Thus $K\in L$.
%----------------------------------------------------------
\par We have proved that for very cube $Q\in L$ there exists a cube $K\in L$ such that
%----------------------------------------------------------
$$
\diam K\ge \tfrac{m}{5}\diam Q.
$$
%----------------------------------------------------------
Hence $\sup\{\diam Q:Q\in L\}=\infty$ which contradicts the condition \rf{S-INF}.
%----------------------------------------------------------
\par Thus inequality \rf{A-GM} is not true which proves the first inequality in \rf{D-MAX1} and the proposition.\bx\medskip
%----------------------------------------------------------
%@@@@@@@@@@@@@@@@@@@@@@@@@@@@@@@@@@@@@@@@@@@@@@@@@@@@@@@@@@
%@@@@@@@@@@@@@@@@@@@@@@@@@@@@@@@@@@@@@@@@@@@@@@@@@@@@@@@@@@
%@@@@@@@@@@@@@@@@@@@@@@@@@@@@@@@@@@@@@@@@@@@@@@@@@@@@@@@@@@
%@@@@@@@@@@@@@@@@@@@@@@@@@@@@@@@@@@@@@@@@@@@@@@@@@@@@@@@@@@
%----------------------------------------------------------
\par Let $L\in\LE$ be a lacuna. By $U_L$ we denote the union of all cubes which belong to the lacuna:
%----------------------------------------------------------
$$
U_L:=\cup\{Q:Q\in L\}.
$$
%----------------------------------------------------------
By $\diam L$ we denote the diameter of the set $U_L$:
%----------------------------------------------------------
$$
\diam L:=\diam U_L=\sup\{\|a-b\|:a,b\in U_L\}.
$$
%----------------------------------------------------------
\par We say that $L$ is bounded if $\diam L<\infty$. If $\diam L=\infty$ we say that $L$ is an unbounded lacuna.
%----------------------------------------------------------
\begin{proposition} (i). For every lacuna $L\in\LE$
%----------------------------------------------------------
\bel{M-DM3}
\diam L\sim \sup\{\diam Q: Q\in L\}\sim \dist(V_L,E\setminus V_L)
\ee
%----------------------------------------------------------
with absolute constants in the equivalences;
%----------------------------------------------------------
%@@@@@@@@@@@@@@@@@@@@@@@@@@@@@@@@@@@@@@@@@@@@@@@@@@@@@@@@@@
\par (ii). If $E$ is an unbounded set then  every lacuna $L\in\LE$ is bounded;
%----------------------------------------------------------
\par (iii). If $E$ is bounded, there exists the unique unbounded lacuna $L^{\max}\in\LE$. The lacuna $L^{\max}$ is a true lacuna for which $V_{L^{\max}}=E$.
%----------------------------------------------------------
%@@@@@@@@@@@@@@@@@@@@@@@@@@@@@@@@@@@@@@@@@@@@@@@@@@@@@@@@@@
%----------------------------------------------------------
\end{proposition}
%@@@@@@@@@@@@@@@@@@@@@@@@@@@@@@@@@@@@@@@@@@@@@@@@@@@@@@@@@@
%@@@@@@@@@@@@@@@@@@@@@@@@@@@@@@@@@@@@@@@@@@@@@@@@@@@@@@@@@@
%@@@@@@@@@@@@@@@@@@@@@@@@@@@@@@@@@@@@@@@@@@@@@@@@@@@@@@@@@@
%@@@@@@@@@@@@@@@@@@@@@@@@@@@@@@@@@@@@@@@@@@@@@@@@@@@@@@@@@@
%----------------------------------------------------------
\par {\it Proof.} (i). The second equivalence in \rf{M-DM3} follows from Proposition \reff{D-MAX}.
%----------------------------------------------------------
\par Clearly $\sup\{\diam Q: Q\in L\}\le \diam L$. Prove  that
%----------------------------------------------------------
\bel{DL-30}
\diam L\le 30\sup\{\diam Q:Q\in L\}.
\ee
%----------------------------------------------------------
\par Of course this inequality is trivial for elementary lacunae. Consider a true lacuna $L\in\tL$. Let cubes $K,K'\in L$ and let $a\in K$ and $a'\in K'$. Suppose that $\diam K'\le\diam K$. Since $L$ is a true lacuna,
%----------------------------------------------------------
$$
V_L=(10K)\cap E=(10K')\cap E
$$
%----------------------------------------------------------
so that  $(10K)\cap (10K')\ne\emp.$ Since $\diam K'\le\diam K$, we have
%----------------------------------------------------------
$$
K'\subset 10K'\subset 30K.
$$
%----------------------------------------------------------
Hence $a,a'\in 30K$ so that
%----------------------------------------------------------
$$
\|a-a'\|\le \diam (30K)=30\diam K
$$
%----------------------------------------------------------
proving \rf{DL-30} and part (i) of the proposition.\bigskip
%----------------------------------------------------------
%@@@@@@@@@@@@@@@@@@@@@@@@@@@@@@@@@@@@@@@@@@@@@@@@@@@@@@@@@@
%@@@@@@@@@@@@@@@@@@@@@@@@@@@@@@@@@@@@@@@@@@@@@@@@@@@@@@@@@@
%@@@@@@@@@@@@@@@@@@@@@@@@@@@@@@@@@@@@@@@@@@@@@@@@@@@@@@@@@@
%----------------------------------------------------------
\par (ii). Let $E$ be an unbounded set. Of course, in this case every elementary lacuna of $E$ is bounded. If $L$ is a true lacuna, by part (i) of the proposition,
%----------------------------------------------------------
\bel{DTL}
\diam L\sim \dist(V_L,E\setminus V_L).
\ee
%----------------------------------------------------------
Since $V_L$ is a bounded set and $E$ is an unbounded set,
$\dist(V_L,E\setminus V_L)<\infty$ proving that $\diam L<\infty$.\bigskip
%----------------------------------------------------------
\par (iii). Let $E$ be a bounded set and let $L\in\tL$ be a true lacuna. If $V_L\ne E$, then
%----------------------------------------------------------
$$
\dist(V_L,E\setminus V_L)<\infty
$$
%----------------------------------------------------------
so that, by \rf{DTL},  $\diam L<\infty$.
%----------------------------------------------------------
\par On the other hand, the unique true lacuna $L^{\max}$ such that $V_{L^{\max}}=E$ contains cubes of arbitrarily big sizes so that $\diam L^{\max}=\infty$.\bx\medskip
%----------------------------------------------------------
%@@@@@@@@@@@@@@@@@@@@@@@@@@@@@@@@@@@@@@@@@@@@@@@@@@@@@@@@@@
%@@@@@@@@@@@@@@@@@@@@@@@@@@@@@@@@@@@@@@@@@@@@@@@@@@@@@@@@@@
%@@@@@@@@@@@@@@@@@@@@@@@@@@@@@@@@@@@@@@@@@@@@@@@@@@@@@@@@@@
%@@@@@@@@@@@@@@@@@@@@@@@@@@@@@@@@@@@@@@@@@@@@@@@@@@@@@@@@@@
%----------------------------------------------------------
%@@@@@@@@@@@@@@@@@@@@@@@@@@@@@@@@@@@@@@@@@@@@@@@@@@@@@@@@@@
%----------------------------------------------------------
\begin{proposition}\lbl{CL-MAX} Let $L\in \LE$ be a bounded lacuna. Then there exists a cube  $\QL\in L$ such that
%----------------------------------------------------------
$$
\diam \QL=\max\{\diam K: K\in L\}.
$$
%----------------------------------------------------------
Furthermore,
%----------------------------------------------------------
\bel{CL-D}
\diam \QL\sim \diam L\sim\dist(V_L,E\setminus V_L),
\ee
%----------------------------------------------------------
and
%----------------------------------------------------------
\bel{MQ-I}
V_L\cup U_L\subset \gamma \QL.
\ee
%----------------------------------------------------------
Here the constant $\gamma$ and constants in the equivalences of \rf{CL-D} are absolute .
%@@@@@@@@@@@@@@@@@@@@@@@@@@@@@@@@@@@@@@@@@@@@@@@@@@@@@@@@@@
%----------------------------------------------------------
\end{proposition}
%@@@@@@@@@@@@@@@@@@@@@@@@@@@@@@@@@@@@@@@@@@@@@@@@@@@@@@@@@@
%@@@@@@@@@@@@@@@@@@@@@@@@@@@@@@@@@@@@@@@@@@@@@@@@@@@@@@@@@@
%@@@@@@@@@@@@@@@@@@@@@@@@@@@@@@@@@@@@@@@@@@@@@@@@@@@@@@@@@@
%@@@@@@@@@@@@@@@@@@@@@@@@@@@@@@@@@@@@@@@@@@@@@@@@@@@@@@@@@@
%----------------------------------------------------------
\par {\it Proof.} Clearly equivalences \rf{CL-D} follow from \rf{M-DM3}, so we turn to the proof of the existence
of the cube $\QL$. Fix a cube $\tQ\in L$ and introduce a family of cubes
%----------------------------------------------------------
$$
I_{\tQ}=\{K\in L: \diam K\ge\diam \tQ\}.
$$
%----------------------------------------------------------
Note that the diameters of the cubes from $I_{\tQ}$ are at least $\diam \tQ$, these cubes are non-overlapping and lie in a bounded set so that the family $I_{\tQ}$ is finite. Therefore there exists a cube $\QL\in I_{\tQ}$ such that $\diam \QL=\max\{\diam K: K\in I_{\tQ}\}$. But
%----------------------------------------------------------
$$
\max\{\diam K: K\in I_{\tQ}\}=\max\{\diam K: K\in L\}
$$
%----------------------------------------------------------
proving that
%----------------------------------------------------------
$$
\diam \QL=\max\{\diam K: K\in L\}.
$$
%----------------------------------------------------------
\par It remains to prove inclusion \rf{MQ-I}.
%----------------------------------------------------------
\par Let $\gamma=270$. We know that $V_L=(90 Q)\cap E$ for every cube $Q\in L$. In particular, $V_L\subset 90\QL\subset \gamma\QL$.
%----------------------------------------------------------
\par Since $V_L\subset(90 Q)\cap (90\QL)$, we conclude that  $(90 Q)\cap (90\QL)\ne\emp.$ But $\diam Q\le\diam \QL$ so that $Q\subset 270 \QL$ proving that
%----------------------------------------------------------
$$
U_L=\cup\{Q: Q\in L\}\subset \gamma \QL.\BX
$$
%----------------------------------------------------------
%@@@@@@@@@@@@@@@@@@@@@@@@@@@@@@@@@@@@@@@@@@@@@@@@@@@@@@@@@@
%@@@@@@@@@@@@@@@@@@@@@@@@@@@@@@@@@@@@@@@@@@@@@@@@@@@@@@@@@@
%@@@@@@@@@@@@@@@@@@@@@@@@@@@@@@@@@@@@@@@@@@@@@@@@@@@@@@@@@@
%@@@@@@@@@@@@@@@@@@@@@@@@@@@@@@@@@@@@@@@@@@@@@@@@@@@@@@@@@@
%@@@@@@@@@@@@@@@@@@@@@@@@@@@@@@@@@@@@@@@@@@@@@@@@@@@@@@@@@@
%@@@@@@@@@@@@@@@@@@@@@@@@@@@@@@@@@@@@@@@@@@@@@@@@@@@@@@@@@@
%@@@@@@@@@@@@@@@@@@@@@@@@@@@@@@@@@@@@@@@@@@@@@@@@@@@@@@@@@@
%----------------------------------------------------------
\begin{proposition}\lbl{M-LAC} Let $L\in\LE$ be a lacuna and let $Q\in L$. Suppose that there exist a lacuna $L'\in \LE$, $L\ne L'$, and a cube $Q'\in L'$ such that $Q\cap Q'\ne\emp$. Then:
%----------------------------------------------------------
\par (i). If $L$ is a true lacuna, then $L'$ is an elementary lacuna, i.e., $L'\in \hL_E=\LE\setminus\tL_E$;
%----------------------------------------------------------
\par (ii). Either
%----------------------------------------------------------
\bel{EQV-1}
\diam Q\sim \diam V_L\sim \diam Q_L
\ee
%----------------------------------------------------------
or
%----------------------------------------------------------
\bel{EQV-2}
\diam Q\sim \dist(V_L,E\setminus V_L)\sim \diam\QL
\ee
%----------------------------------------------------------
with absolute constants in the equivalences.
%@@@@@@@@@@@@@@@@@@@@@@@@@@@@@@@@@@@@@@@@@@@@@@@@@@@@@@@@@@
%----------------------------------------------------------
\end{proposition}
%@@@@@@@@@@@@@@@@@@@@@@@@@@@@@@@@@@@@@@@@@@@@@@@@@@@@@@@@@@
%@@@@@@@@@@@@@@@@@@@@@@@@@@@@@@@@@@@@@@@@@@@@@@@@@@@@@@@@@@
%@@@@@@@@@@@@@@@@@@@@@@@@@@@@@@@@@@@@@@@@@@@@@@@@@@@@@@@@@@
%@@@@@@@@@@@@@@@@@@@@@@@@@@@@@@@@@@@@@@@@@@@@@@@@@@@@@@@@@@
\par {\it Proof.} (i). Since $L\in\tL_E$ is a true lacuna, for every $Q\in L$ we have
%----------------------------------------------------------
$$
V_L=(10 Q)\cap E=(\q Q)\cap E.
$$
%----------------------------------------------------------
\par Prove that $L'\in\hL=\Lc\setminus \tL_E$. In fact, if $L'\in\tL_E$, then
%----------------------------------------------------------
$$
V_{L'}=(10 Q')\cap E=(\q Q')\cap E.
$$
%----------------------------------------------------------
Since $L\ne L'$, we have $V_L\ne V_{L'}$. We know that $Q,Q'\in W_E$ and $Q\cap Q'\ne\emp$ so that, by Lemma \reff{Wadd},
%----------------------------------------------------------
\bel{KTW}
1/4\diam Q\le \diam Q'\le 4\diam Q.
\ee
%----------------------------------------------------------
Hence $Q'\subset 9Q$ so that $10 Q'\subset \q Q$ proving that
%----------------------------------------------------------
$$
V_{L'}=(10 Q')\cap E\subset (\q Q)\cap E=V_L.
$$
%----------------------------------------------------------
In the same way we prove that $V_L\subset V_{L'}$. Hence $V_L=V_{L'}$, a contradiction.\medskip
%----------------------------------------------------------
%@@@@@@@@@@@@@@@@@@@@@@@@@@@@@@@@@@@@@@@@@@@@@@@@@@@@@@@@@@
%@@@@@@@@@@@@@@@@@@@@@@@@@@@@@@@@@@@@@@@@@@@@@@@@@@@@@@@@@@
%@@@@@@@@@@@@@@@@@@@@@@@@@@@@@@@@@@@@@@@@@@@@@@@@@@@@@@@@@@
%@@@@@@@@@@@@@@@@@@@@@@@@@@@@@@@@@@@@@@@@@@@@@@@@@@@@@@@@@@
%----------------------------------------------------------
\par (ii). Note that the second equivalence in \rf{EQV-1} follows from \rf{VLT}. In turn the second equivalence in \rf{EQV-2} follows from \rf{CL-D}. Thus we have to prove that either the first equivalence in \rf{EQV-1} or the first equivalence in \rf{EQV-2} holds.
%----------------------------------------------------------
\par Note that if $L\in\hL$, i.e., $L$ is an elementary lacuna, then, by \rf{A-DL}, equivalence \rf{EQV-1} holds. Thus we may assume that $L$ is a true lacuna.
%----------------------------------------------------------
\par In part (i) we have proved that in this case $L'\in\hL=\LE\setminus\tL_E$ so that the cube $Q'$ is the unique cube which belongs to the lacuna $L'$. Note that, by \rf{A-DL},
%----------------------------------------------------------
$$
\diam Q'\le 2\diam ((\q Q')\cap E).
$$
%----------------------------------------------------------
Also recall that $Q'\subset 9Q$, so that
$10 Q'\subset \q Q$. Combining this with \rf{KTW}, we obtain %----------------------------------------------------------
\bel{DQ-A}
\diam Q\le 4\diam Q'\le 8\diam ((\q Q')\cap E)\le
8\diam ((\gamma_1 Q)\cap E)
\ee
%----------------------------------------------------------
with $\gamma_1:=810$.  Now consider two cases.
%----------------------------------------------------------
\par {\it The first case:}
%----------------------------------------------------------
$$
\diam Q\le 8\diam V_L.
$$
%----------------------------------------------------------
In this case the first equivalence in \rf{EQV-1} holds. In fact, since $V_L=(10Q)\cap E$, we have
$\diam V_L\le 10\diam Q$ so that
%----------------------------------------------------------
$$
\tfrac18 \diam Q\le \diam V_L\le 10 \diam Q.
$$
%----------------------------------------------------------
\par {\it The second case:}
%----------------------------------------------------------
\bel{VL-Q}
\diam V_L<\tfrac18 \diam Q.
\ee
%----------------------------------------------------------
Prove that in this case equivalence \rf{EQV-2} holds. In fact, by \rf{DQ-A}, there exist points $a,b\in (\gamma_1 Q)\cap E$ such that
%----------------------------------------------------------
$$
\tfrac18\diam Q\le \|a-b\|
$$
%----------------------------------------------------------
so that, by \rf{VL-Q}, either $a$ or $b$ does not belong to $V_L$. Assume that $a\notin V_L$. Then
%----------------------------------------------------------
$$
\dist(V_L,E\setminus V_L)\le \dist(a, V_L).
$$
%----------------------------------------------------------
Since $V_L\subset 10 Q\subset \gamma_1 Q$ and $a\in\gamma_1 Q$, we conclude that
%----------------------------------------------------------
$$
\dist(a,V_L)\le \diam(\gamma_1 Q)=\gamma_1\diam Q
$$
%----------------------------------------------------------
proving that
%----------------------------------------------------------
\bel{AEV}
\dist(V_L,E\setminus V_L)\le \gamma_1\diam Q.
\ee
%----------------------------------------------------------
\par On the other hand, by \rf{M-DM3},
%----------------------------------------------------------
$$
\diam Q\le \sup_{K\in L}\diam K\le\gamma_2 \dist(V_L,E\setminus V_L)
$$
%----------------------------------------------------------
with an absolute constant $\gamma_2$. This inequality and inequality \rf{AEV} prove equivalence \rf{EQV-2} and finish the proof of the proposition.\bx\medskip
%----------------------------------------------------------
%@@@@@@@@@@@@@@@@@@@@@@@@@@@@@@@@@@@@@@@@@@@@@@@@@@@@@@@@@@
%@@@@@@@@@@@@@@@@@@@@@@@@@@@@@@@@@@@@@@@@@@@@@@@@@@@@@@@@@@
%@@@@@@@@@@@@@@@@@@@@@@@@@@@@@@@@@@@@@@@@@@@@@@@@@@@@@@@@@@
%@@@@@@@@@@@@@@@@@@@@@@@@@@@@@@@@@@@@@@@@@@@@@@@@@@@@@@@@@@
%----------------------------------------------------------
\begin{proposition}\lbl{INT-L1} Let $L\in\LE$ be a lacuna and let
%----------------------------------------------------------
$$
\Ic_L:=\{K\in W_E\setminus L: \exists~Q\in L~~\text{such that}~~K\cap Q\ne\emp\}.
$$
%----------------------------------------------------------
Then $\card\Ic_L\le\gamma(n).$
%----------------------------------------------------------
%@@@@@@@@@@@@@@@@@@@@@@@@@@@@@@@@@@@@@@@@@@@@@@@@@@@@@@@@@@
%----------------------------------------------------------
\end{proposition}
%@@@@@@@@@@@@@@@@@@@@@@@@@@@@@@@@@@@@@@@@@@@@@@@@@@@@@@@@@@
%@@@@@@@@@@@@@@@@@@@@@@@@@@@@@@@@@@@@@@@@@@@@@@@@@@@@@@@@@@
%@@@@@@@@@@@@@@@@@@@@@@@@@@@@@@@@@@@@@@@@@@@@@@@@@@@@@@@@@@
%@@@@@@@@@@@@@@@@@@@@@@@@@@@@@@@@@@@@@@@@@@@@@@@@@@@@@@@@@@
\par {\it Proof.} Let $K\in\Ic_L$, i.e., $K$ belongs to a lacuna $L'$, $L'\ne L$, and $K\cap Q_K\ne\emp$ for some cube $Q_K\in L$. Then, by Proposition \reff{M-LAC},
either $\diam Q_K\sim\diam V_L$ or
%----------------------------------------------------------
$$
\diam Q_K\sim\dist (V_L,E\setminus V_L).
$$
%----------------------------------------------------------
Since
%----------------------------------------------------------
\bel{AJ}
1/4\diam Q_K\le \diam K\le 4\diam Q_K,
\ee
%----------------------------------------------------------
see Lemma \reff{Wadd}, we conclude that either
%----------------------------------------------------------
\bel{O1}
\diam K\sim\diam V_L
\ee
%----------------------------------------------------------
or
%----------------------------------------------------------
\bel{O2}
\diam K\sim\dist (V_L,E\setminus V_L)
\ee
%----------------------------------------------------------
(with absolute constants in the equivalences).
%----------------------------------------------------------
\par Let us denote by $\Ic_L^{(1)}$ a subfamily of $\Ic_L$ consisting of those cubes $K$ which satisfy $\rf{O1}$. By $\Ic_L^{(2)}$ we denote those cubes $K\in\Ic_L$ which satisfy inequality $\rf{O2}$. Then $\Ic_L=\Ic_L^{(1)}\cup\Ic_L^{(2)}$ so that
%----------------------------------------------------------
$$
\card\Ic_L\le\card\Ic_L^{(1)}+\card\Ic_L^{(2)}.
$$
%----------------------------------------------------------
\par Prove that $\card\Ic_L^{(1)}\le \gamma(n).$ To this end we put $R_1:=\diam V_L$. Then
%----------------------------------------------------------
\bel{DI1}
\diam K\sim R_1~~\text{for all}~~K\in\Ic_L^{(1)}.
\ee
%----------------------------------------------------------
\par Since $Q_K\in L$,
%----------------------------------------------------------
$$
V_L=(10 Q_K)\cap E=(\q Q_K)\cap E.
$$
%----------------------------------------------------------
Let us fix a point $x_0\in V_L$. Then for every $y\in K$ we have
%----------------------------------------------------------
$$
\|x_0-y\|\le \|x_0-c_{Q_K}\|+\|c_{Q_K}-c_K\|\le 10 r_{Q_K}+r_{Q_K}+r_K.
$$
%----------------------------------------------------------
(Recall that $c_K$ denotes the center of $K$ and  $r_K=\tfrac12\diam K$.) Hence, by \rf{AJ} and \rf{DI1},
%----------------------------------------------------------
$$
\|x_0-y\|\le 40r_K+4r_K+r_K=45 r_K\le \gamma_1 R_1
$$
%----------------------------------------------------------
with some absolute constant $\gamma_1$. This proves that
%----------------------------------------------------------
\bel{KTQ}
K\subset  Q(x_0,\gamma_1 R_1).
\ee
%----------------------------------------------------------
Since $\Ic_L^{(1)}\subset W_E$, the cubes of the family
$\Ic_L^{(1)}$ are non-overlapping. By \rf{KTQ}, all these cubes are contained in the cube $\tQ_1:=\gamma_1 Q(x_0,R_1)$, and, by \rf{DI1}, the diameter of each such a cube is equivalent to  $\diam \tQ_1\sim R_1$. This proves that the number of cubes in $\Ic_L^{(1)}$ does not exceed a constant depending only on $n$.
%----------------------------------------------------------
\par In the same way we prove that $\card \Ic_L^{(2)}\le \gamma(n).$ In fact, we put $R_2:=\dist(V_L,E\setminus V_L)$ so that $\diam K\sim  R_2$ for all $K\in \Ic_L^{(2)}$. Using the same approach we show that every cube $K\in \Ic_L^{(2)}$ is contained in a cube $\tQ_2:=\gamma_2 Q(x_0,R_2)$ where $\gamma_2$ is an  absolute constant . This implies the required inequality $\card \Ic_L^{(2)}\le \gamma(n).$   %----------------------------------------------------------
\par The proposition is proved.\bx
%----------------------------------------------------------
%@@@@@@@@@@@@@@@@@@@@@@@@@@@@@@@@@@@@@@@@@@@@@@@@@@@@@@@@@@
%@@@@@@@@@@@@@@@@@@@@@@@@@@@@@@@@@@@@@@@@@@@@@@@@@@@@@@@@@@
%@@@@@@@@@@@@@@@@@@@@@@@@@@@@@@@@@@@@@@@@@@@@@@@@@@@@@@@@@@
%@@@@@@@@@@@@@@@@@@@@@@@@@@@@@@@@@@@@@@@@@@@@@@@@@@@@@@@@@@
%----------------------------------------------------------
%@@@@@@@@@@@@@@@@@@@@@@@@@@@@@@@@@@@@@@@@@@@@@@@@@@@@@@@@@@
%@@@@@@@@@@@@@@@@@@@@@@@@@@@@@@@@@@@@@@@@@@@@@@@@@@@@@@@@@@
%@@@@@@@@@@@@@@@@@@@@@@@@@@@@@@@@@@@@@@@@@@@@@@@@@@@@@@@@@@
%@@@@@@@@@@@@@@@@@@@@@@@@@      @@@@@@@@@@@@@@@@@@@@@@@@@@@
%@@@@@@@@@@@@@@@@@@@@@@@          @@@@@@@@@@@@@@@@@@@@@@@@@
%@@@@@@@@@@@@@@@@@@@@@              @@@@@@@@@@@@@@@@@@@@@@@
%@@@@@@@@@@@@@@@@@@@     SECTION 5    @@@@@@@@@@@@@@@@@@@@@
%@@@@@@@@@@@@@@@@@@@@@              @@@@@@@@@@@@@@@@@@@@@@@
%@@@@@@@@@@@@@@@@@@@@@@@          @@@@@@@@@@@@@@@@@@@@@@@@@
%@@@@@@@@@@@@@@@@@@@@@@@@@      @@@@@@@@@@@@@@@@@@@@@@@@@@@
%@@@@@@@@@@@@@@@@@@@@@@@@@@@@@@@@@@@@@@@@@@@@@@@@@@@@@@@@@@
%@@@@@@@@@@@@@@@@@@@@@@@@@@@@@@@@@@@@@@@@@@@@@@@@@@@@@@@@@@
%@@@@@@@@@@@@@@@@@@@@@@@@@@@@@@@@@@@@@@@@@@@@@@@@@@@@@@@@@@
%----------------------------------------------------------
\SECT{5. ``Projections'' of lacunae and interior bridges.} {5}
%----------------------------------------------------------
\addtocontents{toc}{5. ``Projections'' of lacunae and interior bridges. \hfill \thepage\par}
%----------------------------------------------------------
\indent
%@@@@@@@@@@@@@@@@@@@@@@@@@@@@@@@@@@@@@@@@@@@@@@@@@@@@@@@@@@
%----------------------------------------------------------
%@@@@@@@@@@@@@@@@@@@@@@@@@@@@@@@@@@@@@@@@@@@@@@@@@@@@@@@@@@
%@@@@@@@@@@@@@@@@@@@@@@@@@@@@@@@@@@@@@@@@@@@@@@@@@@@@@@@@@@
%@@@@@@@@@@@@@@@@@@@@@@@@@@@@@@@@@@@@@@@@@@@@@@@@@@@@@@@@@@
%@@@@@@@@@@@@@@@@@@@@@@@@@@@@@@@@@@@@@@@@@@@@@@@@@@@@@@@@@@
%@@@@@@@@@@@@@@@@@@@@@@@@@@@@@@@@@@@@@@@@@@@@@@@@@@@@@@@@@@
\par {\bf 5.1. ``Projections'' of the lacunae.}
%----------------------------------------------------------
\addtocontents{toc}{~~~~5.1. ``Projections'' of the lacunae. \hfill \thepage\par}
%----------------------------------------------------------
In this subsection we construct a mapping $\LE\ni L\mapsto \PRL(L)\in E$ which we have mentioned in Section 2. Let us recall its main properties:
%-----------------------------------------------------------
\par (i) For each lacuna $L$ the point $\PRL(L)$ lies in a fixed dilation of the minimal cube of the lacuna, and
%-----------------------------------------------------------
\par (ii) every point $A\in E$ has at most $C(n)$ ``sources", i.e., lacunae $L'\in\Lc_E$ such that $\PRL(L')=A$.\smallskip
%-----------------------------------------------------------
\par We refer to the mapping $\PRL$ as a ``projection'' of $\Lc_E$ into the set $E$.\smallskip
%----------------------------------------------------------
\par The existence of the mapping $\PRL$ is proven in Proposition \reff{L-PE}. This result relies on Lemma \reff{NET-S} below. For its formulation we need the following notions: Let $S$ be a closed subset of $\RN$, and let $\ve>0$. As usual, a set $A\subset S$ is said to be an {\it $\ve$-net in $S$} if for each $x\in S$ there exists a point $a_x\in A$ such that $\|a_x-x\|\le \ve$. We say that points $x,y\in A$ are {\it $\ve$-separated} if $\|x-y\|\ge\ve$.
%----------------------------------------------------------
%@@@@@@@@@@@@@@@@@@@@@@@@@@@@@@@@@@@@@@@@@@@@@@@@@@@@@@@@@@
%@@@@@@@@@@@@@@@@@@@@@@@@@@@@@@@@@@@@@@@@@@@@@@@@@@@@@@@@@@
%@@@@@@@@@@@@@@@@@@@@@@@@@@@@@@@@@@@@@@@@@@@@@@@@@@@@@@@@@@
%@@@@@@@@@@@@@@@@@@@@@@@@@@@@@@@@@@@@@@@@@@@@@@@@@@@@@@@@@@
%----------------------------------------------------------
\begin{lemma}\lbl{NET-S} For every closed set $S\subset\RN$ there exists a decreasing sequence of non-empty closed sets $\{S_i\}_{i\in\mZ}$, $S_{i+1}\subset S_i\subset S$, $i\in\mZ,$ such that for every $i\in\mZ$ the following conditions are satisfied:
%----------------------------------------------------------
\par (i). The points of the set $S_i$ are $2^i$-separated;
%----------------------------------------------------------
%@@@@@@@@@@@@@@@@@@@@@@@@@@@@@@@@@@@@@@@@@@@@@@@@@@@@@@@@@@
%----------------------------------------------------------
\par (ii). $S_i$ is a $2^{i+1}$-net in $S$.
%----------------------------------------------------------
\end{lemma}
%----------------------------------------------------------
%@@@@@@@@@@@@@@@@@@@@@@@@@@@@@@@@@@@@@@@@@@@@@@@@@@@@@@@@@@
%@@@@@@@@@@@@@@@@@@@@@@@@@@@@@@@@@@@@@@@@@@@@@@@@@@@@@@@@@@
%@@@@@@@@@@@@@@@@@@@@@@@@@@@@@@@@@@@@@@@@@@@@@@@@@@@@@@@@@@
\par {\it Proof.} First let us construct the sets $S_i$ for $i\ge 0$.
%@@@@@@@@@@@@@@@@@@@@@@@@@@@@@@@@@@@@@@@@@@@@@@@@@@@@@@@@@@
\par We let $S_0$ denote a {\it maximal} $1$-net in $S$; thus $S_0$ is a $1$-net in $S$ whose points are $1$-separable.
%@@@@@@@@@@@@@@@@@@@@@@@@@@@@@@@@@@@@@@@@@@@@@@@@@@@@@@@@@@
\par By $S_1$ we denote a maximal $2$-net in $S_0$ so that $S_1$ is a $2$-net in $S_0$ whose points are $2$-separable. We continue this procedure and at the $m$-th step we have subsets
%----------------------------------------------------------
$$
S\supset S_0\supset S_1\supset S_2\supset...\supset S_m
$$
%----------------------------------------------------------
such that each set $S_i$ is a $2^i$-net in $S_{i-1}$, and the points of $S_i$ are $2^i$-separable. We let $S_{m+1}$ denote a maximal $2^{m+1}$-net in $S_m$ so that $S_{m+1}$ is a $2^{m+1}$-net in $S_m$, and the points of $S_{m+1}$ are $2^{m+1}$-separable.
%@@@@@@@@@@@@@@@@@@@@@@@@@@@@@@@@@@@@@@@@@@@@@@@@@@@@@@@@@@
\par Prove that the set $S_i$ is a $2^{i+1}$-net in $S$ for every integer $i\ge 0$. In fact, since $S_0$ is a $1$-net in $S$, for every $x\in S$ there exists $x_0\in S_0$ such that $\|x-x_0\|\le 1$. Since $S_1$ is a $2$-net in $S_0$, there exists $x_1\in S_1$ such that $\|x_0-x_1\|\le 2$. Continuing this process we obtain points $x_j\in S_j$, $j=0,1,...,i,$ such that $\|x_j-x_{j+1}\|\le 2^{j+1},$
$j=0,1,...,i-1.$
%@@@@@@@@@@@@@@@@@@@@@@@@@@@@@@@@@@@@@@@@@@@@@@@@@@@@@@@@@@
\par Hence
%----------------------------------------------------------
$$
\|x-x_i\|\le \|x-x_0\|+\sum_{j=0}^{i-1}\|x_j-x_{j+1}\|
\le\sum_{j=0}^{i}2^j=2^{i+1}-1\le 2^{i+1}
$$
%----------------------------------------------------------
proving that $S_i$ is a $2^{i+1}$-net in $S$ whenever $i\ge 0$.
%@@@@@@@@@@@@@@@@@@@@@@@@@@@@@@@@@@@@@@@@@@@@@@@@@@@@@@@@@@
\par Let us construct sets $S_i$ for $i<0$. To this end we let $S_{-1}$ denote a {\it maximal $2^{-1}$-net in $S$ containing $S_0$.} Then the set $S_{-1}$ is a $2^{-1}$-net in $S$ and its points are $2^{-1}$-separable. At the same way we construct a set $S_{-2}$ as a maximal $2^{-2}$-net in $S$ containing $S_1$. We continue this inductive procedure and in this way we obtain the required sequence of sets $\{S_i\}$ whenever $i<0.$
%@@@@@@@@@@@@@@@@@@@@@@@@@@@@@@@@@@@@@@@@@@@@@@@@@@@@@@@@@@
\par The lemma is proved.\bx\medskip
%----------------------------------------------------------
%@@@@@@@@@@@@@@@@@@@@@@@@@@@@@@@@@@@@@@@@@@@@@@@@@@@@@@@@@@
%@@@@@@@@@@@@@@@@@@@@@@@@@@@@@@@@@@@@@@@@@@@@@@@@@@@@@@@@@@
%@@@@@@@@@@@@@@@@@@@@@@@@@@@@@@@@@@@@@@@@@@@@@@@@@@@@@@@@@@
%@@@@@@@@@@@@@@@@@@@@@@@@@@@@@@@@@@@@@@@@@@@@@@@@@@@@@@@@@@
%@@@@@@@@@@@@@@@@@@@@@@@@@@@@@@@@@@@@@@@@@@@@@@@@@@@@@@@@@@
%@@@@@@@@@@@@@@@@@@@@@@@@@@@@@@@@@@@@@@@@@@@@@@@@@@@@@@@@@@
%@@@@@@@@@@@@@@@@@@@@@@@@@@@@@@@@@@@@@@@@@@@@@@@@@@@@@@@@@@
\begin{proposition}\lbl{WN} Let $S$ be a closed subset of $\RN$ and let $\Qc$ be a family of non-overlapping cubes in $\RN$. Suppose that there exist constants $\eta,\theta\ge 1$ such that for every cube $Q\in\Qc$ the following inequality
%----------------------------------------------------------
\bel{DQ}
\diam Q\le \eta\diam(\theta Q\cap S)
\ee
%----------------------------------------------------------
holds.
%----------------------------------------------------------
\par Then there exists a mapping $\Qc\ni Q\mapsto w_Q\in S$ such that
%----------------------------------------------------------
\par (i). $w_Q\in (2\theta)Q\cap S$ for every $Q\in\Qc$;
%----------------------------------------------------------
%@@@@@@@@@@@@@@@@@@@@@@@@@@@@@@@@@@@@@@@@@@@@@@@@@@@@@@@@@@
%----------------------------------------------------------
\par (ii). for every $a\in S$ we have
%----------------------------------------------------------
$$
\card\{Q\in\Qc: w_Q=a\}\le\gamma.
$$
%----------------------------------------------------------
Here $\gamma$ is a constant depending only on $n,\eta$ and $\theta$.

%----------------------------------------------------------
\end{proposition}
%----------------------------------------------------------
%@@@@@@@@@@@@@@@@@@@@@@@@@@@@@@@@@@@@@@@@@@@@@@@@@@@@@@@@@@
%@@@@@@@@@@@@@@@@@@@@@@@@@@@@@@@@@@@@@@@@@@@@@@@@@@@@@@@@@@
%@@@@@@@@@@@@@@@@@@@@@@@@@@@@@@@@@@@@@@@@@@@@@@@@@@@@@@@@@@
\par {\it Proof.} Let $Q\in\Qc$ and let $i_Q\in\mZ$ be such an integer that
%----------------------------------------------------------
\bel{IQ}
2^{i_Q}<\diam(\theta Q\cap S) \le 2^{i_Q+1}.
\ee
%----------------------------------------------------------
\par Let $\{S_i\}_{i\in\mZ}$ be a family of subsets of $S$ satisfying conditions of Lemma \reff{NET-S}. Let $a_Q,b_Q$ be points on $S$ such that $a_Q,b_Q\in (\theta Q)\cap S$ and
%----------------------------------------------------------
$$
\|a_Q-b_Q\|=\diam(\theta Q\cap S).
$$
%----------------------------------------------------------
Then, by \rf{IQ},
%----------------------------------------------------------
\bel{ABQ}
2^{i_Q}<\|a_Q-b_Q\|\le 2^{i_Q+1}.
\ee
%----------------------------------------------------------
By Lemma \reff{NET-S}, the set $S_{i_Q-2}$ is a $2^{i_Q-1}$-net in $S$, so that there exist points $A_Q,B_Q\in S_{i_Q-2}$ such that
%----------------------------------------------------------
\bel{ABab}
\|A_Q-a_Q\|\le 2^{i_Q-1}~~~\text{and}~~~
\|B_Q-b_Q\|\le 2^{i_Q-1}.
\ee
%----------------------------------------------------------
Prove that $A_Q\ne B_Q$ and
%----------------------------------------------------------
\bel{AB}
\{A_Q,B_Q\}\cap (S_{i_Q-2}\setminus S_{i_Q+2})\ne \emp.
\ee
%----------------------------------------------------------
\par In fact, by \rf{ABab},
%----------------------------------------------------------
$$
\|a_Q-b_Q\|\le \|a_Q-A_Q\|+\|A_Q-B_Q\|+\|B_Q-b_Q\|\le 2^{i_Q-1}+\|A_Q-B_Q\|+2^{i_Q-1}
$$
%----------------------------------------------------------
so that
%----------------------------------------------------------
$$
\|A_Q-B_Q\|\ge \|a_Q-b_Q\|- 2^{i_Q}.
$$
%----------------------------------------------------------
By \rf{IQ}, $\|a_Q-b_Q\|> 2^{i_Q}$ so that $\|A_Q-B_Q\|>0$ proving that $A_Q\ne B_Q$.\medskip
%----------------------------------------------------------
\par Prove that $\{A_Q,B_Q\} \nsubseteq S_{i_Q+2}$. If $\{A_Q,B_Q\} \subset S_{i_Q+2}$, then, by property (i) of Lemma \reff{NET-S}, $A_Q$ and $B_Q$ are $2^{i_Q+2}$-separable so that $\|A_Q-B_Q\|\ge 2^{i_Q+2}$.
On the other hand, by \rf{ABQ} and \rf{ABab},
%----------------------------------------------------------
$$
\|A_Q-B_Q\|\le \|A_Q-a_Q\|+\|a_Q-b_Q\|+\|b_Q-B_Q\|\le 2^{i_Q-1}+2^{i_Q+1}+2^{i_Q-1}=3\cdot 2^{i_Q}<2^{i_Q+2},
$$
%----------------------------------------------------------
a contradiction. Since $A_Q,B_Q\in S_{i_Q-2}$, the property \rf{AB} follows.
%----------------------------------------------------------
\par Now we are in a position to define a mapping $w:\Qc\to S$ satisfying conditions (i) and (ii) of the proposition.
Given a cube $Q\in\Qc$ by $w_Q$ we denote a point from the set $\{A_Q,B_Q\}$  which belongs to the set $S_{i_Q-2}\setminus S_{i_Q+2}$ (thus either $w_Q=A_Q$ or $w_Q=B_Q$). By \rf{AB}, such a point exists.
%----------------------------------------------------------
\par Thus $w_Q\in S_{i_Q-2}\setminus S_{i_Q+2}$. Prove that $w_Q\in (2\theta) Q\cap S$, i.e., condition (i) of the proposition is satisfied. In fact, by \rf{IQ},
%----------------------------------------------------------
$$
2^{i_Q}\le\diam(\theta Q\cap S)\le \diam(\theta Q)=\theta\diam Q.
$$
%----------------------------------------------------------
Combining this inequality with \rf{ABab}, we obtain
%----------------------------------------------------------
$$
\|A_Q-a_Q\|\le 2^{i_Q-1}\le \tfrac{1}{2}\theta\diam Q.
$$
%----------------------------------------------------------
Since $a_Q\in\theta Q$, this implies that
$A_Q\in (2\theta) Q$. In the same way we prove that
$B_Q\in (2\theta) Q$.
%----------------------------------------------------------
\par Since $w_Q$ coincides either with $A_Q$ or with $B_Q$, we conclude that $w_Q\in (2\theta) Q$ as well. Since $w_Q\in S$, the property (i) of the proposition follows.\medskip
%----------------------------------------------------------
\par Prove the property (ii) of the proposition. Let $a\in S$. Suppose that $Q,K\in\Qc$ and $w_Q=w_K=a$. Recall that
%----------------------------------------------------------
$$
a=w_Q\in S_{i_Q-2}\setminus S_{i_Q+2}~~~\text{and}~~~a=w_K\in S_{i_K-2}\setminus S_{i_K+2}.
$$
%----------------------------------------------------------
Since $a\notin S_{i_Q+2}$, $a\in S_{i_K-2}$ and $\{S_i\}$ is a decreasing sequence of sets, we conclude that
%----------------------------------------------------------
$$
i_K-2<i_Q+2
$$
%----------------------------------------------------------
so that $i_K<i_Q+4.$ In the same fashion we prove that $i_Q<i_K+4$ so that
%----------------------------------------------------------
\bel{IKQ}
i_Q-4<i_K<i_Q+4.
\ee
%----------------------------------------------------------
But, by \rf{IQ} and \rf{DQ},
%----------------------------------------------------------
$$
\diam Q\le \eta\diam(\theta Q\cap S) \le \eta 2^{i_Q+1}
$$
%----------------------------------------------------------
and
%----------------------------------------------------------
$$
2^{i_Q}\le \diam(\theta Q\cap S) \le \theta\diam Q.
$$
%----------------------------------------------------------
Hence
%----------------------------------------------------------
$$
(2\eta)^{-1}\diam Q\le 2^{i_Q}\le \theta\diam Q.
$$
%----------------------------------------------------------
The same is true for the cube $K$, i.e.,
%----------------------------------------------------------
$$
(2\eta)^{-1}\diam K\le 2^{i_K}\le \theta\diam K.
$$
%----------------------------------------------------------
Combining these inequalities with \rf{IKQ}, we obtain
%----------------------------------------------------------
$$
(2\eta)^{-1}\diam K\le 2^{i_K}\le 16\cdot\theta\diam Q
$$
%----------------------------------------------------------
proving that
%----------------------------------------------------------
\bel{DKQ}
\diam K\le \gamma_1\diam Q
\ee
%----------------------------------------------------------
with $\gamma_1:=32\,\eta\,\theta$. In the same way we show that %----------------------------------------------------------
\bel{DQK}
\diam Q\le \gamma_1\diam K.
\ee
%----------------------------------------------------------
\par Let
%----------------------------------------------------------
$$
\Qc_a:=\{Q\in\Qc:w_Q=a\}.
$$
%----------------------------------------------------------
By \rf{DKQ} and \rf{DQK}, $|Q'|\sim|Q''|$ for every $Q',Q''\in\Qc_a$ and constants of this equivalence depend only on $\eta$ and $\theta$. Furthermore, $(2\theta)Q\ni a$ for each $Q\in\Qc_a$. Since the cubes of the family $\Qc_a$ are non-overlapping, we conclude that
$\#\Qc_a\le C(n,\eta,\theta)$.
%----------------------------------------------------------
\par The proposition is completely proved.\bx\medskip
%----------------------------------------------------------
%@@@@@@@@@@@@@@@@@@@@@@@@@@@@@@@@@@@@@@@@@@@@@@@@@@@@@@@@@@
%@@@@@@@@@@@@@@@@@@@@@@@@@@@@@@@@@@@@@@@@@@@@@@@@@@@@@@@@@@
%@@@@@@@@@@@@@@@@@@@@@@@@@@@@@@@@@@@@@@@@@@@@@@@@@@@@@@@@@@
%@@@@@@@@@@@@@@@@@@@@@@@@@@@@@@@@@@@@@@@@@@@@@@@@@@@@@@@@@@
%@@@@@@@@@@@@@@@@@@@@@@@@@@@@@@@@@@@@@@@@@@@@@@@@@@@@@@@@@@
%@@@@@@@@@@@@@@@@@@@@@@@@@@@@@@@@@@@@@@@@@@@@@@@@@@@@@@@@@@
%----------------------------------------------------------
\begin{proposition}\lbl{L-PE} There exist an absolute constant $\gamma>0$ and a mapping
%----------------------------------------------------------
$$
\LE\ni L~~\longrightarrow~~\PRL(L)\in E
$$
%----------------------------------------------------------
such that:
%----------------------------------------------------------
\par (i). For every lacuna $L\in\LE$ we have
%----------------------------------------------------------
$$
\PRL(L)\in (\gamma\,Q_L)\cap E~;
$$
%----------------------------------------------------------
\par (ii). For every $a\in E$
%----------------------------------------------------------
$$
\card\{L\in\LE:\PRL(L)=a\}\le C(n).
$$
%----------------------------------------------------------
%@@@@@@@@@@@@@@@@@@@@@@@@@@@@@@@@@@@@@@@@@@@@@@@@@@@@@@@@@@
%----------------------------------------------------------
\end{proposition}
%@@@@@@@@@@@@@@@@@@@@@@@@@@@@@@@@@@@@@@@@@@@@@@@@@@@@@@@@@@
%@@@@@@@@@@@@@@@@@@@@@@@@@@@@@@@@@@@@@@@@@@@@@@@@@@@@@@@@@@
\par {\it Proof.} We let $\Qc_E:=\{Q_L: L\in\LE\}$ denote the family of all ``generalized" cubes $Q_L$ of the set $E$. We use the word ``generalized" to emphasize that if $\diam Q_L>0$, then $Q_L\in W_E$ is a Whitney cube, while $Q_L$ is a point of $E$ provided $\diam Q_L=0$. Thus $\Qc_E$ is a subset of the set $W_E\cup E$.
%----------------------------------------------------------
\par We put
%----------------------------------------------------------
\bel{DQZ}
\PRL(L):=Q_L~~~~\text{if}~~~~\diam Q_L=0.
\ee
%----------------------------------------------------------
%@@@@@@@@@@@@@@@@@@@@@@@@@@@@@@@@@@@@@@@@@@@@@@@@@@@@@@@@@@
\par Let
%----------------------------------------------------------
$$
\widetilde{\Qc}_E:=\{Q_L: L\in \LE,~\diam Q_L>0\}.
$$
%----------------------------------------------------------
Since $\widetilde{\Qc}_E\subset W_E$, this family consists of non-overlapping cubes. By Corollary \reff{D-QL}, every cube $Q\in\widetilde{\Qc}_E$ satisfies inequality \rf{DQ} with $\theta=90.$  Hence, by Proposition \reff{WN}, there exists a mapping
%----------------------------------------------------------
$$
\widetilde{\Qc}_E\ni Q~\longrightarrow~~w_Q\in E
$$
%----------------------------------------------------------
such that
%----------------------------------------------------------
$$
w_Q\in (\gamma\,Q)\cap E,~~~Q\in\widetilde{\Qc}_E,
$$
%----------------------------------------------------------
with $\gamma=2\,\theta=180$. Furthermore, for every $a\in E$
%----------------------------------------------------------
$$
\card\{Q\in\widetilde{\Qc}_E: w_Q=a\}\le C(n).
$$
%----------------------------------------------------------
\par We define $\PRL(L)$ by letting $\PRL(L):=w_{Q_L}$. Then, by the latter inequality and \rf{DQZ},
%----------------------------------------------------------
$$
\card\{L\in\LE: \PRL(L)=a\}\le C(n)+1,
$$
%----------------------------------------------------------
and the proposition follows.\bx\medskip
%----------------------------------------------------------
%@@@@@@@@@@@@@@@@@@@@@@@@@@@@@@@@@@@@@@@@@@@@@@@@@@@@@@@@@@
%@@@@@@@@@@@@@@@@@@@@@@@@@@@@@@@@@@@@@@@@@@@@@@@@@@@@@@@@@@
%@@@@@@@@@@@@@@@@@@@@@@@@@@@@@@@@@@@@@@@@@@@@@@@@@@@@@@@@@@
%@@@@@@@@@@@@@@@@@@@@@@@@@@@@@@@@@@@@@@@@@@@@@@@@@@@@@@@@@@
%@@@@@@@@@@@@@@@@@@@@@@@@@@@@@@@@@@@@@@@@@@@@@@@@@@@@@@@@@@
%@@@@@@@@@@@@@@@@@@@@@@@@@@@@@@@@@@@@@@@@@@@@@@@@@@@@@@@@@@
%----------------------------------------------------------
\begin{corollary}\lbl{CR-PE} Let $E$ be a finite subset of $\RN$. Then the number of its lacunae does not exceed $C(n)\# E$.
%----------------------------------------------------------
%@@@@@@@@@@@@@@@@@@@@@@@@@@@@@@@@@@@@@@@@@@@@@@@@@@@@@@@@@@
%----------------------------------------------------------
\end{corollary}
%@@@@@@@@@@@@@@@@@@@@@@@@@@@@@@@@@@@@@@@@@@@@@@@@@@@@@@@@@@
%@@@@@@@@@@@@@@@@@@@@@@@@@@@@@@@@@@@@@@@@@@@@@@@@@@@@@@@@@@
\par {\it Proof.} By part (ii) of Proposition \reff{L-PE}, for each $a\in E$ we have
%----------------------------------------------------------
$$
\card\{L\in\LE:\PRL(L)=a\}\le C(n).
$$
%----------------------------------------------------------
Hence
%----------------------------------------------------------
$$
\#\LE\le \sum_{a\in E}\card\{L\in\LE:\PRL(L)=a\}\le C(n)\# E.\BX
$$
%----------------------------------------------------------
\medskip
%----------------------------------------------------------
%@@@@@@@@@@@@@@@@@@@@@@@@@@@@@@@@@@@@@@@@@@@@@@@@@@@@@@@@@@
%@@@@@@@@@@@@@@@@@@@@@@@@@@@@@@@@@@@@@@@@@@@@@@@@@@@@@@@@@@
%@@@@@@@@@@@@@@@@@@@@@@@@@@@@@@@@@@@@@@@@@@@@@@@@@@@@@@@@@@
%@@@@@@@@@@@@@@@@@@@@@@@@@@@@@@@@@@@@@@@@@@@@@@@@@@@@@@@@@@
%@@@@@@@@@@@@@@@@@@@@@@@@@@@@@@@@@@@@@@@@@@@@@@@@@@@@@@@@@@
%@@@@@@@@@@@@@@@@@@@@@@@@@@@@@@@@@@@@@@@@@@@@@@@@@@@@@@@@@@
%@@@@@@@@@@@@@@@@@@@@@@@@@@@@@@@@@@@@@@@@@@@@@@@@@@@@@@@@@@
\par {\bf 5.2. A graph of lacunae and interior bridges.}
%----------------------------------------------------------
\addtocontents{toc}{~~~~5.2. A graph of lacunae and interior bridges. \hfill \thepage\\\par}
%----------------------------------------------------------
Let $L\in\LE$ be a lacuna. Recall that
%----------------------------------------------------------
$$
U_L=\cup\{Q: Q\in L\}.
$$
%----------------------------------------------------------
%@@@@@@@@@@@@@@@@@@@@@@@@@@@@@@@@@@@@@@@@@@@@@@@@@@@@@@@@@@
%@@@@@@@@@@@@@@@@@@@@@@@@@@@@@@@@@@@@@@@@@@@@@@@@@@@@@@@@@@
%@@@@@@@@@@@@@@@@@@@@@@@@@@@@@@@@@@@@@@@@@@@@@@@@@@@@@@@@@@
%@@@@@@@@@@@@@@@@@@@@@@@@@@@@@@@@@@@@@@@@@@@@@@@@@@@@@@@@@@
%----------------------------------------------------------
\begin{definition} \lbl{D-CNL}{\em Let $L,L'\in\LE$ be lacunae.  We say that $L$ and $L'$ are {\it contacting
lacunae} if $U_L\cap U_{L'}\ne\emp.$ In this case we write $L\leftrightarrow L'$.}
%----------------------------------------------------------
\end{definition}
%----------------------------------------------------------
%@@@@@@@@@@@@@@@@@@@@@@@@@@@@@@@@@@@@@@@@@@@@@@@@@@@@@@@@@@
%@@@@@@@@@@@@@@@@@@@@@@@@@@@@@@@@@@@@@@@@@@@@@@@@@@@@@@@@@@
%----------------------------------------------------------
\par Thus $L\leftrightarrow L'$ whenever there exist cubes $Q\in L$ and $Q'\in L'$ such that $Q\cap Q'\ne\emp$. We refer to the pair of such cubes as {\it contacting cubes}. Let us present several properties of contacting lacunae and contacting cubes which directly follow from results of the previous subsections.
%----------------------------------------------------------
%@@@@@@@@@@@@@@@@@@@@@@@@@@@@@@@@@@@@@@@@@@@@@@@@@@@@@@@@@@
%@@@@@@@@@@@@@@@@@@@@@@@@@@@@@@@@@@@@@@@@@@@@@@@@@@@@@@@@@@
%@@@@@@@@@@@@@@@@@@@@@@@@@@@@@@@@@@@@@@@@@@@@@@@@@@@@@@@@@@
%@@@@@@@@@@@@@@@@@@@@@@@@@@@@@@@@@@@@@@@@@@@@@@@@@@@@@@@@@@
%----------------------------------------------------------
\begin{proposition}\lbl{PRL-1} (i). Every lacuna $L\in\LE$ contacts with at most $C(n)$ lacunae, i.e.,
%----------------------------------------------------------
$$
\card\{L'\in\LE: L'\lr L\}\le C(n);
$$
%----------------------------------------------------------
\par (ii). Every true lacuna contacts only with elementary lacunae.
%----------------------------------------------------------
\end{proposition}
%----------------------------------------------------------
%@@@@@@@@@@@@@@@@@@@@@@@@@@@@@@@@@@@@@@@@@@@@@@@@@@@@@@@@@@
%@@@@@@@@@@@@@@@@@@@@@@@@@@@@@@@@@@@@@@@@@@@@@@@@@@@@@@@@@@
%----------------------------------------------------------
\par Clearly, part (i) of the proposition follows from Proposition \reff{INT-L1}, and part (ii) follows from part (i) of Proposition \reff{M-LAC}.
%----------------------------------------------------------
\par In turn, part (ii) of Proposition \reff{M-LAC} imply the following
%----------------------------------------------------------
%@@@@@@@@@@@@@@@@@@@@@@@@@@@@@@@@@@@@@@@@@@@@@@@@@@@@@@@@@@
%@@@@@@@@@@@@@@@@@@@@@@@@@@@@@@@@@@@@@@@@@@@@@@@@@@@@@@@@@@
%@@@@@@@@@@@@@@@@@@@@@@@@@@@@@@@@@@@@@@@@@@@@@@@@@@@@@@@@@@
%@@@@@@@@@@@@@@@@@@@@@@@@@@@@@@@@@@@@@@@@@@@@@@@@@@@@@@@@@@
%----------------------------------------------------------
\begin{proposition}\lbl{PRL-2} Let $L\in\LE$ be a lacuna and let $Q\in L$ be a contacting cube. (I.e., there exist a lacuna $L'\in\LE$ and a cube $Q'\in L'$ such that $Q\cap Q'\ne\emp$.) Then either
%----------------------------------------------------------
$$
\diam Q\sim \diam V_L\sim
\min\{\diam K: K\in L\}=\diam Q_L
$$
%----------------------------------------------------------
or
%----------------------------------------------------------
$$
\diam Q\sim \dist(V_L,E\setminus V_L)\sim
\max\{\diam K: K\in L\}=\diam \QL
$$
%----------------------------------------------------------
with absolute constants in the equivalences.
%@@@@@@@@@@@@@@@@@@@@@@@@@@@@@@@@@@@@@@@@@@@@@@@@@@@@@@@@@@
%----------------------------------------------------------
\end{proposition}
%----------------------------------------------------------
%@@@@@@@@@@@@@@@@@@@@@@@@@@@@@@@@@@@@@@@@@@@@@@@@@@@@@@@@@@
%@@@@@@@@@@@@@@@@@@@@@@@@@@@@@@@@@@@@@@@@@@@@@@@@@@@@@@@@@@
%----------------------------------------------------------
\par The relation $\lr$ determines a certain graph structure on $\LE$. We denote this graph by $\GE$. Thus two vertices of this graph, i.e., two lacunae $L,L'\in\LE$, are joined by an edge if they are contacting lacunae $\Longleftrightarrow$ there exist Whitney cubes $Q\in L$ and  $Q'\in L'$ such that $Q\cap Q'\ne\emp$.
%----------------------------------------------------------
\par By part (i) of Proposition \reff{PRL-1}, the degree of every vertex of the graph $\GE$ is bounded by a constant $C=C(n)$. In turn, by Corollary \reff{CR-PE}, for every finite set $E\subset\RN$ the number of vertices of $\GE$
does not exceed the number of points of $E$ (up to a multiplicative constant $C=C(n)$).
%----------------------------------------------------------
\par We turn to definition of an interior bridge of a lacuna.
%----------------------------------------------------------
%@@@@@@@@@@@@@@@@@@@@@@@@@@@@@@@@@@@@@@@@@@@@@@@@@@@@@@@@@@
%@@@@@@@@@@@@@@@@@@@@@@@@@@@@@@@@@@@@@@@@@@@@@@@@@@@@@@@@@@
%@@@@@@@@@@@@@@@@@@@@@@@@@@@@@@@@@@@@@@@@@@@@@@@@@@@@@@@@@@
%@@@@@@@@@@@@@@@@@@@@@@@@@@@@@@@@@@@@@@@@@@@@@@@@@@@@@@@@@@
%----------------------------------------------------------
\begin{definition}\lbl{FR-L} {\em Let $L\in\LE$ be a lacuna. We define two points $A_L,B_L\in E$ as follows:
%----------------------------------------------------------
{\it If $\diam V_L>0$, then we choose $A_L,B_L\in V_L$
to be any pair of distinct points in $V_L$ such that
%----------------------------------------------------------
$$
\|A_L-B_L\|=\diam V_L.
$$
%----------------------------------------------------------
\par If $\diam V_L=0$, i.e., $V_L$ is a single point, we set $A_L=V_L$. We choose $B_L$ to be some point in $E\setminus \{A_L\}$ whose distance from $A_L$  is minimal.}
%----------------------------------------------------------
\par We refer to the ordered couple $\Br(L)=(A_L,B_L)$ as an {\it interior bridge} of the lacuna $L$.}
%----------------------------------------------------------
\end{definition}
%----------------------------------------------------------
\medskip
%@@@@@@@@@@@@@@@@@@@@@@@@@@@@@@@@@@@@@@@@@@@@@@@@@@@@@@@@@@
%----------------------------------------------------------
\par Given lacunae $L,L'\in\LE$ we put
%----------------------------------------------------------
$$
d(L,L'):=\diam\{A_L,B_L,A_{L'},B_{L'}\}.
$$
%----------------------------------------------------------
%@@@@@@@@@@@@@@@@@@@@@@@@@@@@@@@@@@@@@@@@@@@@@@@@@@@@@@@@@@
%@@@@@@@@@@@@@@@@@@@@@@@@@@@@@@@@@@@@@@@@@@@@@@@@@@@@@@@@@@
%@@@@@@@@@@@@@@@@@@@@@@@@@@@@@@@@@@@@@@@@@@@@@@@@@@@@@@@@@@
%@@@@@@@@@@@@@@@@@@@@@@@@@@@@@@@@@@@@@@@@@@@@@@@@@@@@@@@@@@
%----------------------------------------------------------
\begin{proposition}\lbl{CL-Q} Let $L,L'\in\LE$ be contacting lacunae and let $Q\in L$ and $Q'\in L'$ be contacting cubes (i.e., $Q\cap Q'\ne\emp$). Then
%----------------------------------------------------------
$$
d(L,L')\sim \diam Q\sim \diam Q'
$$
%----------------------------------------------------------
with absolute constants in the equivalences. Furthermore,
%----------------------------------------------------------
$$
\{A_L,B_L,A_{L'},B_{L'}\}\subset (\gamma Q)\cap (\gamma Q')
$$
%----------------------------------------------------------
where $\gamma>0$ is an absolute constant.
%@@@@@@@@@@@@@@@@@@@@@@@@@@@@@@@@@@@@@@@@@@@@@@@@@@@@@@@@@@
%----------------------------------------------------------
\end{proposition}
%----------------------------------------------------------
%@@@@@@@@@@@@@@@@@@@@@@@@@@@@@@@@@@@@@@@@@@@@@@@@@@@@@@@@@@
%@@@@@@@@@@@@@@@@@@@@@@@@@@@@@@@@@@@@@@@@@@@@@@@@@@@@@@@@@@
%----------------------------------------------------------
\par {\it Proof.} By Lemma \reff{Wadd}, $\diam Q\sim \diam Q'$. Since $Q\cap Q'\ne\emp$, for some absolute constant $\gamma_1>0$ we have
%----------------------------------------------------------
\bel{Q-PG1}
Q\subset \gamma_1 Q'~~~\text{and}~~~Q'\subset \gamma_1 Q.
\ee
%----------------------------------------------------------
\par Recall that $V_L=(\q K)\cap E$ for each cube $K\in L$, see \rf{D-VL}, so that
%----------------------------------------------------------
$$
V_L\subset \q Q~~~\text{and}~~~V_{L'}\subset \q Q'.
$$
%----------------------------------------------------------
Combining this with \rf{Q-PG1} we obtain
%----------------------------------------------------------
\bel{U-LP}
V_L\cup V_{L'}\subset (\gamma_2 Q)\cap(\gamma_2 Q').
\ee
%----------------------------------------------------------
\par Prove that
%----------------------------------------------------------
\bel{IN-LF}
\{A_L,B_L,A_{L'},B_{L'}\}\subset
(\gamma_3 Q)\cap(\gamma_3 Q').
\ee
%----------------------------------------------------------
\par First we note that, by Proposition \reff{M-LAC}, either $L$ or $L'$ is an elementary lacuna. Suppose that $L'$ is such a lacuna, i.e., $L'\in\hL$. Then $L'=\{Q'\}$
and, by \rf{A-DL},
%----------------------------------------------------------
$$
\diam V_{L'}=\diam ((\q Q')\cap E)\sim \diam Q'.
$$
%----------------------------------------------------------
Since $\diam Q'>0$, by Definition \reff{FR-L}, $A_{L'},B_{L'}\in V_{L'}$ and
%----------------------------------------------------------
$$
\|A_{L'}-B_{L'}\|=\diam V_{L'}
$$
%----------------------------------------------------------
so that
%----------------------------------------------------------
\bel{DF-LP}
\|A_{L'}-B_{L'}\|\sim \diam Q'.
\ee
%----------------------------------------------------------
\par Now let us consider the lacuna $L$. If $\diam V_{L}>0$, then again, by Definition \reff{FR-L},
$A_{L},B_{L}\in V_{L}$ so that, by \rf{U-LP}, the inclusion \rf{IN-LF} is satisfied.
%----------------------------------------------------------
\par Suppose that $\diam V_{L}=0$, i.e., $V_L=\{A_L\}$. In this case, by Definition \reff{FR-L}, $B_L$ is a point nearest to $A_L$ on the set $E\setminus \{A_L\}$. Thus
%----------------------------------------------------------
\bel{E24}
\|A_{L}-B_{L}\|=\dist(A_L,E\setminus \{A_L\})=\dist(V_L,E\setminus V_L).
\ee
%----------------------------------------------------------
\par Recall that $Q$ is a contacting cube so that, by Proposition \reff{PRL-2}, either $ \diam Q\sim\diam V_{L}$ or
%----------------------------------------------------------
\bel{SO}
\diam Q\sim\dist(V_L,E\setminus V_L).
\ee
%----------------------------------------------------------
Since $\diam V_{L}=0$, we conclude that equivalence \rf{SO} holds. Combining this with equality \rf{E24}, we obtain
$\|A_{L}-B_{L}\|\sim \diam Q.$
Since $A_L\subset \q Q$, this implies that $B_L\subset \gamma_3 Q$ for a certain absolute constant $\gamma_3\ge\gamma_2$. This proves inclusion \rf{IN-LF}.
%----------------------------------------------------------
\par In particular, by this inclusion,
%----------------------------------------------------------
$$
d(L,L'):=\diam\{A_L,B_L,A_{L'},B_{L'}\}\le\gamma_3 \diam Q
$$
%----------------------------------------------------------
and $d(L,L')\le\gamma_3 \diam Q'$.
%----------------------------------------------------------
\par On the other hand, by \rf{DF-LP},
%----------------------------------------------------------
$$
\diam Q'\sim\|A_{L'}-B_{L'}\|\le \diam\{A_L,B_L,A_{L'},B_{L'}\}=d(L,L')
$$
%----------------------------------------------------------
proving that $d(L,L')\sim\diam Q\sim \diam Q'$.
%----------------------------------------------------------
\par The proof of the  proposition is complete.\bx
%----------------------------------------------------------
%@@@@@@@@@@@@@@@@@@@@@@@@@@@@@@@@@@@@@@@@@@@@@@@@@@@@@@@@@@
%@@@@@@@@@@@@@@@@@@@@@@@@@@@@@@@@@@@@@@@@@@@@@@@@@@@@@@@@@@
%@@@@@@@@@@@@@@@@@@@@@@@@@@@@@@@@@@@@@@@@@@@@@@@@@@@@@@@@@@
%@@@@@@@@@@@@@@@@@@@@@@@@@@@@@@@@@@@@@@@@@@@@@@@@@@@@@@@@@@
%----------------------------------------------------------
%@@@@@@@@@@@@@@@@@@@@@@@@@@@@@@@@@@@@@@@@@@@@@@@@@@@@@@@@@@
%@@@@@@@@@@@@@@@@@@@@@@@@@@@@@@@@@@@@@@@@@@@@@@@@@@@@@@@@@@
%@@@@@@@@@@@@@@@@@@@@@@@@@@@@@@@@@@@@@@@@@@@@@@@@@@@@@@@@@@
%@@@@@@@@@@@@@@@@@@@@@@@@@      @@@@@@@@@@@@@@@@@@@@@@@@@@@
%@@@@@@@@@@@@@@@@@@@@@@@          @@@@@@@@@@@@@@@@@@@@@@@@@
%@@@@@@@@@@@@@@@@@@@@@              @@@@@@@@@@@@@@@@@@@@@@@
%@@@@@@@@@@@@@@@@@@@     SECTION 6    @@@@@@@@@@@@@@@@@@@@@
%@@@@@@@@@@@@@@@@@@@@@              @@@@@@@@@@@@@@@@@@@@@@@
%@@@@@@@@@@@@@@@@@@@@@@@          @@@@@@@@@@@@@@@@@@@@@@@@@
%@@@@@@@@@@@@@@@@@@@@@@@@@      @@@@@@@@@@@@@@@@@@@@@@@@@@@
%@@@@@@@@@@@@@@@@@@@@@@@@@@@@@@@@@@@@@@@@@@@@@@@@@@@@@@@@@@
%@@@@@@@@@@@@@@@@@@@@@@@@@@@@@@@@@@@@@@@@@@@@@@@@@@@@@@@@@@
%@@@@@@@@@@@@@@@@@@@@@@@@@@@@@@@@@@@@@@@@@@@@@@@@@@@@@@@@@@
%----------------------------------------------------------
\SECT{6. Bridges between lacunae.}{6}
%----------------------------------------------------------
\addtocontents{toc}{6. Bridges between lacunae.\hfill \thepage\\\par}
%----------------------------------------------------------
\indent
%@@@@@@@@@@@@@@@@@@@@@@@@@@@@@@@@@@@@@@@@@@@@@@@@@@@@@@@@@@
%----------------------------------------------------------
%@@@@@@@@@@@@@@@@@@@@@@@@@@@@@@@@@@@@@@@@@@@@@@@@@@@@@@@@@@
%@@@@@@@@@@@@@@@@@@@@@@@@@@@@@@@@@@@@@@@@@@@@@@@@@@@@@@@@@@
%@@@@@@@@@@@@@@@@@@@@@@@@@@@@@@@@@@@@@@@@@@@@@@@@@@@@@@@@@@
%@@@@@@@@@@@@@@@@@@@@@@@@@@@@@@@@@@@@@@@@@@@@@@@@@@@@@@@@@@
%----------------------------------------------------------
\par  In this section we study {\it bridges between contacting lacunae} which we have briefly described in Section 2. This important notion provides a useful tool for the study of a lacunary mo\-dification of the Whitney extension method which we introduce in the next section.
%----------------------------------------------------------
\par We define bridges in several steps. First of all let us slightly generalize the notion of the ``triangle'' introduced in the first section. We have defined a triangle as a subset $\Delta=\{z_1,z_2,z_3\}\subset\RT$ consisting of {three non-collinear points.} In what follows we refer to such a subset as a {\it ``true'' triangle}.
%----------------------------------------------------------
\par Let us also consider a subset $\Delta=\{z_1,z_2,z_3\}\subset\RT$ consisting of {\it three collinear points.} In this case we call $\Delta$ a {\it ``\dg'' triangle.} Thus a degenerate triangle is the "triangle" formed by three collinear points.
%----------------------------------------------------------
\par We refer to the points $z_1,z_2,z_3$ as the vertices of the triangle $\Delta=\{z_1,z_2,z_3\}$. Also given three points $A,B,C\in \RT$ we let  $\Delta\{A,B,C\}$ denote the triangle with vertices in these points.
%@@@@@@@@@@@@@@@@@@@@@@@@@@@@@@@@@@@@@@@@@@@@@@@@@@@@@@@@@@
%----------------------------------------------------------
\par We say that a side $[A,B]$ is {\it the smallest side} of the triangle $\Delta\{A,B,C\}$ if $$\|A-B\|<\min\{\|A-C\|,\|B-C\|\}.$$
%@@@@@@@@@@@@@@@@@@@@@@@@@@@@@@@@@@@@@@@@@@@@@@@@@@@@@@@@@@
%@@@@@@@@@@@@@@@@@@@@@@@@@@@@@@@@@@@@@@@@@@@@@@@@@@@@@@@@@@
%@@@@@@@@@@@@@@@@@@@@@@@@@@@@@@@@@@@@@@@@@@@@@@@@@@@@@@@@@@
%@@@@@@@@@@@@@@@@@@@@@@@@@@@@@@@@@@@@@@@@@@@@@@@@@@@@@@@@@@
%----------------------------------------------------------
\begin{definition}\lbl{BR-L} {\em Let $L,L'\in\LE$ be  contacting lacunae such that
%----------------------------------------------------------
\bel{DL-1}
\{A_L,B_L\}\cap\{A_{L'},B_{L'}\}=\emp.
\ee
%---------------------------------------------------------- %@@@@@@@@@@@@@@@@@@@@@@@@@@@@@@@@@@@@@@@@@@@@@@@@@@@@@@@@@@
\par We say that the (non-ordered) pair
%---------------------------------------------------------- %@@@@@@@@@@@@@@@@@@@@@@@@@@@@@@@@@@@@@@@@@@@@@@@@@@@@@@@@@@
$$
\Br(L,L')=\{C(L,L'), C(L',L)\}
$$
%---------------------------------------------------------- %@@@@@@@@@@@@@@@@@@@@@@@@@@@@@@@@@@@@@@@@@@@@@@@@@@@@@@@@@@
where $C(L,L')\in\{A_L,B_L\}$ and  $C(L',L)\in\{A_{L'},B_{L'}\}$, is an {\it exterior bridge} between lacunae $L$ and $L'$ if the following conditions are satisfied:\medskip
%----------------------------------------------------------
\par {\it In the triangle $\Delta\{A_L,B_L,C(L',L)\}$ the side of the triangle which is opposite to the vertex $C(L,L')$ is not the smallest side of the triangle.
%----------------------------------------------------------
\par A similar condition holds for the triangle $\Delta\{A_{L'},B_{L'},C(L,L')\}$: the side of the triangle which is opposite to the vertex $C(L',L)$ is not the smallest side of the triangle.}\medskip
%----------------------------------------------------------
\par We also say that the bridge $\Br(L,L')$ is {\it connected} to the bridge $\Br(L)$ and to the bridge $\Br(L')$.}
%----------------------------------------------------------
\end{definition}
%----------------------------------------------------------
%@@@@@@@@@@@@@@@@@@@@@@@@@@@@@@@@@@@@@@@@@@@@@@@@@@@@@@@@@@
%@@@@@@@@@@@@@@@@@@@@@@@@@@@@@@@@@@@@@@@@@@@@@@@@@@@@@@@@@@
%@@@@@@@@@@@@@@@@@@@@@@@@@@@@@@@@@@@@@@@@@@@@@@@@@@@@@@@@@@
%----------------------------------------------------------
\begin{remark} {\em In the sequel we identify the exterior bridges $\Br(L,L')$ and $\Br(L',L)$ between contacting lacunae $L$ and $L'$; thus  $\Br(L,L')=\Br(L',L)$ for every $L,L'\in \LE$ such that $L\lr L'.$\rbx}
%----------------------------------------------------------
\end{remark}
%----------------------------------------------------------
%@@@@@@@@@@@@@@@@@@@@@@@@@@@@@@@@@@@@@@@@@@@@@@@@@@@@@@@@@@
%@@@@@@@@@@@@@@@@@@@@@@@@@@@@@@@@@@@@@@@@@@@@@@@@@@@@@@@@@@
%@@@@@@@@@@@@@@@@@@@@@@@@@@@@@@@@@@@@@@@@@@@@@@@@@@@@@@@@@@
%@@@@@@@@@@@@@@@@@@@@@@@@@@@@@@@@@@@@@@@@@@@@@@@@@@@@@@@@@@
%----------------------------------------------------------
\begin{proposition} Let $L,L'\in\LE$ be contacting lacunae satisfying condition \rf{DL-1}. Then there exists an exterior bridge $\Br(L,L')=\{C(L,L'), C(L',L)\}$ between these lacunae. Furthermore,
%----------------------------------------------------------
$$
\diam \Delta\{A_L,B_L,C(L',L)\}\sim \diam \Delta\{A_{L'},B_{L'},C(L,L')\}
\sim \diam\{ A_L,B_L,A_{L'},B_{L'}\}
$$
%----------------------------------------------------------
with absolute constants in the equivalences.
%@@@@@@@@@@@@@@@@@@@@@@@@@@@@@@@@@@@@@@@@@@@@@@@@@@@@@@@@@@
%----------------------------------------------------------
\end{proposition}
%----------------------------------------------------------
%@@@@@@@@@@@@@@@@@@@@@@@@@@@@@@@@@@@@@@@@@@@@@@@@@@@@@@@@@@
%@@@@@@@@@@@@@@@@@@@@@@@@@@@@@@@@@@@@@@@@@@@@@@@@@@@@@@@@@@
%----------------------------------------------------------
\par {\it Proof.} Consider three cases.
%----------------------------------------------------------
\par {\it The first case.} Suppose that
%----------------------------------------------------------
\bel{1-C}
\|A_L-B_L\|=\min\{\|a-b\|:a,b\in \{A_{L},B_{L},A_{L'},B_{L'}\}, a\ne b\}.
\ee
%---------------------------------------------------------- %@@@@@@@@@@@@@@@@@@@@@@@@@@@@@@@@@@@@@@@@@@@@@@@@@@@@@@@@@@
We put $C(L,L')=A_L$. Let $C(L',L)\in \{A_{L'},B_{L'}\}$ be such a point that
%----------------------------------------------------------
\bel{2-B}
\|C(L',L)-A_L\|=\min\{\|A_L-A_{L'}\|,\|A_L-B_{L'}\|\}.
\ee
%------------------------------------------------------
Prove that the couple
%------------------------------------------------------
\bel{QD}
\Br(L,L')=\{C(L,L'),C(L',L)\}
\ee
%------------------------------------------------------
provides an exterior bridge between the lacunae $L$ and $L'$.
%----------------------------------------------------------
\par In fact, by \rf{1-C}, the side $[A_L,B_L]$ is a minimal side of the triangle $\Delta\{A_L,B_L,C(L',L)\}$ so that the side of this triangle which is opposite to the vertex $C(L,L')=A_L$ is not the smallest side.
%----------------------------------------------------------
\par On the other hand, by equality \rf{2-B}, the side of  the triangle $\Delta\{A_{L'},B_{L'},C(L,L')\}$ which is opposite to the vertex $C(L',L)$ is not the smallest because it is not smaller than the side $[C(L,L'),C(L',L)]=[A_L,C(L',L)]$.
%----------------------------------------------------------
\par These observations shows that the couple \rf{QD}
satisfies the conditions of Definition \reff{BR-L} providing an exterior bridge between these lacunae.
%----------------------------------------------------------
\par {\it The second case.} Let
%----------------------------------------------------------
\bel{2-C}
\|A_{L'}-B_{L'}\|=\min\{\|a-b\|:a,b\in \{A_{L},B_{L},A_{L'},B_{L'}\}, a\ne b\}.
\ee
%----------------------------------------------------------
We prove the existence of the corresponding bridge in the same fashion as in the previous case.
%----------------------------------------------------------
\par {\it The third case.} We assume that both condition \rf{1-C} and condition \rf{2-C} are not satisfied. In this case there exist points $C(L,L')\in\{A_L,B_L\}$ and  $C(L',L)\in\{A_{L'},B_{L'}\}$ such that
%----------------------------------------------------------
\bel{3-C}
\|C(L,L')-C(L',L)\|=\min\{\|a-b\|: a,b\in \{A_{L},B_{L},A_{L'},B_{L'}\}, a\ne b\}.
\ee
%----------------------------------------------------------
It can be readily seen that in this case the couple \rf{QD} is a an exterior bridge between $L$ and $L'$. In fact, consider the triangle $\Delta=\Delta\{A_L,B_L,C(L',L)\}$ and the vertex $C(L,L')\in\{A_L,B_L\}$. By \rf{3-C}, the side $[C(L,L'),C(L',L)]$ is a side of minimal length in $\Delta$. Therefore the side of $\Delta$ which is opposite to $C(L,L')$ is not the smallest side in this triangle.
%----------------------------------------------------------
\par In the same way we prove  that the side of the triangle $\Delta'=\Delta\{A_{L'},B_{L'},C(L,L')\}$ which is opposite to the vertex $C(L',L)$ is not the smallest side of $\Delta'$. The proof of the proposition is finished.\bx
%----------------------------------------------------------
%@@@@@@@@@@@@@@@@@@@@@@@@@@@@@@@@@@@@@@@@@@@@@@@@@@@@@@@@@@
%@@@@@@@@@@@@@@@@@@@@@@@@@@@@@@@@@@@@@@@@@@@@@@@@@@@@@@@@@@
%@@@@@@@@@@@@@@@@@@@@@@@@@@@@@@@@@@@@@@@@@@@@@@@@@@@@@@@@@@
%@@@@@@@@@@@@@@@@@@@@@@@@@@@@@@@@@@@@@@@@@@@@@@@@@@@@@@@@@@
%----------------------------------------------------------
%@@@@@@@@@@@@@@@@@@@@@@@@@@@@@@@@@@@@@@@@@@@@@@@@@@@@@@@@@@
%----------------------------------------------------------
%@@@@@@@@@@@@@@@@@@@@@@@@@@@@@@@@@@@@@@@@@@@@@@@@@@@@@@@@@@
%@@@@@@@@@@@@@@@@@@@@@@@@@@@@@@@@@@@@@@@@@@@@@@@@@@@@@@@@@@
%@@@@@@@@@@@@@@@@@@@@@@@@@@@@@@@@@@@@@@@@@@@@@@@@@@@@@@@@@@
%@@@@@@@@@@@@@@@@@@@@@@@@@@@@@@@@@@@@@@@@@@@@@@@@@@@@@@@@@@
\medskip
%----------------------------------------------------------
\par Let us consider two contacting lacunae
$L,L'\in\LE$, $L\lr L'$, and their interior bridges $\Br(L)=(A_L,B_L)$ and $\Br(L')=(A_{L'},B_{L'})$. Suppose that $\{A_L,B_L\}\ne\{A_{L'},B_{L'}\}$ but
%----------------------------------------------------------
\bel{CASE-2}
\{A_L,B_L\}\cap\{A_{L'},B_{L'}\}\ne\emp.
\ee
%----------------------------------------------------------
In other words two point sets $\{A_L,B_L\}$ and $\{A_{L'},B_{L'}\}$ are different and have a unique common point. We denote this point by $D(L,L')$. By $E(L)$ we denote the remaining point of the set $\{A_L,B_L\}$; thus $\{A_L,B_L\}=\{D(L,L'),E(L)\}$. In a similar way we define a point $E(L')$; thus $\{A_{L'},B_{L'}\}=\{D(L,L'),E(L')\}$.
%----------------------------------------------------------
\par Let us define external bridges for this case. Consider two cases.\bigskip
%----------------------------------------------------------
\par {\it The first case:} The side $[E(L),E(L')]$ is {\it the smallest} side of a triangle $\Delta$ with vertices in $E(L),E(L'),$ and $D(L',L))$. In this case we put
%----------------------------------------------------------
$$
C(L,L')=E(L),~~~C(L',L)=E(L')
$$
%----------------------------------------------------------
and define an exterior bridge
$\Br(L,L')=\{C(L,L'),C(L',L)\}$ between the lacunae $L$  and $L'$. Note that the bridge $\Br(L,L')$ satisfies the conditions of Definition \reff{BR-L}.
%----------------------------------------------------------
\par As before we say that {\it the bridge $\Br(L,L')$ is connected to the bridge $\Br(L)$ and to the bridge $\Br(L')$}.\bigskip
%----------------------------------------------------------
\par {\it The second case:} The side $[E(L),E(L')]$ is not {\it the smallest} side of the triangle $\Delta$. In this case {\it we do not introduce any additional bridges between $L$ and $L'$. We only say that the (interior) bridges $\Br(L)$ and $\Br(L')$ are connected.}\bigskip
%----------------------------------------------------------
\par It remains to consider the last case of a pair of contacting lacunae $L,L'\in\LE$, $L\lr L',$ such that  $\{A_L,B_L\}=\{A_{L'},B_{L'}\}$. Similar to the previous case we do not introduce any additional bridges between $L$ and $L'$, and refer to (interior) bridges $\Br(L)$ and $\Br(L')$ as {\it connected bridges}.\bigskip
%----------------------------------------------------------
\par We have defined all types of bridges (interior and exterior). We have also introduced the notion of connected bridges.
%----------------------------------------------------------
\par By $\BRE$ we denote the family of all bridges constructed for the set $E$. Given a bridge $\Br\in\BRE$ by $\A{T}$ and $\B{T}$ we denote its ends. Note that $\A{T}\ne \B{T}$ for every  bridge $\Br\in\BRE$. In particular, if $\Br=\Br(L)$ is an interior bridge of a lacuna $L\in\LE$, these points coincide with the points $A_L$ and $B_L$, see Definition \reff{FR-L}.
%----------------------------------------------------------
\par We have also introduced the notion of {\it connected bridges.} We denote connected bridges by the sign $\bcn$. %----------------------------------------------------------
\par Let us note several general properties of bridges.\bigskip
%----------------------------------------------------------
\par \textbullet~If two bridges $\Br$ and $\Br'$ are connected ($\Br\bcn \Br'$), then {\it one of them is an interior bridge.}
%----------------------------------------------------------
\par \textbullet~Interior bridges $\Br,\Br'\in\BRE$ are connected $(\Br\bcn \Br')$ provided they have the same ends, i.e., $\{\A{T},\B{T}\}=\{\A{T'},\B{T'}\}$.
%----------------------------------------------------------
\par \textbullet~If the ends of two connected bridges $\Br$ and $\Br'$ are different ($(\A{T},\B{T})\ne(\A{T'},\B{T'})$) the bridges have a unique common end which we denote by $D=D(\Br,\Br')$. Thus in this case
%----------------------------------------------------------
$$
\#\,\{\A{T},\B{T},\A{T'},\B{T'}\}=3
$$
%----------------------------------------------------------
so that these points form a triangle $\Delta$. This triangle possesses the following property: {\it the side of $\Delta$ which is opposite to the vertex $D$ is not the smallest side of this triangle.}\smallskip
%----------------------------------------------------------
\par In turn, this property implies the following one: Let $\alpha(D)$ be the angle corresponding to the vertex $D$. Then
%----------------------------------------------------------
\bel{SIN}
|\sin \alpha(D)|\sim\tfrac{1}{R_\Delta}\,\diam\Delta
=\cu_\Delta\,\diam\Delta
\ee
%----------------------------------------------------------
with absolute constants in the equivalences. In fact, using a formula from the elementary geometry, we obtain
%----------------------------------------------------------
$$
|\sin \alpha(D)|=\tfrac12\, \cu_\Delta\,\ell
$$
%----------------------------------------------------------
where $\ell$ is the Euclidean length of the side opposite to the vertex $D$. Since this side is not the smallest in $\Delta$, we have $\diam \Delta\sim \ell$ proving \rf{SIN}.\bigskip
%----------------------------------------------------------
\par Finally, by $\GBE$ we denote a graph whose vertices is the family $\BRE$ of all bridges (interior and exterior). In this graph two vertices (bridges) $\Br$ and $\Br'$ are joined by an edge if $\Br\bcn \Br'$, i.e., the bridges $\Br$ and $\Br'$ are connected. This graph possesses the following property: Let $L,L'\in \LE, L\lr L',$ be contacting lacunae. Then either their interior bridges
$\Br(L)$ and $\Br(L')$ are connected, or their exists an (exterior) bridge which is connected to both the bridge $\Br(L)$ and to the bridge $\Br(L')$.
%----------------------------------------------------------
\par The next proposition enables us to every couple of connected bridges to assign a certain Whitney cube which is ``close'' to the set of ends of these bridges.
%----------------------------------------------------------
\begin{proposition}\lbl{L-QLM} Let $\Br,\Br'\in\BRE$, $\Br\bcn\Br'$, be a pair of connected bridges with the ends at the points $\{\A{T},\B{T}\}$ and $\{\A{T'},\B{T'}\}$. Let the bridge $\Br$ be the interior bridge of a lacuna $L\in\LE$. Then there exists a cube $\hQ(\Br,\Br')\in W_E$ which coincides either with the cube $Q_L$ or with the cube $\QL$ such that
%----------------------------------------------------------
\bel{END-1}
\diam \{\A{T},\B{T},\A{T'},\B{T'}\} \sim \diam \hQ(\Br,\Br')
\ee
%----------------------------------------------------------
and
%----------------------------------------------------------
\bel{END-2}
\{\A{T},\B{T},\A{T'},\B{T'}\}\subset \gamma \hQ(\Br,\Br').
\ee
%----------------------------------------------------------
Here the constant $\gamma$ and the constants of the equivalence in \rf{END-1} are absolute.
%----------------------------------------------------------
%@@@@@@@@@@@@@@@@@@@@@@@@@@@@@@@@@@@@@@@@@@@@@@@@@@@@@@@@@@
%----------------------------------------------------------
\end{proposition}
%----------------------------------------------------------
%@@@@@@@@@@@@@@@@@@@@@@@@@@@@@@@@@@@@@@@@@@@@@@@@@@@@@@@@@@
%@@@@@@@@@@@@@@@@@@@@@@@@@@@@@@@@@@@@@@@@@@@@@@@@@@@@@@@@@@
%----------------------------------------------------------
\par {\it Proof.} Since $\Br$ is the interior bridge of the lacuna $L$, we have $A^{[T]}=A_L,B^{[T]}=B_L$. Consider two case:
%----------------------------------------------------------
\par {\it The first case.} Suppose that the bridge $\Br'$ is the interior bridge of a lacuna $L'\in\LE$ which contacts to $L$ ($L\lr L'$). In this case $\A{T'}=A_{L'},\B{T'}=B_{L'}$ and the sets of the ends $\{A_{L},B_{L}\}$ and $\{A_{L'},B_{L'}\}$ have a common point.
%----------------------------------------------------------
\par Let $Q\in L$ and $Q'\in L'$ be contacting cubes. Then, by Proposition \reff{PRL-2}, either
%----------------------------------------------------------
\bel{V-1}
\diam Q\sim \diam Q_L\,,
\ee
%----------------------------------------------------------
or
%----------------------------------------------------------
$$
\diam Q\sim\diam \QL
$$
%----------------------------------------------------------
with absolute constants in these equivalences.
%----------------------------------------------------------
\par On the other hand, by Proposition \reff{CL-Q},
%----------------------------------------------------------
\bel{7-R}
\diam \{\A{T},\B{T},\A{T'},\B{T'}\} \sim \diam Q
\ee
%----------------------------------------------------------
and %----------------------------------------------------------
\bel{8-R}
\{\A{T},\B{T},\A{T'},\B{T'}\}  \subset \gamma_1 Q
\ee
%----------------------------------------------------------
where $\gamma_1$ is an absolute constant.
%----------------------------------------------------------
\par Now we put $\hQ(\Br,\Br')=Q_L$ if \rf{V-1} is satisfied; otherwise we put $\hQ(\Br,\Br')=\QL$. In particular, $\hQ(\Br,\Br')\in L$.
%----------------------------------------------------------
\par In both cases
%----------------------------------------------------------
\bel{R-1}
\diam Q\sim \diam \hQ(\Br,\Br').
\ee
%----------------------------------------------------------
Combining this equivalence with \rf{7-R} we conclude that \rf{END-1} holds.
%----------------------------------------------------------
\par Prove \rf{END-2}. Since $Q,\hQ(\Br,\Br')\in L$,
%----------------------------------------------------------
$$
V_L=(90Q)\cap E=(90\hQ(\Br,\Br'))\cap E
$$
%----------------------------------------------------------
so that
%----------------------------------------------------------
$$
(90Q)\cap (90\hQ(\Br,\Br'))\ne\emp.
$$
%----------------------------------------------------------
This property and \rf{R-1} imply the inclusion $Q\subset \gamma_2\hQ(\Br,\Br')$ with a certain absolute constant $\gamma_2$. Combining this inclusion with \rf{8-R} we obtain \rf{END-2} proving the proposition in the case under consideration.\medskip
%----------------------------------------------------------
\par {\it The second case.} Let $\Br'$ be an exterior bridge which connects the lacuna $L$ to a contacting lacuna $L'\in\LE$. Let $Q\in L, Q'\in L'$ be corresponding contacting cubes (thus $Q\cap Q'\ne\emp$).
In this case the ends of the bridge $\Br'$ can be identified with the points $C(L,L')\in\{A_L,B_L\}$ and
$C(L',L)\in\{A_{L'},B_{L'}\}$. Thus
%----------------------------------------------------------
$$
\{\A{T},\B{T},\A{T'},\B{T'}\}=
\{A_{L},B_{L},C(L',L)\} .
$$
%----------------------------------------------------------
\par By Definition \reff{FR-L}, if $\diam V_L>0$, then
%----------------------------------------------------------
\bel{DL-P1}
A_L,B_L\in V_L~~~\text{and}~~~\|A_L-B_L\|=\diam V_L.
\ee
%----------------------------------------------------------
If $\diam V_L=0$, then
%----------------------------------------------------------
\bel{C-D0}
V_L=\{A_L\},~~~B_L\in E\setminus V_L,
\ee
%----------------------------------------------------------
and
%----------------------------------------------------------
\bel{J-D0}
\|A_L-B_L\|=\dist(V_L,E\setminus V_L).
\ee
%----------------------------------------------------------
\par Suppose that $\diam V_L>0$ and $C(L',L)\in V_L$. In this case
%----------------------------------------------------------
$$
\{\A{T},\B{T},\A{T'},\B{T'}\}
=\{A_{L},B_{L},C(L',L)\}\subset V_L
$$
%----------------------------------------------------------
so that, by \rf{DL-P1},
%----------------------------------------------------------
$$
\diam \{\A{T},\B{T},\A{T'},\B{T'}\}
=\diam \{A_{L},B_{L},C(L',L)\}=\diam V_L\sim \diam Q_L.
$$
%----------------------------------------------------------
Furthermore, since $Q_L\in L$, we have $V_L=(90Q_L)\cap E$ proving that $V_L\subset 90Q_L.$ Thus in this case we can put $\hQ(\Br,\Br'):=Q_L$.
%----------------------------------------------------------
\par Now suppose that $\diam V_L>0$ and $C(L',L)\in E\setminus V_L$. Then
%----------------------------------------------------------
\bel{U-7}
\diam \{\A{T},\B{T},\A{T'},\B{T'}\}
=\diam \{A_{L},B_{L},C(L',L)\}\ge \dist(V_L,E\setminus V_L).
\ee
%----------------------------------------------------------
By Proposition \reff{CL-MAX},
%----------------------------------------------------------
$$
\diam \QL\le \gamma \dist(V_L,E\setminus V_L)
$$
%----------------------------------------------------------
so that
%----------------------------------------------------------
$$
\diam \QL\le \gamma \diam \{\A{T},\B{T},\A{T'},\B{T'}\}.
$$
%----------------------------------------------------------
\par On the other hand, by Proposition \reff{CL-Q},
%----------------------------------------------------------
\be
\diam \{\A{T},\B{T},\A{T'},\B{T'}\}
&=&
\diam \{A_{L},B_{L},C(L',L)\}\le  \diam \{A_{L},B_{L},A_{L'},B_{L'}\}\nn\\&\le& \gamma \diam Q\le \gamma \diam \QL.\nn
\ee
%----------------------------------------------------------
(Recall that the cube $\QL$ has the maximal diameter in the lacuna $L$.) Hence
%----------------------------------------------------------
\bel{F-4}
\diam \{\A{T},\B{T},\A{T'},\B{T'}\}\sim \diam \QL.
\ee
%----------------------------------------------------------
Furthermore, by Proposition \reff{CL-Q},
%----------------------------------------------------------
$$
\{\A{T},\B{T},\A{T'},\B{T'}\}\subset \{A_{L},B_{L},A_{L'},B_{L'}\}\subset \gamma_1 Q.
$$
%----------------------------------------------------------
Since $Q,\QL\in L$, we have $(90 Q)\cap (90\QL)\ne\emp$. But $\diam Q\le\diam \QL$ so that $Q\subset \gamma_2 \QL$ for some absolute constant $\gamma_2$. Hence
%----------------------------------------------------------
\bel{F-5}
\{\A{T},\B{T},\A{T'},\B{T'}\}\subset \gamma_3\QL.
\ee
%----------------------------------------------------------
This enables us to put $\hQ(\Br,\Br'):=\QL.$
%----------------------------------------------------------
\par Consider the remaining case $\diam V_L=0$, see \rf{C-D0} and \rf{J-D0}. Since
%----------------------------------------------------------
$$
\|A_L-B_L\|=\dist(V_L,E\setminus V_L),
$$
%----------------------------------------------------------
inequality \rf{U-7} is satisfied. Now repeating the considerations of the previous case we show that equivalence \rf{F-4} and inclusion \rf{F-5} hold. This again enables us to put in this case $\hQ(\Br,\Br'):=\QL.$
%----------------------------------------------------------
\par The proof of the proposition is complete.\bx\medskip
%----------------------------------------------------------
%@@@@@@@@@@@@@@@@@@@@@@@@@@@@@@@@@@@@@@@@@@@@@@@@@@@@@@@@@@
%@@@@@@@@@@@@@@@@@@@@@@@@@@@@@@@@@@@@@@@@@@@@@@@@@@@@@@@@@@
%@@@@@@@@@@@@@@@@@@@@@@@@@@@@@@@@@@@@@@@@@@@@@@@@@@@@@@@@@@
%@@@@@@@@@@@@@@@@@@@@@@@@@@@@@@@@@@@@@@@@@@@@@@@@@@@@@@@@@@
%@@@@@@@@@@@@@@@@@@@@@@@@@@@@@@@@@@@@@@@@@@@@@@@@@@@@@@@@@@
%@@@@@@@@@@@@@@@@@@@@@@@@@@@@@@@@@@@@@@@@@@@@@@@@@@@@@@@@@@
%----------------------------------------------------------
\begin{definition}\lbl{WSP-C} {\em We say that $\Qc$ is a family of well-separated cubes if for every pair of cubes $Q,Q'\in \Qc,Q\ne Q',$ the following inequality
%----------------------------------------------------------
$$
\diam Q+\diam Q'\le \dist(Q,Q')
$$
%----------------------------------------------------------
holds.}
%----------------------------------------------------------
\end{definition}
%----------------------------------------------------------
%@@@@@@@@@@@@@@@@@@@@@@@@@@@@@@@@@@@@@@@@@@@@@@@@@@@@@@@@@@
%@@@@@@@@@@@@@@@@@@@@@@@@@@@@@@@@@@@@@@@@@@@@@@@@@@@@@@@@@@
%@@@@@@@@@@@@@@@@@@@@@@@@@@@@@@@@@@@@@@@@@@@@@@@@@@@@@@@@@@
%@@@@@@@@@@@@@@@@@@@@@@@@@@@@@@@@@@@@@@@@@@@@@@@@@@@@@@@@@@
%@@@@@@@@@@@@@@@@@@@@@@@@@@@@@@@@@@@@@@@@@@@@@@@@@@@@@@@@@@
%@@@@@@@@@@@@@@@@@@@@@@@@@@@@@@@@@@@@@@@@@@@@@@@@@@@@@@@@@@
%----------------------------------------------------------
\begin{proposition}\lbl{BR-TO-Q} There exists a family $\KE$ of well-separated cubes and a one-to-one mapping which to every pair of connected bridges $\Br,\Br'\in\BRE$, $\Br\bcn\Br',$ assigns a cube $K(\Br,\Br')\in\KE$ such that
%----------------------------------------------------------
\bel{MP-1}
\diam \{\A{T},\B{T},\A{T'},\B{T'}\} \sim \diam K(\Br,\Br')
\ee
%----------------------------------------------------------
and
%----------------------------------------------------------
\bel{MP-2}
\{\A{T},\B{T},\A{T'},\B{T'}\}\subset \gamma K(\Br,\Br').
\ee
%----------------------------------------------------------
Here the constant $\gamma$ and the constants of the equivalence in \rf{MP-1} depend only on $n$.
%----------------------------------------------------------
%@@@@@@@@@@@@@@@@@@@@@@@@@@@@@@@@@@@@@@@@@@@@@@@@@@@@@@@@@@
%----------------------------------------------------------
\end{proposition}
%----------------------------------------------------------
%@@@@@@@@@@@@@@@@@@@@@@@@@@@@@@@@@@@@@@@@@@@@@@@@@@@@@@@@@@
%@@@@@@@@@@@@@@@@@@@@@@@@@@@@@@@@@@@@@@@@@@@@@@@@@@@@@@@@@@
%----------------------------------------------------------
\par {\it Proof.} Let
%----------------------------------------------------------
$$
\Mc_E:=\{Q_L,\QL: L\in\LE\}
$$
%----------------------------------------------------------
be a family of all ``minimal'' and ``maximal'' cubes of lacunae of the set $E$. Let
%----------------------------------------------------------
$$
\Cc\Bc_E:=\{(\Br,\Br'): \Br,\Br'\in\BRE,\,\Br\bcn\Br'\}
$$
%----------------------------------------------------------
be the family of all (non-ordered) pairs of connected bridges. In Proposition \reff{L-QLM} we have constructed a mapping
%----------------------------------------------------------
\bel{MAP-1}
\Cc\Bc_E\ni (\Br,\Br')\mapsto \hQ(\Br,\Br')\in\Mc_E
\ee
%----------------------------------------------------------
which to every pair of connected bridges $(\Br,\Br')\in\Cc\Bc_E$ assigns a cube $\hQ(\Br,\Br')\in\Mc_E$ satisfying conditions \rf{END-1} and \rf{END-2}. Of course, by these conditions, this mapping satisfies the conditions \rf{MP-1} and \rf{MP-2} as well. However in general the mapping \rf{MAP-1} is not one-to-one mapping so that we can not put $K(\Br,\Br')=\hQ(\Br,\Br')$.
%----------------------------------------------------------
\par Nevertheless the mapping $(\Br,\Br')\to \hQ(\Br,\Br')$ is ``almost'' one-to-one, i.e., every cube $Q\in \Mc_E$
has at most $N(n)$ origins. This enables us to obtain the required mapping $(\Br,\Br')\to K(\Br,\Br')$ by a slight modification of the mapping \rf{MAP-1}. Namely we fix a family consisting of $N(n)$ equal pairwise disjoint subcubes of the cube $Q$ whose diameters are equivalent to $\diam Q$. Then to each origin of $Q$ we assign in a one-to-one way a subcube from this family.
%----------------------------------------------------------
\par Following this scheme let us first prove that the mapping $(\Br,\Br')\to \hQ(\Br,\Br')$ is ``almost'' one-to-one. In fact, let  $L\in\LE$ be a lacuna. Then, by part (i) of Proposition \reff{PRL-1}, there exist at most $C(n)$ lacunae $L'\in\LE$ which contact to $L~(L\lr L')$.
We recall that $L$ is connected by bridges only to {\it contacting} lacunae so that the number of bridges connected to the lacuna $L$ is bounded by $C(n)$ as well.
%----------------------------------------------------------
\par Let $T_L=(A_L,B_L)$ be the interior bridge of the lacuna $L$. By Proposition \reff{L-QLM}, for every bridge $\Br$ connected to $T_L$ the cube $Q(T_L,\Br)$ coincides either with the cube $Q_L$ (i.e., with the ``minimal'' cube of the lacuna $L$), or with the cube $\QL$ (the ``maximal'' cube of the lacuna $L$). This motivates us to introduce two family of cubes: a family
%----------------------------------------------------------
$$
\Ac_L:=\{T\in\BRE: T\bcn T_L,\,\, Q(T_L,\Br)=Q_L\}
$$
%----------------------------------------------------------
and a family
%----------------------------------------------------------
$$
\Ac^{(L)}:=\{T\in\BRE: T\bcn T_L,\,\, Q(T_L,\Br)=\QL\}.
$$
%----------------------------------------------------------
Clearly, $\Ac_L$ and $\Ac^{(L)}$ is a partition of the family
%----------------------------------------------------------
$$
\Ac(L):=\{T\in\BRE: T\bcn T_L\}
$$
%----------------------------------------------------------
of all bridges connected to $T_L$, i.e., $\Ac(L)=\Ac_L\cup\Ac^{(L)}$ and $\Ac_L\cap\Ac^{(L)}=\emp.$
We know that the family $\Ac(L)$ is finite and $\#\Ac(L)\le C(n).$
%----------------------------------------------------------
\par Let us define the mapping $(\Br,\Br')\to K(\Br,\Br')$ on the pairs $(T_L,\Br)$ where $\Br\in\Ac_L$. Since
%----------------------------------------------------------
$$
m_L:=\#\Ac_L\le C(n)
$$
%----------------------------------------------------------
we can represent $\Ac_L$ in the form
%----------------------------------------------------------
$$
\Ac_L=\{T^{(1)},T^{(2)},...,T^{(m_L)}\}.
$$
%----------------------------------------------------------
Then, by \rf{END-1} and \rf{END-2}, for every $i=1,...,m_L,$
%----------------------------------------------------------
\bel{END-3}
\diam \{A_{L},B_{L},A_{T^{(i)}},B_{T^{(i)}}\}\sim\diam Q_L
\ee
%----------------------------------------------------------
and
%----------------------------------------------------------
\bel{END-4}
\{A_{L},B_{L},A_{T^{(i)}},B_{T^{(i)}}\}\subset \gamma Q_L.
\ee
%----------------------------------------------------------
\par By subdividing each edge of the cube $Q_L$ into $m_L$ equal parts we can partition this cube into a family $\{\tK_1,\tK_2,...,\tK_{M_L}\}$ consisting of $M_L:=m_L^n$ congruent cubes of diameter $\diam Q_L/m_L$. We put
%----------------------------------------------------------
\bel{D8}
K_i=\tfrac18\, \tK_i,~~~i=1,...,M_L.
\ee
%----------------------------------------------------------
Since
%----------------------------------------------------------
$$
\diam K_i=\diam Q_L/(8m_L),~~~i=1,...,m_L,
$$
%----------------------------------------------------------
for some absolute constant $\gamma=\gamma(n)$ we have
%----------------------------------------------------------
$$
Q_L\subset \gamma K_i,~~~i=1,...,m_L.
$$
%----------------------------------------------------------
\par Hence, by \rf{END-3} and \rf{END-4},
for every $i=1,...,m_L,$
%----------------------------------------------------------
\bel{END-5}
\diam \{A_{L},B_{L},A_{T^{(i)}},B_{T^{(i)}}\}\sim\diam K_i
\ee
%----------------------------------------------------------
and
%----------------------------------------------------------
\bel{END-6}
\{A_{L},B_{L},A_{T^{(i)}},B_{T^{(i)}}\}\subset \gamma K_i.
\ee
%----------------------------------------------------------
\par Finally we put
%----------------------------------------------------------
$$
K(T_L,T^{(i)}):=K_i,~~~i=1,...,m_L.
$$
%----------------------------------------------------------
Properties \rf{END-5} and \rf{END-6} show that this formula defines the required  mapping $(\Br,\Br')\to K(\Br,\Br')$ for all pairs $(T_L,\Br)$ where $\Br\in\Ac_L$.
%----------------------------------------------------------
\par In the same way we define the mapping $(\Br,\Br')\to K(\Br,\Br')$ for all pairs $(T_L,\Br)$ with $\Br\in\Ac^{(L)}$.
%----------------------------------------------------------
\par Thus we have defined the required mapping on the family
%----------------------------------------------------------
$$
\{(T_L,T): T\in\BRE,\, T\bcn T_L\}
$$
%----------------------------------------------------------
of bridges connected to the (interior) bridge $T_L$.
%----------------------------------------------------------
\par Let us apply this procedure to every lacuna $L\in\LE$. We obtain a mapping $(T,T')\to K(T,T')$ which to every element of the set
%----------------------------------------------------------
$$
\{(T,T')\in\BRE\times\BRE:T\bcn T', T\ne T'\}
$$
%----------------------------------------------------------
assigns a cube
%----------------------------------------------------------
$$
K(T,T')\in \Kc_E.
$$
%----------------------------------------------------------
Here, by \rf{D8}, $\Kc_E$ is a family of cubes defined as follows:
%----------------------------------------------------------
$$
\Kc_E:=\tfrac18\,\widetilde{\Kc}_E=\{\tfrac18\,\tK: \tK\in \widetilde{\Kc}_E\}
$$
%----------------------------------------------------------
where
%----------------------------------------------------------
$$
\widetilde{\Kc}_E=\bigcup_{L\in\LE}
\{\tK_1,\tK_2,...,\tK_{M_L}\}.
$$
%----------------------------------------------------------
\par By \rf{END-5} and \rf{END-6}, this mapping satisfies the required conditions \rf{MP-1} and \rf{MP-2}. Furthermore, since the Whitney cubes are non-overlapping, the cubes of the family
%----------------------------------------------------------
$$
\widetilde{\Kc}_E=\bigcup_{L\in\LE}
\{\tK_1,\tK_2,...,\tK_{M_L}\}
$$
%----------------------------------------------------------
are non-overlapping as well, so that the cubes of the family $\Kc_E:=\tfrac18\,\widetilde{\Kc}_E$ are {\it well-separated}. (See Definition \reff{WSP-C}.)
%----------------------------------------------------------
\par However, in general we can not guarantee that this mapping is well-defined. In fact, if $T,T'\in\BRE, T\bcn T',$ is a pair of {\it connected interior bridges}, the cubes $K(T,T')$ and $K(T',T)$ may be different. To avoid this situation, in this case to the pair $(T,T')$ {\it we simply assign one of these cubes} (no matter which one). As a result we obtain a {\it well-defined} one-to-one mapping satisfying all the conditions of the proposition.
%----------------------------------------------------------
\par The proposition is proved.\bx
%----------------------------------------------------------
%@@@@@@@@@@@@@@@@@@@@@@@@@@@@@@@@@@@@@@@@@@@@@@@@@@@@@@@@@@
%@@@@@@@@@@@@@@@@@@@@@@@@@@@@@@@@@@@@@@@@@@@@@@@@@@@@@@@@@@
%@@@@@@@@@@@@@@@@@@@@@@@@@@@@@@@@@@@@@@@@@@@@@@@@@@@@@@@@@@
%@@@@@@@@@@@@@@@@@@@@@@@@@@@@@@@@@@@@@@@@@@@@@@@@@@@@@@@@@@
%----------------------------------------------------------
%@@@@@@@@@@@@@@@@@@@@@@@@@@@@@@@@@@@@@@@@@@@@@@@@@@@@@@@@@@
%@@@@@@@@@@@@@@@@@@@@@@@@@@@@@@@@@@@@@@@@@@@@@@@@@@@@@@@@@@
%@@@@@@@@@@@@@@@@@@@@@@@@@@@@@@@@@@@@@@@@@@@@@@@@@@@@@@@@@@
%@@@@@@@@@@@@@@@@@@@@@@@@@      @@@@@@@@@@@@@@@@@@@@@@@@@@@
%@@@@@@@@@@@@@@@@@@@@@@@          @@@@@@@@@@@@@@@@@@@@@@@@@
%@@@@@@@@@@@@@@@@@@@@@              @@@@@@@@@@@@@@@@@@@@@@@
%@@@@@@@@@@@@@@@@@@@     SECTION 7    @@@@@@@@@@@@@@@@@@@@@
%@@@@@@@@@@@@@@@@@@@@@              @@@@@@@@@@@@@@@@@@@@@@@
%@@@@@@@@@@@@@@@@@@@@@@@          @@@@@@@@@@@@@@@@@@@@@@@@@
%@@@@@@@@@@@@@@@@@@@@@@@@@      @@@@@@@@@@@@@@@@@@@@@@@@@@@
%@@@@@@@@@@@@@@@@@@@@@@@@@@@@@@@@@@@@@@@@@@@@@@@@@@@@@@@@@@
%@@@@@@@@@@@@@@@@@@@@@@@@@@@@@@@@@@@@@@@@@@@@@@@@@@@@@@@@@@
%@@@@@@@@@@@@@@@@@@@@@@@@@@@@@@@@@@@@@@@@@@@@@@@@@@@@@@@@@@
%----------------------------------------------------------
\SECT{7. A lacunary Whitney-type extension operator.}{7}
%----------------------------------------------------------
\addtocontents{toc}{7. A lacunary Whitney-type extension operator.\hfill \thepage\\\par}
%----------------------------------------------------------
\indent
%@@@@@@@@@@@@@@@@@@@@@@@@@@@@@@@@@@@@@@@@@@@@@@@@@@@@@@@@@@
%----------------------------------------------------------
\par We return to the proof of the sufficiency part of Theorem \reff{MAIN}. Thus in this and the next section  $E$ is a finite subset of $\RT$ and  $p\in(2,\infty)$.
%----------------------------------------------------------
\par Our aim is to prove a stronger result, Theorem \reff{S-MAIN}, which immediately implies the sufficiency in Theorem \reff{MAIN}.
%----------------------------------------------------------
%@@@@@@@@@@@@@@@@@@@@@@@@@@@@@@@@@@@@@@@@@@@@@@@@@@@@@@@@@@
%@@@@@@@@@@@@@@@@@@@@@@@@@@@@@@@@@@@@@@@@@@@@@@@@@@@@@@@@@@
%@@@@@@@@@@@@@@@@@@@@@@@@@@@@@@@@@@@@@@@@@@@@@@@@@@@@@@@@@@
%@@@@@@@@@@@@@@@@@@@@@@@@@@@@@@@@@@@@@@@@@@@@@@@@@@@@@@@@@@
%----------------------------------------------------------
\begin{theorem}\lbl{S-MAIN} Let $E$ be a finite subset of $\RT$ and let $f$ be a function on $E$. Then
%----------------------------------------------------------
\bel{SM-H}
\|f\|_{\LTP|_E}\le C(p)\lambda^{\frac{1}{p}}
\ee
%----------------------------------------------------------
provided $\lambda$ is a positive constant which satisfies all of the following conditions for a certain absolute positive constant $\gamma$:
%----------------------------------------------------------
\par (a). The condition (i) of Theorem \reff{MAIN} holds;
%----------------------------------------------------------
\par (b). Let $\Qc$ and $\Kc$ be finite families of pairwise disjoint squares. Assume that to each \sq $K\in\Kc$ we have arbitrarily assigned a triangle $\Delta(K)$ in $E$ such that
%----------------------------------------------------------
$$
\Delta(K)\subset\gamma K~~~\text{and}~~~\diam K\le\gamma\diam \Delta(K).
$$
%----------------------------------------------------------
\par Suppose that to each \sq $Q\in\Qc$ we have arbitrarily assigned a pair of \sqs $Q',Q''\in\Qc$ such that
$Q'\cup Q''\subset \gamma Q$ and
%----------------------------------------------------------
\bel{CV-S}
(\diam Q')^{p-2}\sigma_p(Q';\Kc)+
(\diam Q'')^{p-2}\sigma_p(Q'';\Kc)\le 1.
\ee
%----------------------------------------------------------
\par Then the following inequality
%----------------------------------------------------------
\bel{CRT-S}
\smed\limits_{Q\in\Qc}\,
\left(\frac{\diam Q' \diam Q''}{\diam Q} \right)^{p-2} S_p(f:Q',Q'';\Kc)
\le\lambda
\ee
%----------------------------------------------------------
holds.
%----------------------------------------------------------
\end{theorem}
%----------------------------------------------------------
%@@@@@@@@@@@@@@@@@@@@@@@@@@@@@@@@@@@@@@@@@@@@@@@@@@@@@@@@@@
%@@@@@@@@@@@@@@@@@@@@@@@@@@@@@@@@@@@@@@@@@@@@@@@@@@@@@@@@@@
%@@@@@@@@@@@@@@@@@@@@@@@@@@@@@@@@@@@@@@@@@@@@@@@@@@@@@@@@@@
%----------------------------------------------------------
\medskip
%----------------------------------------------------------
\par We recall that the quantities $\sigma_p$ and $S_p$ have been defined in Section 1, see formulation of Theorem \reff{MAIN}.\medskip
%----------------------------------------------------------
%@@@@@@@@@@@@@@@@@@@@@@@@@@@@@@@@@@@@@@@@@@@@@@@@@@@@@@@@@@
%@@@@@@@@@@@@@@@@@@@@@@@@@@@@@@@@@@@@@@@@@@@@@@@@@@@@@@@@@@
%@@@@@@@@@@@@@@@@@@@@@@@@@@@@@@@@@@@@@@@@@@@@@@@@@@@@@@@@@@
%----------------------------------------------------------
\par {\it Proof.} We prove Theorem \reff{S-MAIN} in several stages. At the first stage given a family of affine interpolating polynomials we construct a Whitney-type extension of the function $f$ from $E$ to all of $\RT$. Then we estimate the $L^2_p$-norm of this extension via oscillations of these polynomials on contacting Whitney squares.
%----------------------------------------------------------
\par Let  $$\{P_L\in\PO: L\in\LE\}$$ be a family of affine polynomials such that for every lacuna $L\in\LE$
%----------------------------------------------------------
\bel{INT-P}
P_L(A_L)=f(A_L),~~~P_L(B_L)=f(B_L).
\ee
%----------------------------------------------------------
\par Let $L\in\LE$ be a lacuna of Whitney squares. To every \sq $Q\in L$ we assign an affine polynomial
%----------------------------------------------------------
$$
\PQ:=P_L.
$$
%----------------------------------------------------------
Then we construct an extension of the function $f$ using the classical Whitney formula
%----------------------------------------------------------
\bel{W-EXT}
F(x):=\left \{
%----------------------------------------------------------
\begin{array}{ll}
f(x),& x\in E,\\\\
\sum\limits_{Q\in W_E}
\varphi_Q(x)\,\PQ(x),& x\in\RT\setminus E.
\end{array}
%----------------------------------------------------------
\right.
\ee
%----------------------------------------------------------
Here as usual $\{\varphi_Q:Q\in W_E\}$ denotes a smooth partition of unity  subordinated to the Whitney decomposition $W_E$. Let us recall its main properties, see, e.g. \cite{St}, Ch. 6.
%@@@@@@@@@@@@@@@@@@@@@@@@@@@@@@@@@@@@@@@@@@@@@@@@@@@@@@@@@@
\begin{lemma}\lbl{P-U} The family of functions $\{\varphi_Q:Q\in W_E\}$ has the following properties:
%----------------------------------------------------------
%@@@@@@@@@@@@@@@@@@@@@@@@@@@@@@@@@@@@@@@@@@@@@@@@@@@@@@@@@@
%----------------------------------------------------------
\par (a). $\varphi_Q\in C^\infty(\RT)$ and
$0\le\varphi_Q\le 1$ for every $Q\in W_E$;\smallskip
%----------------------------------------------------------
\par (b). $\supp \varphi_Q\subset Q^*(:=\frac{9}{8}Q),$
$Q\in W_E$;\smallskip
%----------------------------------------------------------
\par (c). $\sum\{\varphi_Q(x):Q\in W_E\}=1$ for every
$x\in\RT\setminus E$;\smallskip
%----------------------------------------------------------
\par (d). $ |D^\beta\varphi_Q(x)| \le C(\diam Q)^{-|\beta|}\,\,$ for every $Q\in W_E$, every $x\in\RT$ and every multiindex $\beta$ of order $|\beta|\le 2$.
%----------------------------------------------------------
%@@@@@@@@@@@@@@@@@@@@@@@@@@@@@@@@@@@@@@@@@@@@@@@@@@@@@@@@@@
\end{lemma}
%@@@@@@@@@@@@@@@@@@@@@@@@@@@@@@@@@@@@@@@@@@@@@@@@@@@@@@@@@@
\smallskip
%----------------------------------------------------------
\par Given lacunae $L,L'\in\LE$ we introduce a \sq
%----------------------------------------------------------
\bel{Q-T1}
Q(L,L'):=Q(A_L,d(L,L')).
\ee
%----------------------------------------------------------
Recall that
%----------------------------------------------------------
$$
d(L,L'):=\diam\{A_L,B_L,A_{L'},B_{L'}\}.
$$
%----------------------------------------------------------
%@@@@@@@@@@@@@@@@@@@@@@@@@@@@@@@@@@@@@@@@@@@@@@@@@@@@@@@@@@
%@@@@@@@@@@@@@@@@@@@@@@@@@@@@@@@@@@@@@@@@@@@@@@@@@@@@@@@@@@
%@@@@@@@@@@@@@@@@@@@@@@@@@@@@@@@@@@@@@@@@@@@@@@@@@@@@@@@@@@
%@@@@@@@@@@@@@@@@@@@@@@@@@@@@@@@@@@@@@@@@@@@@@@@@@@@@@@@@@@
%@@@@@@@@@@@@@@@@@@@@@@@@@@@@@@@@@@@@@@@@@@@@@@@@@@@@@@@@@@
%@@@@@@@@@@@@@@@@@@@@@@@@@@@@@@@@@@@@@@@@@@@@@@@@@@@@@@@@@@
%----------------------------------------------------------
\begin{proposition}\lbl{F-INAF} The $L^2_p$-norm of the extension $F$ satisfies the following inequality:
%----------------------------------------------------------
$$
\|F\|^p_{\LTP}\le C\,\sum \left\{d(L,L')^{2-2p}
\max\limits_{Q(L,L')}|P_L-P_{L'}|^p: L,L'\in\LE, L\lr L'\right\}.
$$
%----------------------------------------------------------
\end{proposition}
%----------------------------------------------------------
%@@@@@@@@@@@@@@@@@@@@@@@@@@@@@@@@@@@@@@@@@@@@@@@@@@@@@@@@@@
%@@@@@@@@@@@@@@@@@@@@@@@@@@@@@@@@@@@@@@@@@@@@@@@@@@@@@@@@@@
%@@@@@@@@@@@@@@@@@@@@@@@@@@@@@@@@@@@@@@@@@@@@@@@@@@@@@@@@@@
%----------------------------------------------------------
\par {\it Proof.} Since every point $z\in E$ is isolated, there exists a (unique) true lacuna $L\in\LE$ such that $V_L=\{z\}$. Hence $A_L=z$ so that, by \rf{INT-P}, $P_L(z)=f(z)$. Since in this case the set $\{\cup Q: Q\in L\}$ contains a certain neighborhood of $z$, by the formula \rf{W-EXT}, $F$ coincides with $P_L$ on this neighborhood.
%----------------------------------------------------------
\par On the other hand, by the same formula, $F\in C^{\infty}(\RT\setminus E)$ so that $F\in C^{\infty}(\RT)$. This shows that distributional partial derivatives of $F$ can be identified with its regular derivatives. Let us estimate such a derivative of $F$ of order two on a Whitney square.  %----------------------------------------------------------
\par Let $K\in W_E$ be a Whitney \sq which belongs to a lacuna $L\in\LE$. Let $x\in K$ and let $\alpha$ be a multiindex, $|\alpha|=2$. Then, by the formula \rf{W-EXT}, Lemma \reff{P-U} and Lemma \reff{Wadd},
%----------------------------------------------------------
\be
|D^\alpha F(x)|&=&|D^\alpha(F(x)-P_L(x))|=
\left|D^\alpha\left(\sum\left\{
\varphi_Q(x)\,(\PQ(x)-P_L(x)):Q\in W_E\right\}\right)\right|\nn\\&=&
\left|D^\alpha\left(\sum\{\varphi_Q(x)\,(\PQ(x)-P_L(x))
:Q\in W_E, Q^*\ni x\}\right)\right|\nn\\&\le&
\sum\{|D^\alpha(\varphi_Q(\PQ-P_L))(x)|:Q\in W_E, Q^*\ni x\}.\nn
%----------------------------------------------------------
\ee
%----------------------------------------------------------
\par Given a \sq $Q\in W_E$ such that $Q^*\ni x$ let us estimate the quantity
%----------------------------------------------------------
$$
I_Q(x):=|D^\alpha(\varphi_Q(\PQ-P_L))(x)|.
$$
%----------------------------------------------------------
We have
%----------------------------------------------------------
$$
I_Q(x)\le C\sum\limits_{|\beta|\le 2} |D^{\alpha-\beta}\varphi_Q(x)|\cdot|D^{\beta}(\PQ-P_L)(x)|
$$
%----------------------------------------------------------
so that, by Lemma \reff{P-U},
%----------------------------------------------------------
$$
I_Q(x)\le C\sum\limits_{|\beta|\le 2} (\diam Q)^{|\beta|-|\alpha|} \max_{Q^*}|D^{\beta}(\PQ-P_L)|.
$$
%----------------------------------------------------------
By Markov's inequality,
%----------------------------------------------------------
$$
\max_{Q^*}|D^{\beta}(\PQ-P_L)|\le C(\diam Q^*)^{-|\beta|} \max_{Q^*}|\PQ-P_L|.
$$
%----------------------------------------------------------
Clearly, for every polynomial $P\in\PO$, every \sq $Q\subset\RT$ and every $\gamma>1$ we have
%----------------------------------------------------------
\bel{IN-M1}
\max_{\gamma Q}|P|\le C(\gamma)\max_{Q}|P|\,.
\ee
%----------------------------------------------------------
Hence
%----------------------------------------------------------
$$
\max_{Q^*}|D^{\beta}(\PQ-P_L)|\le C(\diam Q)^{-|\beta|} \max_{Q}|\PQ-P_L|
$$
%----------------------------------------------------------
so that
%----------------------------------------------------------
$$
I_Q(x)\le C(\diam Q)^{-|\alpha|}\max_Q|\PQ-P_L|=C(\diam Q)^{-2}\max_Q|\PQ-P_L|\,.
$$
%----------------------------------------------------------
\par We obtain
%----------------------------------------------------------
\be
|D^\alpha F(x)|&\le& C\sum\left\{(\diam Q)^{-2}\max\limits_Q|\PQ-P_L|:Q\in W_E, Q^*\ni x\right\}\nn\\
&\le& C\sum\left\{(\diam Q)^{-2}\max\limits_Q|\PQ-P_L|:Q\in W_E, Q^*\cap K\ne\emp\right\}\nn
\ee
%----------------------------------------------------------
so that, by part (3) of Lemma \reff{Wadd},
%----------------------------------------------------------
$$
|D^\alpha F(x)|\le C\sum\left\{(\diam Q)^{-2}\max\limits_Q|\PQ-P_L|:Q\in W_E, Q\cap K\ne\emp\right\}.
$$
%----------------------------------------------------------
\par We let $J_K$ denote a family of Whitney \sqs
%----------------------------------------------------------
$$
J_K:=\{Q\in W_E: Q\cap K\ne\emp\}.
$$
%----------------------------------------------------------
We note that, by part (1) of Lemma \reff{Wadd}, $\diam Q\sim \diam K$ for every \sq $Q\in J_K$. Hence
%----------------------------------------------------------
$$
|D^\alpha F(x)|\le C(\diam K)^{-2}\,\sum_{Q\in J_K} \max\limits_Q|\PQ-P_L|.
$$
%----------------------------------------------------------
By part (2) of Lemma \reff{Wadd}, $\#J_K\le C$ so that
%----------------------------------------------------------
$$
|D^\alpha F(x)|\le C(\diam K)^{-2}\,\max_{Q\in J_K} \max\limits_Q|\PQ-P_L|\,.
$$
%----------------------------------------------------------
\par We let $L^{cont}$ denote a family of {\it contacting} \sqs of the lacuna $L$:
%----------------------------------------------------------
$$
L^{cont}:=\{Q\in L:\exists~L'\in\LE, L'\ne L,\, Q'\in L'~~\text{such that}~~Q'\cap Q\ne\emp\}.
$$
%----------------------------------------------------------
\par Note that {\it if the \sq $K$ is not a contacting square,} i.e., $K\in L\setminus L^{(cont)}$, then $J_K\subset L$ so that  $\PQ=P_L$ for every $Q\in J_K$. Hence
%----------------------------------------------------------
$$
D^\alpha F(x)=0~~~~~\text{for every}~~~K\in L\setminus L^{(cont)}~~~\text{and every}~~~x\in K.
$$
%----------------------------------------------------------
\par Let $K\in L^{(cont)}$ be a contacting square. Recall that $\PQ=P_{L}$ for every \sq $Q\in L$. Therefore there exists a contacting lacuna $L_K\in\LE$, $L_K\lr L$,
$L_K\ne L$, and a \sq $\tK\in L_K$ such that $\tK\in J_K$ and
%----------------------------------------------------------
$$
\max_{Q\in J_K} \max\limits_Q|\PQ-P_L|=\max\limits_{\tK}|P^{(\tK)}-P_L|=
\max\limits_{\tK}|P_{L_K}-P_L|.
$$
%----------------------------------------------------------
(In particular $\tK$ is a contacting \sq for $K$.) Hence
%----------------------------------------------------------
$$
|D^\alpha F(x)|\le C\,(\diam K)^{-2} \max\limits_{\tK}|P_{L_K}-P_L|,~~~x\in K.
$$
%----------------------------------------------------------
\par Integrating this inequality to the power $p$ on the \sq $K$ we obtain
%----------------------------------------------------------
$$
\intl_K|D^\alpha F(x)|^p\,dx\le C\,(\diam K)^{-2p+2} \max\limits_{\tK}|P_{L_K}-P_L|^p.
$$
%----------------------------------------------------------
\par Note that, by Proposition \reff{CL-Q},
$\diam K\sim d(L,L_K).$
%----------------------------------------------------------
Furthermore,
%----------------------------------------------------------
$$
\{A_L,B_L,A_{L_K},B_{L_K}\}\subset \gamma K
$$
%----------------------------------------------------------
so that the \sq $Q(L,L_K)=Q(A_L,d(L,L_K))$ have common points with $K$. Hence $K\subset \gamma_1 Q(L,L_K)$ with some absolute constant $\gamma_1$.
%----------------------------------------------------------
\par We obtain
%----------------------------------------------------------
$$
\intl_K|D^\alpha F(x)|^p\,dx\le C\,d(L,L_K)^{-2p+2} \max\limits_{\gamma_1Q(L,L_K)}|P_{L_K}-P_L|^p.
$$
%----------------------------------------------------------
By \rf{IN-M1},
%----------------------------------------------------------
$$
\max\limits_{\gamma_1Q(L,L_K)}|P_{L_K}-P_L|\le
C\max\limits_{Q(L,L_K)}|P_{L_K}-P_L|
$$
%----------------------------------------------------------
so that
%----------------------------------------------------------
$$
\intl_K|D^\alpha F(x)|^p\,dx\le C\,d(L,L_K)^{-2p+2} \max\limits_{Q(L,L_K)}|P_{L_K}-P_L|^p.
$$
%----------------------------------------------------------
Hence
%----------------------------------------------------------
$$
\smed\limits_{K\in L}\intl_K|D^\alpha F(x)|^p\,dx\le C\,\sum\left\{d(L,L')^{-2p+2} \max\limits_{Q(L,L')}|P_{L'}-P_L|^p: L'\in\LE, L'\lr L\right\}
$$
%----------------------------------------------------------
which implies that
%----------------------------------------------------------
\be
\|F\|^p_{\LTP}&=&\intl_{\RT}|D^\alpha F(x)|^p\,dx\le \sum\limits_{L\in\LE}\sum\limits_{K\in L}\intl_K|D^\alpha F(x)|^p\,dx\nn\\&\le&
C\,\sum\limits_{L\in\LE}
\sum\limits_{L'\lr L} d(L,L')^{-2p+2} \max\limits_{Q(L,L')}|P_{L'}-P_L|^p\nn
\ee
%----------------------------------------------------------
proving the proposition.\bx\medskip
%----------------------------------------------------------
%@@@@@@@@@@@@@@@@@@@@@@@@@@@@@@@@@@@@@@@@@@@@@@@@@@@@@@@@@@
%@@@@@@@@@@@@@@@@@@@@@@@@@@@@@@@@@@@@@@@@@@@@@@@@@@@@@@@@@@
%@@@@@@@@@@@@@@@@@@@@@@@@@@@@@@@@@@@@@@@@@@@@@@@@@@@@@@@@@@
%@@@@@@@@@@@@@@@@@@@@@@@@@@@@@@@@@@@@@@@@@@@@@@@@@@@@@@@@@@
%@@@@@@@@@@@@@@@@@@@@@@@@@@@@@@@@@@@@@@@@@@@@@@@@@@@@@@@@@@
%@@@@@@@@@@@@@@@@@@@@@@@@@@@@@@@@@@@@@@@@@@@@@@@@@@@@@@@@@@
%----------------------------------------------------------
\par We follow to the scheme of the proof described in Section 2. Our next goal is to express the Sobolev norm of the extension $F$ defined by the formula \rf{W-EXT} via {\it the gradients} of interpolating polynomials.
%----------------------------------------------------------
\par Let $g:\BRE\to\RT$ be a mapping which to every bridge $T\in\BRE$ (interior or exterior) assigns a vector in $\RT$. Suppose that for every \br $T\in\BRE$ with the ends at points $\A{T},\B{T}\in E$ the following equality
%----------------------------------------------------------
\bel{M-EQV}
\ip{g(T),\A{T}-\B{T}}=f(\A{T})-f(\B{T})
\ee
%----------------------------------------------------------
holds. (Recall that $\ip{\cdot,\cdot}$ is the standard inner product in $\RT$.)
%----------------------------------------------------------
\par Let $L\in\LE$ be a lacuna and let $P_L\in\PO$ be an affine polynomial defined by the following formula
%----------------------------------------------------------
\bel{PL-G}
P_L(x):=f(A_L)+\ip{g(T_L),x-A_L},~~~x\in\RT.
\ee
%----------------------------------------------------------
Recall that $T_L$ denotes the interior bridge $T_L=(A_L,B_L)\in\BRE$. Clearly, by \rf{M-EQV},
%----------------------------------------------------------
\bel{PL-INT}
P_L(A_L)=f(A_L),~~~P_L(B_L)=f(B_L).
\ee
%----------------------------------------------------------
Let us also note that $g(T_L)=\nabla P_L.$
%----------------------------------------------------------
\par Let $T,T'\in\BRE$ be a pair of \brs with the ends at points $\{\A{T},\B{T}\}$ and $\{\A{T'},\B{T'}\}$ correspondingly. Let
%----------------------------------------------------------
\bel{D-TTP}
D(T,T'):=\diam \{\A{T},\B{T},\A{T'},\B{T'}\}.
\ee
%----------------------------------------------------------
%@@@@@@@@@@@@@@@@@@@@@@@@@@@@@@@@@@@@@@@@@@@@@@@@@@@@@@@@@@
%@@@@@@@@@@@@@@@@@@@@@@@@@@@@@@@@@@@@@@@@@@@@@@@@@@@@@@@@@@
%@@@@@@@@@@@@@@@@@@@@@@@@@@@@@@@@@@@@@@@@@@@@@@@@@@@@@@@@@@
%@@@@@@@@@@@@@@@@@@@@@@@@@@@@@@@@@@@@@@@@@@@@@@@@@@@@@@@@@@
%----------------------------------------------------------
\begin{proposition}\lbl{F-GN} Let $F$ be the extension of the function $f$ defined by the formula \rf{W-EXT} where for every lacuna $L\in\LE$ the polynomial $P_L\in\PO$ is determined by \rf{PL-G}. Then
%----------------------------------------------------------
$$
\|F\|^p_{\LTP}\le C\,\sum \left\{D(T,T')^{2-p}
\|g(T)-g(T')\|^p: T,T'\in\BRE, T\bcn T'\right\}.
$$
%----------------------------------------------------------
\end{proposition}
%----------------------------------------------------------
%@@@@@@@@@@@@@@@@@@@@@@@@@@@@@@@@@@@@@@@@@@@@@@@@@@@@@@@@@@
%@@@@@@@@@@@@@@@@@@@@@@@@@@@@@@@@@@@@@@@@@@@@@@@@@@@@@@@@@@
%@@@@@@@@@@@@@@@@@@@@@@@@@@@@@@@@@@@@@@@@@@@@@@@@@@@@@@@@@@
%----------------------------------------------------------
\par {\it Proof.} Let $L,L'\in\LE$, $L\lr L'$, be contacting lacunae. Let
%----------------------------------------------------------
$$
T_{L}=(A_{L},B_{L})~~~\text{and}~~~
T_{L'}=(A_{L'},B_{L'})
$$
%----------------------------------------------------------
be interior bridges of $L$ and $L'$ with the ends $\{A_{L},B_{L}\}$ and $\{A_{L'},B_{L'}\}$.
%----------------------------------------------------------
\par Suppose that these bridges are {\it connected to each other by an exterior bridge}
%----------------------------------------------------------
$$
T=\{C(L,L'),C(L',L)\}.
$$
%----------------------------------------------------------
Thus $T_L\bcn T$ and $T_{L'}\bcn T$ and
%----------------------------------------------------------
$$
C(L,L')\in \{A_L,B_L\}~~~\text{and}~~~
C(L',L)\in\{A_{L'},B_{L'}\}.
$$
%----------------------------------------------------------
\par Let us estimate the quantity
%----------------------------------------------------------
$$
I(L,L'):=d(L,L')^{2-2p}
\max\limits_{Q(L,L')}|P_L-P_{L'}|^p
$$
%----------------------------------------------------------
which appears in the righthand side of the inequality of Proposition \reff{F-INAF}. Let
%----------------------------------------------------------
\bel{PT-X}
P^{[T]}(x):=f(C(L,L'))+\ip{g(T),x-C(L,L')},~~~x\in\RT.
\ee
%----------------------------------------------------------
Note that, by \rf{M-EQV},
%----------------------------------------------------------
$$
\ip{g(T),C(L,L')-C(L',L)}=f(C(L,L'))-f(C(L',L))
$$
%----------------------------------------------------------
so that
%----------------------------------------------------------
$$
P^{[T]}(C(L,L'))=f(C(L,L'))~~~\text{and}~~~
P^{[T]}(C(L',L))=f(C(L',L)).
$$
%----------------------------------------------------------
Since $C(L,L')\in\{A_L,B_L\}$, by \rf{PL-INT}, $P_L(C(L,L'))=f(C(L,L'))$. Since $P_L\in\PO$, it can be  represented in the form
%----------------------------------------------------------
$$
P_L(x):=f(C(L,L'))+\ip{\nabla P_L,x-C(L,L')},~~~x\in\RT.
$$
%----------------------------------------------------------
Hence
%----------------------------------------------------------
$$
P_L(x):=f(C(L,L'))+\ip{g(T_L),x-C(L,L')}
$$
%----------------------------------------------------------
so that, by \rf{PT-X}, for every $x\in\RT$ we have
%----------------------------------------------------------
$$
|P_L(x)-P^{[T]}(x)|=|\ip{g(T_L)-g(T),x-C(L,L')}|.
$$
%----------------------------------------------------------
This implies the following inequality:
%----------------------------------------------------------
\bel{FR-1}
|P_L(x)-P^{[T]}(x)|\le C\,\|g(T_L)-g(T)\|\, \|x-C(L,L')\|,~~~x\in\RT.
\ee
%----------------------------------------------------------
\par Recall that
%----------------------------------------------------------
\bel{DT-DF}
D(T_L,T_{L'})=\diam\{A_L,B_L,A_{L'},B_{L'}\},
\ee
%----------------------------------------------------------
see \rf{D-TTP}, and
%----------------------------------------------------------
$$
Q(L,L')=Q(A_L,d(L,L'))
$$
%----------------------------------------------------------
where $d(L,L')=\diam\{A_L,B_L,A_{L'},B_{L'}\}$, see \rf{Q-T1}. (Thus $d(L,L')=D(T_L,T_{L'})$.)
%----------------------------------------------------------
\par In particular,
%----------------------------------------------------------
$$
A_L,B_L,A_{L'},B_{L'}\in Q(L,L')
$$
%----------------------------------------------------------
so that $C(L,L'),C(L',L)\in Q(L,L')$ as well. Hence, by \rf{FR-1},
%----------------------------------------------------------
$$
\max\limits_{Q(L,L')}|P_L-P^{[T]}|\le C\diam Q(L,L')\,\|g(T_L)-g(T)\|\le C\,d(L,L')\,\|g(T_L)-g(T)\|.
$$
%----------------------------------------------------------
In the same way we obtain the following inequality
%----------------------------------------------------------
$$
\max\limits_{Q(L,L')}|P_{L'}-P^{[T]}|\le C\,d(L,L')\,\|g(T_{L'})-g(T)\|.
$$
%----------------------------------------------------------
Hence
%----------------------------------------------------------
\be
I(L,L')&=&d(L,L')^{2-2p}
\max\limits_{Q(L,L')}|P_L-P_{L'}|^p\nn\\
&\le& C d(L,L')^{2-p}\,\{\|g(T_L)-g(T)\|^p+\|g(T_{L'})-g(T)\|^p\}.
\nn
\ee
%----------------------------------------------------------
But, by \rf{DT-DF},
%----------------------------------------------------------
$$
D(T_L,T)=\diam\{A_L,B_L,C(L',L)\}\le D(T_L,T_{L'})= d(L,L')
$$
%----------------------------------------------------------
and
%----------------------------------------------------------
$$
D(T_{L'},T)=\diam\{A_{L'},B_{L'},C(L',L)\}\le D(T_L,T_{L'})= d(L,L')
$$
%----------------------------------------------------------
so that
%----------------------------------------------------------
$$
I(L,L')\le C\,\{ D(T_L,T)^{2-p}\,\,\|g(T_L)-g(T)\|^p+
D(T_{L'},T)^{2-p}\,\,\|g(T_{L'})-g(T)\|^p\}.
$$
%----------------------------------------------------------
\par We obtain the same estimate of $I(L,L')$ whenever {\it the interior bridges $T_L$ and $T_{L'}$ are connected.} In this case we can simply put $T:=T_{L'}$ and repeat the same proof as for the previous case where  $T_L$ and $T_{L'}$ are connected by an exterior bridge.
%----------------------------------------------------------
\par It remains to replace in the inequality of Proposition \reff{F-INAF} the quantity $I(L,L')$ by its estimate via the mapping $g$ given above, and the proposition follows.\bx\bigskip
%----------------------------------------------------------
%@@@@@@@@@@@@@@@@@@@@@@@@@@@@@@@@@@@@@@@@@@@@@@@@@@@@@@@@@@
%@@@@@@@@@@@@@@@@@@@@@@@@@@@@@@@@@@@@@@@@@@@@@@@@@@@@@@@@@@
%@@@@@@@@@@@@@@@@@@@@@@@@@@@@@@@@@@@@@@@@@@@@@@@@@@@@@@@@@@
%@@@@@@@@@@@@@@@@@@@@@@@@@@@@@@@@@@@@@@@@@@@@@@@@@@@@@@@@@@
%@@@@@@@@@@@@@@@@@@@@@@@@@@@@@@@@@@@@@@@@@@@@@@@@@@@@@@@@@@
%@@@@@@@@@@@@@@@@@@@@@@@@@@@@@@@@@@@@@@@@@@@@@@@@@@@@@@@@@@
%----------------------------------------------------------
\par  We turn to the next step of the proof of Theorem \reff{S-MAIN}. We shall estimate the $L^2_p$-norm of the extension $F$ of the function $f$  via $L_p$-norm of an additional parameter, a function $h:\RT\to \R_+$. We will see that averages of this function on certain \sqs majorize distances between values of the mapping $g:\BRE\to\RT$ on connected bridges, see Proposition \reff{A-H} below.
%----------------------------------------------------------
\par The proof of this proposition uses on an auxiliary result related to the theory of Muckenhoupt's weights.
We recall, see, e.g. \cite{GR}, that a non-negative function $w\in L_{1,loc}(\RT)$ is said to be $A_1$-weight if there exists $\lambda>0$ such that for every \sq $Q\in \RT$ the following inequality
%----------------------------------------------------------
$$
\frac{1}{|Q|}\intl_Q w(u)du\le \lambda \essinf w
$$
%---------------------------------------------------------
holds. We put $\|w\|_{A_1}=\inf \lambda$.
%---------------------------------------------------------
\par Clearly, a weight $w\in A_1$ if and only if
%----------------------------------------------------------
$$
\Mc[w]\le \lambda w(x)~~~~~a.e.~on~ \RT,
$$
%----------------------------------------------------------
and
%----------------------------------------------------------
$$
\|w\|_{A_1}\sim \esssup\limits_{\RT}\frac{\Mc[w](x)}{w(x)}.
$$
%----------------------------------------------------------
Recall that
%----------------------------------------------------------
$$
\Mc[h](x)=\sup_{Q\ni x}\frac{1}{|Q|}\intl_Q|h|dx,~~~x\in\RT,
$$
%----------------------------------------------------------
denotes the \HLM function of a function $h\in L_{1,\,loc}(\RT)$. We also note the following property of a weight $w\in A_1$: if $K,Q$ are two \sqs in $\RT$ such that $K\subset Q$, then
%----------------------------------------------------------
\bel{4.C}
\frac{1}{|Q|}\intl_Q w(u)du\le \|w\|_{A_1}\frac{1}{|K|}\intl_K w(u)du.
\ee
%---------------------------------------------------------
(In other words, the average of a function $w\in A_1$ is a quasi-monotone non-increasing function of a square.)
%----------------------------------------------------------
%@@@@@@@@@@@@@@@@@@@@@@@@@@@@@@@@@@@@@@@@@@@@@@@@@@@@@@@@@@
%----------------------------------------------------------
\par The following remarkable result of Coifman and Rochberg \cite{CR} presents an important property of $A_1$-weights.
%---------------------------------------------------------
\begin{theorem}\lbl{C-R} Let $w\in L_{1,loc}(\RT)$, and let $\Mc[w](x)<\infty$ a.e. Then
$\Mc[w]^{\theta}\in A_1(\RT)$ for every $0<\theta<1$. Furthermore,
$\|\Mc[w]^{\theta}\|_{A_1}\le C(\theta).$
%---------------------------------------------------------
\end{theorem}
%----------------------------------------------------------
%@@@@@@@@@@@@@@@@@@@@@@@@@@@@@@@@@@@@@@@@@@@@@@@@@@@@@@@@@@
%@@@@@@@@@@@@@@@@@@@@@@@@@@@@@@@@@@@@@@@@@@@@@@@@@@@@@@@@@@
\par Let $q=(p+2)/2$; thus $2<q<p$. Let $T,T'\in\BRE$ be a pair of bridges with the ends at points $\{\A{T},\B{T}\}$ and $\{\A{T'},\B{T'}\}$ respectively. By $Q(T,T')$ we denote a \sq
%----------------------------------------------------------
\bel{SQ-T}
Q(T,T')=Q(\A{T},D(T,T')).
\ee
%----------------------------------------------------------
Recall that $D(T,T'):=\diam \{\A{T},\B{T},\A{T'},\B{T'}\}$, see \rf{D-TTP}.
%----------------------------------------------------------
%@@@@@@@@@@@@@@@@@@@@@@@@@@@@@@@@@@@@@@@@@@@@@@@@@@@@@@@@@@
%@@@@@@@@@@@@@@@@@@@@@@@@@@@@@@@@@@@@@@@@@@@@@@@@@@@@@@@@@@
%@@@@@@@@@@@@@@@@@@@@@@@@@@@@@@@@@@@@@@@@@@@@@@@@@@@@@@@@@@
%@@@@@@@@@@@@@@@@@@@@@@@@@@@@@@@@@@@@@@@@@@@@@@@@@@@@@@@@@@
%----------------------------------------------------------
\begin{proposition}\lbl{A-H} Let $g:\BRE\to\RT$ be a mapping satisfying condition \rf{M-EQV}, and let $F$ be the extension of $f$ defined by \rf{W-EXT} where the polynomials $P_L, L\in\LE,$ are given by the formula \rf{PL-G}.
%----------------------------------------------------------
\par Let $\gamma>1$ and let $h:\RT\to\R_+$ be a non-negative function such that $h\in\LPRT$. Suppose that for every pair of connected bridges $T,T'\in\BRE,$ $T\bcn T'$, the following inequality
%----------------------------------------------------------
\bel{GS-A}
\|g(T)-g(T')\|\le \diam Q(T,T')\,\left(\frac{1}{|\gamma Q(T,T')|}
\intl_{\gamma Q(T,T')}h^q(z)dz\right)^{\frac1q}
\ee
%----------------------------------------------------------
holds. Then $F\in\LTP$ and
%----------------------------------------------------------
$$
\|F\|_{\LTP}\le C(p,\gamma)\|h\|_{\LPRT}.
$$
%----------------------------------------------------------
\end{proposition}
%----------------------------------------------------------
%@@@@@@@@@@@@@@@@@@@@@@@@@@@@@@@@@@@@@@@@@@@@@@@@@@@@@@@@@@
%@@@@@@@@@@@@@@@@@@@@@@@@@@@@@@@@@@@@@@@@@@@@@@@@@@@@@@@@@@
%@@@@@@@@@@@@@@@@@@@@@@@@@@@@@@@@@@@@@@@@@@@@@@@@@@@@@@@@@@
%----------------------------------------------------------
\par {\it Proof.} Let $s:=(q+p)/2$; thus $q<s<p.$ We let $\tth:\RT\to\R_+$ denote a function
%----------------------------------------------------------
$$
\tth(z):=(\Mc[h^{s}](z))^{\frac{1}{s}},~~~z\in\RT.
$$
%----------------------------------------------------------
By the Lebesgue theorem, $h^s\le \Mc[h^{s}]$ a.e. so that
%----------------------------------------------------------
\bel{C-HTQ}
h\le\tth~~~\text{a.e. on}~~~ \RT.
\ee
%----------------------------------------------------------
Since $p/s>1$, by the Hardy-Littlewood maximal theorem,
%----------------------------------------------------------
$$
\|\tth\|_{\LPRT}=\left(\,\,\intl_{\RT}
(\Mc[h^{s}])^{\frac{p}{s}}\,dz\right)^{\frac1p}
\le C(p) \left(\,\,\intl_{\RT}
(h^{s})^{\frac{p}{s}}\,dz\right)^{\frac1p}
$$
%----------------------------------------------------------
proving that
%----------------------------------------------------------
\bel{LP-H}
\|\tth\|_{\LPRT}\le C\|h\|_{\LPRT}.
\ee
%----------------------------------------------------------
Since $\theta:=q/s<1$, by Theorem \reff{C-R}, the function
%----------------------------------------------------------
$$
\tth^q=(\Mc[h^{s}])^{\frac{q}{s}}=(\Mc[h^{s}])^{\theta}
$$
%----------------------------------------------------------
belongs to the class $A_1(\RT)$ and $\|\tth^s\|_{A_1(\RT)}\le C(\theta)=C(p)$. Hence, by \rf{4.C},
%----------------------------------------------------------
\bel{M-HQ}
\frac{1}{|Q|}\intl_Q \tth^s\,dz\le C\,\frac{1}{|K|}\intl_K \tth^s\,dz
\ee
%----------------------------------------------------------
provided $K,Q$ are arbitrary \sqs in $\RT$ such that $K\subset Q$.
%----------------------------------------------------------
\par By \rf{C-HTQ} and \rf{GS-A}, for every pair of connected bridges $T,T'\in\BRE$, $T\bcn T'$, the following inequality
%----------------------------------------------------------
\bel{G-D}
\|g(T)-g(T')\|\le \diam Q(T,T')\,\left(\frac{1}{|\gamma Q(T,T')|}
\intl_{\gamma Q(T,T')}\tth^q\,dz\right)^{\frac1q}
\ee
%----------------------------------------------------------
holds. Let
%----------------------------------------------------------
$$
\BRE\ni T,T'\mapsto K(T,T')\in \KE
$$
%----------------------------------------------------------
be a one-to-one mapping constructed in Proposition \reff{BR-TO-Q}. To every pair of connected bridges $\Br,\Br'\in\BRE$, $\Br\bcn\Br',$ this mapping assigns a \sq $K(\Br,\Br')$ which belongs to a family $\KE$ of pairwise disjoint \sqs satisfying conditions \rf{MP-1} and \rf{MP-2}.
%----------------------------------------------------------
\par Let us compare the \sq $\gamma Q(T,T')=Q(\A{T}, \gamma D(T,T'))$ with the \sq $K(\Br,\Br')$. By \rf{MP-1},
%----------------------------------------------------------
\bel{C-DMK}
~~~~~\diam K(T,T')\sim \diam \{\A{T},\B{T},\A{T'},\B{T'}\}=D(T,T')=\tfrac12\diam Q(T,T').
\ee
%----------------------------------------------------------
In particular, $|Q(T,T')|\sim |K(T,T')|$.
%----------------------------------------------------------
\par Note that, by \rf{MP-2}, the point $\A{T}\in \gamma_1 K(T,T')$ where $\gamma_1$ is an absolute constant. Hence
%----------------------------------------------------------
$$
\gamma Q(T,T')\cap (\gamma_1 K(T,T'))\ne\emp.
$$
%----------------------------------------------------------
Combining this property with equivalence \rf{C-DMK} we conclude that for some absolute constant $\gamma_2=\gamma_2(\gamma)\ge\gamma_1$ the following inclusion
%----------------------------------------------------------
$$
\gamma Q(T,T')\subset \tK:= \gamma_2 K(T,T')
$$
%----------------------------------------------------------
holds. Hence
%----------------------------------------------------------
$$
\frac{1}{|\gamma Q(T,T')|}\intl_{\gamma Q(T,T')} \tth^q\,dz\le C\,\frac{1}{|\tK|}\intl_{\tK} \tth^q\,dz.
$$
%----------------------------------------------------------
Since  $K(T,T')\subset\tK$, by \rf{M-HQ},
%----------------------------------------------------------
$$
\frac{1}{|\tK|}\intl_{\tK} \tth^q\,dz\le C(p)\,\frac{1}{|K(T,T')|}\intl_{K(T,T')} \tth^q\,dz.
$$
%----------------------------------------------------------
Combining this inequality with \rf{G-D} and \rf{C-DMK} we obtain
%----------------------------------------------------------
$$
\|g(T)-g(T')\|\le C\diam K(T,T')\,
\left(\frac{1}{|K(T,T')|}\intl_{K(T,T')} \tth^q\,dz\right)^{\frac1q}
$$
%----------------------------------------------------------
provided $T,T'\in\BRE, T\bcn T'$. Since $2<q<p$, by the H\"{o}lder inequality,
%----------------------------------------------------------
\be
\|g(T)-g(T')\|&\le& C\,\diam K(T,T')\,
\left(\frac{1}{|K(T,T')|}\intl_{K(T,T')} \tth^p\,dz\right)^{\frac1p}\nn\\
&=& C\,(\diam K(T,T'))^{1-2/p}\left(\intl_{K(T,T')} \tth^p\,dz\right)^{\frac1p}.\nn
\ee
%----------------------------------------------------------
Since $D(T,T')\sim \diam K(T,T')$, see \rf{C-DMK}, we obtain %----------------------------------------------------------
$$
\|g(T)-g(T')\|^p\le C\,D(T,T')^{p-2} \intl_{K(T,T')} \tth^p\,dz.
$$
%----------------------------------------------------------
Finally, by Proposition \reff{F-GN},
%----------------------------------------------------------
\be
\|F\|^p_{\LTP}&\le&
C\,\sum \left\{\,\,D(T,T')^{2-p}
\|g(T)-g(T')\|^p: T,T'\in\BRE, T\bcn T'\right\}\nn\\
&\le& C\,
\smed\left\{\intl_{K(T,T')}\tth^p\,dz: T,T'\in\BRE, T\bcn T'\right\}\nn\\
&=& C\smed_{K\in\KE}\,\, \intl_{K}\tth^p\,dz.\nn
\ee
%----------------------------------------------------------
\par Since the family $\KE$ consists of pairwise disjoint squares, we obtain
%----------------------------------------------------------
$$
\|F\|^p_{\LTP}\le
C\,\intl_{\RT}\tth^p\,dz=C\|\tth\|_{\LPRT}
$$
%----------------------------------------------------------
so that, by \rf{LP-H},
$\|F\|^p_{\LTP}\le C\,\|h\|_{\LPRT}^p$.
%----------------------------------------------------------
\par The proposition is proved.\bx
%----------------------------------------------------------
%@@@@@@@@@@@@@@@@@@@@@@@@@@@@@@@@@@@@@@@@@@@@@@@@@@@@@@@@@@
%@@@@@@@@@@@@@@@@@@@@@@@@@@@@@@@@@@@@@@@@@@@@@@@@@@@@@@@@@@
%@@@@@@@@@@@@@@@@@@@@@@@@@@@@@@@@@@@@@@@@@@@@@@@@@@@@@@@@@@
%@@@@@@@@@@@@@@@@@@@@@@@@@@@@@@@@@@@@@@@@@@@@@@@@@@@@@@@@@@
%----------------------------------------------------------
%@@@@@@@@@@@@@@@@@@@@@@@@@@@@@@@@@@@@@@@@@@@@@@@@@@@@@@@@@@
%@@@@@@@@@@@@@@@@@@@@@@@@@@@@@@@@@@@@@@@@@@@@@@@@@@@@@@@@@@
%@@@@@@@@@@@@@@@@@@@@@@@@@@@@@@@@@@@@@@@@@@@@@@@@@@@@@@@@@@
%@@@@@@@@@@@@@@@@@@@@@@@@@      @@@@@@@@@@@@@@@@@@@@@@@@@@@
%@@@@@@@@@@@@@@@@@@@@@@@          @@@@@@@@@@@@@@@@@@@@@@@@@
%@@@@@@@@@@@@@@@@@@@@@              @@@@@@@@@@@@@@@@@@@@@@@
%@@@@@@@@@@@@@@@@@@@     SECTION 8    @@@@@@@@@@@@@@@@@@@@@
%@@@@@@@@@@@@@@@@@@@@@              @@@@@@@@@@@@@@@@@@@@@@@
%@@@@@@@@@@@@@@@@@@@@@@@          @@@@@@@@@@@@@@@@@@@@@@@@@
%@@@@@@@@@@@@@@@@@@@@@@@@@      @@@@@@@@@@@@@@@@@@@@@@@@@@@
%@@@@@@@@@@@@@@@@@@@@@@@@@@@@@@@@@@@@@@@@@@@@@@@@@@@@@@@@@@
%@@@@@@@@@@@@@@@@@@@@@@@@@@@@@@@@@@@@@@@@@@@@@@@@@@@@@@@@@@
%@@@@@@@@@@@@@@@@@@@@@@@@@@@@@@@@@@@@@@@@@@@@@@@@@@@@@@@@@@
%----------------------------------------------------------
\SECTLONG{8. Sobolev-type selections of set-valued mappings and a decomposition }{of the sum $\VS:=\VLOP+\VLPME$.}{8}
%----------------------------------------------------------
\addtocontents{toc}{8. Sobolev-type selections of set-valued mappings and a decomposition
%----------------------------------------------------------
\par \hspace*{5mm}of the sum $\VS:=\VLOP+\VLPME$.\hfill \thepage\par}
%----------------------------------------------------------
\indent
%@@@@@@@@@@@@@@@@@@@@@@@@@@@@@@@@@@@@@@@@@@@@@@@@@@@@@@@@@@
%----------------------------------------------------------
\addtocontents{toc}{~~~~8.1. Sobolev-type selections. \hfill \thepage\par}
%----------------------------------------------------------
\par {\bf 8.1. Sobolev-type selections.} The extension criterion formulated in Proposition \reff{A-H} admits a geometrical reformulation in terms of {\it set-valued mappings and their selections.}
%----------------------------------------------------------
\par Let $f:E\to\R$ be a function on the set $E$ and let $T\in\BRE$ be a bridge (interior or exterior) with the ends at points $\A{T},\B{T}\in E$. We let $G_f(T)$ denote a straight line in $\RT$
%----------------------------------------------------------
\bel{D-GT}
G_f(T):=\{z\in\RT: \ip{z,\A{T}-\B{T}}=f(\A{T})-f(\B{T})\}.
\ee
%----------------------------------------------------------
\par Let $\ART$ be the family of all straight lines in $\RT$. The formula \rf{D-GT} defines a mapping
%----------------------------------------------------------
$$
G_f:\BRE\to \ART
$$
%----------------------------------------------------------
which to every bridge $T\in\BRE$ assigns a {\it subset} of $\RT$, the straight line $G_f(T)$. We refer to $G_f$ as a {\it set-valued mapping}.
%----------------------------------------------------------
\par Let $g:\BRE\to\RT$ be a (regular) mapping such that for every \br $T\in\BRE$ with the ends at points $\A{T},\B{T}\in E$ the following equality
%----------------------------------------------------------
$$
\ip{g(T),\A{T}-\B{T}}=f(\A{T})-f(\B{T})
$$
%----------------------------------------------------------
holds. By \rf{D-GT}, this property can be reformulated in the following way:
%----------------------------------------------------------
$$
g(T)\in G_f(T)~~~\text{for every bridge}~~~~T\in\BRE.
$$
%----------------------------------------------------------
We refer to the mapping $g$ as a {\it selection} of the set-valued mapping $G_f$.
%----------------------------------------------------------
\par We are interested in a special type of selections  which satisfy the inequality \rf{GS-A} from Proposition \reff{A-H}:
%----------------------------------------------------------
\bel{S-C}
\|g(T)-g(T')\|\le \diam Q(T,T')\,\left(\frac{1}{|\gamma Q(T,T')|} \intl_{\gamma Q(T,T')}h^q(z)dz\right)^{\frac1q}
\ee
%----------------------------------------------------------
provided $\gamma\ge 1$ is an absolute constant, $2<q<p$, $h$ is a non-negative $\LPRT$-function, and $T$ and $T'$ are arbitrary connected bridges from $\BRE$. This inequality is a natural analog of the Sobolev-Poincar\'e inequality \rf{S-2}. This analogy  motivates us to refer to a selection $g$ satisfying \rf{GS-A} as a {\it Sobolev-type selection} (with respect to the function $h$) of the set-valued mapping $G_f$.
%----------------------------------------------------------
\par Using this terminology we can formulate the following criterion for the trace space $\LTP|_E$.
%----------------------------------------------------------
%@@@@@@@@@@@@@@@@@@@@@@@@@@@@@@@@@@@@@@@@@@@@@@@@@@@@@@@@@@
%@@@@@@@@@@@@@@@@@@@@@@@@@@@@@@@@@@@@@@@@@@@@@@@@@@@@@@@@@@
%@@@@@@@@@@@@@@@@@@@@@@@@@@@@@@@@@@@@@@@@@@@@@@@@@@@@@@@@@@
%@@@@@@@@@@@@@@@@@@@@@@@@@@@@@@@@@@@@@@@@@@@@@@@@@@@@@@@@@@
%@@@@@@@@@@@@@@@@@@@@@@@@@@@@@@@@@@@@@@@@@@@@@@@@@@@@@@@@@@
%@@@@@@@@@@@@@@@@@@@@@@@@@@@@@@@@@@@@@@@@@@@@@@@@@@@@@@@@@@
%----------------------------------------------------------
\begin{claim}\lbl{CL-CR} Let $2<q<p<\infty$ and let $f$ be a function defined on a finite subset $E\subset\RT$. Then the following equivalence
%----------------------------------------------------------
$$
\|f\|_{\LTP|_E}\sim \inf\|h\|_{\LPRT}
$$
%----------------------------------------------------------
holds with constants depending only on $p$ and $q$. Here the infimum is taken over all non-negative functions $h\in \LPRT$ such that the set-valued mapping $G_f:\BRE\to\ART$ has a Sobolev-type selection with respect to $h$.
%----------------------------------------------------------
\end{claim}
%----------------------------------------------------------
%@@@@@@@@@@@@@@@@@@@@@@@@@@@@@@@@@@@@@@@@@@@@@@@@@@@@@@@@@@
%@@@@@@@@@@@@@@@@@@@@@@@@@@@@@@@@@@@@@@@@@@@@@@@@@@@@@@@@@@
%@@@@@@@@@@@@@@@@@@@@@@@@@@@@@@@@@@@@@@@@@@@@@@@@@@@@@@@@@@
%----------------------------------------------------------
\par {\it Proof.} The inequality
%----------------------------------------------------------
$$
\|f\|_{\LTP|_E}\le C\inf\|h\|_{\LPRT}
$$
%----------------------------------------------------------
directly follows from Proposition \reff{A-H}.
%----------------------------------------------------------
\par Let us prove the converse inequality. Let $\gamma\ge 1$ be an arbitrary constant. Let $F\in\LTP, F|_E=f,$ be an arbitrary extension of $f$ and let $T\in\BRE$ be a bridge with the ends at the points $\A{T},\B{T}\in E$. By the Lagrange theorem there exists a point $Z_T\in[\A{T},\B{T}]$ such that
%----------------------------------------------------------
\bel{CH-S}
\ip{\nabla F(Z_T),\A{T}-\B{T}}=f(\A{T})-f(\B{T}).
\ee
%----------------------------------------------------------
We put
%----------------------------------------------------------
$$
g(T):=\nabla F(Z_T),~~~T\in\BRE.
$$
%----------------------------------------------------------
Then, by \rf{CH-S}, $g$ is a {\it selection} of the set-valued mapping $G_f$, see \rf{D-GT}.
%----------------------------------------------------------
\par Now let us consider an arbitrary pair of connected bridges $T,T'\in\BRE, T\bcn T'.$ Clearly the \sq
%----------------------------------------------------------
$$
Q(T,T')=Q(\A{T},D(T,T'))
$$
%----------------------------------------------------------
where $D(T,T')=\diam \{\A{T},\B{T},\A{T'},\B{T'}\}$, contains the set $\{\A{T},\B{T},\A{T'},\B{T'}\}$. Hence $Z_T,Z_{T'}\in Q(T,T')$ as well so that, by the Sobolev-Poincar\'e inequality \rf{SP-IN},
%@@@@@@@@@@@@@@@@@@@@@@@@@@@@@@@@@@@@@@@@@@@@@@@@@@@@@@@@@@
%----------------------------------------------------------
\be
\|g(T)-g(T')\|&=&\|\nabla F(Z_T)-\nabla F(Z_{T'})\|\nn\\
&\le& C\diam Q(T,T') \left(\frac{1}{|Q(T,T')|} \intl_{Q(T,T')}(\nabla^2F)^q\,dz\right)^{\frac{1}{q}}\nn\\&\le&
C(\gamma)\diam Q(T,T') \left(\frac{1}{|\gamma Q(T,T')|} \intl_{\gamma Q(T,T')} (\nabla^2F)^q\,dz\right)^{\frac{1}{q}}.\nn
\ee
%----------------------------------------------------------
This shows that $g$ satisfies inequality \rf{S-C} with $h:=C(\gamma)\nabla^2 F$. Hence
%----------------------------------------------------------
$$
\|h\|_{\LPRT}=C\|\nabla^2 F\|_{\LPRT}=C\|F\|_{\LTP}.
$$
%----------------------------------------------------------
Taking the infimum in this inequality over all functions $F\in\LTP$ such that $F|_E=f$ we obtain the required inequality
%----------------------------------------------------------
$$
\inf\|h\|_{\LPRT}\le C\|f\|_{\LTP|_E}
$$
%----------------------------------------------------------
with $C=C(q,\gamma)$. The claim is proved.\bx\medskip
%----------------------------------------------------------
%@@@@@@@@@@@@@@@@@@@@@@@@@@@@@@@@@@@@@@@@@@@@@@@@@@@@@@@@@@
%@@@@@@@@@@@@@@@@@@@@@@@@@@@@@@@@@@@@@@@@@@@@@@@@@@@@@@@@@@
%@@@@@@@@@@@@@@@@@@@@@@@@@@@@@@@@@@@@@@@@@@@@@@@@@@@@@@@@@@
%@@@@@@@@@@@@@@@@@@@@@@@@@@@@@@@@@@@@@@@@@@@@@@@@@@@@@@@@@@
%@@@@@@@@@@@@@@@@@@@@@@@@@@@@@@@@@@@@@@@@@@@@@@@@@@@@@@@@@@
%@@@@@@@@@@@@@@@@@@@@@@@@@@@@@@@@@@@@@@@@@@@@@@@@@@@@@@@@@@
%----------------------------------------------------------
\begin{remark} {\em Given $q\in(2,p)$ consider
a function
%----------------------------------------------------------
$$
\delta_h(T,T'):=\diam Q(T,T')\,\left(\frac{1}{|Q(T,T')|}
\intl_{Q(T,T')}h^q(z)dz\right)^{\frac1q},~~~~~T,T'\in\BRE.
$$
%----------------------------------------------------------
Functions of such a kind have been studied in the author's paper \cite{S4}. Using the methods of this work one can show that the function $\delta$ is equivalent to a certain {\it metric} $\tdl_h$ on $\BRE$ provided $h^q\in A_1(\RT)$.
%----------------------------------------------------------
\par Thus in this case inequality \rf{S-C} is equivalent to the inequality
%----------------------------------------------------------
$$
\|g(T)-g(T')\|\le \tdl_h(T,T'),~~~~T,T'\in\BRE.
$$
%----------------------------------------------------------
In other words, the mapping $g:\BRE\to\RT$ is {\it Lipschitz} with respect to the metric $\tdl$ so that
$g$ is a {\it Lipschitz selection} of the set valued mapping $G_f$. Our aim is to {\it find the order of magnitude of the minimal $\LPRT$-norm of the function $h$ provided there exists a Lipschitz selection of $G_f$ with respect to the metric $\tdl_h$.}
%----------------------------------------------------------
\par This type of selection problems is a generalization of the so-called Lipschitz selection problem. We refer the reader to the papers \cite{S-GAN,S-GF,S-JF} and references therein for numerous results related to this problem  and techniques for obtaining them. In particular, our approach to the Sobolev-type selection problem which we develop in the next subsection is a generalization of ideas and methods suggested in these papers for solution to the Lipschitz selection problem.\rbx}
%----------------------------------------------------------
\end{remark}
%----------------------------------------------------------
\bigskip
%----------------------------------------------------------
%@@@@@@@@@@@@@@@@@@@@@@@@@@@@@@@@@@@@@@@@@@@@@@@@@@@@@@@@@@
%@@@@@@@@@@@@@@@@@@@@@@@@@@@@@@@@@@@@@@@@@@@@@@@@@@@@@@@@@@
%@@@@@@@@@@@@@@@@@@@@@@@@@@@@@@@@@@@@@@@@@@@@@@@@@@@@@@@@@@
%@@@@@@@@@@@@@@@@@@@@@@@@@@@@@@@@@@@@@@@@@@@@@@@@@@@@@@@@@@
%@@@@@@@@@@@@@@@@@@@@@@@@@@@@@@@@@@@@@@@@@@@@@@@@@@@@@@@@@@
%@@@@@@@@@@@@@@@@@@@@@@@@@@@@@@@@@@@@@@@@@@@@@@@@@@@@@@@@@@
%@@@@@@@@@@@@@@@@@@@@@@@@@@@@@@@@@@@@@@@@@@@@@@@@@@@@@@@@@@
%@@@@@@@@@@@@@@@@@@@@@@@@@@@@@@@@@@@@@@@@@@@@@@@@@@@@@@@@@@
%@@@@@@@@@@@@@@@@@@@@@@@@@@@@@@@@@@@@@@@@@@@@@@@@@@@@@@@@@@
%@@@@@@@@@@@@@@@@@@@@@@@@@@@@@@@@@@@@@@@@@@@@@@@@@@@@@@@@@@
%----------------------------------------------------------
\par {\bf 8.2\,\, The mapping $\Tc(f)$ and its norm in $\VLOP+\VLPME$.}
%----------------------------------------------------------
\addtocontents{toc}{~~~~8.2\,\, The mapping $\Tc(f)$ and its norm in $\VLOP+\VLPME$. \hfill \thepage\par}
%----------------------------------------------------------
Following the scheme of the proof given in Section 2, in this subsection we show that the mapping $\Tc(f)$ defined in Section 2 belongs to the space $\VS:=\VLOP+\VLPME$ and its norm in this space is bounded by $C(p)\lambda^{\frac1p}$, see Proposition \reff{TAU-L}.
%----------------------------------------------------------
\par We recall that, by Proposition  \reff{BR-TO-Q}, there exist a one-to-one mapping,
%----------------------------------------------------------
\bel{M-ON}
\Br,\Br'\in\BRE, \Br\bcn\Br'~\leftrightarrow~~\text{a square}~K(\Br,\Br')\in\KE,
\ee
%----------------------------------------------------------
between the family of (non-ordered) pairs of connected bridges $\Br,\Br'\in\BRE, \Br\bcn\Br'$ and the family $\KE$ of well-separated \sqs $K(\Br,\Br')\in\KE$ satisfying conditions \rf{MP-1} and \rf{MP-2}.
%----------------------------------------------------------
\par The mapping \rf{M-ON} enables us to change our notation related to pairs of connected bridges. Since this mapping is one-to-one, we may consider the converse mapping which to every \sq $K\in\KE$ assigns a pair of connected bridges $T_K,T'_K\in\BRE, T_K\bcn T'_K$. We denote the ends of these bridges by $A_{T_K},B_{T_K}$ and $A_{T'_K},B_{T'_K}$ respectively. Since $T_K\bcn T'_K$,
we have
%----------------------------------------------------------
$$
\{\A{T_K},\B{T_K}\}\cap \{\A{T'_K},\B{T'_K}\}\ne\emp,
$$
%----------------------------------------------------------
i.e., the bridges $T_K$ and $T'_K$ have a common end. Thus the set of ends
%----------------------------------------------------------
$$
\{\A{T_K},\B{T_K},\A{T'_K},\B{T'_K}\}
$$
%----------------------------------------------------------
consists of either three or two points.
%----------------------------------------------------------
\par Let us consider two cases.\medskip
%----------------------------------------------------------
\par {\it The first case.} $\{\A{T_K},\B{T_K}\}\ne
\{\A{T'_K},\B{T'_K}\},$ i.e.,
%----------------------------------------------------------
$$
\#\{\A{T_K},\B{T_K},\A{T'_K},\B{T'_K}\}=3.
$$
%----------------------------------------------------------
\par In this case there exists the unique common end
%----------------------------------------------------------
\bel{CK-T}
C^{(K)}:=\{\A{T_K},\B{T_K}\}\cap \{\A{T'_K},\B{T'_K}\}.
\ee
%----------------------------------------------------------
\par Let
%----------------------------------------------------------
\bel{F-TRI}
A^{(K)}~~~\text{and}~~~B^{(K)}~~~\text{be two remaining points of the set of ends.}
\ee
%----------------------------------------------------------
Thus this set forms a triangle
$\Delta(K):=\Delta\{A^{(K)},B^{(K)},C^{(K)}\}$ (true or \dg). %----------------------------------------------------------
\par Since $T_K$ and $T'_K$ are connected bridges, according to our definitions the side $[A^{(K)},B^{(K)}]$
of the triangle $\Delta(K)$ (which is opposite to the vertex $C^{(K)}$) {\it is not the smallest side of $\Delta(K)$}.
%----------------------------------------------------------
\par Let $\alpha_K$ be the angle in $\Delta(K)$ corresponding to the vertex $C^{(K)}$. (In particular,  $\alpha_K=0$ whenever $\Delta(K)$ is a \dg triangle.) Then, by \rf{SIN},
%----------------------------------------------------------
\bel{SIN-K}
|\sin \alpha_K|\sim \cu_{\Delta(K)}\,\diam \Delta(K).
\ee
%----------------------------------------------------------
Recall that $\cu_{\Delta(K)}=1/R_{\Delta(K)}$ is the Menger curvature of the triangle $\Delta(K)$, see \rf{M-CURV}.
%----------------------------------------------------------
\par Also we recall that, by \rf{MP-1},
%----------------------------------------------------------
\bel{K-D}
\diam \Delta(K)=\diam \{\A{T_K},\B{T_K}, \A{T'_K},\B{T'_K}\} \sim \diam K.
\ee
%----------------------------------------------------------
In turn, by \rf{MP-2},
%----------------------------------------------------------
\bel{D-INK}
\Delta(K)=\Delta\{A^{(K)},B^{(K)},C^{(K)}\}\subset
\gamma K.
\ee
%----------------------------------------------------------
This inclusion and equivalence \rf{K-D} imply the following useful inclusions
%----------------------------------------------------------
\bel{Q-T-IN}
Q(T_K,T'_K)\subset \gamma K~~~~\text{and}~~~~K\subset \gamma Q(T_K,T'_K).
\ee
%----------------------------------------------------------
In all cases $\gamma>0$ is an absolute constant. Recall that
%----------------------------------------------------------
$$
Q(T_K,T'_K)=Q(\A{T_K},D(T_K,T'_K))
$$
%----------------------------------------------------------
where
%----------------------------------------------------------
$$
D(T_K,T'_K)=\diam \{\A{T_K},\B{T_K},\A{T'_K},\B{T'_K}\}
$$
%----------------------------------------------------------
is defined in \rf{SQ-T}.
%----------------------------------------------------------
\par Note that, by \rf{SIN-K} and \rf{K-D},
%----------------------------------------------------------
\bel{K-S}
|\sin \alpha_K|\sim \cu_{\Delta(K)}\,\diam K=\frac{1}{R_{\Delta(K)}}\diam K\,.
\ee
%----------------------------------------------------------
%@@@@@@@@@@@@@@@@@@@@@@@@@@@@@@@@@@@@@@@@@@@@@@@@@@@@@@@@@@
%@@@@@@@@@@@@@@@@@@@@@@@@@@@@@@@@@@@@@@@@@@@@@@@@@@@@@@@@@@
%@@@@@@@@@@@@@@@@@@@@@@@@@@@@@@@@@@@@@@@@@@@@@@@@@@@@@@@@@@
%@@@@@@@@@@@@@@@@@@@@@@@@@@@@@@@@@@@@@@@@@@@@@@@@@@@@@@@@@@
%----------------------------------------------------------
\medskip
%----------------------------------------------------------
\par {\it The second case.}
%----------------------------------------------------------
$$
\#\{\A{T_K},\B{T_K},\A{T'_K},\B{T'_K}\}=2,
$$
%----------------------------------------------------------
i.e., $\{\A{T_K},\B{T_K}\}=\{\A{T'_K},\B{T'_K}\}.$ In this case we denote the ends of the bridges by  $A^{(K)}$ and $B^{(K)}$ and put $C^{(K)}:=A^{(K)}$. Then, by \rf{MP-1} and \rf{MP-2},
%----------------------------------------------------------
$$
\|A^{(K)}-B^{(K)}\|=\|C^{(K)}-B^{(K)}\|\sim \diam K,
$$
%----------------------------------------------------------
and, as in the first case,
%----------------------------------------------------------
$$
\{A^{(K)},B^{(K)},C^{(K)}\}\subset \gamma K.
$$
%----------------------------------------------------------
\bigskip
\par We turn to definitions of the measure $\mu_E$ on $\RT$ and the mapping $\Tc:\RT\to\RT$. These definitions are similar to those given in Section 3, see \rf{D-MU}, \rf{DGV-1} and \rf{DGV-2}.
%----------------------------------------------------------
\par Let
%----------------------------------------------------------
\bel{CC-E}
\Cc_E:=\{c_K: K\in \Kc_E\}
\ee
%----------------------------------------------------------
be the family of centers of all \sqs from $\Kc_E$. By $\mu_E$ we denote a discrete Borel measure with $\supp \mu_E\subset \Cc_E$ such that
%----------------------------------------------------------
\bel{D-MU-S}
\mu_E(\{c_K\}):=\cu_{\Delta(K)}^p\,|K|~~~~~\text{for every \sq}~~ K\in\Kc_E.
\ee
%----------------------------------------------------------
Thus
%----------------------------------------------------------
\bel{MU-SET}
\mu_E(S)=\sum\left\{\cu_{\Delta(K)}^p\,|K|: K\in\Kc_E,
c_K\in S\right\}
\ee
%----------------------------------------------------------
for every set $S\subset\R^2$. We refer to $\mu_E$ as {\it the Menger curvature measure} generated by the set $E$.
%----------------------------------------------------------
\par Let us define the mapping $\Tc(f):\RT\to\RT$. For every \sq $K\in\Kc_E$ we put
%----------------------------------------------------------
\bel{DGV-1-S}
\Tc(c_K;f):=\left \{
%----------------------------------------------------------
\begin{array}{ll}
\nabla P_{\Delta(K)}[f],&\Delta(K)~~\text{is a true triangle},\\\\
0,&\Delta(K)~~\text{is a degenerate triangle.}
\end{array}
%----------------------------------------------------------
\right.
\ee
%----------------------------------------------------------
For every $x\in \RT\setminus \Cc_E$ we set
%----------------------------------------------------------
\bel{DGV-2-S}
\Tc(x;f):=0.
\ee
%----------------------------------------------------------
\par Let $\VS:=\VLOP+\VLPME$. Our aim is to show that
$\|\Tc(f)\|_{\VS}\le C(p)\lambda$ provided the part (b) of the hypothesis of Theorem \reff{S-MAIN} holds. We present a proof of this statement in Proposition \reff{TAU-L}. The  main ingredient of this proof is the following
%----------------------------------------------------------
%@@@@@@@@@@@@@@@@@@@@@@@@@@@@@@@@@@@@@@@@@@@@@@@@@@@@@@@@@@
%@@@@@@@@@@@@@@@@@@@@@@@@@@@@@@@@@@@@@@@@@@@@@@@@@@@@@@@@@@
%@@@@@@@@@@@@@@@@@@@@@@@@@@@@@@@@@@@@@@@@@@@@@@@@@@@@@@@@@@
%@@@@@@@@@@@@@@@@@@@@@@@@@@@@@@@@@@@@@@@@@@@@@@@@@@@@@@@@@@
%@@@@@@@@@@@@@@@@@@@@@@@@@@@@@@@@@@@@@@@@@@@@@@@@@@@@@@@@@@
%----------------------------------------------------------
\begin{theorem}\lbl{S-CMS} Let $2<p<\infty$ and let $\tau>0$. Let $\mu$ be a non-trivial non-negative Borel measure on $\RT$. A mapping $\Vc:\RT\to\RT$ belongs to the space $\VS=\VLOP+\VLPM$ provided $\Vc\in \vec{L}_{p,loc}(\RT;\mu)$ and there exists a constant $\lambda>0$ such that the following statement is true for a certain absolute positive constant $\gamma$:
%----------------------------------------------------------
\par Let $\Qc$ be a finite family of pairwise disjoint \sqs in $\RT$. Suppose that to each \sq $Q\in\Qc$ we have arbitrarily assigned two \sqs $Q',Q''\in\Qc$ such that
$Q'\cup Q''\subset \gamma Q$ and
%----------------------------------------------------------
\bel{R-W}
(\diam Q')^{p-2}\mu(Q')+(\diam Q'')^{p-2}\mu(Q'')\le \tau.
\ee
%----------------------------------------------------------
Then the following inequality
%----------------------------------------------------------
$$
\sbig_{Q\in\Qc}\,\,
\left(\frac{\diam Q' \diam Q''}{\diam Q}\right)^{p-n} \iint \limits_{Q'\times Q''}
\|\Vc(x)-\Vc(y)\|^p\, d\mu(x)d\mu(y)
\le \lambda
$$
%----------------------------------------------------------
holds.
%----------------------------------------------------------
\par Furthermore,
$\|\Vc\|_{\VS}\le C(p,\tau)\lambda^{\frac{1}{p}}$.
%----------------------------------------------------------
\end{theorem}
%----------------------------------------------------------
%@@@@@@@@@@@@@@@@@@@@@@@@@@@@@@@@@@@@@@@@@@@@@@@@@@@@@@@@@@
%@@@@@@@@@@@@@@@@@@@@@@@@@@@@@@@@@@@@@@@@@@@@@@@@@@@@@@@@@@
%@@@@@@@@@@@@@@@@@@@@@@@@@@@@@@@@@@@@@@@@@@@@@@@@@@@@@@@@@@
%@@@@@@@@@@@@@@@@@@@@@@@@@@@@@@@@@@@@@@@@@@@@@@@@@@@@@@@@@@
%----------------------------------------------------------
\par For scalar mappings, i.e., for functions on $\RT$, the proof of this theorem is given in \cite{S3}, Theorem 6.1.
We obtain the statement of Theorem \reff{S-CMS} by applying this result to every component of the mapping $\Vc$.
%----------------------------------------------------------
\par We also note that it suffices to prove this theorem for a certain fixed value of the parameter $\tau$, say for $\tau=1$. Then the general case follows from this particular case by transition to a new measure $\tilde{\mu}=\tfrac{1}{\tau}\,\mu$.
%----------------------------------------------------------
%@@@@@@@@@@@@@@@@@@@@@@@@@@@@@@@@@@@@@@@@@@@@@@@@@@@@@@@@@@
%@@@@@@@@@@@@@@@@@@@@@@@@@@@@@@@@@@@@@@@@@@@@@@@@@@@@@@@@@@
%@@@@@@@@@@@@@@@@@@@@@@@@@@@@@@@@@@@@@@@@@@@@@@@@@@@@@@@@@@
%@@@@@@@@@@@@@@@@@@@@@@@@@@@@@@@@@@@@@@@@@@@@@@@@@@@@@@@@@@
%----------------------------------------------------------
\begin{proposition}\lbl{TAU-L} Let $\mu_E$ and $\Tc(f)$ be the measure and the mapping defined by the formulas \rf{D-MU-S}, \rf{DGV-1-S} and \rf{DGV-2-S}. The mapping $\Tc(f):\RT\to\RT$ belongs to the space
%----------------------------------------------------------
$$
\VS:=\VLOP+\VLPME
$$
%----------------------------------------------------------
provided there exists a constant $\lambda>0$ such that the condition (b) of Theorem \reff{S-MAIN} is satisfied.
%----------------------------------------------------------
\par Furthermore,
%----------------------------------------------------------
$$
\|\Tc(f)\|_{\VS}\le C(p)\lambda^{\frac1p}.
$$
%----------------------------------------------------------
\end{proposition}
%----------------------------------------------------------
%@@@@@@@@@@@@@@@@@@@@@@@@@@@@@@@@@@@@@@@@@@@@@@@@@@@@@@@@@@
%@@@@@@@@@@@@@@@@@@@@@@@@@@@@@@@@@@@@@@@@@@@@@@@@@@@@@@@@@@
%@@@@@@@@@@@@@@@@@@@@@@@@@@@@@@@@@@@@@@@@@@@@@@@@@@@@@@@@@@
%@@@@@@@@@@@@@@@@@@@@@@@@@@@@@@@@@@@@@@@@@@@@@@@@@@@@@@@@@@
%----------------------------------------------------------
\par {\it Proof.} By Theorem \reff{S-CMS}, it suffices to show that for every finite family $\Qc$ of pairwise disjoint \sqs in $\RT$ and for any choice of \sqs $Q',Q''\in\Qc$ such that $Q'\cup Q''\subset \gamma Q$ and
%----------------------------------------------------------
\bel{M-DM-N}
(\diam Q')^{p-2}\mu_E(Q')+(\diam Q'')^{p-2}\mu_E(Q'')\le 1
\ee
%----------------------------------------------------------
the following inequality
%----------------------------------------------------------
$$
I:=\smed_{Q\in\Qc}\,\,
\left(\frac{\diam Q' \diam Q''}{\diam Q}\right)^{p-n} \iint \limits_{Q'\times Q''}
\|\Tc(x;f)-\Tc(y;f)\|^p\, d\mu(x)d\mu(y)
\le \lambda
$$
%----------------------------------------------------------
holds.
%----------------------------------------------------------
\par We note that, by \rf{MU-SET}, for every \sq $Q\in\Qc$ we have
%----------------------------------------------------------
$$
\mu_E(Q)=\sum\left\{\cu_{\Delta(K)}^p\,|K|: K\in\Kc,
c_K\in Q\right\}
$$
%----------------------------------------------------------
so that $\mu_E(Q)=\sigma_p(Q;\Kc),$ see \rf{SGM}. Hence
%----------------------------------------------------------
\be
&&(\diam Q')^{p-2}\mu_E(Q')+(\diam Q'')^{p-2}\mu_E(Q'')\nn\\
&=&
(\diam Q')^{p-2}\sigma_p(Q';\Kc)+
(\diam Q'')^{p-2}\sigma_p(Q'';\Kc)\nn
\ee
%----------------------------------------------------------
proving that inequality \rf{M-DM-N} is equivalent to inequality \rf{CV-S} from part (b) of Theorem \reff{S-MAIN}.
%----------------------------------------------------------
\par In turn, by definitions \rf{D-MU-S}, \rf{DGV-1-S} and \rf{DGV-2-S},
%----------------------------------------------------------
\be
&&
\iint \limits_{Q'\times Q''}
\|\Tc(x;f)-\Tc(y;f)\|^p\, d\mu(x)d\mu(y)\nn\\
&=&
\sum_{\substack {K'\in\,\Kc\\c_{K'}\in Q'}}\,\,
\sum_{\substack {K''\in\,\Kc\\c_{K''}\in Q''}}
\|\nabla P_{\Delta(K')}[f]-
\nabla P_{\Delta(K'')}[f]\|^p\,
\cu_{\Delta(K')}^p |K'|\,\cu_{\Delta(K'')}^p |K''|\nn\\
&=&S_p(f:Q',Q'';\Kc).\nn
\ee
%----------------------------------------------------------
Hence
%----------------------------------------------------------
$$
I=\smed_{Q\in\Qc}\,\,
\left(\frac{\diam Q' \diam Q''}{\diam Q}\right)^{p-n} S_p(f:Q',Q'';\Kc)
$$
%----------------------------------------------------------
so that, by \rf{CRT-S}, $I\le\lambda$ proving the proposition.\bx
%----------------------------------------------------------
\bigskip
%----------------------------------------------------------
%@@@@@@@@@@@@@@@@@@@@@@@@@@@@@@@@@@@@@@@@@@@@@@@@@@@@@@@@@@
%@@@@@@@@@@@@@@@@@@@@@@@@@@@@@@@@@@@@@@@@@@@@@@@@@@@@@@@@@@
%@@@@@@@@@@@@@@@@@@@@@@@@@@@@@@@@@@@@@@@@@@@@@@@@@@@@@@@@@@
%@@@@@@@@@@@@@@@@@@@@@@@@@@@@@@@@@@@@@@@@@@@@@@@@@@@@@@@@@@
%----------------------------------------------------------
%@@@@@@@@@@@@@@@@@@@@@@@@@@@@@@@@@@@@@@@@@@@@@@@@@@@@@@@@@@
%@@@@@@@@@@@@@@@@@@@@@@@@@@@@@@@@@@@@@@@@@@@@@@@@@@@@@@@@@@
%@@@@@@@@@@@@@@@@@@@@@@@@@@@@@@@@@@@@@@@@@@@@@@@@@@@@@@@@@@
%@@@@@@@@@@@@@@@@@@@@@@@@@@@@@@@@@@@@@@@@@@@@@@@@@@@@@@@@@@
%----------------------------------------------------------
\par {\bf 8.3. Pre-selections and selections of the set-valued mapping $G_f$.}
%----------------------------------------------------------
\addtocontents{toc}{~~~~8.3. Pre-selections and selections of the set-valued mapping $G_f$.\hfill\thepage\\\par}
%----------------------------------------------------------
We turn to the last step of the proof of Theorem \reff{S-MAIN}. Following to the scheme presented in Section 2, at this step we define a pre-selection $\tg(f):\BRE\to\RT$ and a selection
$g(f):\BRE\to\RT$ of the set-valued mapping $G_f:\BRE\to\ART$.
%----------------------------------------------------------
\par By Proposition \reff{TAU-L}, the mapping $\Tc(f):\RT\to\RT$ belongs to the space
%----------------------------------------------------------
$$
\VS:=\VLOP+\VLPME,
$$
%----------------------------------------------------------
and its norm in this space is bounded by $C(p)\lambda$. Thus there exist mappings
%----------------------------------------------------------
$$
\Tc_1(f):\RT\to\RT~~~~~\text{and}~~~~~~\Tc_2(f):\RT\to\RT
$$
%----------------------------------------------------------
such that
%----------------------------------------------------------
$$
\Tc(f)=\Tc_1(f)+\Tc_2(f)
$$
%----------------------------------------------------------
and the following conditions are satisfied:
%----------------------------------------------------------
\par (i). The mapping $\Tc_1(f)\in \VLOP$ and
%----------------------------------------------------------
\bel{N-1}
\|\Tc_1(f)\|_{\VLOP}\le C(p)\lambda^{\frac1p};
\ee
%----------------------------------------------------------
%@@@@@@@@@@@@@@@@@@@@@@@@@@@@@@@@@@@@@@@@@@@@@@@@@@@@@@@@@@
%@@@@@@@@@@@@@@@@@@@@@@@@@@@@@@@@@@@@@@@@@@@@@@@@@@@@@@@@@@
%----------------------------------------------------------
\par (ii). The mapping $\Tc_2(f)\in \VLPME$ and
%----------------------------------------------------------
\bel{N-2}
\|\Tc_2(f)\|_{\VLPME}\le C(p)\lambda^{\frac1p}.
\ee
%----------------------------------------------------------
\par Now by $\tg(f):\BRE\to\RT$ we denote a mapping which to every bridge $T\in\BRE$ with ends at points $\A{T},\B{T}\in E$ assigns a vector
%----------------------------------------------------------
\bel{D-PRS}
\tg(T;f):=\Tc_1(\A{T};f).
\ee
%----------------------------------------------------------
We refer to the mapping $\tg(f):\BRE\to\RT$ as a {\it pre-selection} of the the set-valued mapping $G_f.$
%----------------------------------------------------------
\par Finally, we define the required selection $g(f):\BRE\to\RT$ of the set-valued mapping $G_f:\BRE\to\ART$ by the following formula:
%----------------------------------------------------------
\bel{D-SEL}
g(T;f):=\PR(\tg(T;f);G_f(T)),~~~~T\in\BRE.
\ee
%----------------------------------------------------------
Here given a straight line $\ell\subset\RT$ and a point $x\in\RT$ by $\PR(x;\ell)$ we denote {\it the orthogonal projection of $x$ onto $\ell$}.\medskip
%----------------------------------------------------------
\par Our aim is to show that the selection $g(f)$ is a Sobolev-type mapping with respect to a certain non-negative $\LPRT$-function $h(f)$ such that $\|h(f)_{\LPRT}\|\le C(p)\lambda^{\frac1p}$. Hence, by Claim \reff{CL-CR} (or by Proposition \reff{A-H}), we obtain the required estimate
%----------------------------------------------------------
\bel{FIN-R}
\|f\|_{\LTP|_E}\le C(p)\lambda^{\frac1p}.
\ee
%----------------------------------------------------------
\par We construct the function $h(f)$ as a sum of three non-negative functions
%----------------------------------------------------------
$$
h_1(f),h_2(f),h_3(f)\in\LPRT.
$$
%----------------------------------------------------------
Their definitions are motivated by inequalities \rf{N-1}, \rf{N-2} and by part (i) of the hypo\-the\-sis of Theorem \reff{MAIN} respectively.\medskip
%----------------------------------------------------------
%@@@@@@@@@@@@@@@@@@@@@@@@@@@@@@@@@@@@@@@@@@@@@@@@@@@@@@@@@@
%@@@@@@@@@@@@@@@@@@@@@@@@@@@@@@@@@@@@@@@@@@@@@@@@@@@@@@@@@@
%@@@@@@@@@@@@@@@@@@@@@@@@@@@@@@@@@@@@@@@@@@@@@@@@@@@@@@@@@@
%@@@@@@@@@@@@@@@@@@@@@@@@@@@@@@@@@@@@@@@@@@@@@@@@@@@@@@@@@@
%@@@@@@@@@@@@@@@@@@@@@@@@@@@@@@@@@@@@@@@@@@@@@@@@@@@@@@@@@@
%@@@@@@@@@@@@@@@@@@@@@@@@@@@@@@@@@@@@@@@@@@@@@@@@@@@@@@@@@@
%----------------------------------------------------------
\par We begin with a definition of the function $h_1(f):\RT\to\RT$. We put
%----------------------------------------------------------
\bel{H-1}
h_1(f):=\|\nabla\Tc_1(f)\|.
\ee
%----------------------------------------------------------
Here as usual for a vector function $\Tc_1(f)$ the sign $\nabla$ means differentiation according to coordinates; thus
%----------------------------------------------------------
$$
h_1(f)=\|(\nabla\Tc_1^{(1)}(f),\nabla\Tc_1^{(2)}(f))\|
$$
%----------------------------------------------------------
provided $\Tc_1(f)=(\Tc_1^{(1)}(f),\Tc_1^{(2)}(f))$. By \rf{N-1} and definition \rf{N-VLOP},
%----------------------------------------------------------
\be
\|\Tc_1(f)\|_{\VLOP}&=&
\|\Tc_1^{(1)}(f)\|_{\LOP}+\|\Tc_1^{(2)}(f)\|_{\LOP}\nn\\
&=&
\|\nabla\Tc_1^{(1)}(f)\|_{\LPRT}
+\|\nabla\Tc_1^{(2)}(f)\|_{\LPRT}\nn\\
&\le&C(p)\lambda^{\frac1p}.\nn
\ee
%----------------------------------------------------------
Hence
%----------------------------------------------------------
\bel{H1-LP}
\|h_1(f)\|_{\LPRT}\le C(p)\lambda^{\frac1p}.
\ee
%----------------------------------------------------------
%@@@@@@@@@@@@@@@@@@@@@@@@@@@@@@@@@@@@@@@@@@@@@@@@@@@@@@@@@@
%@@@@@@@@@@@@@@@@@@@@@@@@@@@@@@@@@@@@@@@@@@@@@@@@@@@@@@@@@@
%@@@@@@@@@@@@@@@@@@@@@@@@@@@@@@@@@@@@@@@@@@@@@@@@@@@@@@@@@@
%@@@@@@@@@@@@@@@@@@@@@@@@@@@@@@@@@@@@@@@@@@@@@@@@@@@@@@@@@@
%@@@@@@@@@@@@@@@@@@@@@@@@@@@@@@@@@@@@@@@@@@@@@@@@@@@@@@@@@@
%@@@@@@@@@@@@@@@@@@@@@@@@@@@@@@@@@@@@@@@@@@@@@@@@@@@@@@@@@@
\medskip
%----------------------------------------------------------
\par Let us introduce the function $h_2(f):\RT\to\RT$. By \rf{N-2},
%----------------------------------------------------------
$$
\|\Tc_2(f)\|_{\VLPME}=\|\Tc(f)-\Tc_1(f)\|_{\VLPME}\le C(p)\lambda^{\frac1p}.
$$
%----------------------------------------------------------
By definitions \rf{D-MU-S} and \rf{DGV-1-S},
%----------------------------------------------------------
\be
\|\Tc_2(f)\|^p_{\VLPME}&=&
\intl_{\RT}\|\Tc_2(x,f)\|^p\,d\mu_E(x)
=\intl_{\RT}\|\Tc(x,f)-\Tc_1(x,f)\|^p\,d\mu_E(x)\nn\\
&=&
\sum_{K\in\Kc}
\|\Tc(c_K,f)-\Tc_1(c_K,f)\|^p\,\cu^p_{\Delta(K)}|K|\nn\\
&=&\sum_{K\in\Kc}
\|\nabla P_{\Delta(K)}[f]-\Tc_1(c_K,f)\|^p\,
\cu^p_{\Delta(K)}|K|.\nn
\ee
%----------------------------------------------------------
This inequality motivates us to define the function $h_2(f):\RT\to\RT$ by the following formula
%----------------------------------------------------------
\bel{H-2}
h_2(f):=\sum_{K\in\Kc}
\|\nabla P_{\Delta(K)}[f]-\Tc_1(c_K,f)\|\,
\cu_{\Delta(K)}\,\chi_K.
\ee
%----------------------------------------------------------
Since the \sqs of the collection $\Kc$ are pairwise disjoint,
%----------------------------------------------------------
$$
\|h_2(f)\|_{\LPRT}=\sum_{K\in\Kc}
\|\nabla P_{\Delta(K)}[f]-\Tc_1(c_K,f)\|^p\,
\cu^p_{\Delta(K)}|K|
$$
%----------------------------------------------------------
proving that $\|h_2(f)\|_{\LPRT}=\|\Tc_2(f)\|_{\VLPME}$. Hence, by \rf{N-2},
%----------------------------------------------------------
\bel{H2-LP}
\|h_2(f)\|_{\LPRT}\le C(p)\lambda^{\frac1p}.
\ee
%----------------------------------------------------------
%@@@@@@@@@@@@@@@@@@@@@@@@@@@@@@@@@@@@@@@@@@@@@@@@@@@@@@@@@@
%@@@@@@@@@@@@@@@@@@@@@@@@@@@@@@@@@@@@@@@@@@@@@@@@@@@@@@@@@@
%@@@@@@@@@@@@@@@@@@@@@@@@@@@@@@@@@@@@@@@@@@@@@@@@@@@@@@@@@@
%@@@@@@@@@@@@@@@@@@@@@@@@@@@@@@@@@@@@@@@@@@@@@@@@@@@@@@@@@@
%@@@@@@@@@@@@@@@@@@@@@@@@@@@@@@@@@@@@@@@@@@@@@@@@@@@@@@@@@@
%@@@@@@@@@@@@@@@@@@@@@@@@@@@@@@@@@@@@@@@@@@@@@@@@@@@@@@@@@@
\medskip
%----------------------------------------------------------
\par We turn to a definition of the function $h_3(f):\RT\to\RT$. We let $\KFL$ denote a subfamily of the family $\Kc$ consisting of \sqs $K\in\Kc$ such that the triangle $\Delta(K)=\Delta\{A^{(K)},B^{(K)},C^{(K)}\}$ is {\it degenerate}. Thus {\it $A^{(K)},B^{(K)},C^{(K)}$ are three different collinear points.}
%----------------------------------------------------------
\par Let us denote the vertices $A^{(K)},B^{(K)},C^{(K)}$  by letters $z^{(K)}_1,z^{(K)}_2$ and $z^{(K)}_3$ in such a way that $z^{(K)}_2\in(z^{(K)}_1,z^{(K)}_3).$ Since
$A^{(K)},B^{(K)},C^{(K)}\in\gamma K,$ see \rf{D-INK},
%----------------------------------------------------------
\bel{Z-INC}
z^{(K)}_1,z^{(K)}_2,z^{(K)}_3 \in\gamma K
\ee
%----------------------------------------------------------
as well.
%----------------------------------------------------------
\par We put
%----------------------------------------------------------
\bel{H-3}
h_3(f):=\smed_{K\in\KFL}\,(\diam K)^{-1}
\left|
\frac{f(z^{(K)}_1)-f(z^{(K)}_2)}{\|z^{(K)}_1-z^{(K)}_2\|_2}
-\frac{f(z^{(K)}_2)-f(z^{(K)}_3)}{\|z^{(K)}_2-z^{(K)}_3\|_2}
\right|\,\chi_K.
\ee
%----------------------------------------------------------
Recall that $\|\cdot\|_2$ denotes the Euclidean norm in $\RT$. Since the \sqs of the family $\KFL$ are pairwise disjoint,
%----------------------------------------------------------
$$
\|h_3(f)\|_{\LPRT}^p=\smed_{K\in\KFL}\,
\left|
\frac{f(z^{(K)}_1)-f(z^{(K)}_2)}{\|z^{(K)}_1-z^{(K)}_2\|_2}
-\frac{f(z^{(K)}_2)-f(z^{(K)}_3)}{\|z^{(K)}_2-z^{(K)}_3\|_2}
\right|^p\,(\diam K)^{2-p}
$$
%----------------------------------------------------------
so that, by \rf{Z-INC} and \rf{DDIF},
%----------------------------------------------------------
\bel{H3-LP}
\|h_3(f)\|_{\LPRT}\le C(p)\lambda^{\frac1p}.
\ee
%----------------------------------------------------------
%@@@@@@@@@@@@@@@@@@@@@@@@@@@@@@@@@@@@@@@@@@@@@@@@@@@@@@@@@@
%@@@@@@@@@@@@@@@@@@@@@@@@@@@@@@@@@@@@@@@@@@@@@@@@@@@@@@@@@@
%@@@@@@@@@@@@@@@@@@@@@@@@@@@@@@@@@@@@@@@@@@@@@@@@@@@@@@@@@@
%@@@@@@@@@@@@@@@@@@@@@@@@@@@@@@@@@@@@@@@@@@@@@@@@@@@@@@@@@@
%@@@@@@@@@@@@@@@@@@@@@@@@@@@@@@@@@@@@@@@@@@@@@@@@@@@@@@@@@@
%@@@@@@@@@@@@@@@@@@@@@@@@@@@@@@@@@@@@@@@@@@@@@@@@@@@@@@@@@@
%----------------------------------------------------------
\par Finally we put
%----------------------------------------------------------
$$
h(f):=h_1(f)+h_2(f)+h_3(f).
$$
%----------------------------------------------------------
Then, by \rf{H1-LP}, \rf{H2-LP}, and \rf{H3-LP},
%----------------------------------------------------------
$$
\|h(f)\|_{\LPRT}\le C(p)\lambda^{\frac1p}.
$$
%----------------------------------------------------------
%@@@@@@@@@@@@@@@@@@@@@@@@@@@@@@@@@@@@@@@@@@@@@@@@@@@@@@@@@@
%@@@@@@@@@@@@@@@@@@@@@@@@@@@@@@@@@@@@@@@@@@@@@@@@@@@@@@@@@@
%@@@@@@@@@@@@@@@@@@@@@@@@@@@@@@@@@@@@@@@@@@@@@@@@@@@@@@@@@@
%@@@@@@@@@@@@@@@@@@@@@@@@@@@@@@@@@@@@@@@@@@@@@@@@@@@@@@@@@@
%@@@@@@@@@@@@@@@@@@@@@@@@@@@@@@@@@@@@@@@@@@@@@@@@@@@@@@@@@@
%@@@@@@@@@@@@@@@@@@@@@@@@@@@@@@@@@@@@@@@@@@@@@@@@@@@@@@@@@@
%----------------------------------------------------------
\begin{proposition}\lbl{FINAL} Let $2<p<\infty$ and let $q:=(p+2)/2$. Let $f:E\to\R$ be a function satisfying the hypo\-the\-sis of Theorem \reff{S-MAIN}.
%----------------------------------------------------------
\par The mapping $g(f):\RT\to\RT$ defined by the formula \rf{D-SEL} has the following properties:
%----------------------------------------------------------
\par (i). For every bridge $T\in\BRE$ with ends $\A{T},\B{T}\in E$ we have
%----------------------------------------------------------
$$
\ip{g(T;f),\A{T}-\B{T}}=f(\A{T})-f(\B{T});
$$
%----------------------------------------------------------
%@@@@@@@@@@@@@@@@@@@@@@@@@@@@@@@@@@@@@@@@@@@@@@@@@@@@@@@@@@
%----------------------------------------------------------
\par (ii). For every pair of connected bridges $T,T'\in\BRE,$ $T\bcn T',$ the following inequality
%----------------------------------------------------------
$$
\|g(T;f)-g(T';f)\|\le C(p)\diam Q(T,T')\,\left(
\frac{1}{|\gamma Q(T,T')|}
\intl_{\gamma Q(T,T')} h^q(z;f)\,dz\right)^{\frac1q}
$$
%----------------------------------------------------------
holds. Here $\gamma>0$ is an absolute constant.
%----------------------------------------------------------
\end{proposition}
%----------------------------------------------------------
%@@@@@@@@@@@@@@@@@@@@@@@@@@@@@@@@@@@@@@@@@@@@@@@@@@@@@@@@@@
%@@@@@@@@@@@@@@@@@@@@@@@@@@@@@@@@@@@@@@@@@@@@@@@@@@@@@@@@@@
%@@@@@@@@@@@@@@@@@@@@@@@@@@@@@@@@@@@@@@@@@@@@@@@@@@@@@@@@@@
%@@@@@@@@@@@@@@@@@@@@@@@@@@@@@@@@@@@@@@@@@@@@@@@@@@@@@@@@@@
%----------------------------------------------------------
\par We postpone a proof of the proposition to the end of this section. Here we only note that this proposition and Proposition \reff{A-H} (see also Claim \reff{CL-CR}) immediately imply the required inequality \rf{FIN-R}.
%----------------------------------------------------------
\par The proof of the proposition relies on a series of auxiliary statements. First of them is the classical Sobolev-Poincar\'e inequality for vector functions, see e.g., \cite{M}.
%----------------------------------------------------------
%@@@@@@@@@@@@@@@@@@@@@@@@@@@@@@@@@@@@@@@@@@@@@@@@@@@@@@@@@@
%@@@@@@@@@@@@@@@@@@@@@@@@@@@@@@@@@@@@@@@@@@@@@@@@@@@@@@@@@@
%@@@@@@@@@@@@@@@@@@@@@@@@@@@@@@@@@@@@@@@@@@@@@@@@@@@@@@@@@@
%@@@@@@@@@@@@@@@@@@@@@@@@@@@@@@@@@@@@@@@@@@@@@@@@@@@@@@@@@@
%@@@@@@@@@@@@@@@@@@@@@@@@@@@@@@@@@@@@@@@@@@@@@@@@@@@@@@@@@@
%@@@@@@@@@@@@@@@@@@@@@@@@@@@@@@@@@@@@@@@@@@@@@@@@@@@@@@@@@@
%----------------------------------------------------------
\begin{proposition}\lbl{SP-V} Let $2<q\le p<\infty$ and let a mapping $\vec{F}\in \VLOP$. Then for every \sq $Q\subset\RT$ and every $x,y\in Q$ the following inequality
%----------------------------------------------------------
$$
\|\vec{F}(x)-\vec{F}(y)\|\le C(q)\diam Q\left(
\frac{1}{|Q|}
\intl_{Q} \|\vec{F}(z)\|^q\,dz\right)^{\frac1q}
$$
%----------------------------------------------------------
holds.
%----------------------------------------------------------
\end{proposition}
%----------------------------------------------------------
%@@@@@@@@@@@@@@@@@@@@@@@@@@@@@@@@@@@@@@@@@@@@@@@@@@@@@@@@@@
%@@@@@@@@@@@@@@@@@@@@@@@@@@@@@@@@@@@@@@@@@@@@@@@@@@@@@@@@@@
%@@@@@@@@@@@@@@@@@@@@@@@@@@@@@@@@@@@@@@@@@@@@@@@@@@@@@@@@@@
%@@@@@@@@@@@@@@@@@@@@@@@@@@@@@@@@@@@@@@@@@@@@@@@@@@@@@@@@@@
%----------------------------------------------------------
\par The next auxiliary result is the following
%----------------------------------------------------------
%@@@@@@@@@@@@@@@@@@@@@@@@@@@@@@@@@@@@@@@@@@@@@@@@@@@@@@@@@@
%@@@@@@@@@@@@@@@@@@@@@@@@@@@@@@@@@@@@@@@@@@@@@@@@@@@@@@@@@@
%@@@@@@@@@@@@@@@@@@@@@@@@@@@@@@@@@@@@@@@@@@@@@@@@@@@@@@@@@@
%@@@@@@@@@@@@@@@@@@@@@@@@@@@@@@@@@@@@@@@@@@@@@@@@@@@@@@@@@@
%@@@@@@@@@@@@@@@@@@@@@@@@@@@@@@@@@@@@@@@@@@@@@@@@@@@@@@@@@@
%@@@@@@@@@@@@@@@@@@@@@@@@@@@@@@@@@@@@@@@@@@@@@@@@@@@@@@@@@@
%----------------------------------------------------------
\begin{lemma}\lbl{PR-ES} For every pair of connected bridges $T,T'\in\BRE, T\bcn T'$, the following inequality
%----------------------------------------------------------
$$
\|\tg(T;f)-\tg(T';f)\|\le C(p)\diam Q(T,T')\,\left(
\frac{1}{|Q(T,T')|}
\intl_{Q(T,T')} h_1^q(z;f)\,dz\right)^{\frac1q}
$$
%----------------------------------------------------------
holds.
%----------------------------------------------------------
\end{lemma}
%----------------------------------------------------------
%@@@@@@@@@@@@@@@@@@@@@@@@@@@@@@@@@@@@@@@@@@@@@@@@@@@@@@@@@@
%@@@@@@@@@@@@@@@@@@@@@@@@@@@@@@@@@@@@@@@@@@@@@@@@@@@@@@@@@@
%@@@@@@@@@@@@@@@@@@@@@@@@@@@@@@@@@@@@@@@@@@@@@@@@@@@@@@@@@@
%@@@@@@@@@@@@@@@@@@@@@@@@@@@@@@@@@@@@@@@@@@@@@@@@@@@@@@@@@@
%----------------------------------------------------------
\par {\it Proof.} We recall that
%----------------------------------------------------------
$$
\tg(T;f)=\Tc_1(\A{T};f)~~~~~\text{and}~~~~~
\tg(T';f)=\Tc_1(\A{T'};f)
$$
%----------------------------------------------------------
where $\Tc_1(f):\RT\to\RT$ is a Sobolev $\VLOP$-mapping. We also recall that $h_1(f)=\|\nabla\Tc_1(f)\|$ and $Q(T,T')=Q(\A{T},D(T,T'))$ with $D(T,T')=\diam\{\A{T},\B{T},\A{T'},\B{T'}\}$.
%----------------------------------------------------------
\par Hence $\A{T},\A{T'}\in Q(T,T')$. It remains to apply the Sobolev-Poincar\'e inequality to the mapping $\Tc(f)$, the points $\A{T},\A{T'},$ and the \sq $Q(T,T')$, and the lemma follows.\bx\medskip
%----------------------------------------------------------
%@@@@@@@@@@@@@@@@@@@@@@@@@@@@@@@@@@@@@@@@@@@@@@@@@@@@@@@@@@
%@@@@@@@@@@@@@@@@@@@@@@@@@@@@@@@@@@@@@@@@@@@@@@@@@@@@@@@@@@
%@@@@@@@@@@@@@@@@@@@@@@@@@@@@@@@@@@@@@@@@@@@@@@@@@@@@@@@@@@
%@@@@@@@@@@@@@@@@@@@@@@@@@@@@@@@@@@@@@@@@@@@@@@@@@@@@@@@@@@
%----------------------------------------------------------
\par Let us consider a pair of connected bridges $T,T'\in\BRE, T\bcn T'$, such that their ends $\A{T},\B{T},\A{T'},\B{T'}$ {\it are not collinear points} in $\RT$. Recall that in this case
%----------------------------------------------------------
$$
\#\{\A{T},\B{T},\A{T'},\B{T'}\}=3
$$
%----------------------------------------------------------
so that these points are vertices of a triangle. Let $K=K(T,T')\in\Kc$ be the \sq which we assign to the bridges $T$ and $T'$ so that $\Delta(K)=\Delta\{A^{(K)},B^{(K)},C^{(K)}\}$ is the triangle with vertices in $\{\A{T},\B{T},\A{T'},\B{T'}\}$, see \rf{CK-T} and \rf{F-TRI}. We recall that
$\diam \Delta(K)\sim \diam K$ and
$\Delta(K)\subset \gamma K$ with an absolute constant $\gamma$, see \rf{K-D} and \rf{D-INK}.
%----------------------------------------------------------
\par The next lemma presents additional geometric properties of the triangle $\Delta(K)$ and its connections with the set-valued mapping $G_f$, see \rf{D-GT}.
%----------------------------------------------------------
%@@@@@@@@@@@@@@@@@@@@@@@@@@@@@@@@@@@@@@@@@@@@@@@@@@@@@@@@@@
%@@@@@@@@@@@@@@@@@@@@@@@@@@@@@@@@@@@@@@@@@@@@@@@@@@@@@@@@@@
%@@@@@@@@@@@@@@@@@@@@@@@@@@@@@@@@@@@@@@@@@@@@@@@@@@@@@@@@@@
%@@@@@@@@@@@@@@@@@@@@@@@@@@@@@@@@@@@@@@@@@@@@@@@@@@@@@@@@@@
%@@@@@@@@@@@@@@@@@@@@@@@@@@@@@@@@@@@@@@@@@@@@@@@@@@@@@@@@@@
%@@@@@@@@@@@@@@@@@@@@@@@@@@@@@@@@@@@@@@@@@@@@@@@@@@@@@@@@@@
%----------------------------------------------------------
\begin{lemma}\lbl{G-DK} Let $K\in\Kc$ and let
$\Delta(K)=\Delta\{A^{(K)},B^{(K)},C^{(K)}\}$ be a  triangle formed by a pair of connected bridges $T,T'\in\BRE, T\bcn T'$ (thus $K=K(T,T')$).\medskip
%----------------------------------------------------------
\par (i). If $\Delta(K)$ is a \dg triangle, i.e., $A^{(K)},B^{(K)}$ and $C^{(K)}$ are three different collinear points, then the straight lines $G_f(T)$ and $G_f(T')$ are parallel.\medskip
%----------------------------------------------------------
\par (ii). If $\Delta(K)$ is a true triangle, then  $G_f(T)$ and $G_f(T')$ are not parallel. They intersect each other at the point
%----------------------------------------------------------
$$
G_f(T)\cap G_f(T')=\nabla P_{\Delta(K)}[f].
$$
%----------------------------------------------------------
%@@@@@@@@@@@@@@@@@@@@@@@@@@@@@@@@@@@@@@@@@@@@@@@@@@@@@@@@@@
%----------------------------------------------------------
\par Furthermore, the angle $\alpha\in(0,\pi/2)$ between these straight lines coincides with the angle $\alpha_K$ of the vertex $C^{(K)}=C(T,T')$ in the triangle $\Delta(K)$, and the following equivalences
%----------------------------------------------------------
\bel{M-SIN}
\sin\alpha\sim \cu_{\Delta(K)}\,\diam \Delta(K)\sim \cu_{\Delta(K)}\,\diam Q(T,T')
\ee
%----------------------------------------------------------
hold. The constants of these equivalences are absolute.
%----------------------------------------------------------
\end{lemma}
%----------------------------------------------------------
%@@@@@@@@@@@@@@@@@@@@@@@@@@@@@@@@@@@@@@@@@@@@@@@@@@@@@@@@@@
%@@@@@@@@@@@@@@@@@@@@@@@@@@@@@@@@@@@@@@@@@@@@@@@@@@@@@@@@@@
%@@@@@@@@@@@@@@@@@@@@@@@@@@@@@@@@@@@@@@@@@@@@@@@@@@@@@@@@@@
%@@@@@@@@@@@@@@@@@@@@@@@@@@@@@@@@@@@@@@@@@@@@@@@@@@@@@@@@@@
%----------------------------------------------------------
\par {\it Proof.}  The straight line $G_f(T)$ is defined by the equation
%----------------------------------------------------------
$$
\ip{z,\A{T}-\B{T}}=f(\A{T})-f(\B{T})
$$
%----------------------------------------------------------
so that $\vec{n}_T:=\A{T}-\B{T} \perp G_f(T)$, i.e., the vector $\vec{n}_T$ is orthogonal to  $G_f(T)$.
The same is true for the bridge $T'$, i.e.,
$\vec{n}_{T'}:=\A{T'}-\B{T'} \perp G_f(T')$.
%----------------------------------------------------------
\par Recall that $T$ and $T'$ are connected bridges and $\Delta(K)=\Delta\{A^{(K)},B^{(K)},C^{(K)}\}$ is a  triangle formed by the ends of $T,T'$, i.e., by the points $\A{T},\B{T},\A{T'},\B{T'}$. We also recall that in this case $C^{(K)}$ is the (unique) common point of the sets $\{\A{T},\B{T}\}$ and $\{\A{T'},\B{T'}\}$, see \rf{CK-T}. Without loss of generality we may assume that
%----------------------------------------------------------
$$
\{\A{T},\B{T}\}=\{A^{(K)},C^{(K)}\}~~~\text{and}~~~
\{\A{T'},\B{T'}\}=\{B^{(K)},C^{(K)}\}.
$$
%----------------------------------------------------------
\par Hence $A^{(K)}-C^{(K)}=\pm n_T$ so that
%----------------------------------------------------------
\bel{ORT-1}
A^{(K)}-C^{(K)}\perp G_f(T).
\ee
%----------------------------------------------------------
The same is true for the bridge $T'$, i.e.,
%----------------------------------------------------------
\bel{ORT-2}
B^{(K)}-C^{(K)}\perp G_f(T').
\ee
%----------------------------------------------------------
\par Thus if $\Delta(K)=\Delta\{A^{(K)},B^{(K)},C^{(K)}\}$
is a \dg triangle, i.e.,  $A^{(K)},B^{(K)}$ and $C^{(K)}$ lie on a certain straight line, then $(A^{(K)}-C^{(K)})\parallel (B^{(K)}-C^{(K)})$. Combining this property with \rf{ORT-1} and \rf{ORT-2} we conclude that $G_f(T)$ and $G_f(T')$ are parallel as well. This proves part (i) of the lemma.
%----------------------------------------------------------
\par Prove (ii). Suppose that $\Delta(K)=\Delta\{A^{(K)},B^{(K)},C^{(K)}\}$ is a true triangle. Then
%----------------------------------------------------------
$$
(A^{(K)}-C^{(K)})\nparallel (B^{(K)}-C^{(K)})
$$
%----------------------------------------------------------
proving that $G_f(T)\nparallel G_f(T')$.
%----------------------------------------------------------
\par Let us note that the points on the straight line
%----------------------------------------------------------
\be
G_f(T)&:=&\{z\in\RT: \ip{z,\A{T}-\B{T}}=f(\A{T})-f(\B{T})\}\nn\\
&=&
\{z\in\RT: \ip{z,A^{(K)}-C^{(K)}}
=f(A^{(K)})-f(C^{(K)})\}.\nn
\ee
%----------------------------------------------------------
can be identified with {\it the gradients of polynomials $P\in\PO$ which interpolate $f$ at the points $\A{T}$ and $\B{T}$}, or equivalently at $A^{(K)}$ and $C^{(K)}$). In other words,
%----------------------------------------------------------
$$
G_f(T):=\{\nabla P: P\in\PO, P(A^{(K)})=f(A^{(K)}),  P(C^{(K)})=f(C^{(K)})\}.
$$
%----------------------------------------------------------
Analogously,
%----------------------------------------------------------
$$
G_f(T'):=\{\nabla P: P\in\PO, P(B^{(K)})=f(B^{(K)}),  P(C^{(K)})=f(C^{(K)})\}.
$$
%----------------------------------------------------------
Thus the point $G_f(T)\cap G_f(T')$ can be identified with the gradient of the affine polynomial which interpolates $f$ at the points $A^{(K)},B^{(K)}$ and $C^{(K)}$, i.e., at the vertices of the triangle $\Delta(K)$. This shows that
%----------------------------------------------------------
$$
G_f(T)\cap G_f(T')=\nabla P_{\Delta(K)}[f].
$$
%----------------------------------------------------------
%@@@@@@@@@@@@@@@@@@@@@@@@@@@@@@@@@@@@@@@@@@@@@@@@@@@@@@@@@@
%@@@@@@@@@@@@@@@@@@@@@@@@@@@@@@@@@@@@@@@@@@@@@@@@@@@@@@@@@@
%----------------------------------------------------------
\par Let us prove equivalences in \rf{M-SIN}. By \rf{ORT-1} and \rf{ORT-2}, the angle $\alpha$ between the straight lines $G_f(T)$ and $G_f(T')$ coincides with the angle $\alpha_K$ between the sides $[A^{(K)},C^{(K)}]$ and $[B^{(K)},C^{(K)}]$ of the triangle $\Delta(K)$, i.e., with the angle of the vertex $C^{(K)}$ in $\Delta(K)$.
%----------------------------------------------------------
\par Hence, by \rf{K-S},
%----------------------------------------------------------
$$
\sin \alpha=\sin \alpha_K\sim \cu_{\Delta(K)}\,\diam K.
$$
%----------------------------------------------------------
On the other hand, by \rf{K-D},
$\diam \Delta(K)\sim \diam K,$ and
%----------------------------------------------------------
$$
\diam Q(T,T')=2\diam\{\A{T},\B{T},\A{T'},\B{T'}\}
\sim \diam K.
$$
%----------------------------------------------------------
These equivalences imply the required equivalences in \rf{M-SIN} proving the lemma.\bx\medskip
%----------------------------------------------------------
%@@@@@@@@@@@@@@@@@@@@@@@@@@@@@@@@@@@@@@@@@@@@@@@@@@@@@@@@@@
%@@@@@@@@@@@@@@@@@@@@@@@@@@@@@@@@@@@@@@@@@@@@@@@@@@@@@@@@@@
%@@@@@@@@@@@@@@@@@@@@@@@@@@@@@@@@@@@@@@@@@@@@@@@@@@@@@@@@@@
%@@@@@@@@@@@@@@@@@@@@@@@@@@@@@@@@@@@@@@@@@@@@@@@@@@@@@@@@@@
%----------------------------------------------------------
\begin{lemma}\lbl{D-INT} Let $K\in\Kc$ and let $\Delta(K)$ be a true triangle formed by a pair of connected bridges $T,T'\in\BRE, T\bcn T'$ (thus $K=K(T,T')$). Then
%----------------------------------------------------------
$$
\|\tg(T;f)-\nabla P_{\Delta(K)}[f]\|\cu_{\Delta(K)}\le C\,\left(
\frac{1}{|\gamma Q(T,T')|}
\intl_{\gamma Q(T,T')} (h_1+h_2)^q(z;f)\,dz\right)^{\frac1q}
$$
%----------------------------------------------------------
where $\gamma\ge 1$ is an absolute constant.
%----------------------------------------------------------
\end{lemma}
%----------------------------------------------------------
%@@@@@@@@@@@@@@@@@@@@@@@@@@@@@@@@@@@@@@@@@@@@@@@@@@@@@@@@@@
%@@@@@@@@@@@@@@@@@@@@@@@@@@@@@@@@@@@@@@@@@@@@@@@@@@@@@@@@@@
%@@@@@@@@@@@@@@@@@@@@@@@@@@@@@@@@@@@@@@@@@@@@@@@@@@@@@@@@@@
%@@@@@@@@@@@@@@@@@@@@@@@@@@@@@@@@@@@@@@@@@@@@@@@@@@@@@@@@@@
%----------------------------------------------------------
\par {\it Proof.} Let $K\in\Kc$. Since the \sqs of the family $\Kc$ are pairwise disjoint, by \rf{H-2},
%----------------------------------------------------------
$$
h_2(x;f)=
\|\nabla P_{\Delta(K)}[f]-\Tc_1(c_K,f)\|\,
\cu_{\Delta(K)}~~~\text{for every}~~~x\in K.
$$
%----------------------------------------------------------
Hence,
%----------------------------------------------------------
$$
I_1:=\|\nabla P_{\Delta(K)}[f]-\Tc_1(c_K,f)\|\,
\cu_{\Delta(K)}=\left(\frac{1}{|K|}\intl_{K} h_2^q(z;f)\,dz\right)^{\frac1q}
$$
%----------------------------------------------------------
so that, by \rf{Q-T-IN},
%----------------------------------------------------------
$$
I_1=\left(\frac{1}{|K|}\intl_{K} h_2^q(z;f)\,dz\right)^{\frac1q}\le
C\,\left(\frac{1}{|\gamma Q(T,T')|}
\intl_{\gamma Q(T,T')} h_2^q(z;f)\,dz\right)^{\frac1q}.
$$
%----------------------------------------------------------
\par In turn, by the second inclusion in \rf{Q-T-IN}, the square $K\subset \gamma Q(T,T')$ so that $c_K,\A{T}\in \gamma Q(T,T')$. Applying the Sobolev-Poincar\'e inequality, see Proposition \reff{SP-V}, to the mapping $\Tc_1(f)$, the points  $c_K$ and $\A{T}$, and the \sq $\gamma Q(T,T')$, we obtain
%----------------------------------------------------------
\be
\|\Tc_1(c_K;f)-\tg(T;f)\|&=&
\|\Tc_1(c_K;f)-\Tc_1(\A{T};f)\|\nn\\&\le&
C\,\diam Q(T,T')\left(\frac{1}{|\gamma Q(T,T')|}
\intl_{\gamma Q(T,T')}
\|\nabla\Tc_1(z;f)\|^q\,dz\right)^{\frac1q}.
\nn
\ee
%----------------------------------------------------------
But, by \rf{H-1}, $h_1(f)=\|\nabla\Tc_1(f)\|$ so that
%----------------------------------------------------------
$$
\|\Tc_1(c_K;f)-\tg(T;f)\|\le
C\,\diam Q(T,T')\left(\frac{1}{|\gamma Q(T,T')|}
\intl_{\gamma Q(T,T')}
h_1(z;f)^q\,dz\right)^{\frac1q}.
$$
%----------------------------------------------------------
Hence
%----------------------------------------------------------
\be
I_2&:=&\|\Tc_1(c_K;f)-\tg(T;f)\|\cu_{\Delta(K)}\nn\\
&\le&
C\,\cu_{\Delta(K)}\diam Q(T,T')\left(\frac{1}{|\gamma Q(T,T')|}
\intl_{\gamma Q(T,T')}
h_1(z;f)^q\,dz\right)^{\frac1q}.\nn
\ee
%----------------------------------------------------------
Note that, by \rf{K-D} and \rf{K-S},
%----------------------------------------------------------
$$
\cu_{\Delta(K)}\diam Q(T,T')\le C\,\cu_{\Delta(K)}\diam K\le C\,|\sin \alpha_K|
$$
%----------------------------------------------------------
so that
%----------------------------------------------------------
\be
I_2&\le&
C\,|\sin \alpha_K|\left(\frac{1}{|\gamma Q(T,T')|}
\intl_{\gamma Q(T,T')}
h_1(z;f)^q\,dz\right)^{\frac1q}\nn\\
&\le&
C\,\left(\frac{1}{|\gamma Q(T,T')|}
\intl_{\gamma Q(T,T')}
h_1(z;f)^q\,dz\right)^{\frac1q}.\nn
\ee
%----------------------------------------------------------
Finally, since $h_1(f)$ and $h_2(f)$ are non-negative, we have
%----------------------------------------------------------
\be
\|\tg(T;f)-\nabla P_{\Delta(K)}[f]\|\,\cu_{\Delta(K)}
&\le& I_1+I_2\le C\left(\frac{1}{|\gamma Q(T,T')|}
\intl_{\gamma Q(T,T')}
h_2(z;f)^q\,dz\right)^{\frac1q}\nn\\
&+&
\left(\frac{1}{|\gamma Q(T,T')|}
\intl_{\gamma Q(T,T')}
h_1(z;f)^q\,dz\right)^{\frac1q}\nn\\
&\le&  C\,\left(
\frac{1}{|\gamma Q(T,T')|}
\intl_{\gamma Q(T,T')} (h_1+h_2)^q(z;f)\,dz\right)^{\frac1q}.\nn
\ee
%----------------------------------------------------------
\par The lemma is proved.\bx\medskip
%----------------------------------------------------------
\par Given subsets $A_1,A_2\subset\RT$ we let $\dist_2(A_1,A_2)$ denote the Euclidean distance between $A_1$ and $A_2$:
%----------------------------------------------------------
$$
\dist_2(A_1,A_2)=
\inf\{\|a_1-a_2\|_2: a_1\in A_1, a_2\in A_2\}.
$$
%----------------------------------------------------------
Recall that by $\|\cdot\|_2$ we denote the Euclidean norm in $\RT$.
%----------------------------------------------------------
%@@@@@@@@@@@@@@@@@@@@@@@@@@@@@@@@@@@@@@@@@@@@@@@@@@@@@@@@@@
%@@@@@@@@@@@@@@@@@@@@@@@@@@@@@@@@@@@@@@@@@@@@@@@@@@@@@@@@@@
%@@@@@@@@@@@@@@@@@@@@@@@@@@@@@@@@@@@@@@@@@@@@@@@@@@@@@@@@@@
%@@@@@@@@@@@@@@@@@@@@@@@@@@@@@@@@@@@@@@@@@@@@@@@@@@@@@@@@@@
%----------------------------------------------------------
\begin{lemma}\lbl{FL-T} Let $K\in\Kc$ and let $\Delta(K)$ be a \dg triangle formed by a pair of connected bridges $T,T'\in\BRE, T\bcn T'$. Then the straight lines $G_f(T)$ and $G_f(T')$ are parallel. Furthermore,
%----------------------------------------------------------
$$
\dist_2(G_f(T),G_f(T'))\le
C\,\diam Q(T,T')\,
\left(\frac{1}{|\gamma Q(T,T')|}
\intl_{\gamma Q(T,T')} h_3^q(z;f)\,dz\right)^{\frac1q}.
$$
%----------------------------------------------------------
Here $\gamma\ge 1$ is an absolute constant.
%----------------------------------------------------------
\end{lemma}
%----------------------------------------------------------
%@@@@@@@@@@@@@@@@@@@@@@@@@@@@@@@@@@@@@@@@@@@@@@@@@@@@@@@@@@
%@@@@@@@@@@@@@@@@@@@@@@@@@@@@@@@@@@@@@@@@@@@@@@@@@@@@@@@@@@
%@@@@@@@@@@@@@@@@@@@@@@@@@@@@@@@@@@@@@@@@@@@@@@@@@@@@@@@@@@
%@@@@@@@@@@@@@@@@@@@@@@@@@@@@@@@@@@@@@@@@@@@@@@@@@@@@@@@@@@
%----------------------------------------------------------
\par {\it Proof.}  We recall that $A^{(K)},B^{(K)}$ and $C^{(K)}$ are three different collinear points in $\RT$ provided $\Delta(K)=\Delta\{A^{(K)},B^{(K)},C^{(K)}\}$ is a \dg triangle. By part (i) of Lemma \reff{G-DK}, in this case $G_f(T)$ and $G_f(T')$ are parallel straight lines.
%----------------------------------------------------------
\par Let us denote the vertices $A^{(K)},B^{(K)},C^{(K)}$  by letters $z^{(K)}_1,z^{(K)}_2$ and $z^{(K)}_3$ in such a way that $z^{(K)}_2\in(z^{(K)}_1,z^{(K)}_3).$ Then the straight lines $G_f(T)$ and $G_f(T')$ are determined by the equations
%----------------------------------------------------------
$$
\ip{z,z^{(K)}_1-z^{(K)}_2}=f(z^{(K)}_1)-f(z^{(K)}_2)
$$
%----------------------------------------------------------
and
%----------------------------------------------------------
$$
\ip{z,z^{(K)}_2-z^{(K)}_3}=f(z^{(K)}_2)-f(z^{(K)}_3)
$$
%----------------------------------------------------------
respectively. The Euclidean distance between these parallel straight lines is equal to
%----------------------------------------------------------
$$
\dist_2(G_f(T),G_f(T'))=
\left|
\frac{f(z^{(K)}_1)-f(z^{(K)}_2)}{\|z^{(K)}_1-z^{(K)}_2\|_2}
-\frac{f(z^{(K)}_2)-f(z^{(K)}_3)}{\|z^{(K)}_2-z^{(K)}_3\|_2}
\right|.
$$
%----------------------------------------------------------
\par Hence, by \rf{H-3},
%----------------------------------------------------------
$$
h_3(x,f)=\dist_2(G_f(T),G_f(T'))/\diam K~~~~\text{for every}~~~x\in K.
$$
%----------------------------------------------------------
Raising this equality to the power $q$ and then integrating on the \sq $K$ we obtain
%----------------------------------------------------------
$$
\dist_2(G_f(T),G_f(T'))=
\diam K\,
\left(\frac{1}{|K|}
\intl_{K} h_3^q(z;f)\,dz\right)^{\frac1q}.
$$
%----------------------------------------------------------
It remains to note that, by \rf{Q-T-IN}, $K\subset \gamma Q(T,T')$ and $\diam K\sim \diam Q(T,T')$, so that
%----------------------------------------------------------
$$
\dist(G_f(T),G_f(T'))\le
C\,\diam Q(T,T')\,
\left(\frac{1}{|\gamma Q(T,T')|}
\intl_{\gamma Q(T,T')} h_3^q(z;f)\,dz\right)^{\frac1q}
$$
%----------------------------------------------------------
proving the lemma.\bx\medskip
%----------------------------------------------------------
%@@@@@@@@@@@@@@@@@@@@@@@@@@@@@@@@@@@@@@@@@@@@@@@@@@@@@@@@@@
%@@@@@@@@@@@@@@@@@@@@@@@@@@@@@@@@@@@@@@@@@@@@@@@@@@@@@@@@@@
%@@@@@@@@@@@@@@@@@@@@@@@@@@@@@@@@@@@@@@@@@@@@@@@@@@@@@@@@@@
%@@@@@@@@@@@@@@@@@@@@@@@@@@@@@@@@@@@@@@@@@@@@@@@@@@@@@@@@@@
%----------------------------------------------------------
\begin{lemma}\lbl{G-PR} Let $\ell_1$ and $\ell_2$ be two non-parallel straight lines in $\RT$ intersecting at a point $a\in\RT$. Let $\alpha\in(0,\pi/2)$ be the angle between these straight lines.
%----------------------------------------------------------
\par Then for every $x,y\in\RT$ the following inequality
%----------------------------------------------------------
\bel{D-PRJ}
\|\PR(x;\ell_1)-\PR(y;\ell_2)\|_2\le \|x-y\|_2+2\sqrt{2}\,\sin\alpha\,\|x-a\|_2
\ee
%----------------------------------------------------------
holds.
%----------------------------------------------------------
\end{lemma}
%----------------------------------------------------------
%@@@@@@@@@@@@@@@@@@@@@@@@@@@@@@@@@@@@@@@@@@@@@@@@@@@@@@@@@@
%@@@@@@@@@@@@@@@@@@@@@@@@@@@@@@@@@@@@@@@@@@@@@@@@@@@@@@@@@@
%@@@@@@@@@@@@@@@@@@@@@@@@@@@@@@@@@@@@@@@@@@@@@@@@@@@@@@@@@@
%@@@@@@@@@@@@@@@@@@@@@@@@@@@@@@@@@@@@@@@@@@@@@@@@@@@@@@@@@@
%----------------------------------------------------------
\par {\it Proof.} Without loss of generality we may assume that $a=0$, i.e., $\ell_1$ and $\ell_2$ are one dimensional linear subspaces of $\RT$. We may also assume that $\ell_1=A\ell_2$ where $A:\RT\to\RT$ is {\it a rotation operator by the angle $\alpha$}. Note that
for every $z\in\RT$ and every linear subspace $\ell\subset\RT$, $\dim \ell=1$, we have
$\PR(Az;A\ell)=A(\PR(z;\ell))$ so that
%----------------------------------------------------------
\bel{COM}
\PR(Ax;\ell_1)=A(\PR(x;\ell_2)).
\ee
%----------------------------------------------------------
\par Now we have
%----------------------------------------------------------
$$
\|\PR(x;\ell_1)-\PR(y;\ell_2)\|_2\le
\|\PR(x;\ell_1)-\PR(x;\ell_2)\|_2+
\|\PR(x;\ell_2)-\PR(y;\ell_2)\|_2
$$
%----------------------------------------------------------
so that
\bel{NT-1}
\|\PR(x;\ell_1)-\PR(y;\ell_2)\|_2
\le\|\PR(x;\ell_1)-\PR(x;\ell_2)\|_2+\|x-y\|_2.
\ee
%----------------------------------------------------------
\par Let us estimate the quantity
%----------------------------------------------------------
$$
I:=\|\PR(x;\ell_1)-\PR(x;\ell_2)\|_2.
$$
%----------------------------------------------------------
We have
%----------------------------------------------------------
\be
I&\le& \|\PR(x;\ell_1)-\PR(Ax;\ell_1)\|_2
+\|\PR(Ax;\ell_1)-\PR(x;\ell_2)\|_2\nn\\
&\le&
\|x-Ax\|_2+\|\PR(Ax;\ell_1)-\PR(x;\ell_2)\|_2.\nn
\ee
%----------------------------------------------------------
Hence, by \rf{COM},
%----------------------------------------------------------
\bel{ID-IN}
I\le
\|x-Ax\|_2+\|A(\PR(x;\ell_2))-\PR(x;\ell_2)\|_2\,.
\ee
%----------------------------------------------------------
Let $\Id$ be the identity operator on $\RT$ and let $B:=\Id-A$. Then, by \rf{ID-IN},
%----------------------------------------------------------
$$
I\le\|Bx\|_2+\|B((\PR(x;\ell_2))\|_2.
$$
%----------------------------------------------------------
Since $B$ commutes with rotations and $\|\PR(x;\ell_2)\|_2\le\|x\|_2$, we have
%----------------------------------------------------------
$$
\|B((\PR(x;\ell_2))\|_2\le \|Bx\|_2
$$
%----------------------------------------------------------
proving that $I\le 2\|Bx\|_2$.
%----------------------------------------------------------
\par Also, using a rotation, we obtain that
%----------------------------------------------------------
$$
\|Bx\|_2=\|Be_1\|_2\,\|x\|_2~~~\text{where}~~~ e_1=(1,0).
$$
%----------------------------------------------------------
\par Simple calculation shows that $\|Be_1\|_2=2\sin(\alpha/2).$ But $2\sin(\alpha/2)\le \sqrt{2}\sin\alpha$ for $\alpha\in(0,\pi/2)$ so that
$\|Be_1\|_2\le \sqrt{2}\sin\alpha$. Hence
%----------------------------------------------------------
$$
I\le 2\|Bx\|_2\le 2\sqrt{2}\,\sin\alpha\,\|x\|_2
=2\sqrt{2}\,\sin\alpha\,\|x-a\|_2.
$$
%----------------------------------------------------------
Combining this estimate with \rf{NT-1} we obtain the required inequality \rf{D-PRJ}.\bx
%----------------------------------------------------------
%@@@@@@@@@@@@@@@@@@@@@@@@@@@@@@@@@@@@@@@@@@@@@@@@@@@@@@@@@@
%@@@@@@@@@@@@@@@@@@@@@@@@@@@@@@@@@@@@@@@@@@@@@@@@@@@@@@@@@@
%@@@@@@@@@@@@@@@@@@@@@@@@@@@@@@@@@@@@@@@@@@@@@@@@@@@@@@@@@@
%@@@@@@@@@@@@@@@@@@@@@@@@@@@@@@@@@@@@@@@@@@@@@@@@@@@@@@@@@@
\bigskip
%----------------------------------------------------------
\par We are in a position to prove Proposition \reff{FINAL}.\smallskip
%----------------------------------------------------------
\par {\it Proof of Proposition \reff{FINAL}.} Part (i) of the proposition is trivial because for each $T\in\BRE$ the point $g(T;f)=\PR(\tg(T;f);G_f(T))$ belongs to the straight line $G_f(T)$ which is determined by the equation
$
\ip{z,\A{T}-\B{T}}=f(\A{T})-f(\B{T}).
$
%----------------------------------------------------------
\medskip
\par Let us prove part (ii) of the proposition. Let $T,T'\in\BRE, T\bcn T',$ be a pair of connected bridges. Let $K=K(T,T')$ be the \sq from the family $\Kc$ corresponding to the bridges $T$ and $T'$. Consider three cases.
%----------------------------------------------------------
%@@@@@@@@@@@@@@@@@@@@@@@@@@@@@@@@@@@@@@@@@@@@@@@@@@@@@@@@@@
%@@@@@@@@@@@@@@@@@@@@@@@@@@@@@@@@@@@@@@@@@@@@@@@@@@@@@@@@@@
%@@@@@@@@@@@@@@@@@@@@@@@@@@@@@@@@@@@@@@@@@@@@@@@@@@@@@@@@@@
%@@@@@@@@@@@@@@@@@@@@@@@@@@@@@@@@@@@@@@@@@@@@@@@@@@@@@@@@@@
\medskip
%----------------------------------------------------------
\par {\it The first case: the triangle $\Delta(K)$ formed by the ends of the bridges $T$ and $T'$ is a true triangle.}
%----------------------------------------------------------
\par By Lemma \reff{G-DK}, the straight lines $G_f(T)$ and $G_f(T')$ are non-parallel and
%----------------------------------------------------------
$$
G_f(T)\cap G_f(T')=\nabla P_{\Delta(K)}[f].
$$
%----------------------------------------------------------
Let $\alpha\in(0,\pi/2)$ be the angle between these straight lines. Then, by Lemma \reff{G-PR},
%----------------------------------------------------------
\be
\|g(T;f)-g(T';f)\|&=&
\|\PR(\tg(T;f);G_f(T))-\PR(\tg(T';f);G_f(T'))\|\nn\\
&\le& C\{
\|\tg(T;f)-\tg(T';f)\|+
\|\tg(T;f)-\nabla P_{\Delta(K)}[f]\|\sin\alpha\}\nn\\
&=&C\{I_1+I_2\}.\nn
\ee
%----------------------------------------------------------
\par By Lemma \reff{PR-ES},
%----------------------------------------------------------
$$
I_1:=\|\tg(T;f)-\tg(T';f)\|\le C(p)\diam Q(T,T')\,\left(
\frac{1}{|Q(T,T')|}
\intl_{Q(T,T')} h_1^q(z;f)\,dz\right)^{\frac1q}.
$$
%----------------------------------------------------------
\par Let us estimate the quantity %----------------------------------------------------------
$$
I_2:=
\|\tg(T;f)-\nabla P_{\Delta(K)}[f]\|\sin\alpha.
$$
%----------------------------------------------------------
By equivalence \rf{M-SIN},
%----------------------------------------------------------
$$
\sin\alpha\sim\cu_{\Delta(K)}\,\diam Q(T,T')
$$
%----------------------------------------------------------
so that
%----------------------------------------------------------
$$
I_2
\le C\,\|\tg(T;f)-\nabla P_{\Delta(K)}[f]\|\,\cu_{\Delta(K)}\,\diam Q(T,T').
$$
%----------------------------------------------------------
Hence, by Lemma \reff{D-INT},
%----------------------------------------------------------
$$
I_2
\le C\diam Q(T,T')\,
\left(\frac{1}{|\gamma Q(T,T')|}
\intl_{\gamma Q(T,T')} (h_1+h_2)^q(z;f)\,dz\right)^{\frac1q}
$$
%----------------------------------------------------------
with some absolute constant $\gamma\ge 1$.
%----------------------------------------------------------
\par Summarizing the estimates for $I_1$ and $I_2$ we obtain
%----------------------------------------------------------
\be
\|g(T;f)-g(T';f)\|&\le& C\{I_1+I_2\}\nn\\
&\le& C\diam Q(T,T')\,
\left(\frac{1}{|\gamma Q(T,T')|}
\intl_{\gamma Q(T,T')} (h_1+h_2)^q(z;f)\,dz\right)^{\frac1q}.\nn
\ee
%----------------------------------------------------------
%@@@@@@@@@@@@@@@@@@@@@@@@@@@@@@@@@@@@@@@@@@@@@@@@@@@@@@@@@@
%@@@@@@@@@@@@@@@@@@@@@@@@@@@@@@@@@@@@@@@@@@@@@@@@@@@@@@@@@@
%@@@@@@@@@@@@@@@@@@@@@@@@@@@@@@@@@@@@@@@@@@@@@@@@@@@@@@@@@@
%@@@@@@@@@@@@@@@@@@@@@@@@@@@@@@@@@@@@@@@@@@@@@@@@@@@@@@@@@@
\medskip
%----------------------------------------------------------
\par {\it The second case: the triangle $\Delta(K)$ is a \dg triangle.}
%----------------------------------------------------------
\par By Lemma \reff{FL-T}, $G_f(T)$ and $G_f(T')$ are parallel straight lines so that
%----------------------------------------------------------
\be
\|g(T;f)-g(T';f)\|&\le&\|g(T;f)-g(T';f)\|_2\nn\\&=&
\|\PR(\tg(T;f);G_f(T))-\PR(\tg(T';f);G_f(T'))\|_2\nn\\
&\le& \|\PR(\tg(T;f);G_f(T))-\PR(\tg(T';f);G_f(T))\|_2\nn\\
&+&
\|\PR(\tg(T';f);G_f(T))-\PR(\tg(T';f);G_f(T'))\|_2\nn\\
&\le&
\|\tg(T;f)-\tg(T';f)\|_2+\dist_2(G_f(T),G_f(T')).\nn
\ee
%----------------------------------------------------------
\par By Lemma \reff{PR-ES},
%----------------------------------------------------------
$$
\|\tg(T;f)-\tg(T';f)\|\le C\,\diam Q(T,T')\,\left(
\frac{1}{|Q(T,T')|}
\intl_{Q(T,T')} h_1^q(z;f)\,dz\right)^{\frac1q}.
$$
%----------------------------------------------------------
In turn, by Lemma \reff{FL-T},
%----------------------------------------------------------
$$
\dist(G_f(T),G_f(T'))_2\le
C\,\diam Q(T,T')\,
\left(\frac{1}{|\gamma Q(T,T')|}
\intl_{\gamma Q(T,T')} h_3^q(z;f)\,dz\right)^{\frac1q}.
$$
%----------------------------------------------------------
Hence
%----------------------------------------------------------
$$
\|g(T;f)-g(T';f)\|\le C\diam Q(T,T')\,
\left(\frac{1}{|\gamma Q(T,T')|}
\intl_{\gamma Q(T,T')} (h_1+h_3)^q(z;f)\,dz\right)^{\frac1q}.
$$
%----------------------------------------------------------
%@@@@@@@@@@@@@@@@@@@@@@@@@@@@@@@@@@@@@@@@@@@@@@@@@@@@@@@@@@
%@@@@@@@@@@@@@@@@@@@@@@@@@@@@@@@@@@@@@@@@@@@@@@@@@@@@@@@@@@
%@@@@@@@@@@@@@@@@@@@@@@@@@@@@@@@@@@@@@@@@@@@@@@@@@@@@@@@@@@
%@@@@@@@@@@@@@@@@@@@@@@@@@@@@@@@@@@@@@@@@@@@@@@@@@@@@@@@@@@
\medskip
%----------------------------------------------------------
\par {\it The third  case: the set of the ends of the bridges $T$ and $T'$ is a two point set.}
%----------------------------------------------------------
\par In this case $\{\A{T},\B{T}\}=\{\A{T'},\B{T'}\}$ so that, by definition \rf{D-GT}, $G_f(T)=G_f(T')$. Hence
%----------------------------------------------------------
\be
\|g(T;f)-g(T';f)\|&\le&\|g(T;f)-g(T';f)\|_2\nn\\&=&
\|\PR(\tg(T;f);G_f(T))-\PR(\tg(T';f);G_f(T'))\|_2\nn\\
&\le&
\|\tg(T;f)-\tg(T';f)\|_2\nn
\ee
%----------------------------------------------------------
so that, by Lemma \reff{PR-ES},
%----------------------------------------------------------
$$
\|g(T;f)-g(T';f)\|\le C\diam Q(T,T')\,
\left(\frac{1}{|Q(T,T')|}
\intl_{Q(T,T')} h_1^q(z;f)\,dz\right)^{\frac1q}.
$$
%----------------------------------------------------------
\par We have shown that in each of three possible cases the condition (ii) of the proposition holds. This proves part (ii) and finishes the proof of the proposition.\bx\medskip
%----------------------------------------------------------
%@@@@@@@@@@@@@@@@@@@@@@@@@@@@@@@@@@@@@@@@@@@@@@@@@@@@@@@@@@
%@@@@@@@@@@@@@@@@@@@@@@@@@@@@@@@@@@@@@@@@@@@@@@@@@@@@@@@@@@
%@@@@@@@@@@@@@@@@@@@@@@@@@@@@@@@@@@@@@@@@@@@@@@@@@@@@@@@@@@
%@@@@@@@@@@@@@@@@@@@@@@@@@@@@@@@@@@@@@@@@@@@@@@@@@@@@@@@@@@
%----------------------------------------------------------
\par Theorem \reff{S-MAIN} is completely proved.\bx
%----------------------------------------------------------
%@@@@@@@@@@@@@@@@@@@@@@@@@@@@@@@@@@@@@@@@@@@@@@@@@@@@@@@@@@
%@@@@@@@@@@@@@@@@@@@@@@@@@@@@@@@@@@@@@@@@@@@@@@@@@@@@@@@@@@
%@@@@@@@@@@@@@@@@@@@@@@@@@@@@@@@@@@@@@@@@@@@@@@@@@@@@@@@@@@
%@@@@@@@@@@@@@@@@@@@@@@@@@@@@@@@@@@@@@@@@@@@@@@@@@@@@@@@@@@
%----------------------------------------------------------
%@@@@@@@@@@@@@@@@@@@@@@@@@@@@@@@@@@@@@@@@@@@@@@@@@@@@@@@@@@
%@@@@@@@@@@@@@@@@@@@@@@@@@@@@@@@@@@@@@@@@@@@@@@@@@@@@@@@@@@
%@@@@@@@@@@@@@@@@@@@@@@@@@@@@@@@@@@@@@@@@@@@@@@@@@@@@@@@@@@
%@@@@@@@@@@@@@@@@@@@@@@@@@      @@@@@@@@@@@@@@@@@@@@@@@@@@@
%@@@@@@@@@@@@@@@@@@@@@@@          @@@@@@@@@@@@@@@@@@@@@@@@@
%@@@@@@@@@@@@@@@@@@@@@              @@@@@@@@@@@@@@@@@@@@@@@
%@@@@@@@@@@@@@@@@@@@     SECTION 9    @@@@@@@@@@@@@@@@@@@@@
%@@@@@@@@@@@@@@@@@@@@@              @@@@@@@@@@@@@@@@@@@@@@@
%@@@@@@@@@@@@@@@@@@@@@@@          @@@@@@@@@@@@@@@@@@@@@@@@@
%@@@@@@@@@@@@@@@@@@@@@@@@@      @@@@@@@@@@@@@@@@@@@@@@@@@@@
%@@@@@@@@@@@@@@@@@@@@@@@@@@@@@@@@@@@@@@@@@@@@@@@@@@@@@@@@@@
%@@@@@@@@@@@@@@@@@@@@@@@@@@@@@@@@@@@@@@@@@@@@@@@@@@@@@@@@@@
%@@@@@@@@@@@@@@@@@@@@@@@@@@@@@@@@@@@@@@@@@@@@@@@@@@@@@@@@@@
%----------------------------------------------------------
\SECT{9. Algorithms of an almost optimal decomposition of $\VLOP+\VLPME$.}{9}
%----------------------------------------------------------
\addtocontents{toc}{9. Algorithms of an almost optimal decomposition of $\VLOP+\VLPME$.\hfill \thepage\par}
%----------------------------------------------------------
\indent
%----------------------------------------------------------
\par We recall that in Subsection 8.2 we have defined a discrete Borel measure $\mu_E$  (which we call the Menger curvature measure) and a mapping $\Tc:\RT\to\RT$, see formulas \rf{D-MU-S} - \rf{DGV-2-S}.
%----------------------------------------------------------
\par Basing on Theorem \reff{S-CMS} and Proposition \reff{TAU-L}, in Subsection 8.3 we claim the existence of the mappings
%----------------------------------------------------------
$$
\Tc_1(f):\RT\to\RT~~~~~\text{and}~~~~~~\Tc_2(f):\RT\to\RT
$$
%----------------------------------------------------------
such that $\Tc(f)=\Tc_1(f)+\Tc_2(f)$ with the norms satisfying inequalities \rf{N-1} and \rf{N-2}.
%----------------------------------------------------------
\par Note that the mapping $\Tc_1$ is one of the main ingredients of our extension construction because this operator determines the pre-selection $\tg(f)$ given by formula \rf{D-PRS}.
%----------------------------------------------------------
\par Let $\VS=\VLOP+\VLPME$. In this section we describe a constructive algorithm which to every mapping $\Tc:\RT\to\RT$ assigns a mapping
$\Tc_1:\RT\to\RT$ possessing the following properties: the mapping $\Tc_1$ {\it depends linearly on $\Tc$} and the following inequalities
%----------------------------------------------------------
$$
\|\Tc_1\|_{\VLOP}\le C(p)\|\Tc\|_{\VS}~~~\text{and}~~~
\|\Tc-\Tc_1\|_{\VLPME}\le C(p)\|\Tc)\|_{\VS}
$$
%----------------------------------------------------------
hold. This enables us to put $\Tc_1(f)=(\Tc(f))_1$\,.
%----------------------------------------------------------
\par The algorithm is based on an approach suggested in the author's paper \cite{S3}. In the next two subsections we describe two main parts of this algorithm.
%----------------------------------------------------------
%@@@@@@@@@@@@@@@@@@@@@@@@@@@@@@@@@@@@@@@@@@@@@@@@@@@@@@@@@@
%@@@@@@@@@@@@@@@@@@@@@@@@@@@@@@@@@@@@@@@@@@@@@@@@@@@@@@@@@@
%@@@@@@@@@@@@@@@@@@@@@@@@@@@@@@@@@@@@@@@@@@@@@@@@@@@@@@@@@@
%@@@@@@@@@@@@@@@@@@@@@@@@@@@@@@@@@@@@@@@@@@@@@@@@@@@@@@@@@@
%@@@@@@@@@@@@@@@@@@@@@@@@@@@@@@@@@@@@@@@@@@@@@@@@@@@@@@@@@@
%@@@@@@@@@@@@@@@@@@@@@@@@@@@@@@@@@@@@@@@@@@@@@@@@@@@@@@@@@@
%----------------------------------------------------------
\bigskip
%----------------------------------------------------------
\par {\bf 9.1\,\, ``Important Squares'' and concentration of the Menger curvature measure.}
%----------------------------------------------------------
\addtocontents{toc}{~~~~9.1\,\, ``Important Squares'' and concentration of the  Menger curvature measure. \hfill \thepage\par}
%----------------------------------------------------------
\par Let $\Kc_E$ be the family of well-separated \sqs constructed in Proposition \reff{BR-TO-Q}, and let $\mu_E$ be the Menger curvature measure defined in Subsection 8.2, see \rf{D-MU-S} and \rf{MU-SET}. Recall that
%----------------------------------------------------------
$$
\Cc_E=\{c_K:K\in\Kc_E\}
$$
%----------------------------------------------------------
denotes the family of centers of all \sqs from $\Kc_E$.
%----------------------------------------------------------
\par We also recall that for each \sq $K\in\Kc_E$ the triangle $\Delta(K)$ has the following properties: $\Delta(K)\subset \gamma K$ and $\diam \Delta(K)\sim \diam K$ where $\gamma$ and the constants of the last equivalence are absolute. Hence
%----------------------------------------------------------
$$
\cu_{\Delta(K)}=\frac{1}{R_{\Delta(K)}}\le \frac{C}{\diam \Delta(K)}\le \frac{C}{r_K}
$$
%----------------------------------------------------------
with an absolute constant $C$. Then, by \rf{D-MU-S},
%----------------------------------------------------------
$$
\mu_E(K)=\cu_{\Delta(K)}^p|K|\le \frac{C^p\,(2r_K)^2}{r_K^p}=\,4 C^p\, r_K^{2-p}
$$
%----------------------------------------------------------
so that
%----------------------------------------------------------
\bel{C-DK}
\mu_E(K)\le \,\eta\, r_K^{2-p}
\ee
%----------------------------------------------------------
where $\eta:=\,4 C^p$.
%----------------------------------------------------------
\par Let us consider a function
%----------------------------------------------------------
$$
S_K(r):=\mu_E(Q(c_K,r)),~~~r\ge 0.
$$
%----------------------------------------------------------
This function is non-decreasing non-negative and right continuous on $[0,\infty)$. Furthermore, by \rf{C-DK},
$S_K(r)\le \eta r_K^{2-p}$.
%----------------------------------------------------------
\par Let $V(r):=\eta r^{2-p}, r>0$. Then
%----------------------------------------------------------
\bel{RK}
S_K(r_K)\le V(r_K).
\ee
%----------------------------------------------------------
Since $p>2$, the function $V=V(r)$ is {\it strictly decreasing}. These properties of the functions $S_K(r)$ and $V(r)$ imply the existence of a unique number $R_K\in(0,\infty)$ such that
%----------------------------------------------------------
\bel{BS}
\mu_E(Q(c_K,r))>\eta R_K^{2-p}~~~~\text{for every}~~~r>R_K,
\ee
%----------------------------------------------------------
and
%----------------------------------------------------------
\bel{MS}
\mu_E(Q(c_K,r))< \eta R_K^{2-p}~~~~\text{for every}~~~r<R_K.
\ee
%----------------------------------------------------------
\par Note that, by \rf{RK},
%----------------------------------------------------------
\bel{R-R}
R_K\ge r_K.
\ee
%----------------------------------------------------------
\par Also note that, since the function $S_K$ is right continuous on $(0,\infty)$, by \rf{BS},
%----------------------------------------------------------
$$
\mu_E(Q(c_K,R_K))\ge \eta R_K^{2-p}.
$$
%----------------------------------------------------------
\par Note one more property of the function $K\to R_K$ proven in \cite{S3}, Lemma 3.4: for every $K,K'\in \Kc_E$ the following inequality
%----------------------------------------------------------
$$
|R_K-R_{K'}|\le \|c_K-c_{K'}\|
$$
%----------------------------------------------------------
holds.
%----------------------------------------------------------
\par Given a \sq $K\in\Kc_E$ we let $\tK$ denote a \sq
%----------------------------------------------------------
\bel{D-KW}
\tK:=Q(c_K,R_K).
\ee
%----------------------------------------------------------
Note that, by \rf{R-R}, $K\subset\tK$ for every $K\in\Kc_E$. Let
%----------------------------------------------------------
$$
\widetilde{\Kc}_E:=\{\tK: K\in\Kc_E\}.
$$
%---------------------------------------------------------- %@@@@@@@@@@@@@@@@@@@@@@@@@@@@@@@@@@@@@@@@@@@@@@@@@@@@@@@@@@
%@@@@@@@@@@@@@@@@@@@@@@@@@@@@@@@@@@@@@@@@@@@@@@@@@@@@@@@@@@
%@@@@@@@@@@@@@@@@@@@@@@@@@@@@@@@@@@@@@@@@@@@@@@@@@@@@@@@@@@
%----------------------------------------------------------
\par Our aim at this stage of the algorithm is to construct a subfamily $\Qc_E$ of the family $\Kc_E$ consisting of  well separated \sqs which provides a certain {\it net} in the collection $\Kc_E$. This means that for every \sq $\tK\in\widetilde{\Kc}_E$ its fixed dilation (say,  by a factor of $12$) contains a \sq from the family $\Qc_E$.
%----------------------------------------------------------
\par The existence of such a family $\Qc_E$ is proven in \cite{S3}, Proposition 3.5. See also \cite{S3}, Subsection 6.3.
%---------------------------------------------------------- %@@@@@@@@@@@@@@@@@@@@@@@@@@@@@@@@@@@@@@@@@@@@@@@@@@@@@@@@@@
%@@@@@@@@@@@@@@@@@@@@@@@@@@@@@@@@@@@@@@@@@@@@@@@@@@@@@@@@@@
%@@@@@@@@@@@@@@@@@@@@@@@@@@@@@@@@@@@@@@@@@@@@@@@@@@@@@@@@@@
%----------------------------------------------------------
%@@@@@@@@@@@@@@@@@@@@@@@@@@@@@@@@@@@@@@@@@@@@@@@@@@@@@@@@@@
%@@@@@@@@@@@@@@@@@@@@@@@@@@@@@@@@@@@@@@@@@@@@@@@@@@@@@@@@@@
%@@@@@@@@@@@@@@@@@@@@@@@@@@@@@@@@@@@@@@@@@@@@@@@@@@@@@@@@@@
%@@@@@@@@@@@@@@@@@@@@@@@@@@@@@@@@@@@@@@@@@@@@@@@@@@@@@@@@@@
%----------------------------------------------------------
\begin{proposition}\lbl{NET-E} There exists a subfamily $\Qc_E$ of the family $\widetilde{\Kc}_E$ such that the following conditions are satisfied:\medskip
%----------------------------------------------------------
\par (i). The \sqs of the family
%----------------------------------------------------------
$$
6\Qc_E:=\{6 Q: Q\in\Qc_E\}
$$
%----------------------------------------------------------
are pairwise disjoint;\medskip
%---------------------------------------------------------- %@@@@@@@@@@@@@@@@@@@@@@@@@@@@@@@@@@@@@@@@@@@@@@@@@@@@@@@@@@
%@@@@@@@@@@@@@@@@@@@@@@@@@@@@@@@@@@@@@@@@@@@@@@@@@@@@@@@@@@
%@@@@@@@@@@@@@@@@@@@@@@@@@@@@@@@@@@@@@@@@@@@@@@@@@@@@@@@@@@
%----------------------------------------------------------
\par (ii). For each \sq $K\in\widetilde{\Kc}_E$ there exists a \sq $\hat{K}\in\Qc_E$ such that
%----------------------------------------------------------
$$
\hat{K}\subset 12K.
$$
%----------------------------------------------------------
%@@@@@@@@@@@@@@@@@@@@@@@@@@@@@@@@@@@@@@@@@@@@@@@@@@@@@@@@@@
%----------------------------------------------------------
\end{proposition}
%----------------------------------------------------------
%@@@@@@@@@@@@@@@@@@@@@@@@@@@@@@@@@@@@@@@@@@@@@@@@@@@@@@@@@@
%@@@@@@@@@@@@@@@@@@@@@@@@@@@@@@@@@@@@@@@@@@@@@@@@@@@@@@@@@@
%@@@@@@@@@@@@@@@@@@@@@@@@@@@@@@@@@@@@@@@@@@@@@@@@@@@@@@@@@@
%@@@@@@@@@@@@@@@@@@@@@@@@@@@@@@@@@@@@@@@@@@@@@@@@@@@@@@@@@@
%----------------------------------------------------------
\par {\it Proof.} In \cite{S3} we prove a general result of such a kind for an {\it arbitrary} (not necessarily finite) family of cubes in $\RN$. Here for the case of a {\it finite} family of \sqs $\widetilde{\Kc}_E$ we present a short proof of this result due to V. Dol'nikov.
%----------------------------------------------------------
\par Let $\Ac:=6\,\widetilde{\Kc}_E$. Let $K_1$ be a \sq of the minimal diameter among all the \sqs of the family $\Ac_1:=\Ac$. By $G_1$ we denote all \sqs of $\Ac_1$ which intersect $K_1$.
%----------------------------------------------------------
\par We put $\Ac_2:=\Ac_1\setminus G_1$. If $\Ac_2=\emp$ we stop and put $\Bc=\{K_1\}$.  If  $\Ac_2\ne\emp$, by $K_2$ we denote a \sq of the minimal diameter among all the \sqs of the family $\Ac_2$. We continue this procedure.  Since $\Ac$ is finite, this process will stop on a certain (finite) step $m$.
%----------------------------------------------------------
\par As a result we obtain a finite collection of pairwise disjoint \sqs $\Bc=\{K_1,...,K_m\}$ and a partition $\{G_1,...,G_m\}$ of $\Ac$ such that for each $1\le i\le m$ the following conditions are satisfied: the \sq $K_i\in G_i$,  $K_i\cap Q\ne\emp$, and $\diam K_i\le\diam Q$ for every $Q\in G_i$.
%----------------------------------------------------------
\par Thus $\Bc$ is a collection of pairwise disjoint \sqs possessing the following property: for each \sq $Q\in\Ac$ there exists a \sq $\bar{Q}\in\Bc$ such that $\bar{Q}\cap Q\ne\emp$ and $\diam \bar{Q}\le \diam Q.$ Clearly, $\bar{Q}\subset 2 Q$.
%----------------------------------------------------------
\par It remains to put
%----------------------------------------------------------
$$
\Qc_E:=\tfrac16 \Bc=\{\tfrac16 Q: Q\in\Bc\}.
$$
%----------------------------------------------------------
\par Then $\Qc_E$ is a subfamily $\widetilde{\Kc}_E$ satisfying conditions (i) and (ii) of the proposition. In fact, since the \sqs of the family $\Bc=6\Qc_E$ are pairwise disjoint, the condition (i) holds. To prove the condition (ii) fix a \sq $K\in\widetilde{\Kc}_E$. Then the \sq $Q=6K\in\Ac$ so that there exists a \sq $\bar{Q}\in\Bc$ such that $\bar{Q}\subset 2Q=12 K$.
%----------------------------------------------------------
\par Then the \sq $\hat{K}:=\tfrac16 \bar{Q}\in \tfrac16 \Bc=\Qc_E$, and $\hat{K}\subset \bar{Q}\subset 12 K$ proving the proposition.\bx\medskip
%----------------------------------------------------------
%@@@@@@@@@@@@@@@@@@@@@@@@@@@@@@@@@@@@@@@@@@@@@@@@@@@@@@@@@@
%@@@@@@@@@@@@@@@@@@@@@@@@@@@@@@@@@@@@@@@@@@@@@@@@@@@@@@@@@@
%@@@@@@@@@@@@@@@@@@@@@@@@@@@@@@@@@@@@@@@@@@@@@@@@@@@@@@@@@@
%----------------------------------------------------------
\par We refer to each \sq $Q$ which belong to the family $\Qc_E$ as an {\it ``important'' square}. As we will see in the next subsection the \sqs of the family $\Qc_E$ accumulate all information which is necessary for an almost optimal decomposition of every mapping $\Tc\in\VLOP+\VLPME$.
%----------------------------------------------------------
\par Let us present main properties of the ``important'' squares. For the proof of these properties we refer the reader to the work \cite{S3}, Section 3.\bigskip
%----------------------------------------------------------
%@@@@@@@@@@@@@@@@@@@@@@@@@@@@@@@@@@@@@@@@@@@@@@@@@@@@@@@@@@
%@@@@@@@@@@@@@@@@@@@@@@@@@@@@@@@@@@@@@@@@@@@@@@@@@@@@@@@@@@
%@@@@@@@@@@@@@@@@@@@@@@@@@@@@@@@@@@@@@@@@@@@@@@@@@@@@@@@@@@
%----------------------------------------------------------
\par Let
%----------------------------------------------------------
$$
S_E:=\bigcup_{Q\in\Qc_E} Q\,.
$$
%----------------------------------------------------------
%@@@@@@@@@@@@@@@@@@@@@@@@@@@@@@@@@@@@@@@@@@@@@@@@@@@@@@@@@@
%@@@@@@@@@@@@@@@@@@@@@@@@@@@@@@@@@@@@@@@@@@@@@@@@@@@@@@@@@@
%@@@@@@@@@@@@@@@@@@@@@@@@@@@@@@@@@@@@@@@@@@@@@@@@@@@@@@@@@@
%----------------------------------------------------------
\par \textbullet~ For every $Q\in \Qc_E$
%----------------------------------------------------------
$$
\mu_E(Q)\sim \mu_E(5Q)\sim (\diam Q)^{2-p}
$$
%----------------------------------------------------------
with constants of the equivalences depending only on $p$.
\medskip
%----------------------------------------------------------
%@@@@@@@@@@@@@@@@@@@@@@@@@@@@@@@@@@@@@@@@@@@@@@@@@@@@@@@@@@
%@@@@@@@@@@@@@@@@@@@@@@@@@@@@@@@@@@@@@@@@@@@@@@@@@@@@@@@@@@
%@@@@@@@@@@@@@@@@@@@@@@@@@@@@@@@@@@@@@@@@@@@@@@@@@@@@@@@@@@
%----------------------------------------------------------
\par \textbullet~ Let $x\in\RT$ and let $r(x)$ be a (unique) positive number such that
%@@@@@@@@@@@@@@@@@@@@@@@@@@@@@@@@@@@@@@@@@@@@@@@@@@@@@@@@@@
%----------------------------------------------------------
$$
\mu_E(Q(x,r))>\eta r(x)^{2-p}~~~~\text{for every}~~~r>r(x),
$$
%----------------------------------------------------------
and
%----------------------------------------------------------
$$
\mu_E(Q(x,r))<\eta r(x)^{2-p}~~~~\text{for every}~~~r<r(x).
$$
%----------------------------------------------------------
Clearly, by \rf{BS} and \rf{MS}, $r(c_K)=R_K$ for every \sq $K\in\Kc_E$ so that
$\widetilde{\Kc}_E\subset \widetilde{\Kc}$ where
%----------------------------------------------------------
$$
\widetilde{\Kc}:=\{Q(x,r(x)): x\in\RT\}.
$$
%----------------------------------------------------------
\par In spite of the fact that the family $\widetilde{\Kc}$ contains $\widetilde{\Kc}_E$, {\it the same family $\Qc_E$ provides a similar ``net'' in $\widetilde{\Kc}$ as well.} In other words, there exists a constant $\gamma>0$ such that for every $x\in\RT$ there exists a \sq $K^{(x)}\in\Qc_E$ for which the following inclusion
%----------------------------------------------------------
$$
K^{(x)}\subset \gamma Q(x,r(x))
$$
%----------------------------------------------------------
holds.\medskip
%----------------------------------------------------------
%@@@@@@@@@@@@@@@@@@@@@@@@@@@@@@@@@@@@@@@@@@@@@@@@@@@@@@@@@@
%@@@@@@@@@@@@@@@@@@@@@@@@@@@@@@@@@@@@@@@@@@@@@@@@@@@@@@@@@@
%@@@@@@@@@@@@@@@@@@@@@@@@@@@@@@@@@@@@@@@@@@@@@@@@@@@@@@@@@@
%----------------------------------------------------------
\par \textbullet~{\bf ``Important Squares'' as sets of  concentration of the Menger curvature measure.} Let $Q\in\Qc_E$ and let $\theta>0$. Let $K$ be a \sq in $\RT$ such that $K\supset Q$ and
%----------------------------------------------------------
$$
\diam K\le \theta\dist(K,S_E\setminus Q).
$$
%----------------------------------------------------------
Then
%----------------------------------------------------------
$$
\mu_E(K)\le C\mu_E(Q)
$$
%----------------------------------------------------------
where $C$ is a constant depending only on $p$ and $\theta$.\medskip
%----------------------------------------------------------
%@@@@@@@@@@@@@@@@@@@@@@@@@@@@@@@@@@@@@@@@@@@@@@@@@@@@@@@@@@
%@@@@@@@@@@@@@@@@@@@@@@@@@@@@@@@@@@@@@@@@@@@@@@@@@@@@@@@@@@
%@@@@@@@@@@@@@@@@@@@@@@@@@@@@@@@@@@@@@@@@@@@@@@@@@@@@@@@@@@
%----------------------------------------------------------
\par This property easily follows from the next result proven in \cite{S3}, Lemma 3.7: for every $\theta>0$ and every \sq $Q\subset\RT$ such that $\diam Q\le\theta \dist(Q,S_E)$ the following inequality
%----------------------------------------------------------
$$
\mu_E(Q)\le C(\theta,p)(\diam Q)^{2-p}
$$
%----------------------------------------------------------
holds.
%----------------------------------------------------------
%@@@@@@@@@@@@@@@@@@@@@@@@@@@@@@@@@@@@@@@@@@@@@@@@@@@@@@@@@@
%@@@@@@@@@@@@@@@@@@@@@@@@@@@@@@@@@@@@@@@@@@@@@@@@@@@@@@@@@@
%@@@@@@@@@@@@@@@@@@@@@@@@@@@@@@@@@@@@@@@@@@@@@@@@@@@@@@@@@@
%@@@@@@@@@@@@@@@@@@@@@@@@@@@@@@@@@@@@@@@@@@@@@@@@@@@@@@@@@@
%@@@@@@@@@@@@@@@@@@@@@@@@@@@@@@@@@@@@@@@@@@@@@@@@@@@@@@@@@@
%@@@@@@@@@@@@@@@@@@@@@@@@@@@@@@@@@@@@@@@@@@@@@@@@@@@@@@@@@@
%----------------------------------------------------------
\bigskip
%----------------------------------------------------------
\par {\bf 9.2\,\, Whitney-type extensions from ``important squares''.}
%----------------------------------------------------------
\addtocontents{toc}{~~~~9.2\,\, Whitney-type extensions from ``important squares''. \hfill \thepage\\\par}
%----------------------------------------------------------
\par We are in a position to construct an almost optimal decomposition of a mapping
%----------------------------------------------------------
$$
\Tc\in\VS=\VSUM.
$$
%----------------------------------------------------------
\par We let $\Cc$ denote the set of  centers of all \sqs from the collection $\Qc_E$:
%----------------------------------------------------------
$$
\Cc:=\{c_Q:Q\in\Qc_E\}.
$$
%----------------------------------------------------------
Thus $\Cc\subset\Cc_E$.
%----------------------------------------------------------
\par Given $x\in \Cc$ we let $Q^{(x)}$ denote  the (unique) \sq from $\Qc_E$ with center at the point $x$. Hence $\Qc_E=\{Q^{(x)}:x\in\Cc\}$.
%----------------------------------------------------------
%@@@@@@@@@@@@@@@@@@@@@@@@@@@@@@@@@@@@@@@@@@@@@@@@@@@@@@@@@@
%@@@@@@@@@@@@@@@@@@@@@@@@@@@@@@@@@@@@@@@@@@@@@@@@@@@@@@@@@@
%@@@@@@@@@@@@@@@@@@@@@@@@@@@@@@@@@@@@@@@@@@@@@@@@@@@@@@@@@@
%@@@@@@@@@@@@@@@@@@@@@@@@@@@@@@@@@@@@@@@@@@@@@@@@@@@@@@@@@@
%@@@@@@@@@@@@@@@@@@@@@@@@@@@@@@@@@@@@@@@@@@@@@@@@@@@@@@@@@@
%@@@@@@@@@@@@@@@@@@@@@@@@@@@@@@@@@@@@@@@@@@@@@@@@@@@@@@@@@@
%----------------------------------------------------------
\medskip
%----------------------------------------------------------
\par The $\mu_E$-measure concentration property of ``important \sqs'' motivates us to determine the component $\Tc_1\in\VLOP$ of an almost optimal decomposition of $\Tc$ using the following Whitney-type extension construction:
%----------------------------------------------------------
%@@@@@@@@@@@@@@@@@@@@@@@@@@@@@@@@@@@@@@@@@@@@@@@@@@@@@@@@@@
%@@@@@@@@@@@@@@@@@@@@@@@@@@@@@@@@@@@@@@@@@@@@@@@@@@@@@@@@@@
%@@@@@@@@@@@@@@@@@@@@@@@@@@@@@@@@@@@@@@@@@@@@@@@@@@@@@@@@@@
%@@@@@@@@@@@@@@@@@@@@@@@@@@@@@@@@@@@@@@@@@@@@@@@@@@@@@@@@@@
%@@@@@@@@@@@@@@@@@@@@@@@@@@@@@@@@@@@@@@@@@@@@@@@@@@@@@@@@@@
%@@@@@@@@@@@@@@@@@@@@@@@@@@@@@@@@@@@@@@@@@@@@@@@@@@@@@@@@@@
%----------------------------------------------------------
\par {\it Step 1.} We define a new mapping $\tTc:\Cc\to\RT$ by the formula
%----------------------------------------------------------
$$
\tTc(x)=\frac{1}{\mu_E(Q^{(x)})}
\intl_{Q^{(x)}}\Tc\,d\mu_E\,,~~~~x\in\Cc.
$$
%----------------------------------------------------------
\medskip
%----------------------------------------------------------
%@@@@@@@@@@@@@@@@@@@@@@@@@@@@@@@@@@@@@@@@@@@@@@@@@@@@@@@@@@
%@@@@@@@@@@@@@@@@@@@@@@@@@@@@@@@@@@@@@@@@@@@@@@@@@@@@@@@@@@
%@@@@@@@@@@@@@@@@@@@@@@@@@@@@@@@@@@@@@@@@@@@@@@@@@@@@@@@@@@
%@@@@@@@@@@@@@@@@@@@@@@@@@@@@@@@@@@@@@@@@@@@@@@@@@@@@@@@@@@
%@@@@@@@@@@@@@@@@@@@@@@@@@@@@@@@@@@@@@@@@@@@@@@@@@@@@@@@@@@
%@@@@@@@@@@@@@@@@@@@@@@@@@@@@@@@@@@@@@@@@@@@@@@@@@@@@@@@@@@
%----------------------------------------------------------
\par {\it Step 2.} We extend the mapping $\tTc$ from the set $\Cc$ to all of $\RT$ using the classical Whitney Extension Method. More specifically, let $W_{\Cc}$ be a Whitney covering of the open set $\RT\setminus\Cc$, and let $\{\psi_Q:Q\in W_{\Cc}\}$ be a smooth partition of unity subordinated to $W_{\Cc}$, see Section 4, Theorem \reff{Wcov}, and Section 7, Lemma \reff{P-U}.
%----------------------------------------------------------
\par Let $\eta>1$. Given a \sq $Q\in W_{\Cc}$ we let $A_Q$ denote a point which belongs to the set $(\eta Q)\cap \Cc$. (In particular, one can choose $A_Q$  to be a point nearest to $Q$ on the set $\Cc$.)
%----------------------------------------------------------
\par Then we define the mapping $\Tc_1$ by the Whitney formula:
%----------------------------------------------------------
\bel{W-TAU}
\Tc_1(x):=\left \{
%----------------------------------------------------------
\begin{array}{ll}
\tTc(x),& x\in \Cc,\\\\
\sum\limits_{Q\in W_{\Cc}}
\psi_Q(x)\,\tTc(A_Q),& x\in\RT\setminus \Cc.
\end{array}
%----------------------------------------------------------
\right.
\ee
%----------------------------------------------------------
\par The next theorem is a particular case of a general result proven in \cite{S3}.
%----------------------------------------------------------
%@@@@@@@@@@@@@@@@@@@@@@@@@@@@@@@@@@@@@@@@@@@@@@@@@@@@@@@@@@
%@@@@@@@@@@@@@@@@@@@@@@@@@@@@@@@@@@@@@@@@@@@@@@@@@@@@@@@@@@
%@@@@@@@@@@@@@@@@@@@@@@@@@@@@@@@@@@@@@@@@@@@@@@@@@@@@@@@@@@
%@@@@@@@@@@@@@@@@@@@@@@@@@@@@@@@@@@@@@@@@@@@@@@@@@@@@@@@@@@
%----------------------------------------------------------
\begin{theorem} \lbl{OPT-D} Let $2<p<\infty$ and let a mapping $\Tc\in\VS=\VSUM$. Then the mapping $\Tc_1\in\VLOP$ and the mapping $\Tc-\Tc_1\in\VLPME$. Furthermore, the following inequalities
%----------------------------------------------------------
$$
\|\Tc_1\|_{\VLOP}\le C\|\Tc\|_{\VS}\,\,,~~~~\|\Tc-\Tc_1\|_{\VLPME}\le C\|\Tc\|_{\VS}
$$
%----------------------------------------------------------
hold. Here $C$ is a constant depending only on $p$ and $\eta$.
%----------------------------------------------------------
\end{theorem}
%@@@@@@@@@@@@@@@@@@@@@@@@@@@@@@@@@@@@@@@@@@@@@@@@@@@@@@@@@@
%@@@@@@@@@@@@@@@@@@@@@@@@@@@@@@@@@@@@@@@@@@@@@@@@@@@@@@@@@@
%@@@@@@@@@@@@@@@@@@@@@@@@@@@@@@@@@@@@@@@@@@@@@@@@@@@@@@@@@@
%@@@@@@@@@@@@@@@@@@@@@@@@@@@@@@@@@@@@@@@@@@@@@@@@@@@@@@@@@@
%@@@@@@@@@@@@@@@@@@@@@@@@@@@@@@@@@@@@@@@@@@@@@@@@@@@@@@@@@@
%@@@@@@@@@@@@@@@@@@@@@@@@@@@@@@@@@@@@@@@@@@@@@@@@@@@@@@@@@@
%----------------------------------------------------------
\smallskip
%----------------------------------------------------------
\begin{remark} \lbl{R-NL} {\em In \cite{S3}, Section 6, we extend the mapping $\tTc$ to a mapping $\Tc_1\in\VLOP$ using a {\it lacunary modification of the Whitney extension method.} This modification relies on the same ideas as our modification of the Whitney construction given in Section 7.
%----------------------------------------------------------
\par In particular, let $\Lc_{\Cc}$ be the family of all lacunas of Whitney \sqs of the set $\Cc$, see Section 4. By Proposition \reff{L-PE}, there exists a ``projection''
which to every lacuna $L\in\Lc_{\Cc}$ assigns a point $\PRL(L)\in\Cc$ such that
%----------------------------------------------------------
\bel{P-L-C}
\PRL(L)\in(\gamma Q_L)\cap \Cc.
\ee
%----------------------------------------------------------
Furthermore, the mapping $L\to \PRL(L)$ is an ``almost'' one-to-one mapping, i.e., each point $a\in\RT$ has at most $N$ pre-images:
%----------------------------------------------------------
\bel{O-T-O-C}
\#\{L\in\Lc_{\Cc}: \PRL(L)=a\}\le N\, ~~~\text{for every}~~~~ a\in\RT.
\ee
%----------------------------------------------------------
Here $\gamma$ and $N$ are absolute constants.
%----------------------------------------------------------
\par Given a lacuna $L\in\Lc_{\Cc}$ we put
%----------------------------------------------------------
$$
A_Q:=\PRL(L)~~~\text{for all \sqs}~~~~Q\in L.
$$
%----------------------------------------------------------
Then we define a new mapping $\Tc_1^{(\ell)}$ by the same  formula \rf{W-TAU}.
%----------------------------------------------------------
\par By Proposition \reff{OPT-D},
%----------------------------------------------------------
$$
\|\Tc_1^{(\ell)}\|_{\VLOP}\le C\|\Tc\|_{\VS}\,\,,~~~~\|\Tc-\Tc_1^{(\ell)}\|_{\VLPME}\le C\|\Tc\|_{\VS}.
$$
%----------------------------------------------------------
\par In the forthcoming paper \cite{S-ALG} we will show that this modification provides a continuous linear extension operator with a rather simple and ``nice'' structure.\rbx}
%----------------------------------------------------------
\end{remark}
%----------------------------------------------------------
%@@@@@@@@@@@@@@@@@@@@@@@@@@@@@@@@@@@@@@@@@@@@@@@@@@@@@@@@@@
%@@@@@@@@@@@@@@@@@@@@@@@@@@@@@@@@@@@@@@@@@@@@@@@@@@@@@@@@@@
%----------------------------------------------------------
%@@@@@@@@@@@@@@@@@@@@@@@@@@@@@@@@@@@@@@@@@@@@@@@@@@@@@@@@@@
%@@@@@@@@@@@@@@@@@@@@@@@@@@@@@@@@@@@@@@@@@@@@@@@@@@@@@@@@@@
%@@@@@@@@@@@@@@@@@@@@@@@@@@@@@@@@@@@@@@@@@@@@@@@@@@@@@@@@@@
%@@@@@@@@@@@@@@@@@@@@@@@@@      @@@@@@@@@@@@@@@@@@@@@@@@@@@
%@@@@@@@@@@@@@@@@@@@@@@@          @@@@@@@@@@@@@@@@@@@@@@@@@
%@@@@@@@@@@@@@@@@@@@@@              @@@@@@@@@@@@@@@@@@@@@@@
%@@@@@@@@@@@@@@@@@@@     SECTION 10    @@@@@@@@@@@@@@@@@@@@@
%@@@@@@@@@@@@@@@@@@@@@              @@@@@@@@@@@@@@@@@@@@@@@
%@@@@@@@@@@@@@@@@@@@@@@@          @@@@@@@@@@@@@@@@@@@@@@@@@
%@@@@@@@@@@@@@@@@@@@@@@@@@      @@@@@@@@@@@@@@@@@@@@@@@@@@@
%@@@@@@@@@@@@@@@@@@@@@@@@@@@@@@@@@@@@@@@@@@@@@@@@@@@@@@@@@@
%@@@@@@@@@@@@@@@@@@@@@@@@@@@@@@@@@@@@@@@@@@@@@@@@@@@@@@@@@@
%@@@@@@@@@@@@@@@@@@@@@@@@@@@@@@@@@@@@@@@@@@@@@@@@@@@@@@@@@@
%----------------------------------------------------------
\SECT{10. The extension algorithm and its linearity.}{10}
%----------------------------------------------------------
\addtocontents{toc}{10. The extension algorithm and its linearity. \hfill \thepage\par}
%----------------------------------------------------------
\indent \par In this section we describe main stages of the algorithm which provides an almost optimal extension of every function from the trace space $\LTP|_E $ whenever $E\subset\RT$ is an arbitrary finite set and $p>2$. An optimality of this algorithm follows from the proof of Theorem \reff{MAIN}. We will see that this algorithm provides a {\it continuous linear extension operator} for the trace space $\LTP|_E $.
%----------------------------------------------------------
\par The algorithm consists of two main parts which we call ``Pre-work'' and ``Main Algorithm''.
%----------------------------------------------------------
\par In the part ``Pre-work'' we introduce a series of geometrical objects (lacunae, bridges, families of  triangles, measures, etc.) which we use later on in our extension construction. Note that all these objects are determined only by a geometrical structure of the set $E$ and does not depend on values of functions defined on $E$.
%----------------------------------------------------------
\par In the part ``Main Algorithm'' we fix a function $f:E\to\R$. Its values on the set $E$ we consider as input data of our algorithm. Then we show  how these input data transform into the values of an almost optimal extension $F\in\LTP$ of the function $f$. At every step of the algorithm we control its ``linearity'', i.e., linear dependence of all elements of our construction on the data.
%----------------------------------------------------------
\par In this section our aim is only to present more or less detailed description of the extension algorithm. In this paper we do not estimate its complexity and the order of magnitude of the number of elementary operation which is necessary for computation of the extension $F\in\LTP$. We will study these and other problems related to the optimality of the algorithm in the forthcoming paper \cite {S-ALG}.
%---------------------------------------------------------- %@@@@@@@@@@@@@@@@@@@@@@@@@@@@@@@@@@@@@@@@@@@@@@@@@@@@@@@@@@
%@@@@@@@@@@@@@@@@@@@@@@@@@@@@@@@@@@@@@@@@@@@@@@@@@@@@@@@@@@
%@@@@@@@@@@@@@@@@@@@@@@@@@@@@@@@@@@@@@@@@@@@@@@@@@@@@@@@@@@
%@@@@@@@@@@@@@@@@@@@@@@@@@@@@@@@@@@@@@@@@@@@@@@@@@@@@@@@@@@
%@@@@@@@@@@@@@@@@@@@@@@@@@@@@@@@@@@@@@@@@@@@@@@@@@@@@@@@@@@
%@@@@@@@@@@@@@@@@@@@@@@@@@@@@@@@@@@@@@@@@@@@@@@@@@@@@@@@@@@
%@@@@@@@@@@@@@@@@@@@@@@@@@@@@@@@@@@@@@@@@@@@@@@@@@@@@@@@@@@
%@@@@@@@@@@@@@@@@@@@@@@@@@@@@@@@@@@@@@@@@@@@@@@@@@@@@@@@@@@
%@@@@@@@@@@@@@@@@@@@@@@@@@@@@@@@@@@@@@@@@@@@@@@@@@@@@@@@@@@
%----------------------------------------------------------
\bigskip
\par {\bf 10.1. Pre-work: lacunae, bridges, and measure concentration squares.}
%----------------------------------------------------------
\addtocontents{toc}{~~~~10.1. Pre-work: lacunae, bridges, and measure concentration squares. \hfill \thepage\par}
%----------------------------------------------------------
%@@@@@@@@@@@@@@@@@@@@@@@@@@@@@@@@@@@@@@@@@@@@@@@@@@@@@@@@@@
%@@@@@@@@@@@@@@@@@@@@@@@@@@@@@@@@@@@@@@@@@@@@@@@@@@@@@@@@@@
%@@@@@@@@@@@@@@@@@@@@@@@@@@@@@@@@@@@@@@@@@@@@@@@@@@@@@@@@@@
%@@@@@@@@@@@@@@@@@@@@@@@@@@@@@@@@@@@@@@@@@@@@@@@@@@@@@@@@@@
%@@@@@@@@@@@@@@@@@@@@@@@@@@@@@@@@@@@@@@@@@@@@@@@@@@@@@@@@@@
%@@@@@@@@@@@@@@@@@@@@@@@@@@@@@@@@@@@@@@@@@@@@@@@@@@@@@@@@@@
%----------------------------------------------------------
\medskip
\par {\bf Step 1: Whitney squares and a partition of unity.} We fix a family $W_E$ of Whitney \sqs satisfying conditions of Theorem \reff{Wcov}. We also fix a smooth partition of unity subordinated to the Whitney decomposition $W_E$, see Lemma \reff{P-U}. There are various algorithms for constructing of these two classical elements of Whitney extension methods, see, e.g. Stein \cite{St}.
%---------------------------------------------------------- %@@@@@@@@@@@@@@@@@@@@@@@@@@@@@@@@@@@@@@@@@@@@@@@@@@@@@@@@@@
%@@@@@@@@@@@@@@@@@@@@@@@@@@@@@@@@@@@@@@@@@@@@@@@@@@@@@@@@@@
%@@@@@@@@@@@@@@@@@@@@@@@@@@@@@@@@@@@@@@@@@@@@@@@@@@@@@@@@@@
%@@@@@@@@@@@@@@@@@@@@@@@@@@@@@@@@@@@@@@@@@@@@@@@@@@@@@@@@@@
%@@@@@@@@@@@@@@@@@@@@@@@@@@@@@@@@@@@@@@@@@@@@@@@@@@@@@@@@@@
%@@@@@@@@@@@@@@@@@@@@@@@@@@@@@@@@@@@@@@@@@@@@@@@@@@@@@@@@@@
%----------------------------------------------------------
\medskip
\par {\bf Step 2: Lacunae of Whitney squares, a graph of lacunae, and interior bridges.} At this step of the pre-work we construct the family $\Lc_E$ of lacunae of the set $E$, see Subsection 4.1. For every lacuna $L\in\Lc_E$ we fix all its contacting lacunae $L'\in\Lc_E$, $L'\lr L$. See Subsection 5.2, Definition \reff{D-CNL}.
%----------------------------------------------------------
\par Then we construct an interior bridge $(A_L,B_L)$ of the lacuna $L$. See Definition \reff{FR-L}.
%----------------------------------------------------------
%@@@@@@@@@@@@@@@@@@@@@@@@@@@@@@@@@@@@@@@@@@@@@@@@@@@@@@@@@@
%@@@@@@@@@@@@@@@@@@@@@@@@@@@@@@@@@@@@@@@@@@@@@@@@@@@@@@@@@@
%@@@@@@@@@@@@@@@@@@@@@@@@@@@@@@@@@@@@@@@@@@@@@@@@@@@@@@@@@@
%@@@@@@@@@@@@@@@@@@@@@@@@@@@@@@@@@@@@@@@@@@@@@@@@@@@@@@@@@@
%@@@@@@@@@@@@@@@@@@@@@@@@@@@@@@@@@@@@@@@@@@@@@@@@@@@@@@@@@@
%@@@@@@@@@@@@@@@@@@@@@@@@@@@@@@@@@@@@@@@@@@@@@@@@@@@@@@@@@@
%----------------------------------------------------------
\medskip
\par {\bf Step 3: Bridges between lacunae.} Using the technique described in Section 6 we construct the family $\BRE$ of all bridges between lacunae. Simultaneously with constructing of bridges we determine pairs of {\it connected bridges} $T,T'\in\BRE$, $T\bcn T'$, as it was done in Section 6 for bridges satisfying conditions \rf{DL-1} and \rf{CASE-2}.
%----------------------------------------------------------
\par Then, basing on the algorithm suggested in Proposition \reff{BR-TO-Q}, we construct a one-to-one mapping which to every pair $(T,T')$ of connected bridges assigns a \sq $K(T,T')$ satisfying conditions \rf{MP-1} and \rf{MP-2} of the proposition.
%----------------------------------------------------------
\par As a result we obtain a one-to-one mapping defined on a family of well-separated \sqs $\Kc_E$ which to every \sq $K\in\Kc_E$ assigns a pair of connected bridges $T_K,T'_K\in\BRE$, $T_K\bcn T'_K$. The ends of these bridges form a triangle $\Delta(K)$ (true or degenerate) as it described in Subsection 8.2.
%---------------------------------------------------------- %@@@@@@@@@@@@@@@@@@@@@@@@@@@@@@@@@@@@@@@@@@@@@@@@@@@@@@@@@@
%@@@@@@@@@@@@@@@@@@@@@@@@@@@@@@@@@@@@@@@@@@@@@@@@@@@@@@@@@@
%@@@@@@@@@@@@@@@@@@@@@@@@@@@@@@@@@@@@@@@@@@@@@@@@@@@@@@@@@@
%@@@@@@@@@@@@@@@@@@@@@@@@@@@@@@@@@@@@@@@@@@@@@@@@@@@@@@@@@@
%@@@@@@@@@@@@@@@@@@@@@@@@@@@@@@@@@@@@@@@@@@@@@@@@@@@@@@@@@@
%@@@@@@@@@@@@@@@@@@@@@@@@@@@@@@@@@@@@@@@@@@@@@@@@@@@@@@@@@@
%----------------------------------------------------------
\medskip
\par {\bf Step 4: The Menger curvature measure $\mu_E$ and ``important'' squares.} We fix the set $\Cc_E=\{c_Q:Q\in\Kc_E\}$ of centers of the \sqs from the family $\Kc_E$, see \rf{CC-E}. Then we construct the measure $\mu_E$ by formulas \rf{D-MU-S} and \rf{MU-SET}.
%----------------------------------------------------------
\par Basing on the approach suggested in Subsection 9.1, we determine the family of ``important'' squares generated by the measure $\mu_E$. We do this in two steps. First for each \sq $K\in \Kc_E$ we construct the \sq $\tK$, see \rf{D-KW}, whose ``radius'' $R_{\tK}$ satisfies inequalities \rf{BS}{ and \rf{MS}. (Recall that $\tK\supset K$.)
%----------------------------------------------------------
\par We obtain a new family of \sqs $\WKE=\{\tK:K\in\Kc_E\}$. Then we extract from this family a subfamily $\Qc_E$ of ``important'' squares satisfying  conditions (i) and (ii) of Proposition \reff{NET-E}. We do this basing on a constructive filtering procedure suggested in the proof of this proposition.
%---------------------------------------------------------- %@@@@@@@@@@@@@@@@@@@@@@@@@@@@@@@@@@@@@@@@@@@@@@@@@@@@@@@@@@
%@@@@@@@@@@@@@@@@@@@@@@@@@@@@@@@@@@@@@@@@@@@@@@@@@@@@@@@@@@
%@@@@@@@@@@@@@@@@@@@@@@@@@@@@@@@@@@@@@@@@@@@@@@@@@@@@@@@@@@
%@@@@@@@@@@@@@@@@@@@@@@@@@@@@@@@@@@@@@@@@@@@@@@@@@@@@@@@@@@
%@@@@@@@@@@@@@@@@@@@@@@@@@@@@@@@@@@@@@@@@@@@@@@@@@@@@@@@@@@
%@@@@@@@@@@@@@@@@@@@@@@@@@@@@@@@@@@@@@@@@@@@@@@@@@@@@@@@@@@
%----------------------------------------------------------
\medskip
\par {\bf Step 5: Centers of ``important'' \sqs and lacunary Whitney-type extensions.} We fix the set $\Cc=\{c_Q:Q\in\Qc_E\}$ of centers of the ``important'' \sqs. At this step of the pre-work we construct all necessary ingredients of a {\it lacunary Whitney-type extension} from the set $\Cc$. We have described this modification of the Whitney method in Remark \reff{R-NL}. %----------------------------------------------------------
\par First we fix a Whitney covering $W_{\Cc}$ of the open set $\RT\setminus \Cc$ and a smooth partition of unity $\{\psi_Q:Q\in W_{\Cc}\}$ subordinated to $W_{\Cc}$.
%----------------------------------------------------------
\par Then we determine a family $\Lc_{\Cc}$ of all lacunae of the set $\Cc$. Finally, we construct a ``projection'' of $\Lc_{\Cc}$ into $\Cc$, i.e., an ``almost'' one-to-one mapping $\Lc_{\Cc}\ni L\mapsto \PRL(L)\in\Cc$ satisfying conditions \rf{P-L-C} and \rf{O-T-O-C}.
%----------------------------------------------------------
\par All preparations are finished, and we turn to the second part of the algorithm.
%---------------------------------------------------------- %@@@@@@@@@@@@@@@@@@@@@@@@@@@@@@@@@@@@@@@@@@@@@@@@@@@@@@@@@@
%@@@@@@@@@@@@@@@@@@@@@@@@@@@@@@@@@@@@@@@@@@@@@@@@@@@@@@@@@@
%@@@@@@@@@@@@@@@@@@@@@@@@@@@@@@@@@@@@@@@@@@@@@@@@@@@@@@@@@@
%@@@@@@@@@@@@@@@@@@@@@@@@@@@@@@@@@@@@@@@@@@@@@@@@@@@@@@@@@@
%@@@@@@@@@@@@@@@@@@@@@@@@@@@@@@@@@@@@@@@@@@@@@@@@@@@@@@@@@@
%@@@@@@@@@@@@@@@@@@@@@@@@@@@@@@@@@@@@@@@@@@@@@@@@@@@@@@@@@@
%@@@@@@@@@@@@@@@@@@@@@@@@@@@@@@@@@@@@@@@@@@@@@@@@@@@@@@@@@@
%@@@@@@@@@@@@@@@@@@@@@@@@@@@@@@@@@@@@@@@@@@@@@@@@@@@@@@@@@@
%@@@@@@@@@@@@@@@@@@@@@@@@@@@@@@@@@@@@@@@@@@@@@@@@@@@@@@@@@@
%----------------------------------------------------------
\bigskip
\par {\bf 10.2. Main Algorithm: from values of a function to its almost optimal $\LTP$-extension.}
%----------------------------------------------------------
\addtocontents{toc}{~~~~10.2. Main Algorithm: from values of a function to its almost optimal\\
%----------------------------------------------------------
\hspace*{13mm}~~~~ $\LTP$-extension. \hfill \thepage\\\par}
%----------------------------------------------------------
\medskip
\par In this subsection we present main steps of the extension algorithm which to an arbitrary function $f:E\to\R$ assigns its almost optimal extension $F\in\LTP$, $p>2$. Furthermore, the extension $F$ depends linearly  on the function $f$. This algorithm uses only the values of the function $f$ on $E$ and the geometrical objects which we have constructed at the stage of the pre-work of the algorithm.
%----------------------------------------------------------
%@@@@@@@@@@@@@@@@@@@@@@@@@@@@@@@@@@@@@@@@@@@@@@@@@@@@@@@@@@
%@@@@@@@@@@@@@@@@@@@@@@@@@@@@@@@@@@@@@@@@@@@@@@@@@@@@@@@@@@
%@@@@@@@@@@@@@@@@@@@@@@@@@@@@@@@@@@@@@@@@@@@@@@@@@@@@@@@@@@
%@@@@@@@@@@@@@@@@@@@@@@@@@@@@@@@@@@@@@@@@@@@@@@@@@@@@@@@@@@
%@@@@@@@@@@@@@@@@@@@@@@@@@@@@@@@@@@@@@@@@@@@@@@@@@@@@@@@@@@
%@@@@@@@@@@@@@@@@@@@@@@@@@@@@@@@@@@@@@@@@@@@@@@@@@@@@@@@@@@
%----------------------------------------------------------
\medskip
\par {\bf Step 1: the mapping $\Tc(f)$ and its averages on ``important'' squares.} Using formulas \rf{DGV-1-S} and \rf{DGV-2-S}, we construct the mapping $\Tc(f):\RT\to\RT$; thus
%----------------------------------------------------------
$$
\Tc(c_K;f)=\nabla P_{\Delta(K)}[f]
$$
%----------------------------------------------------------
whenever $K\in\Kc_E$ and $\Delta(K)$  is a true triangle, and $\Tc(c_K;f)=0$ if $\Delta(K)$ is a degenerate triangle. The mapping $\Tc(x;f)=0$ for every point $x\in\RT\setminus\Cc_E$ (i.e., out of the family of centers of the \sqs from $\Kc_E$).
%----------------------------------------------------------
\par Recall that $P_{\Delta(K)}[f]$ is the affine polynomial which interpolates $f$ on the vertices of $\Delta(K)$. Clearly, its gradient depends linearly  on $f$ so that the mapping $\Tc(f)$  depends linearly on $f$ as well.
%----------------------------------------------------------
\par Let $\Cc=\{c_Q:Q\in\Qc_E\}$ be the family of centers of the ``important'' squares, see Step 4 of the pre-work. We define a new mapping $\TT(f):\Cc\to\RT$ such that for every``important'' \sq $Q\in\Qc_E$
%----------------------------------------------------------
$$
\TT(c_Q;f):=\frac{1}{\mu_E(Q)}\intl_Q
\Tc(z;f)\,d\mu_E(z).
$$
%----------------------------------------------------------
\par Note that, by definition of $\Tc(f)$ and the Menger curvature measure, see \rf{MU-SET}, for each $Q\in\Qc_E$ we have
%----------------------------------------------------------
$$
\mu_E(Q)=\sum\left\{\cu_{\Delta(K)}^p\,|K|: K\in\Kc_E,
c_K\in Q\right\}
$$
%----------------------------------------------------------
and
%----------------------------------------------------------
$$
\TT(c_Q;f)=\frac{1}{\mu_E(Q)}\,
\sum \left\{\nabla P_{\Delta(K)}[f]\,\cu^p_{\Delta(K)}\,|K|:
K\in\Kc_E,c_K\in Q\right\}.
$$
%----------------------------------------------------------
\par Since the affine polynomials $P_{\Delta(K)}[f]$ depend linearly on $f$, the mapping $\TT(f)$ depends linearly  on $f$ as well.
%----------------------------------------------------------
%@@@@@@@@@@@@@@@@@@@@@@@@@@@@@@@@@@@@@@@@@@@@@@@@@@@@@@@@@@
%@@@@@@@@@@@@@@@@@@@@@@@@@@@@@@@@@@@@@@@@@@@@@@@@@@@@@@@@@@
%@@@@@@@@@@@@@@@@@@@@@@@@@@@@@@@@@@@@@@@@@@@@@@@@@@@@@@@@@@
%@@@@@@@@@@@@@@@@@@@@@@@@@@@@@@@@@@@@@@@@@@@@@@@@@@@@@@@@@@
%@@@@@@@@@@@@@@@@@@@@@@@@@@@@@@@@@@@@@@@@@@@@@@@@@@@@@@@@@@
%@@@@@@@@@@@@@@@@@@@@@@@@@@@@@@@@@@@@@@@@@@@@@@@@@@@@@@@@@@
%----------------------------------------------------------
\medskip
\par {\bf Step 2: the component $\Tc_1(f)$.} At this step we extend the mapping $\TT(f)$ from the set $\Cc$ to all of $\RT$ using the lacunary modification of the Whitney method described at Step 5 of the pre-work. More specifically, let $W_{\Cc}$ be the Whitney covering of $\RT\setminus\Cc$ introduced at this step. Let $L\in\Lc_{\Cc}$ be a lacuna of Whitney \sqs from $W_{\Cc}$, and let $Q\in L$. Recall that at Step 5 of the pre-work we have also constructed a ``projection'' $\Lc_{\Cc}\ni L\mapsto \PRL(L)\in \Cc$. We put
%----------------------------------------------------------
$$
A_Q:=\PRL(L)~~~\text{for every}~~~Q\in L.
$$
%----------------------------------------------------------
Then we construct a Whitney-type extension $\Tc_1(f):\RT\to\RT$ of the mapping $\TT(f)$:
%----------------------------------------------------------
$$
\Tc_1(x;f):=\left \{
%----------------------------------------------------------
\begin{array}{ll}
\TT(x;f),& x\in \Cc,\\\\
\sum\limits_{Q\in W_{\Cc}}
\psi_Q(x)\,\TT(A_Q;f),& x\in\RT\setminus \Cc.
\end{array}
%----------------------------------------------------------
\right.
$$
%----------------------------------------------------------
Here $\{\psi_Q:Q\in W_{\Cc}\}$  is a smooth partition of unity subordinated to $W_{\Cc}$ which we have determined at Step 5 of the pre-work.
%----------------------------------------------------------
\par Obviously, since $\TT(f)$ depends linearly on $f$ and the Whitney extension operator is a linear operator, the mapping  $\Tc_1(f)$ depends linearly  on $f$.
%----------------------------------------------------------
%@@@@@@@@@@@@@@@@@@@@@@@@@@@@@@@@@@@@@@@@@@@@@@@@@@@@@@@@@@
%@@@@@@@@@@@@@@@@@@@@@@@@@@@@@@@@@@@@@@@@@@@@@@@@@@@@@@@@@@
%@@@@@@@@@@@@@@@@@@@@@@@@@@@@@@@@@@@@@@@@@@@@@@@@@@@@@@@@@@
%@@@@@@@@@@@@@@@@@@@@@@@@@@@@@@@@@@@@@@@@@@@@@@@@@@@@@@@@@@
%@@@@@@@@@@@@@@@@@@@@@@@@@@@@@@@@@@@@@@@@@@@@@@@@@@@@@@@@@@
%@@@@@@@@@@@@@@@@@@@@@@@@@@@@@@@@@@@@@@@@@@@@@@@@@@@@@@@@@@
%----------------------------------------------------------
\medskip
\par {\bf Step 3: pre-selections and selections.} Let $L\in\Lc_E$ be an arbitrary lacuna of the set $E$, and let
$T=T(L)=(A_L,B_L)$ be its interior bridge with the ends at points $A_L,B_L\in E$, see Step 2 of the pre-work. Using the formulas \rf{D-PRS} and \rf{D-SEL} at this step we construct a pre-selection $\tg(T;f)$ and a selection $g(T;f)$ of the set valued mapping $G_f$, see \rf{D-GT}.
%----------------------------------------------------------
\par Thus we put
%----------------------------------------------------------
$$
\tg(T;f)=\Tc_1(A_L;f),
$$
%----------------------------------------------------------
and
%----------------------------------------------------------
$$
g(T;f)=\PR(\tg(T;f);G_f(T)).
$$
%----------------------------------------------------------
Recall that given a straight line $\ell\subset\RT$  by $\PR(x;\ell)$ we denote the orthogonal projection of a point $x\in\RT$ onto $\ell$.
%----------------------------------------------------------
\par Clearly, since $\Tc_1(f)$ depends linearly on $f$, the same is true for the pre-selection $\tg(f)$.
%----------------------------------------------------------
\par Let us see that the selection $g(f)$ depends linearly  on $f$ as well. In fact, let us present the straight line
%----------------------------------------------------------
$$
G_f(T)=\{z\in \RT: \ip{z,A_L-B_L}=f(A_L)-f(B_L)\,\}
$$
%----------------------------------------------------------
in the form
%----------------------------------------------------------
$$
G_f(T)=\{z\in \RT: \ip{z,n_T}=D_T(f)\,\}
$$
%----------------------------------------------------------
where
%----------------------------------------------------------
$$
n_T=\frac{A_L-B_L}{\|A_L-B_L\|_2}~~~\text{and}~~~~
D_T(f)=\frac{f(A_L)-f(B_L)}{\|A_L-B_L\|_2}\,.
$$
%----------------------------------------------------------
By $H_T$ we denote a one dimensional linear subspace of $\RT$
%----------------------------------------------------------
$$
H_T=\{z\in \RT: \ip{z,n_T}=0\,\}.
$$
%----------------------------------------------------------
\par Then
%----------------------------------------------------------
$$
G_f(T)=H_T+ D_T(f)\,n_T
$$
%----------------------------------------------------------
so that for every $x\in\RT$
%----------------------------------------------------------
$$
\Pr(x;G_f(T))=\Pr(x;H_T)+D_T(f)\,n_T \,.
$$
%----------------------------------------------------------
\par In particular,
%----------------------------------------------------------
$$
g(T;f)=\Pr(\tg(T;f);G_f(T))
=\Pr(\tg(T;f);H_T)+D_T(f)\,n_T\,. \,.
$$
%----------------------------------------------------------
Note that $\Pr(\cdot;H_T)$ is a linear operator, the vector function $f\to D_T(f)\,n_T$ depends linearly  on $f$, and the same is true for the mapping $f\to \tg(f)$. Hence we conclude that the selection $g(f)$ depends linearly on $f$.
%----------------------------------------------------------
%@@@@@@@@@@@@@@@@@@@@@@@@@@@@@@@@@@@@@@@@@@@@@@@@@@@@@@@@@@
%@@@@@@@@@@@@@@@@@@@@@@@@@@@@@@@@@@@@@@@@@@@@@@@@@@@@@@@@@@
%@@@@@@@@@@@@@@@@@@@@@@@@@@@@@@@@@@@@@@@@@@@@@@@@@@@@@@@@@@
%@@@@@@@@@@@@@@@@@@@@@@@@@@@@@@@@@@@@@@@@@@@@@@@@@@@@@@@@@@
%@@@@@@@@@@@@@@@@@@@@@@@@@@@@@@@@@@@@@@@@@@@@@@@@@@@@@@@@@@
%@@@@@@@@@@@@@@@@@@@@@@@@@@@@@@@@@@@@@@@@@@@@@@@@@@@@@@@@@@
%----------------------------------------------------------
\medskip
\par {\bf Step 4: ``almost optimal'' affine polynomials.} At this step given an arbitrary lacuna $L\in\Lc_E$ and its interior bridge $T(L)$ with the ends at points $A_L,B_L\in E$ we construct an affine polynomial $P_L(f)\in\PO$ such that
%----------------------------------------------------------
$$
P_L(A_L;f)=f(A_L),~~~P_L(B_L;f)=f(B_L)
$$
%----------------------------------------------------------
and
%----------------------------------------------------------
$$
\nabla P_L(f)=g(T(L);f)
$$
%----------------------------------------------------------
where $T(L)=(A_L,B_L)\in\BRE$ is the interior bridge of the lacuna $L$. Thus
%----------------------------------------------------------
\bel{D-PLF}
P_L(x;f)=f(A_L)+\ip{g(T(L);f),x-A_L},~~~~x\in\RT\,.
\ee
%----------------------------------------------------------
\par As we have proved at the previous step, the vector function $f\to g(T(L);f)$ is a linear function of $f$, so that, by formula \rf{D-PLF}, the affine polynomial $P_L(f)$ depends linearly on $f$.
%----------------------------------------------------------
%@@@@@@@@@@@@@@@@@@@@@@@@@@@@@@@@@@@@@@@@@@@@@@@@@@@@@@@@@@
%@@@@@@@@@@@@@@@@@@@@@@@@@@@@@@@@@@@@@@@@@@@@@@@@@@@@@@@@@@
%@@@@@@@@@@@@@@@@@@@@@@@@@@@@@@@@@@@@@@@@@@@@@@@@@@@@@@@@@@
%@@@@@@@@@@@@@@@@@@@@@@@@@@@@@@@@@@@@@@@@@@@@@@@@@@@@@@@@@@
%@@@@@@@@@@@@@@@@@@@@@@@@@@@@@@@@@@@@@@@@@@@@@@@@@@@@@@@@@@
%@@@@@@@@@@@@@@@@@@@@@@@@@@@@@@@@@@@@@@@@@@@@@@@@@@@@@@@@@@
%----------------------------------------------------------
\medskip
\par {\bf Step 5: the extension operator.} This is the final step of the algorithm. We apply the lacunary Whitney-type extension operator suggested in Section 7 to the family of affine polynomials $\{P_L(f):L\in\Lc_E\}$ from the previous step and construct the required extension of the function $f$.
%----------------------------------------------------------
\par Let $L\in\Lc_E$ be a lacuna.  We put
%----------------------------------------------------------
$$
\PQ(f)=P_L(f)~~~\text{for every square} ~~~ Q\in L
$$
%----------------------------------------------------------
where $P_L(f)$ is the polynomial defined by \rf{D-PLF}.
%----------------------------------------------------------
\par Finally, we construct the extension $F(f):\RT\to\RT$ by the formula
%----------------------------------------------------------
$$
F(x;f):=\left \{
%----------------------------------------------------------
\begin{array}{ll}
f(x),& x\in E,\\\\
\sum\limits_{Q\in W_E}
\varphi_Q(x)\,\PQ(x;f),& x\in\RT\setminus E.
\end{array}
%----------------------------------------------------------
\right.
$$
%----------------------------------------------------------
Here $W_E$ and $\{\varphi_Q:Q\in W_E\}$ are the Whitney covering and the smooth partition of unity subordinated to $W_E$ respectively. See Step 1 of the pre-work.
%----------------------------------------------------------
\par Since the polynomials of the family $\{P_L(f):L\in\Lc_E\}$ depend linearly on $f$, {\it the extension operator $f\to F(f)$ is linear}.
%----------------------------------------------------------
%@@@@@@@@@@@@@@@@@@@@@@@@@@@@@@@@@@@@@@@@@@@@@@@@@@@@@@@@@@
%@@@@@@@@@@@@@@@@@@@@@@@@@@@@@@@@@@@@@@@@@@@@@@@@@@@@@@@@@@
%@@@@@@@@@@@@@@@@@@@@@@@@@@@@@@@@@@@@@@@@@@@@@@@@@@@@@@@@@@
%@@@@@@@@@@@@@@@@@@@@@@@@@      @@@@@@@@@@@@@@@@@@@@@@@@@@@
%@@@@@@@@@@@@@@@@@@@@@@@          @@@@@@@@@@@@@@@@@@@@@@@@@
%@@@@@@@@@@@@@@@@@@@@@              @@@@@@@@@@@@@@@@@@@@@@@
%@@@@@@@@@@@@@@@@@@@     SECTION 11    @@@@@@@@@@@@@@@@@@@@@
%@@@@@@@@@@@@@@@@@@@@@              @@@@@@@@@@@@@@@@@@@@@@@
%@@@@@@@@@@@@@@@@@@@@@@@          @@@@@@@@@@@@@@@@@@@@@@@@@
%@@@@@@@@@@@@@@@@@@@@@@@@@      @@@@@@@@@@@@@@@@@@@@@@@@@@@
%@@@@@@@@@@@@@@@@@@@@@@@@@@@@@@@@@@@@@@@@@@@@@@@@@@@@@@@@@@
%@@@@@@@@@@@@@@@@@@@@@@@@@@@@@@@@@@@@@@@@@@@@@@@@@@@@@@@@@@
%@@@@@@@@@@@@@@@@@@@@@@@@@@@@@@@@@@@@@@@@@@@@@@@@@@@@@@@@@@
\bigskip
%----------------------------------------------------------
\SECT{11. Refinements of the trace criterion: Theorems \reff{REF-MAIN}, \reff{REF-SP} and \reff{NEW-MAIN}.} {11}
%----------------------------------------------------------
\addtocontents{toc}{11. Refinements of the trace criterion: Theorems \reff{REF-MAIN}, \reff{REF-SP} and \reff{NEW-MAIN}.\hfill \thepage\par}
%----------------------------------------------------------
\indent
%----------------------------------------------------------
%@@@@@@@@@@@@@@@@@@@@@@@@@@@@@@@@@@@@@@@@@@@@@@@@@@@@@@@@@@
%@@@@@@@@@@@@@@@@@@@@@@@@@@@@@@@@@@@@@@@@@@@@@@@@@@@@@@@@@@
%@@@@@@@@@@@@@@@@@@@@@@@@@@@@@@@@@@@@@@@@@@@@@@@@@@@@@@@@@@
%@@@@@@@@@@@@@@@@@@@@@@@@@@@@@@@@@@@@@@@@@@@@@@@@@@@@@@@@@@
%@@@@@@@@@@@@@@@@@@@@@@@@@@@@@@@@@@@@@@@@@@@@@@@@@@@@@@@@@@
%@@@@@@@@@@@@@@@@@@@@@@@@@@@@@@@@@@@@@@@@@@@@@@@@@@@@@@@@@@
%----------------------------------------------------------
%@@@@@@@@@@@@@@@@@@@@@@@@@@@@@@@@@@@@@@@@@@@@@@@@@@@@@@@@@@
%@@@@@@@@@@@@@@@@@@@@@@@@@@@@@@@@@@@@@@@@@@@@@@@@@@@@@@@@@@
%@@@@@@@@@@@@@@@@@@@@@@@@@@@@@@@@@@@@@@@@@@@@@@@@@@@@@@@@@@
%\par  PICTURE 2: ILLUSTRATION TO THE THEOREM                                   %@@@@@@@@@@@@@@@@@@@@@@@@@@@@@@@@@@@@@@@@@@@@@@@@@@@@@@@@@@
%@@@@@@@@@@@@@@@@@@@@@@@@@@@@@@@@@@@@@@@@@@@@@@@@@@@@@@@@@@
%@@@@@@@@@@@@@@@@@@@@@@@@@@@@@@@@@@@@@@@@@@@@@@@@@@@@@@@@@@
%@@@@@@@@@@@@@@@@@@@@@@@@@@@@@@@@@@@@@@@@@@@@@@@@@@@@@@@@@@
%----------------------------------------------------------
\par {\bf 11.1. Proof of Theorem \reff{REF-MAIN}.}
%----------------------------------------------------------
\addtocontents{toc}{~~~~11.1. Proof of Theorem \reff{REF-MAIN}.\hfill \thepage\par}
%----------------------------------------------------------
%@@@@@@@@@@@@@@@@@@@@@@@@@@@@@@@@@@@@@@@@@@@@@@@@@@@@@@@@@@
%@@@@@@@@@@@@@@@@@@@@@@@@@@@@@@@@@@@@@@@@@@@@@@@@@@@@@@@@@@
%@@@@@@@@@@@@@@@@@@@@@@@@@@@@@@@@@@@@@@@@@@@@@@@@@@@@@@@@@@
%@@@@@@@@@@@@@@@@@@@@@@@@@@@@@@@@@@@@@@@@@@@@@@@@@@@@@@@@@@
%@@@@@@@@@@@@@@@@@@@@@@@@@@@@@@@@@@@@@@@@@@@@@@@@@@@@@@@@@@
%@@@@@@@@@@@@@@@@@@@@@@@@@@@@@@@@@@@@@@@@@@@@@@@@@@@@@@@@@@
%----------------------------------------------------------
\par First we note that for {\it any} choice of objects in parts (i), (ii) and (iii) of Theorem \reff{REF-MAIN}, the quantity in the right-hand side of the theorem's equivalence does not exceed $C(p)\|f\|_{\LTP|_E}$. This follows from the results of Subsection 3.1.
%----------------------------------------------------------
\par In Section 6 we construct a certain family of squares $\Kc_E$ which plays  an important role in our construction. In Subsection 8.2 we have assigned to $\Kc_E$ a one-to-one mapping $\Delta$ from $\Kc_E$ into the family $\Tri(E)$ of all triangles with vertices in $E$. We know that for each $K\in\Kc_E$ we have $\diam \Delta(K)\sim \diam K$ and  $\Delta(K)\subset \gamma K$.
%----------------------------------------------------------
\par Then we construct a measure $\mu_E$ on $\RT$ with $\supp\mu_E\subset \Cc_E$, see \rf{D-MU-S}. Here $\Cc_E$ is the set of centers of squares from $\Kc_E$, see \rf{CC-E}.
(Recall that we call $\mu_E$ the Menger curvature measure generated by $E$.)
%----------------------------------------------------------
\par Prove that
%----------------------------------------------------------
\bel{N-KE}
\#\Kc_E\le C \,\# E.
\ee
%----------------------------------------------------------
In fact, each lacuna $L\in\Lc_E$ generates a family of squares from $\Kc_E$ using only its {\it contacting} squares of $L$, i.e., the squares which have common points with the squares from other lacunae. In other words, each lacuna $L'\in\Lc_E$ such that $L'\lr L$ generates (together with $L$) a square $K\in\Kc_E$.
%----------------------------------------------------------
\par But, by Proposition \reff{PRL-1}, the number of lacunae $L'$ which contact with $L$ is bounded by an absolute constant $C$. Thus
%----------------------------------------------------------
$$
\#\Kc_E\le C \,\# \Lc_E.
$$
%----------------------------------------------------------
In turn, by Corollary \reff{CR-PE}, $\# \Lc_E\le C\# E$ proving \rf{N-KE}.
%----------------------------------------------------------
\par Let us note that the family $\Kc_E$ may be partitioned in a natural way into two families of squares: the family
%----------------------------------------------------------
$$
\Kc_E^{(tr)}:=\{K\in\Kc_E:\Delta(K)~~\text{is a true triangle}\}
$$
%----------------------------------------------------------
and the family
%----------------------------------------------------------
$$
\Kc_E^{(dg)}:=\{K\in\Kc_E:\Delta(K)~~\text{is a degenerate triangle}\}.
$$
%----------------------------------------------------------
%@@@@@@@@@@@@@@@@@@@@@@@@@@@@@@@@@@@@@@@@@@@@@@@@@@@@@@@@@@
%----------------------------------------------------------
\par Let us consider the family $\Kc_E^{(dg)}$. (Of course it can be empty for some $E$.) For each $K\in\Kc_E^{(dg)}$ the vertices of the degenerate triangle $\Delta(K)$ are three collinear points in $\RT$. Let us denote them by $a_1(K),a_2(K),a_3(K)$. We may assume that $a_2(K)\in(a_2(K),a_3(K))$.
%----------------------------------------------------------
\par Let us enumerate the family $\Kc_E^{(dg)}$:
%----------------------------------------------------------
$$
\Kc_E^{(dg)}=\{Q_1,Q_2,...,Q_m\}.
$$
%----------------------------------------------------------
Thus $m\le\#\Kc_E\le C\# E$. We also put
%----------------------------------------------------------
$$
z_1^{(i)}:=a_1(Q_i),~ z_2^{(i)}:=a_2(Q_i),~ z_3^{(i)}:=a_3(Q_i).
$$
%----------------------------------------------------------
Remark that, since $\Delta(K)\subset\gamma K$, we have
%----------------------------------------------------------
$$
z_1^{(i)},z_2^{(i)},z_3^{(i)}\in E\cap (\gamma Q_i).
$$
%----------------------------------------------------------
\par Examining the method of proof suggested in Section 8,
the reader can readily see that the squares $\{Q_1,Q_2,...,Q_m\}$ and the triples of collinear points  $\{z_1^{(i)},z_2^{(i)},z_3^{(i)}\}$ are the required objects in part (i) of Theorem \reff{REF-MAIN}.\medskip
%----------------------------------------------------------
\par We turn to the family $\Kc_E^{(tr)}$. We construct a non-negative Borel measure $\mu_E$ on $\RT$ using the formulas \rf{D-MU-S} and \rf{MU-SET}. Then we construct a ``gradient'' mapping $\Tc(f)$ using the formulas \rf{DGV-1-S} and \rf{DGV-2-S}.
%----------------------------------------------------------
\par At this point we modify the proof of Theorem \reff{MAIN} by applying to the mapping $\Tc(f)$ and the measure $\mu_E$ the following refinement of Theorem \reff{S-CMS}. (Here we present a vector version of this result).
%----------------------------------------------------------
%@@@@@@@@@@@@@@@@@@@@@@@@@@@@@@@@@@@@@@@@@@@@@@@@@@@@@@@@@@
%@@@@@@@@@@@@@@@@@@@@@@@@@@@@@@@@@@@@@@@@@@@@@@@@@@@@@@@@@@
%@@@@@@@@@@@@@@@@@@@@@@@@@@@@@@@@@@@@@@@@@@@@@@@@@@@@@@@@@@
%@@@@@@@@@@@@@@@@@@@@@@@@@@@@@@@@@@@@@@@@@@@@@@@@@@@@@@@@@@
%@@@@@@@@@@@@@@@@@@@@@@@@@@@@@@@@@@@@@@@@@@@@@@@@@@@@@@@@@@
%@@@@@@@@@@@@@@@@@@@@@@@@@@@@@@@@@@@@@@@@@@@@@@@@@@@@@@@@@@
%----------------------------------------------------------
\begin{theorem}\lbl{REF-S1} (\cite{S3}, Subsection 6.1.) Let $2<p<\infty$ and let $\mu$ be a non-trivial non-negative Borel measure on $\RT$.  There exist absolute constants $\gamma=\gamma>0$ and  $N\in\N$, a family $\Qc$ consisting of pairwise disjoint squares and a family $\tQc$ of squares in $\RT$ with covering multiplicity  $M(\tQc)\le N$, mappings
%----------------------------------------------------------
$$
\Qc\ni Q\mapsto Q'\in\tQc~~~\text{and}~~~\Qc\ni Q\mapsto Q''\in\tQc
$$
%----------------------------------------------------------
satisfying the condition
%----------------------------------------------------------
$$
Q'\cup Q''\subset \gamma Q~~~\text{for all}~~~Q\in\Qc,
$$
%----------------------------------------------------------
such that for every mapping  $\Vc\in \vec{L}_{p,loc}(\RT;\mu)$ from $\RT$ into $\RT$ its norm in the space $\VS=\VLOP+\VLPM$ can be calculated (up to a constant depending only on $p$) as follows:
%----------------------------------------------------------
$$
\|\Vc\|_{\VS}\sim\left(\shuge_{Q\in\Qc}\,\,
\frac{(\diam Q)^{n-p}\iint \limits_{Q'\times Q''}
\|\Vc(x)-\Vc(y)|^p\, d\mu(x)d\mu(y)}
{ \{(\diam Q')^{n-p}+\mu(Q')\} \{(\diam Q'')^{n-p}+\mu(Q'')\}}\right)^{\frac1p}.
$$
%----------------------------------------------------------
%@@@@@@@@@@@@@@@@@@@@@@@@@@@@@@@@@@@@@@@@@@@@@@@@@@@@@@@@@@
%----------------------------------------------------------
\end{theorem}
%----------------------------------------------------------
%@@@@@@@@@@@@@@@@@@@@@@@@@@@@@@@@@@@@@@@@@@@@@@@@@@@@@@@@@@
%@@@@@@@@@@@@@@@@@@@@@@@@@@@@@@@@@@@@@@@@@@@@@@@@@@@@@@@@@@
%@@@@@@@@@@@@@@@@@@@@@@@@@@@@@@@@@@@@@@@@@@@@@@@@@@@@@@@@@@
%@@@@@@@@@@@@@@@@@@@@@@@@@@@@@@@@@@@@@@@@@@@@@@@@@@@@@@@@@@
%@@@@@@@@@@@@@@@@@@@@@@@@@@@@@@@@@@@@@@@@@@@@@@@@@@@@@@@@@@
%----------------------------------------------------------
\par The families of squares  $\Qc$ and $\tQc$ and the mappings $Q\mapsto Q'$ and $Q\mapsto Q''$ from the above theorem provides the objects from part (iii) of Theorem \reff{REF-MAIN}.
%----------------------------------------------------------
\par In turn, the objects of part (ii) of this theorem are the family of squares $\Kc_E^{(tr)}$ and the mapping $K\mapsto\Delta(K)$ defined on this family.
%----------------------------------------------------------
\par The method of proof suggested in Sections 7 and 8 and Theorem \reff{REF-S1} shows that in all our considerations in these sections we may restrict ourself only to these particular objects (i.e., families $\Qc$ and $\Kc$, mappings $Q\mapsto Q'$ and $Q\mapsto Q''$, etc., determined in Theorem \reff{REF-S1}).
%----------------------------------------------------------
\par This enables us to put in the proof of Theorem \reff{S-MAIN} the number $\lambda$ to be equal the right-hand side in the equivalence of Theorem \reff{REF-MAIN}. Then, following the method of proof given in Sections 7 and 8, we obtain the required inequality
%----------------------------------------------------------
$$
\|f\|_{\LTP}|_E\le C\,\lambda^{\frac1p}.
$$
%----------------------------------------------------------
\par This completes the proof of Theorem \reff{REF-MAIN}.\bx
%----------------------------------------------------------
\bigskip
%----------------------------------------------------------
%@@@@@@@@@@@@@@@@@@@@@@@@@@@@@@@@@@@@@@@@@@@@@@@@@@@@@@@@@@
%@@@@@@@@@@@@@@@@@@@@@@@@@@@@@@@@@@@@@@@@@@@@@@@@@@@@@@@@@@
%@@@@@@@@@@@@@@@@@@@@@@@@@@@@@@@@@@@@@@@@@@@@@@@@@@@@@@@@@@
%@@@@@@@@@@@@@@@@@@@@@@@@@@@@@@@@@@@@@@@@@@@@@@@@@@@@@@@@@@
%@@@@@@@@@@@@@@@@@@@@@@@@@@@@@@@@@@@@@@@@@@@@@@@@@@@@@@@@@@
%@@@@@@@@@@@@@@@@@@@@@@@@@@@@@@@@@@@@@@@@@@@@@@@@@@@@@@@@@@
%----------------------------------------------------------
\par {\bf 11.2. Proof of Theorem \reff{REF-SP}: sparsification.}
%----------------------------------------------------------
\addtocontents{toc}{~~~~11.2. Proof of Theorem \reff{REF-SP}: sparsification.\hfill \thepage\par}
%----------------------------------------------------------
\par The proof follows the same scheme as the proof of Theorem \reff{REF-MAIN}. The only difference that instead of Theorem \reff{REF-S1} in this case we use the following %----------------------------------------------------------
%@@@@@@@@@@@@@@@@@@@@@@@@@@@@@@@@@@@@@@@@@@@@@@@@@@@@@@@@@@
%@@@@@@@@@@@@@@@@@@@@@@@@@@@@@@@@@@@@@@@@@@@@@@@@@@@@@@@@@@
%@@@@@@@@@@@@@@@@@@@@@@@@@@@@@@@@@@@@@@@@@@@@@@@@@@@@@@@@@@
%@@@@@@@@@@@@@@@@@@@@@@@@@@@@@@@@@@@@@@@@@@@@@@@@@@@@@@@@@@
%@@@@@@@@@@@@@@@@@@@@@@@@@@@@@@@@@@@@@@@@@@@@@@@@@@@@@@@@@@
%@@@@@@@@@@@@@@@@@@@@@@@@@@@@@@@@@@@@@@@@@@@@@@@@@@@@@@@@@@
%----------------------------------------------------------
\begin{theorem}\lbl{REF-S2} (\cite{S3}, Subsection 6.3) Let $\mu$ be a non-trivial non-negative Borel measure on $\RT$, $2<p<\infty$, and let
%----------------------------------------------------------
$$
\VS=\VLOP+\VLPM.
$$
%----------------------------------------------------------
%@@@@@@@@@@@@@@@@@@@@@@@@@@@@@@@@@@@@@@@@@@@@@@@@@@@@@@@@@@
%----------------------------------------------------------
\par There exist families of closed sets $\{G_1,G_2,...\}$ and $\{H_1,H_2,...\}$ in $\RT$ with covering multiplicity $M(\{G_i\}),M(\{H_i\})\le C$ where $C$ is an absolute constant, and a family $\{\lambda_1,\lambda_2,...\}$ of positive numbers such that for every mapping  $\Vc\in \vec{L}_{p,loc}(\RT;\mu)$ from $\RT$ into $\RT$
the following equivalence
%----------------------------------------------------------
$$
\|\Vc\|_{\VS}^p\sim \sum_{i=1}^\infty \,\,
\lambda_i \iint\limits_{G_i\times H_i}
|\Vc(x)-\Vc(y)|^p\, d\mu(x)\,d\mu(y)
$$
%----------------------------------------------------------
holds. The constants of this equivalence depend only on $p$.
%----------------------------------------------------------
\end{theorem}
%----------------------------------------------------------
%@@@@@@@@@@@@@@@@@@@@@@@@@@@@@@@@@@@@@@@@@@@@@@@@@@@@@@@@@@
%@@@@@@@@@@@@@@@@@@@@@@@@@@@@@@@@@@@@@@@@@@@@@@@@@@@@@@@@@@
%@@@@@@@@@@@@@@@@@@@@@@@@@@@@@@@@@@@@@@@@@@@@@@@@@@@@@@@@@@
%@@@@@@@@@@@@@@@@@@@@@@@@@@@@@@@@@@@@@@@@@@@@@@@@@@@@@@@@@@
%@@@@@@@@@@@@@@@@@@@@@@@@@@@@@@@@@@@@@@@@@@@@@@@@@@@@@@@@@@
%----------------------------------------------------------
\par Combining this results with ideas suggested in the proof of Theorem \reff{REF-MAIN}, we obtain the following trace criterion:
%----------------------------------------------------------
%@@@@@@@@@@@@@@@@@@@@@@@@@@@@@@@@@@@@@@@@@@@@@@@@@@@@@@@@@@
%@@@@@@@@@@@@@@@@@@@@@@@@@@@@@@@@@@@@@@@@@@@@@@@@@@@@@@@@@@
%@@@@@@@@@@@@@@@@@@@@@@@@@@@@@@@@@@@@@@@@@@@@@@@@@@@@@@@@@@
%@@@@@@@@@@@@@@@@@@@@@@@@@@@@@@@@@@@@@@@@@@@@@@@@@@@@@@@@@@
%@@@@@@@@@@@@@@@@@@@@@@@@@@@@@@@@@@@@@@@@@@@@@@@@@@@@@@@@@@
%----------------------------------------------------------
\begin{theorem}\lbl{REF-ML} Let $2<p<\infty$ and let $E$ be a finite subset of $\RT$. There exist absolute constants $C>0$ and $N \in \N$ and\,:
\medskip
%----------------------------------------------------------
\par (i) A family $\{Q_i:i=1,...,\ell\}$, $\ell\le C \,\#E$, of pairwise disjoint \sqs and a family
%----------------------------------------------------------
$$
\{z^{(i)}_{1},z^{(i)}_{2},z^{(i)}_{3}\in E: z^{(i)}_{2}\in(z^{(i)}_{1},z^{(i)}_{3}),\,\, i=1,...,\ell\}
$$
%---------------------------------------------------------
of triples of collinear points; \medskip
%----------------------------------------------------------
\par (ii) A family $\Kc$ of pairwise disjoint \sqs with $\#\Kc\le C\,\# E$, and a mapping
%----------------------------------------------------------
$$
\Kc\ni K\mapsto \Delta(K)\in \Tri(E);
$$
%----------------------------------------------------------
\medskip
%----------------------------------------------------------
%@@@@@@@@@@@@@@@@@@@@@@@@@@@@@@@@@@@@@@@@@@@@@@@@@@@@@@@@@@
%----------------------------------------------------------
\par (iii) Families of closed sets
%----------------------------------------------------------
$$
\Uc=\{U_1,U_2,...,U_m\}~~~~\text{and}~~~~
\Vc=\{V_1,V_2,..,V_m\},~~~~m\le C \,\#E,
$$
%----------------------------------------------------------
in $\RT$ with co\-ve\-ring multi\-pli\-city
$M(\Uc)+M(\Vc)\le N,$ and a family of positive numbers $\{\alpha_1,\alpha_2,...,\alpha_m\}$,\bigskip
%----------------------------------------------------------
\par such that for every function $f:E\to\R$ the following equivalence
%----------------------------------------------------------
\be
\|f\|_{\LTP|_E}^p&\sim&
\sbig_{i=1}^{\ell} \left|
\frac{f(z^{(i)}_{1})-f(z^{(i)}_{2})}
{\|z^{(i)}_{1}-z^{(i)}_{2}\|_2}
-\frac{f(z^{(i)}_{2})-f(z^{(i)}_{3})}
{\|z^{(i)}_{2}-z^{(i)}_{3}\|_2}
\right|^p
(\diam Q_i)^{2-p}\nn\\
&+&
\smed\limits_{i=1}^m\,\,\,\,
\alpha_i\,\,S_p(f:U_i,V_i\,;\Kc)
\nn
\ee
%----------------------------------------------------------
holds. The constants of this equivalence depend only on $p$.
%----------------------------------------------------------
\end{theorem}
%----------------------------------------------------------
%@@@@@@@@@@@@@@@@@@@@@@@@@@@@@@@@@@@@@@@@@@@@@@@@@@@@@@@@@@
%@@@@@@@@@@@@@@@@@@@@@@@@@@@@@@@@@@@@@@@@@@@@@@@@@@@@@@@@@@
%@@@@@@@@@@@@@@@@@@@@@@@@@@@@@@@@@@@@@@@@@@@@@@@@@@@@@@@@@@
%@@@@@@@@@@@@@@@@@@@@@@@@@@@@@@@@@@@@@@@@@@@@@@@@@@@@@@@@@@
%@@@@@@@@@@@@@@@@@@@@@@@@@@@@@@@@@@@@@@@@@@@@@@@@@@@@@@@@@@
%----------------------------------------------------------
\par Recall that given sets $U,V\subset\RT$
%----------------------------------------------------------
$$
S_p(f:U,V;\Kc):=
\sum_{\substack {K'\in\,\Kc\\c_{K'}\in U}}\,\,
\sum_{\substack {K''\in\,\Kc\\c_{K''}\in V}}
\|\nabla P_{\Delta(K')}[f]-
\nabla P_{\Delta(K'')}[f]\|^p\,
\cu_{\Delta(K')}^p |K'|\,\cu_{\Delta(K'')}^p |K''|.
$$
%----------------------------------------------------------
\par We are in a position to finish the proof of Theorem \reff{REF-SP}.
%----------------------------------------------------------
\par First we note that the first sum in the above equivalence consists of at most $C\,\# E$ linear functionals to the power $p$ each depending on $3$ values of $f$ on $E$.
%----------------------------------------------------------
\par Let us determine similar linear functionals for the second sum. (But in this case each of these functionals will depend on at most {\it six} points of $E$.)
%----------------------------------------------------------
\par We shall do this using the so-called {\it spectral sparsification} of the quantity $S_p(f:U_i,V_i\,;\Kc)$, $i=1,2,...m.$ We will be needed a certain special form of matrix sparsification related to the existence of  $p$-sparsifiers of non-negative product matrices.
%----------------------------------------------------------
%@@@@@@@@@@@@@@@@@@@@@@@@@@@@@@@@@@@@@@@@@@@@@@@@@@@@@@@@@@
%@@@@@@@@@@@@@@@@@@@@@@@@@@@@@@@@@@@@@@@@@@@@@@@@@@@@@@@@@@
%@@@@@@@@@@@@@@@@@@@@@@@@@@@@@@@@@@@@@@@@@@@@@@@@@@@@@@@@@@
%@@@@@@@@@@@@@@@@@@@@@@@@@@@@@@@@@@@@@@@@@@@@@@@@@@@@@@@@@@
%@@@@@@@@@@@@@@@@@@@@@@@@@@@@@@@@@@@@@@@@@@@@@@@@@@@@@@@@@@
%----------------------------------------------------------
\begin{theorem}(B. Klartag) \lbl{SPR-P} Let $2<p<\infty.$ Let $G=(g_{i}g_{j})$ be $n\times n$ matrix with $g_{i}\ge 0, i=1,...,n$. Then there exists an $n\times n$ matrix $H=(h_{ij})$, $h_{ij}\ge 0$, such that:
%----------------------------------------------------------
\par (i) $\supp(H)\subseteq\supp (G)$~;
%----------------------------------------------------------
\par (ii) $\#(\supp(H))\le C\, n$;
%----------------------------------------------------------
\par (iii)  For every $x\in\RN$ we have
%----------------------------------------------------------
$$
\sum_{i=1}^n\sum_{j=1}^n g_{i}g_{j}|x_i-x_j|^p\le
\sum_{i=1}^n\sum_{j=1}^n h_{ij}|x_i-x_j|^p
\le C(p)\sum_{i=1}^n\sum_{j=1}^n |x_i-x_j|^p.
$$
%----------------------------------------------------------
%@@@@@@@@@@@@@@@@@@@@@@@@@@@@@@@@@@@@@@@@@@@@@@@@@@@@@@@@@@
%----------------------------------------------------------
\end{theorem}
%----------------------------------------------------------
\par We refer to $H$ as a $p$-sparsifier of $G$.
%----------------------------------------------------------
\par We recall that the remarkable sparsification theorem of J. D. Batson, D. A. Spielman and N. Srivastava \cite{BSS} states that for every (not necessarily product) non-negative matrix $G=(g_ig_j)$ there exists a ``good'' $2$-sparsifier.
%----------------------------------------------------------
\par The proof of Theorem \reff{SPR-P} follows from the next statement proven by B. Klartag: {\it if $p\in(2,\infty)$ and $H$ is a $2$-sparsifier of a non-negative product matrix $G$} (with a certain constant in (iii)), {\it then $H$ is a $p$-sparsifier of $G$} (with a bigger constant depending on $p$). A. Naor \cite{N} showed that if $G$ is not a product non-negative matrix, then in general it is not true.
%----------------------------------------------------------
\par We note that the work \cite{BSS} have been remarked by C. Fefferman, A. Israel and G. K. Luli \cite{FIL} as a useful tool in solution to Problem \reff{RT-DEPTH} and related problems.
%----------------------------------------------------------
\par Let us apply Theorem \reff{SPR-P} to the quantity $S_p(f:U,V\,;\Kc)$ whenever $U$ and $V$ are two arbitrary subsets of $\RT$. We obtain the following
%----------------------------------------------------------
%@@@@@@@@@@@@@@@@@@@@@@@@@@@@@@@@@@@@@@@@@@@@@@@@@@@@@@@@@@
%@@@@@@@@@@@@@@@@@@@@@@@@@@@@@@@@@@@@@@@@@@@@@@@@@@@@@@@@@@
%@@@@@@@@@@@@@@@@@@@@@@@@@@@@@@@@@@@@@@@@@@@@@@@@@@@@@@@@@@
%@@@@@@@@@@@@@@@@@@@@@@@@@@@@@@@@@@@@@@@@@@@@@@@@@@@@@@@@@@
%@@@@@@@@@@@@@@@@@@@@@@@@@@@@@@@@@@@@@@@@@@@@@@@@@@@@@@@@@@
%----------------------------------------------------------
\begin{theorem}\lbl{SPR-SPUV} Let $2<p<\infty.$ There exists a function $h=h(K',K'')\ge 0$ defined on the set
%----------------------------------------------------------
$$
\{K'\in\Kc: c_{K'}\in U\}\times\{K''\in\Kc: c_{K''}\in V\}
$$
%----------------------------------------------------------
which satisfies the following conditions:\medskip
%----------------------------------------------------------
\par (i) $h(K',K'')\ne 0$ at most for $C\#(U\cup V)$ pairs $(K',K'')$;\medskip
%----------------------------------------------------------
\par (ii) The following equivalence
%----------------------------------------------------------
$$
S_p(f:U,V;\Kc)\sim H_p(f:U,V;\Kc)
$$
%----------------------------------------------------------
holds with constants depending only on $p$. Here
%----------------------------------------------------------
$$
H_p(f:U,V;\Kc):=
\sum_{\substack {K'\in\,\Kc\\c_{K'}\in U}}\,\,
\sum_{\substack {K''\in\,\Kc\\c_{K''}\in V}}
h(K',K'')\,\|\nabla P_{\Delta(K')}[f]-
\nabla P_{\Delta(K'')}[f]\|^p.
$$
%----------------------------------------------------------
%@@@@@@@@@@@@@@@@@@@@@@@@@@@@@@@@@@@@@@@@@@@@@@@@@@@@@@@@@@
%----------------------------------------------------------
\end{theorem}
%----------------------------------------------------------
\par Using this theorem we finish the proof of Theorem \reff{REF-SP} as follows.
%----------------------------------------------------------
\par Given a vector $x=(x_1,x_2)\in\RT$ we put $(x)_1:=x_1$ and $(x)_2:=x_2$. Then for each $K',K''\in\Kc$ we have the following equivalence
%----------------------------------------------------------
$$
\|\nabla P_{\Delta(K')}[f]-
\nabla P_{\Delta(K'')}[f]\|^p\sim |\lambda_1(f;K',K'')|^p+
|\lambda_2(f;K',K'')|^p
$$
%----------------------------------------------------------
where
%---------------------------------------------------------- %@@@@@@@@@@@@@@@@@@@@@@@@@@@@@@@@@@@@@@@@@@@@@@@@@@@@@@@@@@
$$
\lambda_i(f;K',K''):=(\nabla P_{\Delta(K')}[f]-
\nabla P_{\Delta(K'')}[f])_j,~~~~j=1,2.
$$
%------------------------------------------------------
Clearly, $\lambda_1(f;K',K'')$ and $\lambda_2(f;K',K'')$ are linear functional each depending on six values of $f$ on $E$.
%----------------------------------------------------------
\par Let
%----------------------------------------------------------
$$
Y_i:=\{c_K:K\in\,\Kc,c_{K}\in U_i\},~~~i=1,...,m,
$$
%----------------------------------------------------------
and
%----------------------------------------------------------
$$
X_i:=\{c_K:K\in\,\Kc,c_{K}\in U_i\},~~~i=1,...,m,
$$
%----------------------------------------------------------
and let $\Yc=\{Y_i\}$ and $\Xc=\{X_i\}$. Since the families $\Uc=\{U_i\}$ and $\Vc=\{V_i\}$ has covering multiplicity $M(\Uc)+M(\Vc)\le N,$ we conclude that
%----------------------------------------------------------
\bel{CM-N}
M(\Yc)+M(\Xc)\le N
\ee
%----------------------------------------------------------
as well.
%----------------------------------------------------------
\par Thus, by Theorem \reff{SPR-SPUV} each quantity $S_p(f:U_i,V_i\,;\Kc)$ can be represented as a sum of the $p$-powers of at most $C(\#Y_i+\#X_i)$ linear functionals, so that the second sum in the equivalence of Theorem \reff{REF-ML} depends on at most
%----------------------------------------------------------
$$
L\le C\sum_{i=1}^m (\#Y_i+\#X_i)
$$
%----------------------------------------------------------
linear functionals. But, by \rf{CM-N},
%----------------------------------------------------------
$$
\sum_{i=1}^m (\#Y_i+\#X_i)\le
N\#\left\{\bigcup_{i=1}^m(Y_i\cup X_i)\right\}\le N\#\Kc\,.
$$
%----------------------------------------------------------
It remains to note that, by part (ii) of Theorem \reff{REF-MAIN}, $\#\Kc\le C\,\# E$ so that the number $L$ of the linear functionals is bounded by
%----------------------------------------------------------
$$
L\le C\,\# E.
$$
%----------------------------------------------------------
\par Theorem \reff{REF-SP} is completely proved.\bx
%----------------------------------------------------------
%@@@@@@@@@@@@@@@@@@@@@@@@@@@@@@@@@@@@@@@@@@@@@@@@@@@@@@@@@@
%@@@@@@@@@@@@@@@@@@@@@@@@@@@@@@@@@@@@@@@@@@@@@@@@@@@@@@@@@@
%@@@@@@@@@@@@@@@@@@@@@@@@@@@@@@@@@@@@@@@@@@@@@@@@@@@@@@@@@@
%@@@@@@@@@@@@@@@@@@@@@@@@@@@@@@@@@@@@@@@@@@@@@@@@@@@@@@@@@@
%@@@@@@@@@@@@@@@@@@@@@@@@@@@@@@@@@@@@@@@@@@@@@@@@@@@@@@@@@@
%@@@@@@@@@@@@@@@@@@@@@@@@@@@@@@@@@@@@@@@@@@@@@@@@@@@@@@@@@@
%----------------------------------------------------------
%@@@@@@@@@@@@@@@@@@@@@@@@@@@@@@@@@@@@@@@@@@@@@@@@@@@@@@@@@@
\bigskip
%@@@@@@@@@@@@@@@@@@@@@@@@@@@@@@@@@@@@@@@@@@@@@@@@@@@@@@@@@@
%----------------------------------------------------------
\par {\bf 11.3. A sketch of the proof of Theorem \reff{NEW-MAIN}.}
%----------------------------------------------------------
\addtocontents{toc}{~~~~11.3. A sketch of the proof of Theorem \reff{NEW-MAIN}.\hfill\thepage\\\par}
%----------------------------------------------------------
We note that, by definitions \rf{C-DS} and  \rf{C-DS}, given a disk $B\subset\R^2$ the condition $\cu_p(B:\Dc,\Delta)\le c_B$ is equivalent to the inequality
%----------------------------------------------------------
$$
\sum\limits_{c_D\in B}\cu^p_{\Delta(D)}|D|\le
|B|/r_B^p=\pi r_B^{2-p}.
$$
%----------------------------------------------------------
We let $\tilde{\sigma}_p(B;\Dc)$ denote the quantity
%----------------------------------------------------------
$$
\tilde{\sigma}_p(B;\Dc):= \sum\limits_{c_D\in B} \cu^p_{\Delta(D)}|D|.
$$
%----------------------------------------------------------
This quantity is an analog of the quantity $\sigma_p(Q;\Kc)$ defined by \rf{SGM}. Thus the condition \rf{B-CV} of the theorem is equivalent to the inequality
%----------------------------------------------------------
$$
(\diam B')^{p-2}\tilde{\sigma}_p(B';\Dc)+
(\diam B'')^{p-2}\tilde{\sigma}_p(B'';\Dc)\le C(p).
$$
%----------------------------------------------------------
Obviously, this inequality is a ``disk'' version of inequality \rf{CV-S}.
%----------------------------------------------------------
\par We note that an analogue of Theorem \reff{NEW-MAIN} for \sqs rather than for disks is also true. Its necessity part follows from the necessity part of Theorem \reff{MAIN}. The sufficiency part of such an analogue follows from Theorem \reff{S-MAIN}.
%----------------------------------------------------------
\par We also note that we can slightly modify the condition \rf{CV-S} in formulation of Theorem \reff{S-MAIN}. More specifically, we can replace this condition by a more general one:
%----------------------------------------------------------
\bel{T-Z}
(\diam Q')^{p-2}\sigma_p(Q';\Kc)+
(\diam Q'')^{p-2}\sigma_p(Q'';\Kc)\le \eta
\ee
%----------------------------------------------------------
where $\eta$ is a positive constant. The result of Theorem \reff{S-MAIN} remains true after such a modification, but with constant $C$ in inequality \rf{SM-H} depending on $p$ and also on $\eta$. In fact, in the proof of this theorem we use inequality \rf{T-Z} only to verify condition \rf{R-W} in Theorem \reff{S-CMS}. But obviously this condition holds with $\tau=\eta$.
%----------------------------------------------------------
\par The proof of Theorem \reff{NEW-MAIN} follows the same scheme: we repeat our considerations for families of disks rather than squares. There are no any technical difficulties in such a generalization of the methods and ideas developed for \sqs to the case of disks. We will only  remark two places in the proof where certain non-trivial changes should be done. %----------------------------------------------------------
\par First of them relates to an analogue of the Whitney covering Theorem \reff{Wcov} for disks. Of course, in this case we can not cover the open set $\RT\setminus E$ by {\it non-overlapping} disks $D$ such that
$\diam D\sim\dist(D,E)$. Nevertheless for our purpose it suffice to cover $\RT\setminus E$ by a family $\widetilde{W}_E$ of disks whose {\it covering multiplicity} is bounded by an absolute constant $N$. In other words, every point $x\in\RT$ is covered at most $N$
disks from the family $\widetilde{W}_E$.
%----------------------------------------------------------
\par The existence of a Whitney-type covering of such a kind, i.e., a covering of an open set by a family of Whitney disks with finite multiplicity, follows from a general result proven by M. Guzman \cite{G}. (In turn, this result is based on the Besicovitch  covering theorem \cite{Be}.)
%----------------------------------------------------------
\par Our second remark relates to Theorem \reff{S-CMS} which is an important ingredient of the proof of Theorem \reff{S-MAIN}. For its analogue for disks (rather than squares) we refer the reader to the paper \cite{S3}. (See there Remark 6.7 and Theorem 6.8.)\bx\medskip
%----------------------------------------------------------
%@@@@@@@@@@@@@@@@@@@@@@@@@@@@@@@@@@@@@@@@@@@@@@@@@@@@@@@@@@
%@@@@@@@@@@@@@@@@@@@@@@@@@@@@@@@@@@@@@@@@@@@@@@@@@@@@@@@@@@
%@@@@@@@@@@@@@@@@@@@@@@@@@@@@@@@@@@@@@@@@@@@@@@@@@@@@@@@@@@
%@@@@@@@@@@@@@@@@@@@@@@@@@@@@@@@@@@@@@@@@@@@@@@@@@@@@@@@@@@
%@@@@@@@@@@@@@@@@@@@@@@@@@@@@@@@@@@@@@@@@@@@@@@@@@@@@@@@@@@
%@@@@@@@@@@@@@@@@@@@@@@@@@@@@@@@@@@@@@@@@@@@@@@@@@@@@@@@@@@
%----------------------------------------------------------
\begin{remark}{\em Analyzing the proof of the sufficiency we note that either the disk $B'$ or the disk $B''$ from part (ii) of Theorem \reff{NEW-MAIN} is an {\it``important''} disk, i.e., either equivalence $\cu_p(B':\Dc,\Delta)\sim \cu_{B'}$ or equivalence $\cu_p(B'':\Dc,\Delta)\sim \cu_{B''}$ holds with absolute constants. This enables us to replace the condition
\rf{B-CV} by a slightly stronger condition
%----------------------------------------------------------
$$
\alpha\le \cu_p(B':\Dc,\Delta)/\cu_{B'}+\cu_p(B'':\Dc,\Delta)/ \cu_{B''}\le 1
$$
%----------------------------------------------------------
where $\alpha$ is an absolute positive constant.\rbx}
%----------------------------------------------------------
\end{remark}
%----------------------------------------------------------
%@@@@@@@@@@@@@@@@@@@@@@@@@@@@@@@@@@@@@@@@@@@@@@@@@@@@@@@@@@
%@@@@@@@@@@@@@@@@@@@@@@@@@@@@@@@@@@@@@@@@@@@@@@@@@@@@@@@@@@
%@@@@@@@@@@@@@@@@@@@@@@@@@@@@@@@@@@@@@@@@@@@@@@@@@@@@@@@@@@
%@@@@@@@@@@@@@@@@@@@@@@@@@@@@@@@@@@@@@@@@@@@@@@@@@@@@@@@@@@
%@@@@@@@@@@@@@@@@@@@@@@@@@@@@@@@@@@@@@@@@@@@@@@@@@@@@@@@@@@
%@@@@@@@@@@@@@@@@@@@@@@@@@@@@@@@@@@@@@@@@@@@@@@@@@@@@@@@@@@
%----------------------------------------------------------
\begin{remark} {\em In part (ii) of Theorem \reff{MAIN} we require that the \sqs $Q'$ and $Q''$ belong to the family $\Qc$. Note that the theorem's result remains true if we replace this condition with the following one: $Q',Q''\in\tQc$ where $\tQc$ is an arbitrary family of pairwise disjoint squares in $\RT$.
%----------------------------------------------------------
\par The necessity of this modification of part (ii) of Theorem \reff{MAIN} directly follows from the method of proof suggested in Section 3. In particular, in formulation of Theorem \reff{SL-DEC} the families $\Ac$ and $\Sc$ may be different. This enables us in the proof of Proposition \reff{M-PL} to replace \rf{ASQ} with the following definition:
%----------------------------------------------------------
$$
\Ac=\Qc,~~~~\Sc=\tQc.
$$
%----------------------------------------------------------
\par In turn, the sufficiency follows from the proof of the sufficiency part of Theorem \reff{MAIN} where we can put $\tQc=\Qc$.\rbx}
%----------------------------------------------------------
\end{remark}
%----------------------------------------------------------
%@@@@@@@@@@@@@@@@@@@@@@@@@@@@@@@@@@@@@@@@@@@@@@@@@@@@@@@@@@
%@@@@@@@@@@@@@@@@@@@@@@@@@@@@@@@@@@@@@@@@@@@@@@@@@@@@@@@@@@
%@@@@@@@@@@@@@@@@@@@@@@@@@@@@@@@@@@@@@@@@@@@@@@@@@@@@@@@@@@
%@@@@@@@@@@@@@@@@@@@@@@@@@@@@@@@@@@@@@@@@@@@@@@@@@@@@@@@@@@
%@@@@@@@@@@@@@@@@@@@@@@@@@@@@@@@@@@@@@@@@@@@@@@@@@@@@@@@@@@
%@@@@@@@@@@@@@@@@@@@@@@@@@@@@@@@@@@@@@@@@@@@@@@@@@@@@@@@@@@
%----------------------------------------------------------
%@@@@@@@@@@@@@@@@@@@@@@@@@@@@@@@@@@@@@@@@@@@@@@@@@@@@@@@@@@
\bigskip
%@@@@@@@@@@@@@@@@@@@@@@@@@@@@@@@@@@@@@@@@@@@@@@@@@@@@@@@@@@
%----------------------------------------------------------
%@@@@@@@@@@@@@@@@@@@@@@@@@@@@@@@@@@@@@@@@@@@@@@@@@@@@@@@@@@
%@@@@@@@@@@@@@@@@@@@@@@@@@@@@@@@@@@@@@@@@@@@@@@@@@@@@@@@@@@
%@@@@@@@@@@@@@@@@@@@@@@@@@@@@@@@@@@@@@@@@@@@@@@@@@@@@@@@@@@
%@@@@@@@@@@@@@@@@@@@@@@@@@@@@@@@@@@@@@@@@@@@@@@@@@@@@@@@@@@
%@@@@@@@@@@@@@@@@@@@@@@@@@@@@@@@@@@@@@@@@@@@@@@@@@@@@@@@@@@
%@@@@@@@@@@@@@@@@@@@@@@@@@@@@@@@@@@@@@@@@@@@@@@@@@@@@@@@@@@
%----------------------------------------------------------
%&&&&&&&&&&&&&&&&&&&&&&&&&&&&&&&&&&&&&&&&&&&&&&&&&&&&&&&&&&
%                                                         &
%                      REFERENCES                         &
%_________________________________________________________&
%&&&&&&&&&&&&&&&&&&&&&&&&&&&&&&&&&&&&&&&&&&&&&&&&&&&&&&&&&&

%@@@@@@@@@@@@@@@@@@@@@@@@@@@@@@@@@@@@@@@@@@@@@@@@@@@@@@@@@@
\end{document}